\documentclass[twoside]{report}
\usepackage{amsmath,amsfonts,amssymb,euscript}
\usepackage[british]{babel}
\usepackage{amsthm}
\usepackage[dvipdfmx]{graphicx}
\usepackage{makeidx}
\usepackage{array}

\usepackage{longtable}
\usepackage{slashbox}
\newcommand{\C}{{\mathbb C}}
\newcommand{\U}{{\mathbb U}}

\let\ncm=\newcommand   \let\rncm=\renewcommand \rncm{\o}{\operatorname}
\ncm{\smb}{\symbol}    \ncm{\ctl}{\centerline} \ncm{\nn}{\noindent}
\rncm{\t}{\text}       \ncm{\q}{\quad}         \ncm{\qq}{\qquad}
\ncm{\lt}{\limits}     

\ncm{\ls}{\leqslant}     \ncm{\gs}{\geqslant}
\ncm{\pt}{\partial}     \ncm{\wh}{\widehat}      \rncm{\ss}{\subset}
\ncm{\sq}{\subseteq}    \ncm{\fr}{\frac}         \ncm{\dfr}{\dfrac}
\ncm{\la}{\langle}      \ncm{\ra}{\rangle}       \ncm{\fy}{\infty}
\ncm{\x}{\times}        \ncm{\nem}{\varnothing}	 \ncm{\ol}{\overline}
\ncm{\ovec}{\overrightarrow} \ncm{\Lra}{\Longrightarrow} \ncm{\ul}{\underline}

\ncm{\al}{\alpha}       \ncm{\bt}{\beta}         \ncm{\gam}{\gamma}
\ncm{\dl}{\delta}       \ncm{\eps}{\varepsilon}   \ncm{\vt}{\vartheta}
\ncm{\lm}{\lambda}      \ncm{\sg}{\sigma}        \ncm{\vfi}{\varphi}
\ncm{\om}{\omega}       \ncm{\ta}{\theta}        \ncm{\vr}{\varrho}
\ncm{\e}{\epsilon}      \ncm{\zt}{\zeta}         \ncm{\vsg}{\varsigma}

\ncm{\Gam}{\Gamma}       \ncm{\Dl}{\Delta}        \ncm{\Tt}{\Theta}
\ncm{\Lm}{\Lambda}      \ncm{\Sg}{\Sigma}        \ncm{\Om}{\Omega}
\ncm{\vFi}{\varPhi}     \ncm{\vTt}{\varTheta}    \ncm{\vSg}{\varSigma}
\ncm{\vOm}{\varOmega}   \ncm{\vP}{\varPsi}       \ncm{\vLm}{\varLambda}
\ncm{\vPi}{\varPi}      \ncm{\vUp}{\varUpsilon}  \ncm{\vGm}{\varGamma}
\ncm{\vDl}{\varDelta}   \ncm{\vXi}{\varXi}       \ncm{\Up}{\Upsilon}
\ncm{\Fi}{\Phi}

\def\cl{\o{cl}}                  
     \def\mes{\o{mes}}         \def\const{\o{const}}
	\def\max{\o{max}}	  \def\min{\o{min}}

\def\sign{\o{sign}}

\ncm{\gA}{\mathfrak A}   \ncm{\gB}{\mathfrak B}    \ncm{\gC}{\mathfrak C}
\ncm{\gD}{\mathfrak D}   \ncm{\gF}{\mathfrak F}    \ncm{\gI}{\mathfrak I}
\ncm{\gJ}{\mathfrak J}   \ncm{\gK}{\mathfrak K}    \ncm{\gL}{\mathfrak L}
\ncm{\gM}{\mathfrak M}   \ncm{\gN}{\mathfrak N}    \ncm{\gO}{\mathfrak O}
\ncm{\gP}{\mathfrak P}   \ncm{\gT}{\mathfrak T}    \ncm{\gU}{\mathfrak U}
\ncm{\gX}{\mathfrak X}   \ncm{\gY}{\mathfrak Y}    \ncm{\gZ}{\mathfrak Z}
\ncm{\ga}{\mathfrak a}   \ncm{\gb}{\mathfrak b}    \ncm{\gc}{\mathfrak c}
\ncm{\gd}{\mathfrak d}   \ncm{\gk}{\mathfrak k}    \ncm{\gl}{\mathfrak l}
\ncm{\gm}{\mathfrak m}   \ncm{\gn}{\mathfrak n}    \ncm{\gq}{\mathfrak q}
\ncm{\gr}{\mathfrak r}   \ncm{\gss}{\mathfrak s}    \ncm{\gx}{\mathfrak x}
\ncm{\gy}{\mathfrak y}   \ncm{\gz}{\mathfrak z}    \ncm{\gS}{\mathfrak S}
\ncm{\gV}{\mathfrak V}   \ncm{\gW}{\mathfrak W}    \ncm{\gt}{\mathfrak t}
\ncm{\gu}{\mathfrak u}   \ncm{\gv}{\mathfrak v}    \ncm{\gw}{\mathfrak w}

\ncm{\A}{\mathbb A}       \ncm{\B}{\mathbb B}      \rncm{\C}{\mathbb C}
\ncm{\aH}{\mathbb H}      \ncm{\I}{\mathbb I}      \ncm{\K}{\mathbb K}
\rncm{\L}{\mathbb L}      \ncm{\N}{\mathbb N}      \ncm{\aO}{\mathbb O}
\ncm{\Q}{\mathbb Q}       \ncm{\R}{\mathbb R}      \rncm{\U}{\mathbb U}
\ncm{\Z}{\mathbb Z}       \ncm{\ak}{\varkappa}     \ncm{\D}{\mathbb D}
\ncm{\E}{\mathbb E}       \ncm{\F}{\mathbb F}      \ncm{\G}{\mathbb G}
\ncm{\M}{\mathbb M}       \rncm{\P}{\mathbb P}     \ncm{\J}{\mathbb J}
\ncm{\aS}{\mathbb S}      \ncm{\T}{\mathbb T}      \ncm{\V}{\mathbb V}
\ncm{\W}{\mathbb W}       \ncm{\X}{\mathbb X}      \ncm{\Y}{\mathbb Y}

\ncm\pA{\mathcal A}       \ncm\pB{\mathcal B}      \ncm\pC{\mathcal C}
\ncm\pD{\mathcal D}       \ncm\pE{\mathcal E}      \ncm\pF{\mathcal F}
\ncm\pG{\mathcal G}       \ncm\pH{\mathcal H}      \ncm\pI{\mathcal I}
\ncm\pK{\mathcal K}       \ncm\pL{\mathcal L}      \ncm\pM{\mathcal M}
\ncm\pN{\mathcal N}       \ncm\pO{\mathcal O}      \ncm\pP{\mathcal P}
\ncm\pR{\mathcal R}       \ncm\pS{\mathcal S}      \ncm\pT{\mathcal T}
\ncm\pU{\mathcal U}       \ncm\pV{\mathcal V}      \ncm\pW{\mathcal W}

\input method.sty

\makeatletter
\def\@seccntformat#1{\csname the#1\endcsname.\space}
\renewcommand{\section}{\@startsection{section}{1}{0pt}{-3.5ex plus -1ex minus -.2ex}%
{2.3ex plus .2ex}{\Large\bf}}

\renewcommand{\subsection}{\@startsection{subsection}{2}{\parindent}{0pt}
{0pt}{\bf}}
\makeatother
\let\poonkt=\subsection

\newtheorem{lem}{\bf{Lemma}}[section]
\newtheorem{teo}{\bf{Theorem}}[section]

\newtheorem{sle}{\bf{Corollary}}[section]

\newtheorem{example}{Example}

\newtheorem{zam}{\bf{Comment}}[section]
\newtheorem{utv}{\bf{Proposition}}[section]

\newcommand{\doc}{Proof.}

\begin{document}


\thispagestyle{empty}

\hphantom{T}

\vspace{30mm}
 
\large
\begin{center}
{\bf V.\,Ya. Derr}
\end{center}

\vspace{20mm}
\normalsize
\begin{center}
{\bf {\Large REGULATED \;\; FUNCTIONS.\\
\vskip 10pt
$*$-INTEGRAL}}
\end{center}
\vspace{10mm}

\vfill
\begin{center}
{Izhevsk 2023}
\end{center}
\newpage

\thispagestyle{empty}


\vspace{20mm}
\begin{center}
{V.\,Ya. Derr}
\end{center}

\vspace{20mm}
\begin{center}
{\bf {\large   REGULATED \;\; FUNCTIONS.\\
\vskip 20pt
$*$-INTEGRAL}}\\
\end{center}

\vspace{13mm}

\vfill
\begin{center}
{Izhevsk 2023}
\end{center}
\newpage
\thispagestyle{empty}

\hphantom{r} 



\vspace{40mm}
{\bf V.\,Ya\,Derr} \\  \quad  Izhevsk 2023

%

\pagestyle{empty}

\newpage

\pagestyle{empty}
\renewcommand{\contentsname}{\Large \bf Content}
\thispagestyle{empty}
\tableofcontents
\thispagestyle{empty}

\newpage

\pagestyle{plain}






\begin{center}
\begin{large}
\section*{Foreword}
\end{large}
\end{center}
\addcontentsline{toc}{section}{Foreword}

The present book gives a systematic overview of function theory and the theory of Stieltjes integral.
Most of the material in this book either appeared only in journal articles or was not previously published, and is the content of a number of
special courses given by the author  for graduate and undergraduate students of Udmurt State University.
The style of the book differs significantly from the other books by the author (see \cite{derr08}, \cite{derr13}, \cite{derr21}), 
which contained, in addition to the theoretical material, a large number of exercices of varying difficulty together with solutions to most of these exercices. 
This book is a crossover between a textbook and a monograph. It is the first part of the ``big book'' dedicated to selected sections of the theory of
functions of real variable and functional analysis.

To make the exposition more systematic and to avoid overloading the reader with references, we included along with the new results a number of classical results (at the same time, references to widely used textbooks are not always given). For the same purpose, 
the book begins with a functional-analytical introduction, which contains what is necessary for reading the main body of the book (and for reading the aforementioned ``big book'').

Our account of the theory of {\it functions of bounded variation} is more detailed than what can be found 
in the classical textbooks (see e.g. \cite{ficht3},\,\cite{nat97}).
In particular, we included some not so well known results for the functions of bounded variation from an article of\linebreak
 A.Yu.\,Levin \;\cite{Lev}.

One of the central concepts in this book is that of a {\it regulated} or {\it proper} function (that is, a  function 
 having finite one-sided limits at each point).  (\cite{dje64},\, \cite{Schw},\, \cite{Tolst84};\, see \cite{derr97}~--- \cite{derr20},\, 
 \cite{DK05},\, \cite{DK06},\,  see also \cite{BarRod22},\,\cite{Rod07}~--- \cite{Rod21}). The theory of proper functions is presented 
 in detail,  including previously unpublished results. We describe dense subsets of the Banach space ${\bf{R}}$ of  proper functions, 
 and discuss some useful subspaces of this space. We introduce and study the space ${\bf{N}}$ of $\sigma$-continuous functions (i.e.\,bounded fonctions having at most countably 
 many discontinuity points) and find its dense subset. 
(One has ${\bf{R}}\subset {\bf{N}}$.)  It turned out that spaces 
 ${\bf{R}}$ and ${\bf{N}}$ play an important role in the theory of the classical Riemann --- Stieltjes integral ($RS$-integral).

We present the most complete, in comparison with the classical textbooks, theory of the $RS$-integral. This includes the necessary and sufficient condition
for the existence of this integral (an analogue of the Lebesgue theorem), as well as some sufficient conditions  (for instance, we prove  
integrability of a $\sigma$-continuous function with respect to a continuous function of bounded variation). We also obtain some new limit theorems. We introduce the notion of an 
exact pair (of spaces) for the existence of the $RS$-integral and find sond previously unknown exact pairs.

The notion of  ``$*$-integral'', which allows one to integrate discontinuous functions (integrable in the sense of Riemann) with 
respect to
discontinuous functions (of bounded variation), plays a central role in this book.
As is well known, the RS-integral does not exist if the integrated and the integrating
functions have even a single common point of discontinuity (see, for example, \cite[sect.584]{ficht3}. Meanwhile, in applications, one often has to integrate
discontinuous functions with respect to discontinuous functions (see \cite{derr88},\, \cite{fil},\, \cite{Lev} and many others). For this purpose, one can use the
Perron---Stieltjes integral (PS-integral). However, for all the complexity of defining and handling it, it essentially requires
that both functions --- the integrated and the integrating --- have finite total variation (see \, \cite{Stv79},\,
see also \cite{Kurz57},\;\cite{Kurz58}). Moreover, for this integral there is no {\it suitable} limit theorem of the form
$$
\;\; \int\limits_{a}^{b} x_n(t)\,dg_n(t)\to \int\limits_{a}^{b} x(t)\,dg(t)\qquad (n\to \infty),\;\;
$$
which is needed for applications (see, for example, \cite{derr88},\,\cite{derr18}).
The Lebesgue--
Stieltjes integral \;(LS-integral)\; has similar disadvantages. A large 
body of literature was devoted to various attempts to overcome these problems (see, for example,
\cite{Atk68}, \,\cite{derr88},\, \cite{DK05},\, \cite{fil},\, \cite{Lev}, \,\cite{Rod07}, \,\cite{ZavSes91} etc).

Strictly speaking, the *-integral is not an integral in the usual sense (because it is not defined using the integral sum).
It is defined as a certain bilinear functional equal to the sum of two ``real'' integrals: 
1) {\it RS-integral} of a $\sigma$-continuous function with respect to the continuous part of a function of bounded variation, and \;
2)\;{\it the sum of a series}, which can be interpreted as an LS-integral of the same \;$\sigma$-continuous function with respect to the discrete part 
of the same functions of bounded variation (this is the sum of products of the values of the integrated function and the jumps of the 
integrating function at the points of discontinuity of the latter). The proposed *-integral retains some ``good'' properties of the RS-integral. 
In particular, for this integral one has Helly's limit thoerem, which does not carry over to the general LS-integral (or PS-integral) 
(see Theorems \ref{ChelliM}--\ref{okpredt}).

The specifics of the proposed definitions allows not to use any measure theory.
When dealing with the *-integral, one only needs the concept of a set of zero measure. 
This circumstance makes possible to have the improved limit theorems (in comparison with the limits theorems for LS- and PS-integrals).
The general notion of (Lebesgue) measure is used only to prove some auxiliary statements of general nature in Sec. {\bf{14.2}},
not directly related to the *-integral theory. Moreover, all the properties of the *-integral, all statements of the theory of this 
integral (see below) are deduced from the above definition.

As a direct application of the *-integral, we present the scheme for constructing solution of 
ordinary linear differential equation with generalized functions in the coefficients
, first published in \cite{derr88}. The scheme consists of replacing the differential equation by a quasi-differential equaton whose coefficients are ordinary functions. The LS-integral used in
\cite{derr88},\, \cite{derr89}, \,\cite{derr18} is now replaced by the *-integral, which allowed us to obtain a new theorem on
$\delta$-correctness of the solution.

Another application of *-integral is given in \S 16. There we show that the general form of a continuous linear functional in spaces of proper or $\sigma$-continuous functions is given by *-integral.

To read the book, it is enough to have the basic understanding of the set theory and linear algebra, differential calculus
for functions of one and several variables, the Riemann integral, the theory of series of real numbers and functions.

The book is intended as a textbook for senior undergraduate and graduate students in mathematics, 
as well as for young teachers. It will be useful for beginning (and not only)
 researchers  in the field of function theory, integration theory, theory of differential equations.

\newpage

\begin{center}
\begin{Large}
{\bf{Chapter I. \; Metric spaces}} 
\end{Large}
\end{center}
\addcontentsline{toc}{section}{Chapter I. \; Metric spaces}

\begin{center}
\begin{large}
\section{Metric spaces}
\end{large}
\end{center}

\poonkt{Definition and examples}\phantom{01234567890123456789} \

Let $\mathfrak M$~ be an arbitrary set. The function $\rho:\mathfrak M\times \mathfrak M \to [0,\infty)$ is called a {\it metric}
(or a {\it metric function}) if it has the following properties:

1) $\rho(x,y)=0$ if and only if $x=y$ (the identity axiom);

2) $\rho(y,x)=\rho(x,y)$ for any $x,y \in \mathfrak M$ (the axiom of symmetry);

3)$\rho(x,y)\leqslant\rho(x,z)+\rho(z,y)$ for any $ x,y,z \in \mathfrak M$ (the triangle inequality).

The pair $(\mathfrak M,\rho)$ is called a {\em metric space}. If the metric $\rho$ is clear form the context, then we simply write
$\mathfrak M$. The elements of a metric space are called {\em points}. The value of $\rho(x,y)$ is called the
{\em distance} between points $x$ and $y$.

From the triangle inequality, we obtain the following {\it quadrilateral \linebreak
inequality}:
{\it for any} $x,x',y,y' \in \mathfrak M$
\begin{equation} \label{met_8}
\bigl| \rho(x,y)-\rho(x',y')\bigr| \leqslant \rho(x,x')+\rho(y,y').
\end{equation}
\doc\,
Let $x,y,x',y'\in\mathfrak M.$ Applying the triangle inequality, we successively obtain
$$
\rho(x,y)\leqslant\rho(x,y')+\rho(y',y)\leqslant\rho(x,y')+\rho(y',x')+\rho( x',y),
$$
$$
\rho(x',y')\leqslant\rho(x',y)+\rho(y,y')\leqslant\rho(x',y)+\rho(y,x)+\rho( x,y'),
$$
where
$$
\rho(x,y)-\rho(x',y')\leqslant \rho(x,y')+\rho(x',y),\; \rho(x',y')-\rho(x,y)\leqslant \rho(x,y')+\rho(x',y).
$$
The resulting inequalities are equivalent to the desired inequality \eqref{met_8}.
\hfill $\square$

{\bf{Examples of metric spaces}}


$1^\circ.$ The set $\mathbb R$ $(\mathbb C)$ of real (complex) numbers endowed with the metric
$\rho_p(x,y)=
|x -y|$
is a metric space.
The axioms of identity and symmetry are obvious. The triangle inequality follows from the properties 
of the absolute value of a real number (of the modulus of a complex number, respectively).

$2^\circ.$ The sets $\mathbb R_p^n$ $(\mathbb C_p^n)$ of $n$-tuples $x=(x_1,x_2,...,x_n)$ of $n\;(n\geqslant 2) $ real
(complex) numbers with the metric
$$
\rho_p(x,y)=\left(\sum\limits_{k=1}^n|x_k -y_k| ^p\right)^{\frac{1}{p}} \qquad (1\leqslant p < \infty).
$$
The  axioms of identity and symmetry are obvious. The triangle inequality follows from the Minkowski inequality
(see the Appendix). Example ${\bf{1}}$ follows from Example ${\bf{2}}$ upon taking $n=1$.

$3^\circ.$ The set of sequences of real or complex numbers ~ $x=(x_1,x_2,\ldots)$ such that the series 
$\sum\limits_{k=1}^{\infty}|x_k|^p \;(1\leqslant p<\infty)$ converges. We define the following metric:
$$
\rho(x,y)=\left(\sum\limits_{k=1}^{\infty}|x_k-y_k|^p\right)^\frac{1}{p}.
$$
The first two axioms of the metric are obvious. The triangle inequality is again a consequence of the Minkowski inequality.
The resulting metric space is denoted by ${\bf{l}}_p$. We denote by ${\bf{l}}_{\infty}={\bf{m}}$ the set of bounded
sequences of real or complex numbers endowed with the metric
\begin{equation} \label{met_3}
\rho(x,y)=\underset{k\in\mathbb N}{\sup}|x_k-y_k|.
\end{equation}
Let us verify the triangle inequality for (\ref{met_3}) (axioms 1 and 2 are obvious).
For each $k\in\mathbb N$, by the properties of the absolute value of a real number (modulus of a complex number), we have
\begin{multline}
\notag
|x_k-y_k|=|x_k-z_k+z_k-y_k| \leqslant |x_k-z_k|+|z_k-y_k|\leqslant \\ \leqslant \underset{k}{\sup}|x_k-z_k|+\underset{k}{\sup}|z_k-y_k|.
\end{multline}
Passing to the supremum in the left side, we obtain the triangle inequality for (\ref{met_3}).

$4^\circ.$ The set of convergent sequences of real or complex numbers equipped with metric (\ref {met_3}) is denoted by ${\bf{c}}$;
the set of sequences of real or complex numbers that converge to zero is
denoted by ${\bf{c}}_0,\; {\bf{c}}_0\subset {\bf{c}}\subset {\bf{m}}$.

$5^\circ.$ The set of arbitrary sequences of real or complex numbers endowed with metric
\begin{equation} \label{met_4}
\rho(x,y)=\sum\limits_{k=1}^{\infty} \frac{1}{2^k} \cdot \frac{|x_k-y_k|}{1+|x_k-y_k |}
\end{equation}
is denoted by ${\bf{s}}$.
The axioms of identity and symmetry for function (\ref{met_4}) are, clearly, satisfied. Let us prove the triangle inequality.
Since function $t \mapsto \frac {t}{1+t}$ is increasing and $|x_k-y_k|\leqslant |x_k-z_k| +|z_k-y_k|,$ we have
\begin{multline}
\notag
\frac{|x_k-y_k|}{1+|x_k-y_k|} \leqslant \frac{|x_k-z_k|+|z_k-y_k|}{1+|x_k-z_k|+|z_k-y_k|} \leqslant \\ \leqslant
\frac{|x_k-z_k|}{1+|x_k-z_k|}+\frac{|z_k-y_k|}{1+|z_k-y_k|}.
\end{multline}
Multiplying the left and the right side of the resulting inequality by $2^{-k}$, summing up in $k$ from $k=1$ to $k=N$ and taking $N \to \infty$,
we arrive at the triangle inequality for (\ref{met_4}).

$6^\circ.$ The set ${\bf{C}}[a,b]$ of continuous functions $x:[a,b]\to \mathbb R$ $(\mathbb C)$ equipped with metric
\begin{equation} \label{met_5}
\rho(x,y)=\underset{t\in [a,b]}{\max} |x(t)-y(t)|.
\end{equation}
We will only verify the
triangle inequality. For each $t\in [a,b]$, using the properties of the absolute value (modulus), we have
\begin{multline}
\notag
|x(t)-y(t)| \leqslant |x(t)-z(t)|+|z(t)-y(t)| \leqslant \\ \leqslant \underset{t\in[a,b]}{\max} |x(t)-z(t)|+\underset{t\in[a,b]}{\max} |z(t)-y(t)|.
\end{multline}
Taking the maximum in $t \in [a,b]$ in the left-hand side of this inequality ends the proof.

$7^\circ.$ The set ${\bf{M}}[a,b]$ of bounded functions $x:[a,b] \to \mathbb R$ $(\mathbb C)$ with metric
\begin{equation} \label{met_6}
\rho(x,y)=\underset{t\in [a,b]}{\sup} |x(t)-y(t)|
\end{equation}
The triangle inequality is proved in exactly the same way as in the previous example.

\poonkt{Basic concepts in metric space}\phantom{01234567890123456789} \

$1^\circ.$ The {\it open ball} of radius $r$ centered at point $a\in \mathfrak M$
is defined to be the set $B(a,r)=
\{x \in \mathfrak M: \rho(x,a)<r\}$. A point $x$ is called an {\it interior} point of the set $A\subset \mathfrak M$
if there exists $\varepsilon >0$ such that $B(x,\varepsilon) \subset A$, that is, if $x$ belongs to the set $A$ together with some ball
centered at the point $x$; the set ${\rm int\,}A$ of interior points is called the {\it interior} of the set $A$; the set $A$ is called
{\it open} if all its points are interior, that is, if $A={\rm int\,}A$. Any open set containing point $x\in\mathfrak M$
is called a {\it neighborhood} of this point.

An open ball in space $\mathbb R$ (example {\bf{1}}) is an open interval: $$B(a,r)=
(a-r,a+r).$$ In space $\mathbb R^2_2$ 
(example { \bf{2}}) an open ball is a disk, i.e.\,$$B(a,r)=
\{(x,y):x^2+y^2<r^2\}.$$

{\it Open ball is an open set.}

Proof.
Put $B\doteq B(x_0,r)\;(x_0\in\mathfrak M,r>0)$. Let $x\in B$ be an arbitrary point. The latter means that $\rho(x,x_0)<r.$
Take $\varepsilon <r-\rho (x,x_0).$ Let $y\in B(x,\varepsilon);$ then $\rho (x_0,y)\leqslant \rho (x_0,x)+\rho (x,y)<\rho (x_0,x)+\varepsilon<r,$
i.e. $y\in B.$ This proves that $B(x,\varepsilon)\subset B.$
\hfill $\square$
 
$2^\circ.$ A {\it closed ball} of radius $r$ centered at a point $a\in \mathfrak M$ is defined to be the set
$B[a,r]=\{x\in \mathfrak M: \rho(x,a)\leqslant r\}$

A closed ball in $\mathbb R$ is a closed interval $B[a,r]=[a-r,a+r];$ a closed ball in $\mathbb R^2_2$ is a (closed) disk 
$B[a,r]=\{(x,y):x^2+y^2\leqslant r^2\}.$

$3^\circ.$ We say that a sequence $\{x_n\}_{n=1}^{\infty}$ of points in a metric space $\mathfrak M$ {\it converges} to the point $x\in \mathfrak M$
(we write: $x_n \to x$ or $x=\lim x_n$) if the sequence of numbers $\{\rho(x_n,x)\}_{n=1}^{\infty}$ converges to zero, $\rho(x_n,x)\to 0.$
The point $x$ is called the {\it limit} of the sequence ${x_n}$.

{\it A sequence $\{x_n\}_{n=1}^{\infty}$ in a metric space $\mathfrak M$ can have only one limit.}

Proof.
Let $x_n\to x$ and let $x_n\to y$. By the triangle inequality, $\rho(x,y)
\leqslant\rho(x,x_n)+\rho(x_n,y)$. Since on the right side of this inequality both sequences
tend to zero, we have $\rho(x,y)=0$,
which yields $x=y.$
\hfill $\square$

$4^\circ.$ A point $x' \in \mathfrak M$ is called a {\it limit point} of the set $A$ if any neighborhood of $x'$ contains
another point of $A$ (different from $x'$). Let $r_n\to 0.$ According to the definition of limit point, the ball $B(x',r_n)$ has a point 
$x_n \in A, x_n\ne x'.$ Since $\rho(x_n,x')<r_n$, we have $x_n \to x'$. Thus, {\it point $x'$ is a limit point
of the set \;$A$ \; if and only if
there exists a sequence of elements of this set that converges to $x'.$}

The set of limit points of the set $A$ is denoted by $A'.$ The set $A$ is called {\it closed} if $A' \subset A.$

The points of the set $A \setminus A'$ are called {\it isolated}.

The set ${\rm cl\,}A=A\bigcup A'$ is called the {\it closure} of the set $A.$

A set $A$ is called {\it perfect} if $A' =A.$ A perfect set, therefore, does not contain isolated points

{\it Closed ball is a closed set.}

Proof.
Set $B\doteq B[x_0,r]\;(x_0\in\mathfrak M, r>0)$ and take an arbitrary point $x\in B$. This means that $\rho(x,x_0)\leqslant r.$
Let $x^*$ be a limit point of $B.$ There is a sequence $\{x_n\}_{n=1}^{\infty}\subset B$ such that $x_n\to x^*. $
The last inclusion means that $\rho (x_n,x)\leqslant r, \; n=1,2,\ldots$ Letting $n\to\infty$, 
we obtain $\rho (x^*,x)\leqslant r,$ i.e.\,$x^*\in B.$ Hence $B$ is a closed set.
\hfill $\square$

{\it The closure of a set is a closed set. The closure of a set is the smallest
closed set containing this set.}

Proof.
By the definition of the closure, $cl\,A=A\cup A'\;(A\subset\mathfrak M)$. Any neighborhood of each limit point of the set $A'$
contains points of this set, and hence points of set $A.$ Therefore, the limit points of set
$A'$ are also limit points of set $A,$, i.e.\,they belong to $A'.$ Thus, $(cl\,A)'\subset A'\subset cl\,A$.

Let $F$ be a closed set that contains $A$. The set $F$ contains all its limit points, and hence the limit points of the set $A.$
This means that $F\supset cl\,A.$
\hfill $\square$

$5^\circ.$ A subset $M$ of a metric space $(\mathfrak M,\rho)$ is called {\it dense} in a set $A\subset \mathfrak M$ if
for any $\varepsilon >0$ and any $x\in A$ there exists $y\in M,$ such that $\rho(x,y)< \varepsilon.$ If $A=\mathfrak M,$ then we
say that $M$ {\it is everywhere dense}.

{\it A set $M$ is dense in $A$ if and only if for every $x\in A$ there exists a sequence
$\{y_n\}_{n=1}^{\infty}\subset M$ such that $y_n\to x$; $M$ is dense in $A$ if and only if $A\subset {\rm cl\,}M$}.

Proof.
Let $M$ be dense in $A$ and let $x\in A$ be an arbitrary point. By the definition of a dense subset, for every $n\in\mathbb N$ there is a
point $y_n\in M$ such that $\rho (x,y_n)<\frac{1}{n}\to 0,$ i.e. $y_n\to x$ for $n\to\infty.$

Assume now that for any $x\in A$ there exists a sequence $\{y_n\}_{n=1}^{\infty}\subset M$ such that $y_n\to x$. Let $\varepsilon >0$ 
be arbitrary. There is $n\in\mathbb N$ such that $\rho (x,y_n)<\varepsilon.$ Clearly, $y_n$ is the required point.

Let $M$ be dense in $A.$ It follows from what was just proved that every point of $A$ is a limit point
of $M,$  i.e.\,it belongs to ${\rm cl\,}M.$ Hence, $A\subset {\rm cl}\,M.$ Conversely, let $A\subset {\rm cl}\,M, \varepsilon >0$ and let $x\in A.$
Then $x\in A\subset {\rm cl}\,M=M\cup M'.$ If $x\in M,$ then we put $y=x$. If $x\in M'$ (that is, $x $ is a limit point $M$), then there is
point $y\in M$ such that $\rho (x,y)<\varepsilon.$ Hence $M$ is dense in $A.$
\hfill $\square$

A subset $A$ of a metric space $(\mathfrak M,\rho)$ is called {\it nowhere dense} if any open ball $B\subset \mathfrak M$
contains an open ball $\widetilde B\subset B$ such that $\widetilde B\bigcap A=\varnothing.$

$6^\circ.$ A set $M$ of a metric space $(\mathfrak M,\rho)$ is called {\it bounded} if there exists a ball containing
this set, i.e. if there are $x\in \mathfrak M$ and $r>0$ such that $M\subset B[x,r].$

{\it If $M\!\subset\!\mathfrak M$ can be covered the union of finitely many balls, i.e.\,
$M\subset \bigcup\limits_{k=1}^n B[x_k,r_k]$
$(n\in \mathbb N),$ then $M$ is a bounded set.}

Proof.
Let
$$
M\subset \bigcup\limits_{k=1}^n B[x_k,r_k] \;(n\in \mathbb N),\quad m=\underset{2\leqslant k\leqslant n}{\max} \rho (x_1,x_k),\quad
r=\underset{1\leqslant k\leqslant n}{\max}r_k
$$
and let $x\in M\subset \bigcup\limits_{k=1}^n B[x_k,r_k].$ Then $x\in B[x_k,r_k] $ for some $k.$ Hence
$$
\rho (x_1,x)\leqslant \rho (x_1,x_k) +\rho (x_k,x)\leqslant m+r\Longrightarrow x\in B[x_1,m+r].
$$
Thus, we have proved that $M\subset B[x_1,m+r],$ as needed.
\hfill $\square$

{\it A convergent sequence is bounded.}

Proof.
Let $x_n\to x.$ For $\varepsilon =1$, there is a positive integer $N_1$ such that for all $n>N_1\;\rho (x_n,x)<1.$
Denote $r=\max\{1,\underset{1\leqslant n\leqslant N_1}{\max}\rho (x,x_n)\}.$ Then we have, for all $n\in\mathbb N$,
$\rho (x,x_n)<r,$ i.e. $\{x_n\}_{n=1}^{\infty}\subset B[x,r].$
\hfill $\square$

$7^\circ.$ Let $(\mathfrak M,\rho)$ be a metric space, $\mathfrak N\subset \mathfrak M$, $\rho'$ is the restriction of $\rho$ on $\mathfrak N$.
Then ($\mathfrak N,\rho')$ is also a metric space. It is called a {\it subspace} of the space $(\mathfrak M,\rho).$
For example, the line segment $[a,b]$ (closed interval) is a subspace of the metric space $\mathbb R.$ Since
${\bf{c}}_o\subset {\bf{c}}\subset {\bf{l}}_{\infty}$, and the metric
in these spaces is defined by the same formula (\ref{met_3}), we obtain that ${\bf{c}}_0$ and ${\bf{c}}$ are subspaces of
${\bf{l}}_{\infty},$ and ${\bf{c}}_0$ is a subspace of ${\bf{c}}.$
Similarly, ${\bf{C}}[a,b]\subset {\bf{M}}[a,b],$ and the metric (\ref{met_6}), after we restrict oursleves to the space of continuous
functions, becomes metric (\ref{met_5}).
Thus, ${\bf{C}}[a,b]$ is a subspace of ${\bf{M}}[a,b].$ At the same time, despite the inclusion
${\bf{l}}_p \subset {\bf{c}}_0,$ ${\bf{l}}_p$ is not a subspace of ${\bf{c}}_o,$ since
metric (\ref{met_3}) does not coincide with the metric in ${\bf{l}}_p$ upon restricting to ${\bf{l}}_p.$

$8^\circ.$ Let $(\mathfrak M,\rho)$ and $(\mathfrak N,d)$ be metric spaces. We equip direct product $\mathfrak M \times \mathfrak N$ with the metric
$$h\bigl((x,y),(x',y')\bigr)=\rho(x,x')+d(y,y').$$ \quad The resulting metric space is called the {\it direct product} of metric spaces $(\mathfrak M,\rho)$ and $(\mathfrak N,d)$.

{\it $h$ satisfies the axioms of a metric. Convergence of $(x_n ,y_n)\to (x,y)$ in $\mathfrak M\times \mathfrak N$
means that $x_n \to x$ \, in \,$\mathfrak M$ \, and \, $y_n \to y$ \, in \, $\mathfrak N,$ and convergence
$x_n \to x$  to  $\mathfrak M$ and $y_n \to y$ \, to \, $\mathfrak N$ imply convergence
$(x_n ,y_n)\to (x,y)$ in $\mathfrak M\times \mathfrak N$}.

Proof.
Let $h((x,y),(x',y'))=0.$ Then $\rho(x,x')=0,\;d(y,y')=0.$ Hence, $x=x',y=y',$ i.e. $(x,y)=(x',y').$ The axiom of symmetry is obvious. Let us verify the triangle inequality:
\begin{multline*}
h((x,y),(x',y'))=\rho (x,x')+d(y,y')\leqslant \rho (x,x'')+\rho (x'',x')+d(y,y'')+\\
+d(y'',y')=h((x,y),\,(x'',y''))+h((x'',y''),\,(x', y')).
\end{multline*}

Let $(x_n,y_n)\to (x,y)$ be in $\mathfrak M \times \mathfrak N.$ This means that
\begin{multline*}
h((x_n,y_n),(x,y))\to 0 \Longrightarrow \rho(x_n,x)+d(y_n,y)\to 0 \Longleftrightarrow \\ \Longleftrightarrow \rho(x_n,x)\to 0, d(y_n,y)\to 0.
\end{multline*}
\hfill $\square$

$9^\circ.$ Let $(\mathfrak M,\rho)$ and $(\mathfrak N,d)$ be metric spaces. A mapping $f: \mathfrak M\to \mathfrak N$ is said to be {\it continuous} at a point $x$ if, for any sequence $x_n\to x$ in $\mathfrak M$,
the sequence $f(x_n)\to f(x)$ in $\mathfrak N.$ The mapping is said to be {\it continuous on set} $M \subset \mathfrak M,$
if it is continuous at every point of $M.$

{\it Metric $\rho$ is continuous} ({\it as a mapping from $\mathfrak M \times \mathfrak M$ to $\mathbb R$}).
{\it If $x^*\in \mathfrak M,$ then the function $\rho (\cdot ,x^*): \mathfrak M\to \mathbb R$ is continuous.}

Proof.
Let $x_n\to x,\;y_n\to y.$ We obtain from \ref{met_8} that
$$
|\rho (x_n,y_n)-\rho (x,y)|\leqslant \rho(x_n,x)+\rho (y_n,y)\to 0 \Longrightarrow \rho(x_n,y_n)\to \rho (x,y),
$$
hence function $\rho (\cdot,\cdot)$ is continuous.

Putting  $y_n=y=x^*$ in \ref{met_8}, we arrive at 
$$\bigl(x_n\to x\bigr)\Longrightarrow \bigl(\rho (x_n,x^*)\to \rho (x,x^*)\bigr)\;(n\to \infty),
$$
which gives the required.
\hfill $\square$

$10^\circ.$ Let $(\mathfrak M,\rho)$ and $(\mathfrak N,d)$ be two metric spaces. A mapping 
$f: \mathfrak M\to \mathfrak N$ is called an
{\it isometry} if $d\bigl(f(x),f(y)\bigr)=\rho(x,y)$ for all $x,y \in \mathfrak M$. In other words, $f$ is an isometry if it preserves distance between points.
Metric spaces $(\mathfrak M,\rho)$ and $(\mathfrak N,d)$ are said to be {\it isometric} if there exists a bijective
(that is, one-to-one)
isometry $f: \mathfrak M\to \mathfrak N.$ If $\mathfrak M$ is isometric to $\mathfrak N,$ then we write $\mathfrak M \cong \mathfrak N.$
Two isometric spaces are identified since they cannot be distinguished within the framework of the theory of metric spaces.

\bigskip

\poonkt
{Analysis of convergence in concrete spaces} \

$1^\circ.$
 Spaces $\mathbb R_p^n \;(\mathbb C_p^n).$ In the course of Analysis one learns that the convergence in space
 $\mathbb R_2^n$ is the coordinate-wise convergence.
In the same way, one learns that the convergence in  space $\mathbb R_p^n \; (\mathbb C_p^n)$, for every $p\geqslant 1$, is the coordinate-wise convergence.

$2^\circ.$ Spaces $C[a,b]$ and $M[a,b].$ Let $x_n \to x$ in $C[a,b].$ By definition, the convergence in space $C[ a,b]$ means that
for any $\varepsilon >0$ there is a number $N$ such that for all $n>N$ one has inequality
$$
\underset{t\in [a,b]}{\max} |x_n(t)-x(t)| <\varepsilon.
$$
The latter means that $|x_n(t)-x(t)|<\varepsilon$ for all  $t\in [a,b]$ at the same time, that is, the sequence
$\{x_n\}_{n=1}^{\infty}$ converges uniformly, $x_n(t)\rightrightarrows x(t).$ Arguing in exactly the same way, one obtains that the convergence
in the space of bounded functions $M[a,b]$ is the uniform convergence.

This means that for every
$\varepsilon>0$ there exists a number $N,$ such that for all $n>N$ one has
$$
\left(\sum\limits_{k=1}^{\infty}|x_k^{(n)}-x_k|^p\right)^{\frac{1}{p}}< \varepsilon.
$$
Therefore, for all $k\in \mathbb N$ and for
all $n>N$ one has $|x_k^{(n)}-x_k|<~\varepsilon,$ which in turn yields the uniform coordinate-wise
convergence. Thus, convergence in space ${\bf l}_p$ implies the uniform coordinate-wise convergence. However, the converse is not true:
the coordinate-wise convergence (even if uniform) does not yet imply the convergence in ${\bf l}_p.$
Indeed, let $p\in\mathbb N,$, 
$$
x^{(n)}=\left(\underbrace{\frac{1}{n},\ldots,\frac{1}{n}}_{ n^p},0,\ldots,\right).$$ Here $x_k^{(n)}=\frac{1}{n}$
for $k=1,2,\ldots,n^p,\;x_k^{(n)}=0$ for $k>n^p;\;x_k^{(n)}\leqslant \frac{1 }{n}$ for all $k,$ hence $x_k^{(n)}\to 0$
uniformly with respect to $k\in\mathbb N$. However, $\rho (x^{(n)},0)=(n^p\cdot\frac{1}{n^p})^{\frac{ 1}{p}}=1,$
i.e. $x^{(n)}\nrightarrow 0.$

$4^\circ.$ Spaces ${\bf l}_{\infty},\; {\bf c}, \;{\bf c}_0.$ The metric in these spaces is defined by the formula (\ref{met_3}).
So, we see immediately that $x^{(n)}\to x$ if and only if $x_k^{(n)}\to x_k\;(n\to \infty)$ uniformly with respect to $ k.$
Thus, the convergence in these spaces is the uniform coordinate-wise convergence.

$5^\circ.$ Space ${\bf s}.$ Let $x_k^{(n)}\to x$ for $n \to \infty$ in $s.$ This means that for any
$\varepsilon>0$ there exists $N$ such that for all $n>N$ one has
\begin{equation}\label{met_9}
\rho(x^{(n)},x)=\sum\limits_{k=1}^{\infty} \frac{1}{2^k}\cdot\frac{|x_k^{(n) }-x_k|}{1+|x_k^{(n)}-x_k|}<\varepsilon.
\end{equation}
This implies that for all $n>N$
$$
\frac{1}{2^k}\cdot \frac{|x_k^{(n)}-x_k|}{1+|x_k^{(n)}-x_k|}<\varepsilon,
$$
that is, for every $k\in \mathbb N \;x_k^{(n)}\to x_k$ for $n\to \infty,$ and so convergence in space ${\bf s}$ implies the
coordinate-wise convergence. Conversely, suppose that we have the coordinate-wise convergence: $x_k^{(n)}\to x_k$ as
$n\to \infty$, for every $k\in \mathbb N.$ Since the series in (\ref{met_9}) converges uniformly in $n$
(indeed, it is majorized by the convergent series $\sum\limits_{k=1}^{\infty}\frac{1}{2^k}$), 
we can pass to the limit $n\to \infty$ under the sum. Thus, we obtain that $x^{(n)}\to x,$ and so convergence in space ${\bf s}$ 
is the coordinate-wise convergence.

\newpage


\bigskip

\poonkt
{Complete metric spaces}\phantom{01234567890123456789} \

Let ($\mathfrak M,\rho)$~ be a metric space. 
A sequence $\{x_n\}$ of  elements of $M$ is called {\it fundamental $($or converging in itself, or a 
Cauchy sequence$)$} if for any $\varepsilon >0$ there exists  $N$ such that for all $m,n>N$ one has
$\rho(x_n,x_m) <\varepsilon$.

{\it Every convergent sequence is fundamental.}

Proof.
Let $x_n\to x.$ This means that for every $\varepsilon >0$ there is a positive integer $N$ such that for all
$n>N$ one has $\rho (x_n,x)<\frac{\varepsilon}{2}.$ Let $m>N.$ Then $\rho (x_n,x_m)\leqslant
\rho (x_n,x)+\rho (x,x_m)<\frac{\varepsilon}{2}+\frac{\varepsilon}{2}=\varepsilon,$ so $\{x_n\}_{n=1}^{\infty}$ is fundamental.
\hfill $\square$

{\it Every fundamental sequence is bounded.}

Proof.
Let $\{x_n\}_{n=1}^{\infty}$ be a fundamental sequence. There is an $N$ such that $\rho (x_m,x_N)<1$
for $m>N.$ Denote $M
\doteq\underset{1\leqslant n\leqslant N}{\max}\rho (x_n,x_N).$ Then $\{x_n\}_{n=1}^ {\infty}\subset B[x_N,\,M+1],$
which, by definition, means that $\{x_n\}_{n=1}^{\infty}$ is bounded.
\hfill $\square$

{\it If a subsequence of a fundamental sequence converges, then the entire sequence converges}
({\it to the same limit}).

Proof.
Let $\{x_n\}_{n=1}^{\infty}$ be a fundamental sequence such that its subsequence 
$\{x_{n_k}\}_{k=1}^{\infty}\subset \{x_n\}_ {n=1}^{\infty}$ converges:
$x_{n_k}\to x$ as $k\to\infty.$
We have
$\rho (x_n,x)\leqslant \rho (x_n,x_{n_k})+\rho (x_{n_k},x)\to 0.$ 
The first term in the right-hand side tends to zero as $n\to\infty$ and $k\to \infty$ since
 $\{x_n\}_{n=1}^{\infty}$ is fundamental. The second term tends to zero since $x_{n_k}\to x\;(k\to\infty).$ ~
\hfill $\square$

A metric space is called {\it complete} if any of its fundamental sequences has a limit in this space, i.e. for
any sequence
$\{x_n\} \subset \mathfrak M$ such that
\begin{equation}\label{met_10}
(\forall \varepsilon>0)\;(\exists N\in \mathbb N)\;(\forall n,m>N)\;\bigl(\rho(x_n,x_m)<\varepsilon\bigr)
\end{equation}
one has
$
(\exists x\in \mathfrak M,\;x_n\to x).
$

In the course of Analysis one proves completeness 
of the space of real numbers $\mathbb R$, as well as
of the spaces $\mathbb R^n_p,\,\mathbb C^n_p\;(1\leqslant p<+\infty,
n\geqslant 1)$ (see the Cauchy convergence principle \cite[c. 83]{ficht1}).


{\it Space ${\bf{M}}[a,b]$ of bounded functions is complete.}

Proof.
Let $\{x_n\}_{n=1}^{\infty}$  be a fundamental sequence in ${\bf{M}}[a,b],$
that is, according to (\ref{met_10}) and (\ref{met_6})
$$
(\forall\, \varepsilon>0)\;\;\;(\exists\, N \in \mathbb N)\;\;\;(\forall\, n,m>N) \quad
\left(\underset{t\in [a,b]}{\sup}\; |x_n(t)-x_m(t)| <\frac{\varepsilon}{2}\right).
$$
Hence, for each $t\in [a,b]$, for all $n,m>N$, we have
\begin{equation}\label{met_11}
\bigl|x_n(t)-x_m(t) \bigr|< \frac {\varepsilon}{2}.
\end{equation}
This means that for each $t\in [a,b]$ a sequence of real numbers
$\{x _n(t)\}_{n=1}^{\infty}$ is fundamental in $\mathbb R,$ i.e. for every $t\in [a,b]$ there is the limit
$
x(t)\doteq\underset{n\to \infty}{\lim} x_n(t).
$
Let us show that $x(\cdot)$ is a bounded function. Let us pass to the limit in (\ref{met_11}) as $m \to \infty:$
$|x_n(t)-x(t)|\leqslant \frac{\varepsilon}{2}$ for all $t\in [a,b]$ and $n>N$. Then
\begin{equation}\label{met_12}
\underset{t\in [a,b]}{\sup} |x_n(t)-x(t)| \leqslant \frac {\varepsilon}{2}<\varepsilon
\end{equation}
Fixing $n>N,$ we find $K$ such that $|x_n(t)| \leqslant K$ for all $t\in [a,b].$ 
Taking $\varepsilon<1$, we have
$
|x(t)| \leqslant |x(t)-x_n(t)|+|x_n(t)| <1+K,
$
that is, $x(\cdot)\in {\bf{M}}[a,b].$ It follows from (\ref{met_12})  that $x_n\to x$ in ${\bf{M} }[a,b].$ 
This proves completeness of space ${\bf{M}}[a,b]$.
\hfill $\square$

The following assertion is quite useful in the proofs of completeness of metric spaces.

\begin{teo}\label{mprth1}
{\it Let $\mathfrak N$ be a subspace of a metric space $\mathfrak M.$ The following is true:}

1)\; {\it if $\mathfrak N$ is a complete space, then $\mathfrak N$ is closed in $\mathfrak M;$}

2)\; {\it if $\mathfrak M$ is a complete space and $\mathfrak N$ is closed in $\mathfrak M,$ then $\mathfrak N$~ is a complete space}.
\end{teo}

Proof.
1. Let $\mathfrak N$ be a complete subspace of $\mathfrak M.$ Let
$x^*\in \mathfrak M$ be its limit point.
The latter means that there is a sequence $\{x_n\}_{n=1}^{\infty}\subset \mathfrak N, \;x_n \to x^*$ (in $\mathfrak M$). Since this sequence
is fundamental and $\mathfrak N$ is complete, we have $x^*\in \mathfrak N.$ Therefore, $\mathfrak N$ is closed.

2. Let $\mathfrak M$ be complete, let $\mathfrak N$ be closed, and let $\{x_n\}_{n=1}^{\infty}\subset \mathfrak N$~ be an arbitrary fundamental
sequence. Since $\mathfrak M$ is complete, there exists $x^*\in \mathfrak M$,
$x^*=\lim x_n.$ Therefore, $x^*$ is a  limit point of $\mathfrak N.$
Since $\mathfrak N$ is closed, $x^*\in \mathfrak N.$ Thus, $\mathfrak N$ contains the limits of its fundamental sequences,
that is, $\mathfrak N$ is complete.
\hfill $\square$

Using this theorem, we now prove completeness of the space of 
continuous functions ${\bf{C}}[a,b].$

Since ${\bf{C}}[a,b]$~ is a subspace of ${\bf{M}}[a,b],$ 
it suffices to prove that ${\bf{C}}[a,b]$ is closed in ${\bf{M}}[a,b].$ Let $x^*$~ be a limit point of ${\bf{ C}}[a,b].$ There is a sequence
$\{x_n\}_{n=1}^{\infty}\subset {\bf{C}}[a,b]$ converging to $x^*.$ This means that $x_n(t) \to x^*(t)$ uniformly on $[a,b].$
By one of the basic results in mathematical analysis, the limit function $x^*(\cdot)$ is continuous, i.e.\,it belongs to ${\bf{C}}[a,b].$ Therefore, ${\bf{C}}[a,b]$ is closed.
By Theorem \ref{mprth1}, ${\bf{C}}[a,b]$~ is a complete space.

Proofs of the completeness of spaces ${\bf{s}}, \;{\bf{l}}_p \; (1\leqslant p< \infty), {\bf{l}}_{\infty},\;{\bf{c}}_0, \;{\bf{c}}$
can be found, for example, in \cite[exerc. 1.24, 1.25]{derr13}.

\newpage


\poonkt
{Completion theorem}\;\phantom{01234567890123456789} \

Let $(\mathfrak M,\rho)$ be a metric space. The smallest complete metric space $(\widehat{\mathfrak M}, \hat {\rho})$ that contains
$(\mathfrak M,\rho)$ as a subspace is called a {\it completion} of $(\mathfrak M,\rho)$. 
For example, $\mathbb R$~ is a completion of $\mathbb Q.$

The following Completion theorem is non-trivial only when applied to incomplete spaces.

\begin{teo}\label{mprth2}
{\it Every metric space has a completion which is unique up to isometry.
Moreover, the original metric space is dense in its completion:}
$cl\, \mathfrak M=\widehat{\mathfrak M}.$
\end{teo}
\doc\,
 To shorten the notation, we will write $\{x_n\}$ instead of $\{x_n\}_{n=1}^{\infty}$. Two
fundamental sequences $\{x_n\}$ and $\{y_n\}$ from $\mathfrak M$ are said to be equivalent (written as $\{x_n\} \sim \{y_n\}$),
if 
$\rho(x_n,y_n)\to 0\; (n\to \infty).$ The relation $\sim$ is reflexive, symmetric and transitive, that is, it is an equivalence relation
on the set of all fundamental sequences in $\mathfrak M.$ This relation allows us to represent $\mathfrak M$ as the union of disjoint
classes of equivalent fundamental sequences. Let us denote the set of these equivalence classes by $\widehat{\mathfrak M}.$

Let $\xi$ and $\eta$ be two classes from $\widehat{\mathfrak M}.$ Select arbitrary fundamental sequences $\{x_n\}\in \xi$ and $\{y_n\}\in \eta$.
Due to the quadrilateral inequality (\ref{met_8}),
$$
\bigl|\rho(x_m,y_m)-\rho(x_n,y_n)\bigr| \leqslant \rho(x_n,x_m)+\rho(y_n,y_m) \to 0 \quad (n\to \infty),
$$
the sequence $\bigl\{\rho(x_n,y_n)\bigr\}$ is fundamental in $\mathbb R.$ Since $\mathbb R$ is a complete metric space, this sequence
converges. Put $\hat{\rho}(\xi,\eta)=\linebreak=
\underset{n\to \infty}{\lim}\rho(x_n,y_n).$
If $\{x'_n\}$ and $\{y'_n\}$ are some other representatives of the classes $\xi$ and $\eta$, respectively, then applying again the 
inequality (\ref{met_8}), we obtain
$$
\bigl|\rho(x'_n,y'_n)-\rho(x_n,y_n)\bigr| \leqslant \rho(x'_n,x_n)+\rho(y'_n,y_n) \to 0\quad (n\to \infty),
$$
that is, $\lim\rho(x'_n,y'_n)=\lim\rho(x_n,y_n).$ This means that the function $\hat{\rho}$ does not depend on the choice 
of representatives $\xi$ and $\eta$.

Let us show that $\hat{\rho}$ is a metric. If $\xi=\eta$ then, obviously, $\hat{\rho}(\xi,\eta)=0.$ Let $\hat{\rho}(\xi,\eta)=0.$
Then $\lim \rho(x_n,y_n)=0$, i.e.\,$\{x_n\}\sim \{y_n\}$, which means $\xi=\eta$. 
The  symmetry of $\hat{\rho}$ follows immediately  from the symmetry of $\rho$.
Let $\xi,\eta,\zeta \in\widehat{\mathfrak M},\; \{x_n\},\{y_n\},\{z_n\}$ be representatives of these classes. 
The triangle inequality for
$\mathfrak M$ is $\rho(x_n,y_n)\leqslant \rho(x_n,z_n)+\linebreak+
\rho(z_n,y_n).$ Passing to the limit as $n\to\infty$, we obtain
$$
\hat{\rho}(\xi,\eta)\leqslant \hat{\rho}(\xi,\zeta)+\hat{\rho}(\zeta,\eta),
$$
that is, the triangle inequality is valid for $\hat{\rho}$.

Thus, $(\widehat{\mathfrak M},\hat{\rho})$ is a metric space.

Let us show that $\mathfrak M$ is isometrically embedded in $\widehat{\mathfrak M}.$ For a given $x\in \mathfrak M$, let us denote by
$\xi_x\in \widehat{\mathfrak M}$ the class containing the sequence 
$\hat{x}=(x,x,\ldots,x,\ldots).$ In fact, the class $\xi_x$ contains all sequences converging to $x:$
$$
(x_n\to x)\Rightarrow (\rho(x_n,x)\to 0)\Rightarrow (\{x_n\}\sim\hat{x})\Rightarrow \bigl(\{x_n\}\in \xi_x \bigr).
$$
The converse is also true: if $\{x_n\}\in \xi_x,$ then $\{x_n\} \sim \hat{x},$ i.e., $\rho(x_n,x)\to 0.$ Thus, since $\hat{\rho}(\xi_x,\xi_y)=\rho(x,y),$
identifying the class $\xi_x \in \widehat{\mathfrak M}$ with $x\in \mathfrak M$, we embed $\mathfrak M$ isometrically into $\widehat{\mathfrak M}.$
Thus, up to an isometry, $\mathfrak M\subset \widehat{\mathfrak M},$, i.e.\,$\mathfrak M$~ is a subspace of $\widehat{\mathfrak M}.$

Let $\xi \in\widehat{\mathfrak M}$ and $\{x_n\}$~ be an arbitrary representative of $\xi.$ For an arbitrarily chosen $\varepsilon>0$, let us 
fix a positive integer $N$ such that for all $m,n>N$ $ \rho(x_n,x_m)<\varepsilon.$ Consider the corresponding sequence of classes of equivalence
$\{\xi_{x_n}\}$ (where $\xi_{x_n}$ is identified with $x_n$). Let $n>N.$ Since $\hat{\rho}(\xi_{x_n},\xi)=\underset{m\to \infty}{\lim}
\rho(x_n,x_m)\leqslant \varepsilon$, we have $\xi_{x_n}(\equiv x_n)\to \xi$ in $\widehat{\mathfrak M}.$ Thus, $\mathfrak M$ is dense in
$\widehat{\mathfrak M}.$

 Let us show that $\widehat{\mathfrak M}$ is complete. Let $\{\xi^{(k)}\}$~ be a fundamental sequence in $\widehat{\mathfrak M}$.
Since $\mathfrak M$ is dense in $\widehat{\mathfrak M}$, for every $k$ there exists $x^{(k)}\in \mathfrak M,$ such that
$\hat{\rho}(\xi_{x^{(k)}},\xi^{(k)})<\frac{1}{k}.$ We apply the triangle inequality:
\begin{multline*}
\rho(x^{(k)},x^{(m)})=\hat{\rho}(\xi_{x^{(k)}},\xi_{x^{(m)}}) )\leqslant \hat{\rho}(\xi_{x^{(k)}},\xi^{(k)})+\hat{\rho}(\xi^{(k)},\xi ^{(m)})+\\
+\hat{\rho}(\xi^{(m)},\xi_{x^{(m)}})<
\frac{1}{k}+\frac{1}{m}+\hat{\rho}(\xi^{(k)},\xi^{(m)}) \to 0 \quad (m ,k\to \infty).
\end{multline*}
Thus, sequence $\{x_n\}$ is fundamental in $\mathfrak M.$ Hence there is a class $\xi \in \widehat{\mathfrak M}$
containing this sequence. We have
$$
\hat{\rho}(\xi^{(k)},\xi)\le \hat{\rho}(\xi^{(k)},\xi_{x_n})+\hat{\rho} (\xi_{x_n},\xi)\to 0 \quad (k\to \infty),
$$
i.e.\,$\xi^{(k)}\to \xi$ in $\widehat{\mathfrak M}$. The completeness of $\widehat{\mathfrak M}$ is proved.

Let $\widetilde{\mathfrak M}$ be another complete metric space containing $\mathfrak M.$ The space $\widetilde{\mathfrak M}$
must contain the limits of fundamental sequences of elements of $\mathfrak M.$ Identifying each such limit
$x*$ with the class of equivalent fundamental sequences converging to $x*,$ we obtain $\widehat{\mathfrak M}\subset \widetilde{\mathfrak M}$.
Thus, $\widehat{\mathfrak M}$ is the desired completion $\mathfrak M.$
\hfill $\square$

Note that a completion is unique up to an isometry. That is, various completions can differ only
by the descriptions of their elements. It is often possible to describe a completion departing directly from the definition 
(as, for example, in the case of $\mathbb Q$ and $\mathbb R).$ The following assertion is convenient in applications.

{\it Let $\mathfrak N$ be a subspace of a complete metric space $\mathfrak M.$ Then the completion $\widehat{\mathfrak N}$ of the metric space
${\mathfrak N}$ is its closure $cl\,\mathfrak N.$}


\poonkt
{Nested balls principle}\phantom{01234567890123456789} \

A sequence $\{B_n\}_{n=1}^{\infty}$ of balls in $\mathfrak M$ is called {\it nested} if
$B_{n+1}\subset B_n,$ 
$n=1,2,\ldots$

\begin{teo}[Nested balls principle]
\label{mprth3}
{\it The following statements are equivalent:}

1) {\it the space $\mathfrak M$ is complete;}

2) {\it any sequence of closed nested balls whose radii tend to zero has a non-empty intersection.}
\end{teo}
\doc\,
1. Let space $\mathfrak M$ be complete and let $\{B[x_n,r_n]\}_{n=1}^{\infty}$~ be an arbitrary sequence of closed
nested balls. The sequence of centers $\{x_n\}_{n=1}^{\infty}$ is fundamental, since
$\;
\rho(x_n,x_m) \leqslant \max\{r_n,r_m\}\to 0.
$\;
Since $\mathfrak M$ is complete, there exists $x=\lim\,x_n \in\mathfrak M.$ Since each ball contains all centers,
except perhaps for a finite number of them, we obtain that $x$ is a limit point of each ball. Since the balls are closed, we have
$x\in B[x_n,r_n]$ for all $n\in \mathbb N,$ i.e. $\bigcap\limits_{n=1}^{\infty} B[x_n,r_n]\ne \varnothing.$

2. Assume that assertion 2) holds and let $\{x_n\}_{n=1}^{\infty}$~ be an arbitrary fundamental sequence.
Find a number $n_1$ such that $\rho(x_{n_1},x_n)<\frac{1}{2}$ for all $n>n_1.$ Let the numbers
$n_k >n_{k-1}>\ldots >n_2>n_1$ be such that
$\rho(x_{n_{i}},\,x_n)<\frac{1}{2^i}$ for all $n>n_i,\; i=1,\ldots,k.$
We choose $n_{k+1}>n_k$ so that
$\rho(x_{n_{k+1}}, \,x_n)< \frac{1}{2^{k+1}}$ for all $n>n_{k+1}.$ So we get a subsequence
$\{x_{n_{k}}\}_{k=1}^{\infty}\subset \{x_n\}_{n=1}^{\infty}$ having the following property:
\begin{equation}\label{met_13}
\rho(x_{n_{k}},x_n)<\frac{1}{2^k}\quad \text{for all}\; n>n_k.
\end{equation}
Consider the closed balls $B_k \doteq B[x_{n_k},\,\frac{1}{2^{k-1}}],\; k=1,2,\ldots$ (Note that the radii of the balls have doubled compared to the right side of
(\ref{met_13}).)

Let us show that the sequence of balls $\{B_k\}_{k=1}^{\infty}$ is nested. Let $x\in B_{k+1}.$ This means that
$\rho(x_{n_{k+1}},x)\leqslant \frac{1}{2^k},$
so by (\ref{met_13}) and the triangle inequality
$$
 \rho(x_{n_k},x)\leqslant \rho(x_{n_k},x_{n_{k+1}})+\rho(x_{n_{k+1}},x)\leqslant \frac {1}{2^k}+\frac{1}{2^k}=\frac{1}{2^{k-1}},
$$
so $x\in B_k.$ Thus, $B_{k+1} \subset B_k,\; k=1,2,\ldots,$ as claimed.

By our assumption, there exists a point $x \in \mathfrak M$ belonging to all balls $B_k.$ Since 
$\rho(x_{n_{k}},x)\leqslant \frac{1}{2^{k-1 }} \to 0,$
we have $x= \underset{k\to\infty}{\lim}x_{n_k}.$ According to the third property of fundamental sequences, $x=\underset{n\to\infty}{\lim}x_n.$
Thus, an arbitrary fundamental  sequence in $\mathfrak M$ has limit in $\mathfrak M,$, i.e.\,space $\mathfrak M$ is complete.
\hfill $\square$


\poonkt
{Separable metric spaces}\phantom{01234567890123456789} \

A metric space $\mathfrak M$ is said to be {\it separable} if it contains a countable
everywhere dense set $D.$
For example, the space $\mathbb R$ of real numbers with $\rho(x,y)= |x-y|$ is separable, since one can take as a countable everywhere dense subset the set of rational numbers
$D=\mathbb Q$. The space $\mathbb R_p^n$ is separable, $D=\mathbb Q_p^n$ is the set of vectors
with rational coordinates.

The space ${\bf{C}}[a,b]$ of continuous functions is separable, one can take as a countable everywhere dense subset $D$ the set of polynomials with rational
coefficients.

Indeed, according to the Weierstrass Approximation Theorem, for every $\varepsilon>0$ and any continuous function $x(\cdot)$ on $[a,b]$
there exists a polynomial $p(\cdot)$ such that $\rho(x,p)<\frac{\varepsilon}{2}.$ Further, there is a polynomial $q(\cdot)$
with rational coefficients such that $\rho(p,q)<\frac{\varepsilon}{2}.$ Thus, $\rho(x,q)\le\rho(x,p)+\rho(q ,p)<\varepsilon.$

It is easy to verify that the set of {\it finite} sequences of rational numbers (``finite'' means that the sequence contains only a finite number of non-zero terms) constitutes a countable everywhere dense set in the spaces of sequence 
${\bf{l}}_p\;(1\leqslant p<\infty),\;{\bf{c}},\;{\bf{c}}_0,\;{\bf{s}}$, so these spaces are separable.

Let us present a criterion for non-separability of a metric space.
\begin{teo}
\label{mprth6}
{\it If there exists an uncountable set 
$S\subset \mathfrak M$ and a number $\alpha>0,$ such that for any distinct elements $x,y\in S$ one has $\rho(x,y)\geqslant \alpha,$ then the space $\mathfrak M$ is not separable.}
\end{teo}

\doc\,
Let $D$~ be an arbitrary dense set in $\mathfrak M$. Let us construct a mapping $f:S\to D$ as follows.
For every $x\in S$ we assign an element $f(x)\in D$ such that $\rho\bigl(x,f(x)\bigr)<\frac{\alpha}{3}.$ Since $D$ is everywhere dense in $\mathfrak M,$
there is such $f(x)$. Let us show that the mapping $f$ is injective. Let $x\ne y$. Suppose that $f(x)=f(y).$ Then
$$
\alpha\leqslant \rho(x,y)\leqslant \rho\bigl(x,f(x)\bigr)+\rho\bigl(f(x),f(y)\bigr)+\rho\bigl (f(y),y\bigr)<\frac{\alpha}{3}+0+
\frac{\alpha}{3}=\frac{2}{3}\alpha.
$$
This is a contraction, which proves the injectivity of $f.$ Thus, an\linebreak 
uncountable set $S$ is in one-to-one correspondence with a part of the set $D.$ Hence $D$ is uncountable. Consequently,
the space $\mathfrak M$ is not separable.
\hfill $\square$

{\it The space of bounded functions ${\bf{M}}[a,b]$ is not separable.}

\doc\,
Let $S=\bigcup\limits_{c\in[a,b]}\{x_c(\cdot)\},$ where $x_c(t)=0$ for $t\ne c$,
$x_c(c)=1.$ Obviously, $S$ is uncountable, $S\subset {\bf{M}}[a,b],$ and for
any distinct $x$ and $y\quad \rho(x,y)=1.$ By Theorem \ref{mprth6}, ${\bf{M}}[a,b]$ is not separable.
\hfill $\square$

{\it The space of bounded sequences ${\bf{m}}={\bf{l}}_{\infty}$ is not separable.}

\doc\,
According to Theorem \ref{mprth6}, one needs to specify an uncountable set $S\subset {\bf{l}}_{\infty}$ such that the distance between any two points in $S$ is greater or equal to some positive number. One can take as such set the set of sequences
whose components are zeros and ones. One can associate to every sequence $x\in S$ a
set $X\subset \mathbb N$ whose elements are the indices of the entries of $x$ that are equal to 1. And vice versa, to each subset
$X\subset \mathbb N$ one can associate the sequence $x\in S$ such that $x_i=1$ for $i\in X,\;x_i=0$ for
$i\notin X$. This establishes a one-to-one correspondence between the set $S$ and the set $\mathcal P(\mathbb N)$ of all subsets of the set
natural numbers. It is known that $\mathcal P(\mathbb N)$ is uncountable. Clearly, given any two elements $x,y\in M,\;x\ne y$ one has $\rho (x,y)=1.$
\hfill $\square$


\poonkt
{Compact sets\;\;}\phantom{01234567890123456789} \

\textbf{1.\:Definition. Hausdorff's Theorem}\phantom{01234567890123456789} \

Let $(\mathfrak M,\rho)$ be a metric space. A set $A\subset \mathfrak M$ is called {\it compact}
if any of its infinite subsets $B\subset A$ has a limit point in $A$. This definition is equivalent to
the following: a set $A$ is compact if any sequence $\{x_n\}\subset A$ contains  a 
subsequence that converges to an element of $A$.

Since a finite set cannot contain an infinite subset, all finite sets
are compact by definition.

A compact set is closed.
Indeed, let $A$ be compact and let $x_0$ be is its limit point. The set $A$ contains a sequence $\{x_n\}$ converging to $x_0.$
By the definition of compactness, $x_0 \in A.$ Hence $A$ is closed.

A set $A\subset \mathfrak M$ is called {\it relatively compact} if its closure $cl\,A$ is compact.

By the Bolzano-Weierstrass theorem, a bounded 
(bounded and closed)
subset of $\mathbb R_p^n$ is relatively compact (compact).

Let $A\subset \mathfrak M.$ We say that $B \subset \mathfrak M$ is an $\varepsilon$-net for $A$ if for every $a \in A$ 
there exists $b\in B$ such that $\rho(a,b)<\varepsilon.$

A set $A\subset \mathfrak M$ is called {\it totally bounded} if $\mathfrak M$ contains, for every $\varepsilon>0$, a
finite $\varepsilon$-net for $A.$

Let us show that a totally bounded set is bounded.

Indeed, let $A\subset \mathfrak M$ be totally bounded and let $B=
\{b_1,b_2,\ldots , b_n\}$ be an 1-net for $A.$ For every $x\in A$
one can find $b_k$ such that $\rho(x,b_k)<1,$ that is, $x\in B[b_k,1].$ Thus $A\subset \bigcup\limits_{k=1}^{n}B[b_k,1].$
Clearly, a finite union of  balls of radius $1$ is contained in a ball of finite radius. Therefore, $A$~ is a bounded set.

The concept of compactness can also be applied to $\mathfrak M$ itself. A compact metric space is called a {\it compact}.

The following Hausdorff theorem is a criterion for  relative 
compactness of a subset of a metric space.

\begin{teo}
\label{hausd}
{\it If $A\subset \mathfrak M$ is relatively compact, then it is totally bounded. Moreover, if $\mathfrak M$ is complete, 
then a totally bounded set $A \subset \mathfrak M$ is relatively compact.}
\end{teo}
\doc\,
 {\it Necessity.} Let $A$ be relatively compact. Suppose that $A$ is not totally bounded.
This means that for some $\varepsilon_0$ there is no finite $\varepsilon_0$-net for $A$.
Therefore, for every $a_1\in A$, there is an element $a_2\in A$ such that $\rho(a_1,a_2)\!\!\geqslant \!\! \varepsilon_0$ 
(otherwise $\{a_1\}$~ is an
$\varepsilon_0$-net for $A$ consisting of a single element, which contradicts to out hypothesis that $A$ is not totally bounded). 
Suppose that we already selected 
$a_1,a_2,\ldots,a_k$ such that $\rho(a_i,a_j)\geqslant \varepsilon_0$ for $i\ne j$. Then there is an element $a_{k+1}$ 
such that $\rho(a_{i},a_{k+1})\geqslant
 \varepsilon_0$ for $i=1,2,\ldots ,k$ (if not, then $\{a_1,a_2,\ldots ,a_{k}\}$ is a $\varepsilon_0$-net for $ A$). Thus, we obtain
an infinite sequence $\{a_k\}_{k=1}^{\infty}$ which has no limit point since $\rho(a_i,a_k)\geqslant \varepsilon_0$ for $i\ne k.$
This contradicts to the relative compactness of $A.$ Hence $A$ is totally bounded.

{\it Sufficiency.} Let $\mathfrak M$ be complete, let $A\subset \mathfrak M$ be 
totally
bounded, Let $C$~ be an arbitrary
infinite subset of $A$, let $C_1=\{x_1,x_2,\ldots,x_{n_1}\}$ be a finite 1-net for $A.$
Consider balls $B[x_k,1],\; k=1,2,\ldots,n_1$. These balls cover all $A.$
Since there are finitely many of them, at least one of these balls contains an infinite set of points from $C$,
say, the ball $B[y_1,1]$
$ (y_1=x_k\; \mbox {for some}\; k).$ We can find a finite $\frac{1}{2}$-net for $B[y_1,1]$. 
For every element of this $\frac{1}{2}$-net, let us fix a ball of radius $\frac{1}{2}$ containing this element.
One of these balls contains an infinite set of points from $C$. Denote it by $B[y_2,\frac{1}{2}].$ 
Continuing this process, at the $k$-th step we obtain a ball $B_k=B[y_k,\frac{1}{2^{k-1}}]$ 
containing an infinite set of points from $C$ ($k=1,2, \ldots$).
The balls $\hat B_k=
B[y_k,\frac{1}{2^{k-2}}]$ have radii  twice as large. Each of them also contain an infinite set of points from
$C$. Moreover, these new balls are
nested: if $y\in\hat B_{k+1},$ then $\rho(y,y_{k+1})\leqslant
\frac{1}{2^{k-1}},$ so
$$
\rho(y,y_k)\leqslant \rho(y,y_{k+1})+\rho(y_{k+1},y_k)\leqslant \frac{1}{2^{k-1}} +\frac{1}{2^{k-1}}=\frac{1}{2^{k-2}},
$$
i.e.\,$y\in\hat B_k,\;\hat B_{k+1}\subset \hat B_k.$ Due to the Principle of nested balls
 ($\mathfrak M$ is complete!) there 
exists $x$ that belongs to all $\hat B_k,\; k=1,2,\ldots$ 
Since each of the balls contains an infinite set of points from $C,$ $x$ is a limit point of $C.$ 
Thus, an arbitrary infinite subset of $A$ has a limit point. Hence $A$ is relatively compact.
\hfill $\square$

\begin{sle}
\label{mprcolhaus1}
{\it A compact space is separable.}
\end{sle}

\doc\,
Let $\mathfrak M$ be a compact set and let $C_n$ be a $\frac{1}{n}$-net for $\mathfrak M.$ 
Then $D\doteq\bigcup\limits_{n= 1}^{\infty}C_n$ is a
countable everywhere dense set in $\mathfrak M.$
\hfill $\square$

\begin{sle}
\label{mprcolhaus2}
{\it Consider a subset of $\mathbb R_p^n$. It is relatively compact} ({\it compact}), {\it if and only if 
it is  bounded}
({\it bounded and closed}).
\end{sle}

\doc\,
Sufficiency is a consequence of the Bolzano-Weierstrass theorem mentioned above.
Necessity follows from the Hausdorff theorem, since a 
totally
bounded set is bounded.
\hfill $\square$

\begin{sle}
{\it For a set $A\subset \mathfrak M$ to be relatively compact, it is necessary, and if $\mathfrak M$ is complete, 
then it is sufficient, that $\mathfrak M$ contains,
for every $\varepsilon>0$, a relatively compact $\varepsilon$-net for $A.$}
\end{sle}

\doc\,
Let $B$ be a relatively compact $\frac{\varepsilon}{2}$-net for $A$, let $C$ be a finite
$\frac{\varepsilon}{2}$-net for $B.$ Then $C$ is a finite $\varepsilon$-net for $A.$
\hfill $\square$

\medskip

\textbf{2.\;Criterion for compactness in ${\bf{C}}[a,b]$}\phantom{01234567890123456789}

A set $\Phi \subset {\bf{C}}[a,b]$ is called {\it uniformly bounded} if there exists a constant $K$ such that
$|x(t)|\leqslant K$ for all $t\in[a,b]$ and for all $x\in\Phi.$ A uniformly bounded set lies in the ball $B[0,K],$ and, consequently,
is bounded in ${\bf{C}}[a,b].$

A set $\Phi$ is called {\it equicontinuous} if for any $\varepsilon>0$ there exists $\delta>0$
such that for all $t',t''\in[a,b]$, $|t'-t''|<\delta$ and for all $x\in\Phi$ one has $| x(t')-x(t'')|<\varepsilon.$
That is,
\begin{equation*}
(\forall \varepsilon>0)\,(\exists \delta>0)\,(\forall t',t''\in[a,b]\!:\!\,|t'-t'' | <\delta)\,(\forall x\in \Phi)\,
\bigl(|x(t')-x(t'')| <\varepsilon\bigr).
\end{equation*}
The negation of the equicontinuity looks as follows:
\begin{multline}
\label{met_30}
(\exists \varepsilon_0>0)\;(\forall \delta>0)\;(\exists t',t''\in[a,b]\!: \;|t'-t''| < \delta)\;(\exists x_0\in \Phi) \\
\bigl(|x_0(t')-x_0(t'')|\geqslant \varepsilon_0\bigr).
\end{multline}

Let us consider an example. Let $\Phi=\{\sin\alpha t\}_{\alpha \in\mathbb R}\subset C[0,\pi].$
The set $\Phi$ is uniformly bounded ($K=1$), but it is not equicontinuous.
Let $\varepsilon_0=1,\;\delta>0$ be arbitrary, $t'=0,\; t''=\frac{\delta}{2},\; x_0(t)=\sin\alpha_0 t, \alpha_0=\frac{\pi}{\delta}.$
Then $|t'-t''|\!=\!\frac{\delta}{2}<\delta,$
$|x_0(t')-x_0(t'')|\! =\!\left|0-\sin \frac{\pi}{\delta}\cdot \frac{\delta}{2}\right|\! =\!\sin\frac{\pi}{2}\!=\!1.$
According to (\ref{met_30}), $\Phi$ is not equicontinuous.

On the contrary, the set $\Phi=\{\sin\alpha t\}_{\alpha \in [-100,100]}$ is equicontinuous: given an arbitrary $ \varepsilon>0$,
we set $\delta=\frac{\varepsilon}{100}.$ Then, if $|t'-t''| <\delta,$ 
we obtain using the Lagrange finite increments formula 
$$
|\sin\alpha t'-\sin\alpha t''| =|\alpha \cos\xi(t'-t'')| \leqslant 100\cdot \frac{\varepsilon}{100}=\varepsilon.
$$
Hence the set $\Phi$ is equicontinuous.

The following theorem due to Arzela is a criterion for relative compactness of a set in the space $C[a,b].$
\begin{teo}\label{mprth8}
{\it A set $\Phi \subset {\bf{C}}[a,b]$ is relatively compact if and only if it is uniformly
bounded and equicontinuous.}
\end{teo}
\doc\,
{\it Necessity.} Let $\Phi \subset {\bf{C}}[a,b]$ be relatively compact.
Then, according to the Hausdorff theorem, it is 
totally
bounded, and hence uniformly bounded. Let us prove its equicontinuity.
Let $\varepsilon>0$ be arbitrary, and let $A=\{a_1,a_2,\ldots,a_n\}$ be a $\frac{\varepsilon}{3}$-net for $\Phi.$ Each function $a_k(\cdot)$
is uniformly continuous on $[a,b],$ so there are $\delta_k$ such that if $| t'-t''| <\delta_k,$ then
\begin{equation}\label{met_31}
\bigl|a_k(t')-a_k(t'')\bigr| <\frac{\varepsilon}{3}, \quad k=1,2,\ldots,n.
\end{equation}
We set $\delta=\underset{1\leqslant k\leqslant n}{\min}\delta_k.$ Then, if $|t'-t''|<\delta,$ all inequalities
(\ref{met_31}) are satisfied.

Let $x(\cdot)$ be an arbitrary function from $\Phi.$ By the definition of an $\varepsilon$-net, there exists $a_{k_0}(\cdot)\in A$ such that
$\rho\bigl(x,a_{k_0}\bigr)<\frac{\varepsilon}{3}.$ This means that for all $t\in[a,b]$ one has inequality
\begin{equation}
\label{met_32}
\bigl|x(t)-a_{k_0}(t)\bigr|<\frac{\varepsilon}{3}.
\end{equation}
Taking into account (\ref{met_31}), (\ref{met_32}), we obtain under the condition \linebreak
$|t'-t''| <\delta:\;$
$\;
\bigl|x(t')-x(t'')\bigr| \leqslant \bigl|x(t')-a_{k_0}(t')\bigr|+\bigl|a_{k_0}(t')-a_{k_0}(t'')\bigr| +\linebreak+
\bigl|a_{k_0}(t'')-x(t'')\bigr|<\frac{\varepsilon}{3}+\frac{\varepsilon}{3}+\frac{\varepsilon}{3} =\varepsilon.\;$
This proves that $\Phi$ is equicontinuous.

{\it Sufficiency.} Let $\Phi$ be uniformly bounded and equicontinuous, let $\varepsilon>0$ be arbitrary, and let $K$ and $\delta$
be from the definitions of the uniform boundedness and the equicontinuity.

Let us divide the closed interval $[a,b]$ into subintervals $[t_k,t_{k+1}]$ whose length is less than $\delta$.
Next, we divide the closed interval $[-K,K]$ of the $y$-axis into subintervals whose length is less than
$\varepsilon$. Let us draw the corresponding vertical and horizontal lines. Thus, we obtain a grid. Consider the set of piecewise linear curves
with nodes at 
the nodes of the grid.
There is a finite number of such curves.
For an arbitrary function $x(\cdot)\in\Phi$, consider a continuous piecewise linear function $y(\cdot)$
(its graph is one of the indicated piecewise linear curves) whose graph is situated, 
at the nodes of the grid, at the distance no greater than $\varepsilon$  from the graph of  $x(\cdot)$. The inequalities
$$
\bigl|x(t_k)-x(t_{k+1})\bigr|<\varepsilon,\quad \bigl|y (t_k)-x(t_k)\bigr|<\varepsilon, \quad \bigl| y(t_{k+1})-x(t_{k+1})\bigr|<\varepsilon
$$
yield
$$
\bigl|y(t_k)-y(t_{k+1})\bigr|\!\leqslant
\!\bigl|y(t_k)-x(t_k)\bigr|\!+\!\bigl|x (t_k)-x(t_{k+1})\bigr|\!+\!|y(t_{k+1})\!-\!x(t_{k+1})|\!<3\varepsilon.
$$
Since $y(\cdot)$ is linear on the interval $[t_k,t_{k+1}],$ we have  $|y(t )-y(t_k)|<3\varepsilon$ for all $t\in [t_k,t_{k+1}]$.

Let $t\in[a,b]$ and let $t_k$ be the grid node closest to $t$ on the left. Then
$$
\bigl|x(t)-y(t)\bigr|\!\leqslant \!\bigl|x(t)\!-\!x(t_k)\bigr|\! +\!\bigl|x(t_k)\!-\!y(t_k)\bigr|\!+\!\bigl|y(t_k)-y(t)\bigr|\! 
<\! \varepsilon+\varepsilon+3\varepsilon=5\varepsilon.
$$
Thus, this set of piecewise linear functions
forms a finite $5\varepsilon$-net for $\Phi$. According to Hausdorff's theorem,
 $\Phi$ is relatively compact.
\hfill $\square$

If the set $\Phi \subset {\bf{C}}[a,b]$ consists of differentiable functions, then from the Arzela theorem we obtain the 
following sufficient condition
for relative compactness.

\begin{teo}\label{mprth9}
{\it Let $\Phi \subset {\bf{C[a,b]}}$.  Assume that all functions from $\Phi$ are differentiable and there exist constants
$M_0$ and $M_1$ such that $| x(a)| \leqslant M_0,\;| x'(t)|\leqslant M_1$
for all $x\in\Phi$ and $t\in [a,b].$ Then $\Phi$ is relatively compact.}
\end{teo}
\doc\,
Every function from $\Phi$ admits representation
$$
x(t)=x(a)+\int\limits_a^tx'(s)\,ds, \quad t\in[a,b].
$$
It follows from this representation that
$$
\bigl|x(t)\bigr|=\left|x(a)+\int\limits_a^tx'(s)\,ds\right|\leqslant M_0+M_1(b-a), \quad t\in[ a,b].
$$
This proves the uniform boundedness of $\Phi.$ Further,
$$
\bigl|x(t')-x(t'')\bigr|=\left|\int\limits_{t'}^{t''}\bigl|x'(s)\bigr|\,ds \right|\leqslant M_1|t'-t''|, \quad t',t''\in[a,b].
$$
Taking, for an arbitrary $\varepsilon>0$, $\delta:=\frac{\varepsilon}{M_1},$ we obtain: if $|t'-t''|<\delta,$ then $|x(t')- x(t'')| <\varepsilon.$
Hence, set $\Phi$ is equicontinuous. By the Arzela theorem, $\Phi$ is relatively compact.
\hfill $\square$

\textbf{3.\;Equivalent definition of compactness}\phantom{01234567890123456789}

Let $(\mathfrak M,\,\rho)$~ be a metric space. We say that a family of sets
$\{S_\alpha\}\linebreak 
\Bigl(S_\alpha\subset\mathfrak M\Bigr)$ is a
{\it cover} of the set $A\subset\mathfrak M$ if $A\subset \bigcup\limits_\alpha S_\alpha.$ If the sets $S_\alpha$ are open,
then the cover is said to be open. If sets $S_\alpha\subset A$, then $A=\bigcup\limits_\alpha S_\alpha.$
A cover is said to be countable if it is a countable collection of sets: $\{S\}_{n=1}^{\infty}$.

We say that a family of sets $\{S_\alpha\}\;\;\;\Bigl(S_\alpha\subset\mathfrak M\Bigr)$ is centered if any finite intersection
$\bigcap\limits_{k=1}^n S_k$ is non-empty. A centered family can be countable.
\begin{lem}
\label{equivkomp}
{\it The following statements are equivalent:

a$)\;$ Any open $($open countable$)$ cover of set $A\subset\mathfrak M$ contains a finite subcover.

b$\;)$ Every $($countable$)$ centered family of closed subsets of set $A$ has non-empty intersection.}
\end{lem}
\doc\,
Assume that a) is true. Let $\{G_\alpha\}$ be an open cover of set $A$, i.e.\,$A=\bigcup\limits_\alpha G_\alpha.$
Then the sets $F_\alpha=A\setminus G_\alpha$ are closed and, by the duality laws, $\bigcap\limits_\alpha F_\alpha=\varnothing.$
This means that the family $\{F_\alpha\}$ cannot be centered, that is, there is a finite subfamily $\{F_1,\,F_2,\ldots,F_n\}$
such that $\bigcap\limits_{k=1}^n F_k=\varnothing.$ By the duality laws,
$\bigcup\limits_{k=1}^n G_k=A,$ i.e. $G_k=A\setminus F_k,\; k=1,2,\ldots,n$ is
a finite subcover for $A.$ Therefore, property b) holds. We argue in exactly the same way if the original cover is countable.

Now, assume that b) holds, and let $\{F_\alpha\}$ be a centered family of closed subsets of the set $A.$
Assume that $\bigcap\limits_{\alpha} F_\alpha=\varnothing.$ By the duality laws,
$\bigcup\limits_{\alpha} G_\alpha=A,\;\;G_\alpha=A\setminus F_\alpha,$ i.e. $\{G_\alpha\}$ forms an open cover
of set $A.$ By our assumption, there exists a finite subcover $G_1,G_2,\ldots,G_n,\;\;A=\bigcup\limits_{k=1}^n G_k.$
But this means that $\bigcap\limits_{k=1}^n F_k=\varnothing,\;\;F_k=\T\setminus G_k=1,2,\ldots,n,$ which contradicts to the fact that $\{F_\alpha\}$ is centered. Hence $\bigcap\limits_{\alpha} F_\alpha\ne\varnothing.$ 
We argue in exactly the same way if the centered family is countable.
\hfill $\square$
\begin{teo}
\label{equivkompb}
{\it A set $A\subset\mathfrak M$ is compact if and only if condition a$)$ or condition b$)$
in Lemma \ref{equivkomp} is satisfied.}
\end{teo}
\doc\,
{\it Sufficiency.\;}
Let one of the properties a) or b) hold
(then, by what was already proved, both of these properties hold). Let $B\subset A$ be an arbitrary infinite subset with
no limit points. Then the countable set $C\doteq \{c_1,c_2,\ldots\}\subset B$ has no limit points either.
The sets $C_n=\{c_n,c_{n+1},\ldots\},\;n=1,2,\ldots$ are closed (the sets of their limit points are empty), each of their finite intersections
is not empty. Therefore, the family $\{C_n\}$~ is centered. But $\bigcap\limits_{n=1}^\infty C_n=\varnothing.$ This contradicts to property a), which means that set $B$ has at least one limit point, that is, the set $A$ is
compact.

{\it Necessity.\;}
(We will carry out the proof for countable covers only.)\;
Let $A$ be compact and let $\{G_n\}_{n=1}^\infty$~ be an open cover
for $A$, i.e.
$A=\bigcup\limits_{n=1}^\infty G_n.$ 
Suppose that no finite subfamily $\{G_1,G_2,\ldots,G_n\}$ covers set
$A$, that is, $\bigcup\limits_{k=1}^n G_k\ne A.$ Let $F_n\doteq A\setminus G_n.$ Then, by the duality laws,
$\bigcap\limits_{k=1}^n F_k\ne \varnothing,$ i.e. the family of closed sets $\{F_n\}_{n=1}^\infty$ is centered (and countable).
We must prove that the intersection $\bigcap\limits_{n=1}^\infty F_n\ne\varnothing.$

Consider sets $\Phi_n\doteq\bigcap\limits_{k=1}^n F_k.$ Since
$\{F_n\}_{n=1}^\infty$ is centered, we have $\Phi_n\ne\varnothing,\;n=1,2,\ldots$ By the definition,
\begin{equation}
\label{exprt}
\Phi_1\supset\Phi_2\supset\ldots, \qquad  \bigcap\limits_{n=1}^\infty F_n=\bigcap\limits_{n=1}^\infty\Phi_n.
\end{equation}
We have to distinguish between two possibilities:

1)\; Starting from some number $n_0$, $\Phi_n$'s coincide, i.e.\,$\Phi_{n_0}=\Phi_{n_0+1},\ldots$
Then $\bigcap\limits_{n=1}^\infty\Phi_n=\Phi_{n_0}\ne\varnothing$ and (\ref{exprt}) proves the assertion.

2)\; There are infinitely many different $\Phi_n$'s. Without loss of generality, we may assume that all $\Phi_n$'s are distinct.
Let $x_n\in\Phi_n\setminus\Phi_{n+1},\linebreak
n=1,2,\ldots$ Since this sequence is an infinite set of different
points from $A,$ then, since $A$ is compact, this sequence has at least one limit point. Let us denote this point by $x^*.$
Since every
$\Phi_n$ is closed, we have $x^*\in\Phi_n,\;n=1,2,\ldots$ According to (\ref{exprt}), $x^*\in\bigcap\limits_{n=1} ^\infty F_n,$ so this intersection is not empty.
\hfill $\square$

Theorem \ref{equivkompb} means that we could also take either property a) or property b) as the definition of compactness.
In the class of \textit{topological spaces}, which is wider than the class of metric spaces, a set is called
compact if any open cover has a finite subcover. As Theorem \ref{equivkompb} states, for metric
spaces, both of these definitions are equivalent. We can even assume that the original cover is countable.

\medskip

{\bf{4.\;Criteria for compactness in spaces of
sequences } } \ \phantom{01234567890123456789}

$1^{\circ}.\;$ {\textbf{Space}} ${\bf{l}}_p\;(1\leqslant p<\infty)\;$

\begin{utv}
\label{compl_1}
{\it A set $M\subset {\bf{l}}_p$ is relatively compact if and only if
 
1$)$ it is bounded, and

2$)$ for any $\varepsilon >0$ there exists a number $N$ such that the following inequality holds simultaneously for all $x\in M$}:
$\sum\limits_{k=N+1}^{\infty} | x_k| ^p<\varepsilon^p$ 
(the remainder of the series $\sum\limits_{k=1}^{\infty} | x_k| ^p$  tends to zero uniformly with respect to
$x\in M).$
\end{utv}
\doc\,
{\it Sufficiency.} 
Let conditions 1) , 2) be satisfied, let $\varepsilon >0$ be arbitrary and suppose that $N\in\mathbb N$ corresponds to $\frac{\varepsilon}{2}$ in
condition 2). Denote by ${\bf{l}}_p^{(0,N)}$ the set of finite sequences with $x_k=0$ for $k=N+1,N+2,\ldots$
Consider the bijective mapping $F:\;{\bf{l}}_p^{(0,N)}\to\mathbb R_p^N$, which associates to sequence
$x=(x_1,\ldots,x_N,0,\ldots)\in l_p^{(0,N)}$  the vector $(x_1,\ldots,x_N)\in\mathbb R_p^N.$
Denote by $\tilde F$ the extension of
$F$ to ${\bf{l}}_p$ by continuity (this extension is no longer injective). Since set $\hat M\doteq
\tilde F(M)\subset\mathbb R_p^N$
is bounded, it is relatively compact, so it has a finite $\frac{\varepsilon}{2}$-net
$\hat A\doteq \{\hat a^{(i)}\}_{i=1}^{m} \subset \mathbb R_p^N$. 
Let $a^{(i)}=F^{-1}(\hat a^{(i)})\;(i=
1,2,\ldots,m)$. Let us show that
$A\doteq\{a^{(i)}\}_{i=1}^m$~ is an $\varepsilon$-net for $M.$

Let $x\in M$. Then $\hat x=\tilde F(x)\in\mathbb R_p^N.$ There exists $i_0\;(1\leqslant i_0\leqslant m)$ such that
$\rho_{\mathbb R_p^N}(\hat x,\hat a^{(i_0)})<\frac{\varepsilon}{2}.$ Therefore,
$$
\left(\rho (x,a^{(i_0)})\right)^p\!=\!\sum\limits_{k=1}^N|x_k-a_k^{(i)}|^p \!+\!\sum\limits_{k=N+1}^{\infty}|x_k|^p\!<\!\left(\frac{\varepsilon}{2}\right)^p\! +\!
\left(\frac{\varepsilon}{2}\right)^p\!=\!\frac{\varepsilon ^p}{2^{p-1}}\!\leqslant \!\varepsilon ^p,
$$
i.e. $\rho (x,a^{(i_0)})<\varepsilon$. Hence $A=\{a^{(i)}\}_{i=1}^m$~is an $\varepsilon$-net for $M.$
By the Hausdorff theorem, $M$ is relatively compact in ${\bf{l}}_p.$

{\it Necessity.}
Let $M$ be relatively compact in ${\bf{l}}_p$, let $\varepsilon >0$ be arbitrary. Then $M$ is bounded (condition 1))
and 
there exists a finite $\frac{\varepsilon}{2}$-net $A\doteq\{a^{(i)}\}_{i=1}^m \subset l_p$ for $M.$ For each
$i\;(i\!=\!1,2,\ldots,m)$ there exists
natural numbers $N_i$ such that $\sum\limits_{k=N_i}^{\infty}|a^{(i)}_k|^p\!<\!~\left(\frac{\varepsilon}{ 2}\right)^p.$ Let be
$N\doteq \max\,\{N_1,\ldots,N_m\}.$ Then the inequality
$\sum\limits_{k=N}^{\infty}|a^{(i)}_k|^p<\left(\frac{\varepsilon}{2}\right)^p$
holds for all $i=1,2,\ldots,m$. Taking this inequality into account, we have for every $x\in M$ 
\begin{multline*}
\sum\limits_{k=N+1}^{\infty}|x_k|^p\leqslant \sum\limits_{k=N+1}^{\infty}|x_k-a^{(i)}_k |^p+\sum\limits_{k=N+1}^{\infty}|a^{(i)}_k|^p<\\
<\rho (x,a^{(i)})+\left(\frac{\varepsilon}{2}\right)^p<\left(\frac{\varepsilon}{2}\right)^p+ \left(\frac{\varepsilon}{2}\right)^p\leqslant \varepsilon ^p,
\end{multline*}
that is, condition 2) is satisfied.
\hfill $\square$

\medskip

$2^{\circ}.\;${\bf{Spaces}} ${\bf{c}}_0,\;{\bf{c}}\;$
\begin{utv}
\label{comp c_0}
{\it A set $M\subset {\bf{c}}_0$ is relatively compact if and only if
 
$1)$ it is bounded, and

$2)$ for any $\varepsilon >0$ there exists $N$ such that inequality $| x_k|<\varepsilon$ holds  for all $x\in M$.}
 
\end{utv}
{\it Sufficiency.} 
We use notations introduced in Proposition \ref{compl_1} for $p=\infty.$
Let conditions 1) , 2) be satisfied, let $\varepsilon >0$ be arbitrary and suppose that $N\in\mathbb N$ corresponds to $\varepsilon$ in
condition 2). Since set $\hat M\doteq\tilde F(M)\subset\mathbb R_{\infty}^N$
is bounded, it is relatively compact, so it admits a finite $\varepsilon$-net in $\mathbb R_{\infty}^N$,
$\hat A\doteq \{\hat a^{(i)}\}_{i=1}^{m}.$ Let $a^{(i)}=F^{-1}(\hat a^{(i)})\;
(i=1,2,\ldots,m).$ Let us show that
$A\doteq\{a^{(i)}\}_{i=1}^m$~ is an $\varepsilon$-net for $M.$

Let $x\in M$. Then $\hat x=\tilde F(x)\in\mathbb R_{\infty}^N.$ There exists $i_0\;(1\leqslant i_0\leqslant m)$ such that 
$\rho_{\mathbb R^N_{\infty}}(\hat x,\hat a^{(i_0)})<\varepsilon.$
$$
\rho_{{\bf{c}}_0} (x,a^{(i_0)})=\max\,\{\underset{1\leqslant k\leqslant N}{\max}|x_k-a^ {(i_0)}_k|,\;\underset{k>N}{\sup}|x_k|\}<\varepsilon,
$$
i.e. $A=\{a^{(i)}\}_{i=1}^m$ is an $\varepsilon$-net for $M.$

By the Hausdorff theorem, $M$ is relatively compact in ${\bf{c}}_0.$

{\it Necessity.}
Let  $M$ be relatively compact in ${\bf{c}}_0$ and let $\varepsilon >0$ be arbitrary. Then $M$ is bounded (condition 1))
and 
there exists a finite $\frac{\varepsilon}{2}$-net $A\doteq\{a^{(i)}\}_{i=1}^m \subset {\bf{c}}_0$ for $M.$ For each $i
(i=1,2,\ldots,m)$ there exists $N_i$ such that $\underset{k>N_i}{\sup}|a^{(i)}_k|<\frac{\varepsilon}{2}.$ Let $N\doteq 
\max\,\{N_1,\ldots,N_m\}.$ Then inequality $\underset{k>N}{\sup}|a^{(i)}_k|<\frac{\varepsilon}{2}$
holds for all $i=1,2,\ldots,m.$ Taking this inequality into account, we have for all $x\in M$ and $k>N$ 
$$
|x_k|\leqslant |x_k-a^{(i)}_k|+|a^{(i)}_k|\leqslant \underset{k>N}{\sup}|x_k-a^{(i) }_k|+\underset{k>N}{\sup}|a^{(i)}_k|<
\frac{\varepsilon}{2}+\frac{\varepsilon}{2}=\varepsilon,
$$
which yields condition 2).
\hfill $\square$

\medskip

One proves the following result in exactly the same way as the previous one.

\begin{utv}
\label{comp c}
{\it A set $M\subset {\bf{c}}$ is relatively compact if and only if
 
$1)$ it is bounded, and

$2)$ for any $\varepsilon >0$ there exists $N$ such that inequality }
 $| x_k-\underset{k\to \infty}{\lim}x_k|<\varepsilon$ holds for all $x\in M$.
\end{utv}

\medskip

$3^{\circ}.$\;{\bf{Space}} ${\bf{l}}_\infty\;$

\medskip

{\bf{A}}.\;Let $\mathcal A\subset \mathbb R,\;\;a_n:\mathcal A\to\mathbb R, \;n=1,2,\ldots,$
$$
M=\{a(\alpha)=\bigl(a_1(\alpha),\,a_2(\alpha),\ldots\bigr),\quad \alpha\in\mathcal A\}\subset {\bf{ l}}_\infty.
$$
\begin{teo}
\label{priznkM0}
{\it Let $\mathcal A$ be compact and let the function sequence $a(\cdot)$ be continuous on $\mathcal A.$ Then the set $M$
is compact.}
\end{teo}
\doc\,
An arbitrary infinite sequence from $M$ has form
\begin{equation}
\label{seqvM}
\bigl\{a\bigl(\alpha_k\bigr)\bigr\}_{k=1}^\infty,
\end{equation}
where
$
\{\alpha_k\}_{k=1}^\infty
$
is an arbitrary infinite sequence from $\mathcal A.$ Since $\mathcal A$ is compact, there exists a convergent subsequence
$$
\{\alpha_{k_j}\}_{j=1}^\infty\subset\{\alpha_k\}_{k=1}^\infty,\;\;\underset{j\to\infty}{\lim}\,\alpha_{k_j}\doteq\alpha^*\in\mathcal A.
$$
Since $a(\cdot)$ is continuous, $\underset{j\to\infty}{\lim}\,a(\alpha_{k_j})=a(\alpha^*)$
(meaning the limit in the space ${\bf{l}}_\infty$).
By the definition of $M$, we have $a(\alpha^*)\in M,$ so $M$ is closed. So, an arbitrary sequence (\ref{seqvM})
has a limit point in $M$, that is, $M$ is compact.
\hfill $\square$
\begin{zam}
\label{notnescont}
{\it The continuity condition for $a(\cdot)$ in the theorem is not necessary.}
\end{zam}
Indeed, let $a_n(\alpha)\equiv 0$ for $n>m$ for all $\alpha$ ($m$ is a fixed natural number),
$a_n(\cdot)\;\;(n\leqslant m)$ are characteristic functions of various subsets of $\mathcal A.$ Then $a(\cdot)$ can have
an infinite set of discontinuity points, while $M$ is compact since it contains a finite number of  points  $(2^m)$.
\hfill $\square$

\medskip

{\bf{B}}.\;The following is a proposition from \cite[IV.5.6]{dansch}
\begin{utv}
\label{DS282}
{\it A bounded set $M\subset {\bf{l}}_\infty$ is relatively compact if and only if
if for each $\varepsilon>0$ there exist a finite number of pairwise disjoint sets
$\mathcal N_1,\,\mathcal N_2,\,\ldots,\mathcal N_m$ from $\mathbb N,\;\bigcup\limits_{k=1}^m\mathcal N_k=\mathbb N$, and there are 
$n_i\in \mathcal N_i,\; i=1,2,\ldots,m,$} such that
$$
\underset{n\in\mathcal N_i}{\sup}\bigl|a_{n_i}-a_n\bigr|<\varepsilon,\quad a\in M,\quad i=1,2,\ldots,m.
$$
\end{utv}

\poonkt
{Open and closed sets} \

In subsections {\bf{1}}--{\bf{3}} below
$(\mathfrak M,\rho)$ is an arbitrary metric space.\, The distance $\rho(t,\,A)$ from a point $t\in \mathfrak M$ to a set
$A\subset\mathfrak M$ is defined as $\rho(t,\,A)\doteq \underset{t^{'}\in A}{\inf}\rho(t,t^{'})$.
The distance $\rho(A,\,B)$ between sets $A, B\subset\mathfrak M$ is defined as
$$\rho(A,\,B)\doteq \underset{t^{'}\in A,\,t^{''}\in B}{\inf}\rho(t^{'},t^ {''}).$$

{\bf{1.\;Open sets.\;}} 
$1^{\circ}.$\; {\it Let $A\subset\mathfrak M$ and let $d>0.$ 
Then $D\doteq \{t:\rho(t,A)<d\}$ is an open set and $A\subset D.$}

\doc\,
Since for the points from $A$ one has $\rho(t,A)=0,$ we have $A\subset D.$ Let $t_0\in D$. Since $\rho(t_0,A)<d ,$
there is a point $t^{'}\in A$ such that $\rho(t^{'},t_0)<d.$
Let $\delta\doteq d-\rho(t_0,\,t^{'})$ and let $s$ be an arbitrary point in the ball $B(t_0,\delta).$ Then
$$
\rho(t^{'},s)\leqslant \rho(t^{'},\,t_0)+\rho(t_0,\,s)=d-\delta+\rho(t_0,\,s) <d-\delta+\delta=d.
$$
Since $\rho(t^{'},\,s)<d$, we have furthermore
$\rho(s,A)<d,$ i.e. $s\in D$. Therefore, $B(t_0,\,\delta)\subset D,$ i.e. $D$~ is an open set.
\hfill $\square$

\medskip

$2^{\circ}.$\; \;{\it Let $A, B\subset\mathfrak M,\;A\ne\varnothing,\;B\ne\varnothing,$ where $\rho(A,B)=r>0.$ Let, furthermore,
$$A_1=\{t:\rho(t,A)<\frac{r}{2}\},\quad B_1=\{t:\rho(t,B)<\frac{r}{ 2}\}.$$ Then $A_1\bigcap B_1=\varnothing.$}

\doc\,
Suppose that $A_1\bigcap B_1\ne\varnothing$ and $z\in A_1\bigcap B_1.$ Then $\rho(z,A)<\frac{r}{2}$, $
\rho(z,B)<\frac{r}{2}$ 
and there are points $t_1\in A,\;\;t_2\in B$ such that $\rho(z,t_1)<\frac{r}{2},\linebreak
\rho(z,t_2)<\frac{ r}{2}.$ By the triangle inequality, this yields
$\rho(t_1, t_2)<r,$ and hence $\rho(A,B)<r,$ which contradicts to the definition of $r.$ Hence $A_1\bigcap B_1=\varnothing.$
\hfill $\square$

\medskip

$3^{\circ}.$\; {\it The union of any number of open sets is an open set. The intersection of a finite number of open sets is an open set.
The intersection of an infinite number of open sets may not be open.}

\doc\,
Let $G_{\alpha}\;(\alpha \in A, \;A$~ is an infinite set) be open sets, $G=\bigcup\limits_{\alpha \in A}G_{ \alpha}$ and $x\in G.$
Then $x\in G_{\alpha '}$ for some $\alpha'\in A$. Since $ G_{\alpha '}$ is open, there exists $\varepsilon>0$ such that
$B(x,\varepsilon)\subset G_{\alpha '}.$ Hence, $B(x,\varepsilon)\subset G,$ that is, $G$~ is an open set.

Let $G=\bigcap\limits_{k=1}^n G_k$, where $G_k$~ are open sets and $x\in G.$ Then $x\in G_k$ for $k=1,2 ,\ldots,n.$
There are $\varepsilon_k>0,\;k=1,2,\ldots,n$ such that $B(x,\varepsilon_k)\subset G_k$ for $k=1,2,\ldots,n.$ Let
$\varepsilon\doteq \underset{1\leqslant k\leqslant n}{\min}\varepsilon_k.$ Then $B(x,\varepsilon)\subset G_k$ for $k=1,2,\ldots,n, $ that is,
$B(x,\varepsilon)\subset G$. This means that $G$~ is an open set.
\hfill $\square$

\medskip

An infinite intersection of open sets may not be open. This is seen from the following example $\mathfrak M=\mathbb R:$
$\bigcap\limits_{k=1}^{\infty}(-\frac{1}{k},1+\frac{1}{k})=[0,1].$
\hfill $\square$

\medskip

Let $G_k\;(k=1,2,\ldots)$ be open sets. Then the set $\bigcap\limits_{k=1}^{\infty}G_k$ is called {\it a set of type} $G_{\delta}.$ Obviously, a finite or countable intersection of sets of type $G_{\delta}$
is a set of type $G_{\delta}.$

\medskip

{\bf{2.\;Closed sets.}}\;
$1^{\circ}.$\; {\it A closed set is a set of type $G_{\delta}.$}

\doc\,
Let $F$~ be a closed set. Put 
$G_n\doteq\{t:\rho(t,F)<\frac{1}{n}\};$
By the results of subsection {\bf{1.}1}, 
$G_n$ are open sets. The assertion now follows from
 the representation $F=\bigcup\limits_{n=1}^\infty G_n.$
\hfill $\square$
 
 $2^{\circ}.$\; {\it The intersection of any number of closed sets is a closed set. 
 The union of a finite number of closed sets is a closed set. The union of an infinite collection of closed sets may not be a closed set.}

\doc\,
Let $F_{\alpha}\;(\alpha \in A, \;A$~ is an infinite set)~ be closed sets, $F=\bigcap\limits_{\alpha \in A}F_{ \alpha}$ and $x*$~ be a
limit point of $F.$ There is a sequence $\{x_n\}_{n=1}^{\infty}\subset F$ such that $x_n\to x^*.$ By
the definition of intersection, $\{x_n\}_{n=1}^{\infty}\subset F_{\alpha}$ for all $\alpha \in A.$ Therefore, $x^*$~ is a limit point
of all sets $F_{\alpha}$ and, since they are closed, $x^*\in F_{\alpha}$ for all $\alpha \in A.$ This means that $x^*\in F,$
that is, $F$ is closed.

Let $F=\bigcup\limits_{k=1}^n F_k,$ where $F_k$~are closed sets and $x^*$~is a limit point of $F.$
There is a sequence $\{x_n\}_{n=1}^{\infty}\subset F$ such that $x_n\to x^*.$ This means that $\{x_n\}_{n=1 }^{\infty}\subset F_{k_0}$
for some $k_0\;(1\leqslant k_0\leqslant n).$ Therefore, $x^*$~is a limit point of $F_{k_0}.$
Since $F_{k_0}$~ is a closed set, we have $x^*\in F_{k_0}.$ Hence $x^*\in F,$ i.e. $F$ is closed.

An infinite union of closed sets may not be closed. This is seen from the example $\mathfrak M=\mathbb R:$
$\bigcup\limits_{k=1}^{\infty}[\frac{1}{k},2-\frac{1}{k}]=(0,2).$
\hfill $\square$

Let $F_k\;(k=1,2,\ldots)$~ be closed sets. Then the set $\bigcup\limits_{k=1}^{\infty}F_k$ is said to be
{\it a set of type} $F_{\sigma}.$ Clearly, a finite or countable union of sets of type $F_{\sigma}$
is a set of type $F_{\sigma}.$

 $3^{\circ}.$\; {\it The complement of an open} ({\it closed}) {\it set is a closed} ({\it open}) {\it set.}

\doc\,
Let $G\subset\mathfrak M$ be an open set, and let $x^*$ be a limit point of the complement $\overline{G}\doteq\mathfrak M\setminus G$
Assume that $x^*\ne \overline{G}.$ Hence $x^*\in G.$ Since $G$~ is an open set, there exists $\varepsilon >0,$ such that ball
$B(x^*,\varepsilon)\subset G.$ Consequently, the ball $B(x^*,\varepsilon)$ does not contain points in the set $\overline{G},$, which contradicts
definition of the limit point. Hence $x^*\in\overline{G},$ i.e. $\overline{G}$~ is closed.

Let $F\subset\mathfrak M$~ be a closed set and let $x\in\overline{F}.$ Assume that the set $\overline{F}$ is not open.
This means that for any $n\in\mathbb N$ set
$B(x,\frac{1}{n})\setminus\overline{F}\ne\varnothing.$ Let $x_n\in B(x,\frac{1}{n})\setminus\overline{ F}$
for $n=1,2,\ldots$ It is obvious that $\{x_n\}_{n=1}^{\infty}\subset F$ and $x_n\to x.$ This means that $x$ ~is a  limit point of $F$,  
and that $x\in F,$ which contradicts to $x\in\overline{F}.$ Hence the set $\overline{F}$ is open.
\hfill $\square$

\medskip

{\bf{3.\;Sets of type $F_{\sigma}$.}}\;~
$1^{\circ}.$\; {\it The intersection of sets of type $F_{\sigma}$ is a set of type $F_{\sigma}.$}

\doc\,
Let $$A=\left(\bigcup\limits_{k=1}^\infty A_k\right)\bigcap\left(\bigcup\limits_{j=1}^\infty B_j\right),$$ where $A_k$ and $B_i$ are closed sets.
Then $A=\bigcup\limits_{k,j=1}^\infty\left(A_k\bigcap B_j\right).$
\hfill $\square$

 $2^{\circ}.$\; {\it The complement of a set of type $F_{\sigma}\;\;\Bigl(G_{\delta}\Bigr)$ is a set of type} 
 $G_{\delta}\;\;\Bigl(F_ {\sigma}\Bigr).$

\doc\,
Let $A$~ be a set of type $F_{\sigma},$ i.e. $A=\bigcup\limits_{n=1}^\infty F_n,$ where $F_n$~ are closed sets. Then
$\overline{A}=\bigcap\limits_{n=1}^\infty \overline{F_n}$. According to subsection {\bf{3.}}1, the sets $\overline{F_n}$~ are open.
This proves the first assertion. The second assertion is proved similarly.
\hfill $\square$

\medskip

 $3^{\circ}.$\;From the results of  subsections {\bf{3.}}2, {\bf{2.}}1 and {\bf{2.}}3 we obtain the following statement:

\medskip

{\it An open set is a set of type $F_{\sigma}.$}

\bigskip

 $4^{\circ}.$\;{\it The difference of closed sets is a set of type $F_{\sigma}.$}

 \doc\,
 Since $A\setminus B=A\bigcap\overline B,$ and set $\overline B$ is open, see subsection {\bf{2.}}3, we have, in view of subsection {\bf{3.}}3, that the  set $A\setminus B$ is of type $F_{\sigma}.$
\hfill $\square$

\bigskip

 $5^{\circ}.$\;{\it Let $$M=\bigcup\limits_{k=1}^n A_k,$$ where $A_k$ are some sets of type $F_{\sigma}.$
Then there is another representation $M=\bigcup\limits_{k=1}^n B_k,$ where $B_k$ are also sets of type $F_{\sigma},$ moreover}
$$
B_k\subset A_k\;(k=1,2,\ldots,n),\;\;B_k\bigcap B_j=\varnothing\;\;\text{{\it at}}\;\; j\ne k.
$$

\doc\,
In the end $M=\bigcup\limits_{k=1}^n F_k,$ where $F_k$~ are closed sets, each of which is entirely contained in one of the sets
$A_i.$ We assume (see subsection {\bf{3.}}4)
$$
S_1=F_1,\quad S_2=F_2\setminus F_1,\quad \ldots \quad S_k=F_k\setminus \bigcup\limits_{i=1}^{k-1} F_i.
$$
By construction, $S_k$ are sets of type $F_{\sigma}$, they are pairwise disjoint and $M=\bigcup\limits_{k=1}^n S_k.$

Let us divide the collection $T$ of all $S_k$ into $n$ parts $T_1,T_2,\ldots,T_n$. We do it as follows. In $T_1$ we include those sets $S_k$ which
are contained in $A_1,$ in $T_2$ we include those $S_k\in T\setminus T_1$ that are contained in $A_2,$, etc.
Putting $B_i\doteq \bigcup\limits_{S_k\in T_i} S_k,$ we get the required representation.
\hfill $\square$

\bigskip

{\bf{4.\;The structure of open and closed sets in $\mathbb R.$}}\;
 $1^{\circ}.$\; {\it A non-empty bounded open set in $\mathbb R$ 
 can be represented the union of at most countably many non-intersecting ($\equiv$ pairwise disjoint) open intervals.}

\doc\,
Let~$G$~ be a non-empty open bounded set and let $t_0$ be an arbitrary point in $G$. Consider the set 
$F=[t_0,\,+\fy) \cap \ol G,$ where $\ol G=\R \setminus G$~ is the complement~$G.$ Since the complement of an 
open set is closed, $F$ is the intersection of closed sets, and hence itself is a  closed set. Since~$G$ is bounded, $F$ is non-empty.
It is also obvious that~$F$ is bounded from below (for example, by $t_0$).
Let $\bt=\inf F$. Due to the closedness~$F$, we have $\bt\in F,$ $\bt \gs t_0.$
Since $t_0\in G$ (so $t_0 \not\in \ol G$), $t_0\ne\bt,$
i.e. $t_0<\bt$ and $\bt \not\in G$ (because $\bt\in F \subset \ol G$).

\parindent=1cm
Suppose that there is $s\in [t_0,\,\bt),$ $s\not\in G.$ Then~$s\in F$
and~$s<\bt,$ which contradicts to the definition of $\bt.$ Hence,
$[t_0,\,\bt)\subset G.$

\parindent=1cm
In the same way, we prove the existence of $\al< t_0,$ $\al \not\in G,$
$(\al,\,t_0] \subset G.$ Thus, each point of set~$G$
is contained in some interval whose endpoints are not contained in~$G.$ This
means that~$G$ is an open interval or a union of
non-intersecting open intervals.

\parindent=1cm
Moreover, there can be only finitely many or countably many intervals that constitute $G$.
Indeed, let us fix in each of the constituent intervals a rational
point. This puts the set of intervals that constitute $G$ in a one-to-one correspondence with a subset of the set of rational numbers 
$\Q$, which yields the result.
\hfill $\square$

 $2^{\circ}.$\; {\it A non-empty closed bounded set in~$\R$ is either
a closed interval, or is obtained from a closed interval by deleting a finite or countable
collection of non-intersecting open intervals.}

\doc\,
Let~$F$~ be a non-empty bounded closed set. It follows from the boundedness of $F$ that there is a closed interval $[a,\,b] \supset F,$ $a=\inf F,$ $b=\sup F.$
If ${F\ne[a,\,b]},$
then $G=[a,\,b] \setminus F$ is a non-empty open bounded set.
By the previous result, $G=\bigcup\lt_k (\al_k,\,\bt_k),$
where intervals $(\al_k,\,\bt_k)$ are pairwise disjoint and their are only finitely many or countably many of them. Therefore,
$
F=[a,\,b] \setminus G=[a,\,b] \setminus \bigcup\lt_k (\al_k,\,\bt_k).
$
\hfill $\square$

 $3^{\circ}.$\; {\it A perfect bounded set in~$\R$ is either a closed interval or
is obtained from a closed interval by removing finitely or countably many disjoint
open intervals that do not have common ends and do not have common ends
with the original closed interval.}

\doc\,
Let~$C$~ be a bounded perfect set in~$\R.$ Since~$C$ is closed, it has form
$
C=[a,\,b]$ or $C=[a,\,b] \setminus \bigcup\lt_k (\al_k,\,\bt_k)
$
(see the previous statement). If any of the two open intervals would have a common endpoint,
or the endpoints of one or two intervals are points~$a$ or~$b,$
then all these points would be isolated. But this contradicts to the hypothesis that $C$ is a perfect set.
\hfill $\square$

\medskip

{\bf{5.\;Cantor open and perfect sets.}}\; 
Let us divide interval $J=[0,\,1]$ into three equal parts by points $\fracs 13$
and $\fracs 23$ and remove the middle interval from it~$G_1=(\fracs 13,\,\fracs 23).$
We divide each of the remaining closed intervals $[0,\,\fracs 13]$ and~$[\fracs 23,\,1]$
 into three equal parts by points $\fracs 19$ and~$\fracs 29$~(first)
and~$\fracs 79$ and~$\fracs 89$ (second) and remove the middle intervals, i.e.
the set $G_2=(\fracs 19,\,\fracs 29)\cup (\fracs 79,\,\fracs 89).$
In the third step, we will remove the middle parts (intervals) of the four remaining
closed intervals, that is, the set $G_3=(\fracs{1}{27},\,\fracs{2}{27})\cup
(\fracs{7}{27},\,\fracs{8}{27})\cup(\fracs{19}{27},\,\fracs{20}{27})\cup
(\fracs{25}{27},\,\fracs{26}{27}).$ We continue this process indefinitely:
at the $k$-th step we remove~$2^{k-1}$ intervals, the middle parts of the closed interval,
remaining from the previous steps, that is, the set
${G_k=(\fracs{1}{3^k},\,\fracs{2}{3^k})\cup\ldots\cup(\fracs{3^k-2}{3^k} ,\,
\fracs{3^k-1}{3^k})}.$

As a result, we will remove the open set $G_0=\bigcup\lt_{k=1}^\fy G_k$.
The sets~$G_0$ and~${\cal K}=J\setminus G_0$ are called Cantor sets.
The intervals that make up~$G_k$ are called adjacent intervals of the $k$-th rank
(there are only~$2^{k-1}$ of them). The endpoints of adjacent intervals are called the points of the first
kind. The remaining points~${\cal K}$ are called the points of the second kind. It is clear that 0 and 1  are points of the second kind. We also denote by $\Gamma\doteq\bigcup\limits_{k=1}^\infty \gamma_k$ the set of points of the first kind, numbered
in ascending order of rank and from left to right within the same rank.

Let us prove the following statements.

a)~${\cal K}$~ is a perfect nowhere dense set;

b)~${\cal K}$ has the cardinality of continuum;

c)~$\Gamma$ is everywhere dense in~${\cal K},$ i.e., $cl\,\Gamma={\cal K};$

d)~$G_0$ is everywhere dense in~$[0,1],$ i.e., $cl\,G_0=[0,1];$

e) the sum of the lengths of all adjacent intervals (`length' of set $G_0$) is equal to 1.

\doc\,
a) Since $\pK=[0,\,1] \setminus \bigcup\lt_{k=1}^{\fy} G_k,$ and
the intervals constituting $G_k$ do not have common ends, then,   
by statement $3^{\circ}$ (see subsection {\bf{1.9.4}}), $\pK$ is perfect.

Let $(\al,\,\bt) \subset [0,\,1]$ be an arbitrary interval. There is
a number $N$ such that $\fracs{1}{3^k}<\bt-\al$ for $k>N$ ($N=\Big[\log_3 \frac{1}{\bt-\al} \Big]$). For such~$k$,  $(\al,\,\bt)$ contains an interval to be removed, that is, an interval that does not contain points of the set~$\pK.$
Thus, $\pK$ is nowhere dense.

b) Each point~$t$ of  interval $(0,\,1)$ can be represented as
$$
t=\frac{a_1}{3}+\frac{a_2}{3^2}+\frac{a_3}{3^3}+\ldots,
$$
where~$a_i$ take one of the values 0,\,1,\,2. This representation
is written as an infinite ternary fraction 
$t=0,\,a_1,\,a_2,\,a_3,\ldots$

For example,
$$
\aligned
\frac 12 &= \frac 13+\frac 19+\frac{1}{27} +\ldots=0,\,111\ldots, \\
\frac 14 &= \frac 29+\frac{2}{81}+\frac{2}{729}+\ldots=0,\,020202\ldots, \\
\frac 34 &= \frac 23+\frac{2}{27}+\frac{2}{243}+\ldots=0,\,20202\ldots
\endaligned
$$
In this case, the ends of the removed intervals are ternary rational
numbers. They admit two representations as a ternary fraction:
$$
\aligned
{} & \frac 13=0,1000\ldots=0,0222\ldots;\\
{} & \frac 23=0,2000\ldots=0,1222\ldots; \\
{} & \frac 19=0,0100\ldots=0,00222\ldots\\
{} & \frac 23=0,02000\ldots=0,001222\ldots; \\
{} & \frac 79=0,2100\ldots=0,20222\ldots; \\
{} & \frac 89=0,022000\ldots=0,21222\ldots
\endaligned
$$

For the points of the interval~\;$\big( 0,\,\fracs 13 \big)$ one has $\ {a_1=0}$. For the points
of interval
$\big( \fracs 13,\,\fracs 23 \big)$\;\;one has ${a_1=1}$, and for the points of 
interval~\;$\big( \fracs 23,\,1 \big)\;\;$ $\ {a_1=2}.$

At the first step of the process that gave us $\pK,$ i.e.\,when we removed
interval~$\big( \fracs 13,\,\fracs 23 \big),$ we have in fact removed from~$(0,\,1)$
the numbers whose ternary
representation has $a_1=1$. 
At the second step, 
when we removed intervals~$\big( \fracs 19,\,\fracs 29 \big)$ and
$\big( \fracs 79,\,\fracs 89 \big)$, we actually removed from~$(0,\,1)$ the numbers whose ternary representation has $a_2=1$. At the $k$-th step, when
we removed $2^{k-1}$ intervals of length~$\fracs{1}{3^k}$, we removed from~$(0,\,1)$
those numbers that have $a_k=1$ in the ternary representation. Thus, the set~$G_0$ consists of numbers,
whose ternary representation 
necessarily contains ones, and~$\pK$
consists of numbers having a ternary representation that does not contain 1,
while $0=0,\,00\ldots,$ $1=0,\,2222\ldots$
(obviously, $0\in \pK$ and~$1\in \pK$).

Note now that each point of the interval~$(0,\,1)$ can also be
represented as a binary fraction:
$$
t=\frac{b_1}{2}+\frac{b_2}{2^2}+\frac{b_3}{2^3}+\ldots=0,\,b_1b_2b_3\ldots,
$$
where~$b_i$ take values 0 or 1 (with $1=0,\,111\ldots$).

Let $\vfi: \pK \to[0,\,1]$ be defined as follows: for
$t=0,\,a_1a_2a_3\ldots \in \pK$ ($a_i=0$ or $a_i=2$) we set
$$
{\vfi(t)=0,b_1b_2b_3\ldots \in [0,\,1]}
$$ ($b_i=0$ or $b_i=1$),
where~$b_i=0,$ if~$a_i=0,$ $b_i=1,$ if~$a_i=2.$ Clearly,
$\vfi(\cdot)$ is a bijection between $\pK$ and $[0,\,1]$, so $\pK$ has the cardinality of continuum.

c) By construction, any neighborhood of each point~$t\in \pK$ contains an infinite set of points from~ $\Gamma$. So, $t$
is a limit point for~$\Gamma,$ i.e.\,$cl\,\Gamma={\cal K}.$

d)~Arguing similarly, we obtain that $\cl G_0=[0,\,1],$ i.e. $G_0$ is dense in~$[0,\,1].$

e) The sum of the lengths of adjacent (removed) intervals is equal to
$$
\frac{1}{3}+2\cdot\frac{1}{9}+4\cdot\frac{1}{27}+\ldots+2^{k-1}\cdot\frac{1} {3^k}+\ldots=
\sum\limits_{k=1}^\infty\frac{2^{k-1}}{3^k}=\dfrac{\frac{1}{3}}{1-\frac{2}{ 3}}=1.\quad \hfill \square
$$


\poonkt
{Continuous functions. Interchanging limits} \ 
\phantom{12345678901234567890}

{\bf{1.\;Continuous functions.\;}}

$1^\circ.$
{\it Let $(\mathfrak M,\,\rho)$~ be a metric space, let 
$F_1, F_2,\ldots,F_n\subset\mathfrak M$ be closed, pairwise disjoint sets, and let $F=\bigcup\limits_{k=1}^n F_k.$
Then, if function $x:F\to\mathbb R$ is constant on each of the sets $F_k,$ it is continuous on set $F.$}

\doc\;
Let $t_0\in F$ and let $t_m\to t_0\;(m\to\infty).$ There exists $k_0$ such that $t_0\in F_{k_0}.$
Since sets $F_k$ are pairwise disjoint, $t_0\notin F_k$ for $k\ne k_0,$ and $t_0$ is not a limit point of $F_k$
(since it is closed).
Therefore, sequence $\{t_m\}_{m=1}^\infty$ can contain only a finite number of points belonging to $F_k$ for $k\ne k_0.$
Let $t_{m_0}$ be the last term of this sequence included in one of the sets
$
F_1,\ldots,F_{k_0-1},F_{k_0+1},\ldots,F_n.
$
Then, for $m>m_0$, $x(t_m)=x(t_0)$. Therefore, $x(\cdot)$ is continuous at the point $t_0.$
\hfill $\square$

$2^\circ.$
{\it Let $(\mathfrak M,\,\rho)$~ be a metric space. Given
$n$ pairwise disjoint closed sets 
$F_1,F_2,\ldots,F_n\;\bigl(\subset\mathfrak M\bigr)$ and
$n$ real numbers $y_1,y_2,\ldots,y_n,$
there exists a continuous function $x:\mathfrak M\to \mathbb R$ such that $x(t)=y_k$ for $t\in F_k\;\;(k=\linebreak=
1,2,\ldots,n)$
and}
$$
\underset{t\in \mathfrak M}{\sup}\;\;|x(t)|\leqslant\underset{1\leqslant k\leqslant n}{\max}\;\;|y_k|.
$$

\doc\;
According to what was said in subsection {\bf{1.2}},$9^{\circ}$, 
the function $g(t)\doteq \rho\,(t, A)\;\big(A\subset \mathfrak M\bigr)$ 
is continuous on the whole space $\mathfrak M.$ Now it is easy to verify that the function
$$
x(t)\doteq\sum\limits_{k=1}^n
\dfrac{y_k\cdot \rho\left(t,\,\bigcup\limits_{i\ne k} F_i\right)}{\rho\bigl(t,\,F_k\bigr)+\rho\left( t,\,\bigcup\limits_{i\ne k} F_i\right)}
$$
satisfies our requirements.
\hfill $\square$

\bigskip

$3^\circ.$
\textbf{Continuity theorem for the sum of a series} \phantom{012345678901234}
(cf.\,\cite[sect. 431]{ficht2}.) {\it Let $(\mathfrak M,\,\rho)$~ be a metric space, let functions 
$$a_n:\mathfrak M\to \mathbb R,
\quad n=1,2,\ldots,$$ be continuous on $\mathfrak M,$, assume that the series $\sum\limits_{n=1}^\infty a_n(\cdot)$ converges uniformly on $\mathfrak M.$ 
Then the sum  of the series
$$
S(t)\doteq \sum\limits_{n=1}^\infty a_n(t)
$$
is continuous on} $\mathfrak M$

\doc\;
Let $t_0\in \mathfrak M$ be an arbitrary point. Due to the uniform convergence of the series, for every $\varepsilon>0$
there exists a positive integer $n$ such that the remainder of the series $r_n(t)\doteq \sum\limits_{k=n+1}^{\infty} a_k(t)$ satisfies  inequality
$|r_n(t)|<\dfrac{\varepsilon}{3}$ for all $t\in \mathfrak M$. In particular, $|r_n(t_0)|<\dfrac{\varepsilon}{3}. $ For this $n$, the finite sum
$S_n(t)\doteq \sum\limits_{k=1}^{n} a_k(t)$ of continuous functions is a continuous function. Therefore, there exists $\delta>0$ such that
for all $t\in \mathfrak M$ satisfying  inequality $\rho(t,t_0)<\delta,$ the inequality
$|S_n(t)-S_n(t_0)|<\dfrac{\varepsilon}{3}$ 
is satisfied. Since
$$
S(t)=S_n(t)+r_n(t),\qquad S(t_0)=S_n(t_0)+r_n(t_0),
$$
subtracting from the first representation the second one, we obtain
$$
\bigl|S(t)-S(t_0)\bigr|\leqslant \bigl|S_n(t)-S_n(t_0)\bigr|+\bigl|r_n(t)\bigr|+\bigl|r_n(t_0 )\bigr|<\varepsilon,
$$
i.e.\,function $S(\cdot)$ is continuous at point $t_0.$
\hfill $\square$

\bigskip

$4^\circ.$
\textbf{Passing to limit under the series sign}
(\cite[sect. 433]{ficht2}) 
 
Taking in the previous subsection $\mathfrak M=\mathbb N,\;\rho(x,y)=|x-y|$, we obtain the following statement.

\smallskip

{\it Let\;\;$a:\mathbb N\times\mathbb N\to\mathbb R,\;\;a_{nm}\to a_{n\,\bullet}\;(m\to\infty)$ for every $n\in\mathbb N,$
assume that the series \;\;$\sum\limits_{n=1}^\infty a_{nm}$ converges uniformly with respect to $m\in\mathbb N.$
Then

1)\; the series converges\;$\sum\limits_{n=1}^\infty a_{n\,\bullet},$ and}\;

2)\;$\underset{m\to\infty}{lim}\;\sum\limits_{n=1}^\infty a_{nm}=\sum\limits_{n=1}^\infty a_{ n\,\bullet}.$

\bigskip

$5^\circ.$
\textbf{Theorem on extension by continuity} (cf.\;\cite[I.5.3]{dansch})
  
{\it Let $F$~ be a closed set in a metric space $(\mathfrak M,\,\rho),\,$ let $x:F\to\mathbb R$ be a continuous function on $F.$
 Then there exists a function $y:\mathfrak M\to\mathbb R$ with the following properties:
 
 1)\;$y(\cdot)$ is continuous on $\mathfrak M;$
 
 2)\;$y(t)=x(t)$ for} $t\in F;$
 
 3)\;$\underset{t\in\mathfrak M}{\sup}\;|y(t)|=\underset{t\in F}{\sup}|x(t)|.$

\doc\;(cf.\,\cite[I.5.3]{dansch}.)
Without loss of generality, we can assume that $x(\cdot)$ is not identically zero (the assertion is trivial for such $x(\cdot)$).
Put
\begin{multline*}
x_0(t)\doteq x(t),\;m_0\doteq \underset{t\in F}{\sup}\,|x_0(t)|, 
\quad F_{0-}\doteq\{t\in F:\;x_0(t)\leqslant -\dfrac{m_0}{3}\},\\ 
F_{0+}\doteq\{t\in F:\;x_0(t)\geqslant \dfrac{m_0}{3}\}.
\end{multline*}
The sets $F_{0-}$, $F_{0+}$ 
are closed and have empty intersection. By the statement in $\mathbf{3}^\circ$, 
there exists a continuous function
$y_0:\mathfrak M\to \mathbb R$ such that
\begin{multline*}
y_0(t)=-\dfrac{m_0}{3}\;\;(t\in F_{0-}),\quad y_0(t)=\dfrac{m_0}{3}\;\;(t \in F_{0+}),\\
-\dfrac{m_0}{3}\leqslant y_0(t)\leqslant \dfrac{m_0}{3}\;\;(t\in \mathfrak M).
\end{multline*}
We put $x_1(t)=x_0(t)-y_0(t)\;\;(t\in F).$ The function $x_1(\cdot)$ is continuous on $F$ and
$m_1\doteq \underset{t\in F}{\sup}\,|x_1(t)|\leqslant \dfrac{2}{3}m_0.$

The same argument, applied to $x_1$ and $m_1$ instead of $x_0$ and $m_0,$ gives us function $y_1(\cdot)$ 
defined and continuous everywhere on 
$\mathfrak M.$ Continuing this process, we obtain a sequence $\{y_n\}_{n=0}^\infty$ of continuous functions,
defined on the entire space $\mathfrak M$ and having properties
\begin{equation}
\label{ind1}
\left|x(t)-\sum\limits_{k=0}^n y_k(t)\right|\leqslant \left(\dfrac{2}{3}\right)^{n+1}m_0, \qquad t\in F
\end{equation}
and
\begin{equation}
\label{ind2}
\underset{t\in\mathfrak M}{\sup}\;\;|y_n(t)|\leqslant \dfrac{1}{3}\cdot \left(\dfrac{2}{3}\right)^ {n}m_0.
\end{equation}
Inequalities (\ref{ind2}) mean that the series
\begin{equation}
\label{ind}
y(t)\doteq\sum\limits_{n=0}^\infty y_k(t)
\end{equation}
converges uniformly on $\mathfrak M,$ hence, according to the result of paragraph {\bf{3}},
its sum 
$y(\cdot)$ is continuous on $\mathfrak M$.
Inequalities (\ref{ind1}) imply that for $t\in F$ 
$y(t)=x(t).\;$
Condition 3) follows from the proposition in paragraph {\bf{2}}.
\hfill $\square$

\medskip

{\bf{2}}.\;{\bf{Interchanging limits.\;}} 
Let $(\mathfrak M,\,\rho)$~ be a complete \textit{vector} metric space. Here we will discuss functions
\begin{equation}
\label{functip}
x:J\to \mathfrak M,\;\;\text{where}\;\;J\subset \mathbb R^n \text{ is a domain (possibly closed) }\;(n\geqslant 1)
\end{equation}
For the functions of this type, one can define the basic concepts of
mathematical analysis (limit, continuity, uniform continuity, derivative, Riemann integral, etc.).

So, for example, a point $a\in\mathfrak M$ is called the limit of function
$x:J\to \mathfrak M$ at $t_0$ where $t_0$ is a limit point of $J,$
if, for any sequence $\{t_n\}_{n=1}^\infty\subset J,\;\;t_n\to t_0$ one has $x_n(t)\to a$ in $\mathfrak M.$
The analytic concepts are transferred similarly. Most of the facts of the differential and integral calculus of a scalar function of a scalar argument
remain valid for the such functions. (Except the mean value theorems!)

Let $(\mathfrak M,\,\rho)$~ be a complete vector metric space, let
\begin{equation}
\label{poslnm}
\{x_{nm}\}_{n,m=1}^\infty
\end{equation}
be a sequence of its elements
(essentially, a function of type (\ref{functip}) with $J=\mathbb N\times\mathbb N\subset \mathbb R^2).$

$1^\circ.$ \;The following results serves as the basis for changing the order of limits
(Moore-Shatunovsky Lemma, see \cite[I.6.6]{dansch}).
\begin{teo}
\label{SM}
{\it Assume that the following conditions are satisfied:

\smallskip

a$)$\;for every $m\in\mathbb N$ there is limit $x_{\bullet m}\doteq \underset{n\to\infty}{lim}\;x_{nm};$

b$)$\;limit $x_{n}
\doteq\underset{m\to\infty}{\lim}\;x_{nm}$ exists uniformly with respect to $n\in\mathbb N.$

\medskip

Then the repeated limits
\begin{equation}
\label{SMpovt}
\underset{m\to\infty}{lim}\;x_{\bullet m}\quad\text{and}\quad \underset{n\to\infty}{lim}\;x_{n\bullet}
\end{equation}
and the double limit
\begin{equation}
\label{SMdvo}
\underset{n,m\to\infty}{lim}\;x_{nm}
\end{equation}
exist and are equal.}
\end{teo}
\doc\,
Let $\varepsilon>0$~ be arbitrary. According to condition b), there is a natural number $m_0$ such that, for
$m\geqslant m_0\;$
$\;
\;\rho(x_{nm},\,x_{n\bullet})<\dfrac{\varepsilon}{8}
$\;
for all $n\in\mathbb N.$ It follows that
$$
\rho(x_{nm},\,x_{n,m_0})\leqslant \rho(x_{nm},\,x_{n\bullet})+\rho(x_{n\bullet},\,x_ {n,m_0})<\dfrac{\varepsilon}{4}
$$
for $n\in\mathbb N,\;m\geqslant m_0.$ Due to condition a), there exists a positive integer $n_0$ such that for $n\geqslant n_0$
$\;
\rho(x_{n,m_0},\,x_{\bullet m_0})<\dfrac{\varepsilon}{8},
$\;
therefore
$$
\rho(x_{n,m_0},\,x_{n_0,m_0})\leqslant \rho(x_{n,m_0},\,x_{\bullet m_0})+\rho(x_{\bullet m_0} ,\,x_{n_0,m_0})<\dfrac{\varepsilon}{4},
$$
and
$$
\rho\bigl(x_{nm},\,x_{n_0,m_0}\bigr)\leqslant \rho\bigl(x_{nm},\,x_{n,m_0}\bigr)+\rho\bigl( x_{n,m_0}\,x_{n_0,m_0}\bigr)
<\dfrac{\varepsilon}{2}.
$$
Thus, if $n,\,n^{'}\geqslant n_0,\;\;m,\,m^{'}\geqslant m_0,$ then
$$
\rho\bigl(x_{nm},\,x_{n^{'},m^{'}}\bigr)\leqslant \rho\bigl(x_{nm},\,x_{n_0,m_0}\bigr)+\rho\bigl(x_{n_0,m_0},\,x_{n^{'},m^{'}}\bigr)<\varepsilon,
$$
that is, the sequence $\{x_{nm}\}_{n,m=1}^\infty$ is fundamental. This means that there exists the double limit (\ref{SMdvo}):\;
$p\doteq \underset{n,m\to\infty}{lim}\;x_{nm}.$

Passing to the limit $n,m\to\infty$ in inequality $\rho\bigl(x_{nm},\,x_{n^{'},m^{'}}\bigr)<\varepsilon$  gives us
$$
\rho\bigl(p,\,x_{n^{'},\,m^{'}}\bigr)=\rho\bigl(\underset{n,m\to\infty}{\lim}\; x_{nm},\,x_{n^{'},\,m^{'}}\bigr)=
\underset{n,m\to\infty}{\lim}\;\rho\bigl(x_{nm},\,x_{n^{'},\,m^{'}}\bigr)\leqslant \varepsilon
$$
for $\;n^{'}\geqslant n_0,\;\;m^{'}\geqslant m_0.$
Consequently,
$$
\rho\bigl(p,\,x_{\bullet m^{'}}\bigr)=\rho\bigl(p,\,\underset{n^{'}\to\infty}{lim}\; x_{n^{'},\,m^{'}}\bigr)\leqslant \varepsilon
$$
at\;$m^{'}\geqslant m_0.$
Therefore, $\underset{m\to\infty}{\lim}\;x_{\bullet m}=p.$
One proved in the same way $\underset{n\to\infty}{lim}\;x_{n\bullet}=p.$
\hfill $\square$

\medskip

$2^\circ.$ We will show below (subsection  {\bf{5.5.3}}) that the {\it uniform convergence} of  sequence 
(\ref{poslnm}) in one of its indices
is only a sufficient condition for changing the order of limits. 
The basis for obtaining the necessary and sufficient conditions
for changing the order of limit is the Arzela theorem (see 
\cite[IV.6.11]{dansch}, \cite[sect. 432]{ficht2}; see also \cite{DK02}). Let us state it here
for functions $\mathfrak M\to\mathbb R.$

Let us say that the sequence $\{x_n(\cdot)\}_{n=1}^\infty$ of functions $\mathfrak M\to\mathbb R$ converges to a function $x(\cdot)$ of the same type
{\it quasi-uniformly} on $\mathfrak M$
$$
(\text{write:}\qquad\;\; x_n(t)\to x(t)\;\;(QU)\;\;\text{on}\;\; \mathfrak M),
$$
if
$$
x_n(t)\to x(t)\;\;\text{for all}\;\;t\in\mathfrak M
$$
and
$$
\bigl(\forall\,\varepsilon>0\bigr)\;\;\bigl(\forall\,n_0\in\mathbb N\bigr)\;\;
\bigl(\exists\,p\in\mathbb N\bigr)\;\;\bigl(\exists\,n_1,\,\ldots,\,n_p\,(\geqslant n_0)\bigr)\;\;\bigl(\forall\,t\in\mathfrak M\bigr)\;\;
$$
$$
\underset{i=1,\,\ldots,\,p}{min}\;|x_{n_i}(t)-x(t)|<\varepsilon.
$$

In what follows, it will be convenient sometimes  to  state this definition in the language of the series.

Let us say that the  series
$$
\sum\limits_{k=1}^\infty a_{k}(t)\;\;(t\in \mathfrak M)\;\;\;\left(\sum\limits_{k=1}^ \infty a_{kn}\;\;\;(n\in\mathbb N)\right)
$$
converges quasi-uniformly with respect to $t\in \mathfrak M,\;\;\bigl(n\in\mathbb N\bigr)$
if it converges for all $t\in [a,b]\;\;\bigl(n\in\mathbb N\bigr)$ and
$$
\bigl(\forall\,\varepsilon>0\bigr)\;\bigl(\forall\,n_0\in\mathbb N\bigr)\;
\bigl(\exists\,p\in\mathbb N\bigr)\;\bigl(\exists\,n_1,\,\ldots,\,n_p\,(\geqslant n_0)\bigr)\;\bigl(\forall\,t\in\mathfrak M\bigr)\;
\bigl(\forall\,n\in\mathbb N\bigr)
$$
$$
\underset{i=1,\,\ldots,\,p}{min}\;\left|\sum\limits_{k=n_i+1}^\infty a_{k}(t)\right|<\varepsilon\quad
\left(\underset{i=1,\,\ldots,\,p}{min}\;\left|\sum\limits_{k=n_i+1}^\infty a_{kn}\right|<\varepsilon\right).
$$
\begin{teo}
\label{Arzela}
{\it Let $\mathfrak M$~ be a metric compact, $x_n(\cdot)\in {\bf{C}}(\mathfrak M),\linebreak
n=1,2,\ldots,$ and
$$
x_n(t)\to x(t)\;\;\text{for all}\;\;t\in \mathfrak M.
$$
Then for the continuity of $x(\cdot)$ on $\mathfrak M$ it is necessary and sufficient that 
$ x_n(t)\to x(t)\;(QU)\;\;\text{on}\;\; \mathfrak M.$ }
\end{teo}
\doc\,
{\it{Necessity.}}\, Let $x(\cdot)$ be continuous on $\mathfrak M.$
For arbitrary  $n_0,\,\varepsilon>0$ and arbitrary $t\in\mathfrak M$, there exists $n(t)\geqslant n_0,$ such that
\begin{equation}
\label{Arz1}
|x_{n(t)}(t)-x(t)|<\varepsilon
\end{equation}
Denote
$
G(t)\doteq \{s:\,|x_{n(t)}(s)-x(s)|<\varepsilon\}.
$
Since $x(\cdot)$ is continuous, $G(t)$~ is an open set and $t\in G(t).$ Since $\mathfrak M$ is compact, it is covered by a finite
 number of sets $G(t_1),\,\ldots,\,G(t_p).$ The numbers $n_i\doteq n(t_i),\;i=1,\,\ldots\,p$ satisfy the definition of 
the quasi-uniform convergence. Hence $ x_n(t)\to x(t)$ $(QU)$ on $ \mathfrak M.$

{\it{Sufficiency.}}\,
Let $ x_n(t)\to x(t)$ $(QU)$ on $\mathfrak M.$
Given arbitrary $t\in\mathfrak M$ and $\varepsilon>0$, there exists $n_0$ such that inequality (\ref{Arz1}) holds for $n\geqslant n_0$.
We choose $n_i\geqslant n_0\;\;(i=1,\,\ldots,\,p)$ from the definition of quasi-uniform convergence and denote
$G_i(t)\doteq \{s:\,|x_{n_i}(s)-x(t)|<\varepsilon\},$\linebreak
$i=1,\,\ldots,\,p.$
Since $x_n(\cdot)$ is continuous, the sets $G_i(t)$ are open. The set
$G(t)=\bigcap\limits_{i=1}^p G_i(t)\;(\ni t)$ is also open. Let $s\in G(t)$ be arbitrary and let number $i$ be such that
$|x_{n_i}(s)-x(s)|<\varepsilon.$ Then
$$
|x(s)-x(t)|\leqslant |x(s)-x_{n_i}(s)|+|x_{n_i}(s)-x_{n_i}(t)|+|x_{n_i }(t)-x(t)|<3\,\varepsilon.
$$
Hence $x(\cdot)$ is continuous at point $t.$
\hfill $\square$


\begin{center}
\begin{large}
\section{Linear normed spaces}
\end{large}
\end{center}

\poonkt
{Definition and examples}\phantom{01234567890123456789} \

Let $\mathbb R$ denote, as usual, the field of real numbers, 
let $\mathcal X$~ be a linear (vector) space over the field $\mathbb R.$
A function $\|\cdot\|:\mathcal X\to\mathbb R_+\equiv [0,+\infty)$ is called a {\it norm},
if it has the properties:

1) \;$\|x\|=0\Longleftrightarrow x=0$ (we denote by $0$ both the number zero and the zero element of the space $\mathcal X$);

2) \;$\|\lambda x\|=|\lambda|\cdot \|x\|\;(\lambda\in\mathbb R,\;x\in \mathcal X);$

3)\; $\|x+y\|\leqslant \|x\|+\|y\|\;(x,y\in \mathcal X).$

Properties 1) -- 3) are called the axioms of norm. Namely: 1) is the axiom of identity, 2) is the axiom of positive homogeneity,
3) is the triangle inequality. A linear space endowed with a norm is called a {\it linear normed space} (LNS).
If we need to emphasize that we are talking about the norm in space $\mathcal X,$ then we write $\|x\|_{\mathcal X}.$

$1^{\circ}.$ {\it For any $x,y\in \mathcal X$,\;$\bigl|\|x\|-\|y\|\bigr |\leqslant \|x-y\|.$}

\doc\,
We apply the triangle inequality and use the positive\linebreak
homogeneity of the norm:
$$
\|x\|=\|(x-y)+y\|\leqslant \|x-y\|+\|y\| \Rightarrow \|x\|-\|y\|\leqslant \|x-y\|,
$$
$$
\|y\|=\|(y-x)+x\|\leqslant \|y-x\|+\|x\| \Rightarrow \|y\|-\|x\|\leqslant \|y-x\|.
$$
These two inequalities combined give us the required.
\hfill $\square$

\medskip

$2^{\circ}.$
{\it Every LNS is a metric space}

\doc\,
Let $\rho (x,y)\doteq \|x-y\|.$ The identity axiom for the metric follows from the identity axiom for the norm,
the axiom of symmetry follows from the positive homogeneity of the norm. The triangle inequality for the metric is obtained
from the triangle inequality for the norm as follows:
$$
\rho (x,y)=\|x-y\|=\|(x-z)+(z-y)\|\leqslant \|x-z\|+\|z-y\|=\rho (x,z)+\rho ( z, y).
$$
\hfill $\square$

The latter implies that $\|x\|=\rho(x,0)$ (the norm of an element is its distance to zero). This means that all concepts 
and statements that were introduced and proved for metric spaces extend to LNS.

\medskip

$3^{\circ}.$
{\it The norm is a continuous function on $\mathcal X.$}

\doc\,
Let $x_n\to x \;(n\to\infty).$ This means that $\|x_n-x\|\to 0.$ From the statement {\bf{1}} we see that $\bigl | \|x_n\|-\|x\|\bigr |\leqslant \|x_n-x\|\to 0\\
(n\to\infty);$ hence $\|x_n\|\to \|x\|,$ i.e.\,norm is a continuous function.
\hfill $\square$

\medskip

A complete LNS is called a {\it Banach space} ($B$-space) in honor of the Polish mathematician S. Banach, who introduced this concept.

As some examples of Banach spaces we mention the following spaces that were introduced earlier:
$\mathbb R_p^n,\,{\bf{l}}_p,\,{\bf{c}}_0,\,{\bf{c}},\,{\bf{m}},\ ,{\bf{C}}[a,b],\,{\bf{M}}[a,b].$ 
The norms in all these LNSs are defined by the equality $\|x\|=\rho(x,0).$

\medskip

$4^{\circ}.$
{\it Linear operations in LNS are continuous.}

\doc\,
Let $x,y,x_n,y_n\;(n\in\mathbb N)$ be LNS elements, let $\alpha, \beta$ be real numbers, and let $x_n\to x,\;y_n\to y\;(n\to\infty).$
Then 
\begin{align*}
\|(\alpha x_n+\beta y_n)-(\|\alpha x+\beta y)\| & =\|\alpha (x_n-x)+\beta (y_n-y)\| \\
& \leqslant |\alpha |\|x_n-x\|+|\beta |\|y_n-y\|\to 0. 
\end{align*}
This
means that $\alpha x_n+\beta y_n\to \alpha x+\beta y\;n\to\infty$. Therefore, algebraic operations are continuous.
\hfill $\square$

Let us consider another example.
Denote by ${\bf{C}}^{(n)}[a,b]$ the set of $n$ times ($n\in\mathbb N$) continuously differentiable functions
endowed with the norm
$$
\|x\|=\sum\limits_{k=0}^n\underset{t\in [a,b]}{\max}\bigl|x^{(k)}(t)\bigr|= \sum\limits_{k=0}^n\|x^{(k)}\|_{{\bf{C}}[a,b]}.
$$

$5^{\circ}.$
{\it ${\bf{C}}^{(n)}[a,b]$~ is a Banach space.}

\doc\,
Since one can write $\|x\|_{{\bf{C}}^{(n)}}=\sum\limits_{k=0}^n \|x^{(k)}\|_ {{\bf{C}}},$ the axiom of positive homogeneity holds and
the triangle inequality is obvious. Let $\|x\|_{{\bf{C}}^{(n)}}=0$. Then $\|x^{(k)}\|_{{\bf{C}}}=0$ for $k=0,1,\dots ,n,$ which is possible only if $x$  is identically equal to zero, i.e.\,we have $\|x\|_{{\bf{C}}^{(n)}}\!=\!0\Longrightarrow \!x\!=\!0$ holds. The reverse implication
is obvious. Thus, we have proved that ${\bf{C}}^{(n)}[a,b]$~ is an LNS.

Let $x_m\to x$ in ${\bf{C}}^{(n)}[a,b],$ i.e. $\|x_m-x\|_{{\bf{C}}^{(n)}}\to 0.$ This means that $\|x_m^{(k)}-x^{(k) }\|_{{\bf{C}}}\to 0$ or
$x_m^{(k)}(t)\rightrightarrows x^{(k)}(t)\;(k=0,1,\dots ,n)$. Thus, the convergence of a sequence of functions in the space ${\bf{C}}^{(n)}[a,b]$ is the
uniform convergence of all derivatives up to (including) order $n$.

We will prove the completeness of this space by the induction in $n$. For $n=0\;\; {\bf{C}}^{(n)}[a,b]={\bf{C}}[a,b]$ is a complete space. Suppose that
the space ${\bf{C}}^{(n-1)}[a,b]$ is complete. Let $\{x_m\}_{m=1}^{\infty}\subset {\bf{C}}^{(n)}[a,b]$ be an arbitrary
fundamental sequence. This means that $\|x_m-x_p\|_{{\bf{C}}^{(n)}}\to 0$ for $m,p\to \infty$, so
\begin{equation}\label{C^n}
\sum\limits_{k=0}^{n-1}\|x_m^{(k)}-x_p^{(k)}\|_{{\bf{C}}}\to 0,\quad \|x_m^{(n)}-x_p^{(n)}\|\to 0,\quad(m,p\to \infty).
\end{equation}
The first convergence in (\ref{C^n}) means that the original sequence is fundamental in space ${\bf{C}}^{(n-1)}[a,b]$, which, by our hypothesis, is complete. Therefore, there exists $x\in {\bf{C}}^{(n-1)}[a,b]$ such that $x_m\to x$ in ${\bf{C}}^{( n-1)}[a,b].$
The second limit in (\ref{C^n}) means that there exists a continuous function $y$ such that
$x_m^{(n)}(t)\rightrightarrows y(t)\;(m\to\infty).$
Now we recall the representation of a continuously differentiable function in terms of its derivative:
\begin{equation}\label{NewLeib}
x_m^{(n-1)}(t)=x_m^{(n-1)}(a)+\int\limits_a^t x_m^{(n)}(s)\,ds.
\end{equation}
In view of the uniform convergence of $x_m^{(n)}(t)\rightrightarrows y(t)$, we can pass to the limit
$m\to\infty$ in this representation, obtaining
$$
x^{(n-1)}(t)=x^{(n-1)}(a)+\int\limits_a^t y(s)\,ds.
$$
It follows that $x^{(n)}=y\in {\bf{C}}[a,b],$ or
$x\in {\bf{C}}^{(n)}[a,b]$ and, since \linebreak
$x_m^{(k)}(t)\rightrightarrows x^{(k)}\,
(k =0,1,\dots ,n,\;m\to\infty),$ we have $x_m\to x$ in ${\bf{C}}^{(n)}[a,b].$
By induction, the space ${\bf{C}}^{(n)}[a,b]$ is complete for every $n\in\mathbb N.$
\hfill $\square$



\bigskip

\poonkt
{Finite-dimensional LNS} \ \phantom{01234567890123456789}

Let $\mathcal X$~ be a linear space equipped with two functions $\|\cdot\|_1$ and
$\|\cdot\|_2$ satisfying the axioms of the norm.
We say that the norms $\|\cdot\|_1$ and $\|\cdot\|_2$ are {\it equivalent}
(written as $\|\cdot\|_1 \thicksim\|\cdot\|_2$) if
there are positive numbers $\alpha$ and $\beta$ such that the inequalities
\begin{equation}\label{equi1}
\alpha\|x\|_2\leqslant \|x\|_1\leqslant \beta\|x\|_2
\end{equation}
hold for all $x \in \mathcal X$.

$1^{\circ}.$
{\it The relation $\thicksim$ on the set of all norms that are defined on the same linear space is an equivalence relation, that is, it is reflexive, symmetric and transitive.}

\doc\,
Reflexivity: for any norm $\|x\|\leqslant
\|x\|\leqslant \|x\|$. Symmetry: if $\alpha \|x\|_2\leqslant \|x\|_1\leqslant \beta\|x\|_2,$ where
$\alpha>0,\beta >0,$ then
$$
\frac{1}{\beta}\|x\|_1\leqslant \|x\|_2\leqslant \frac{1}{\alpha}\|x\|_1.
$$
Transitivity: if
$$
\alpha \|x\|_2\leqslant \|x\|_1\leqslant \beta\|x\|_2\quad\mbox{and}\quad \gamma \|x\|_3\leqslant \|x\| _2\leqslant \delta \|x\|_3
$$
$(\alpha ,\beta ,\gamma ,\delta >0)$, then $\alpha \gamma \|x\|_3\leqslant \|x\|_1\leqslant \delta\beta\|x\|_3 .$
\hfill $\square$

$2^{\circ}.$
{\it If a sequence of elements of a linear space
$\mathcal X$ converges with respect to one of the two equivalent norms,
then it converges with respect to both norms.}

\doc\,
Let 
$
\alpha \|x\|_2\leqslant \|x\|_1\leqslant \beta\|x\|_2\;(\alpha>0,\beta >0)
$
 and $\|x_n-x\|_2\to 0$ $(n\to\infty)$. It follows from inequality $\|x\|_1\leqslant \beta\|x\|_2$ that $\|x_n-x\|_1\to 0$. If $\|x_n-x\|_1\to 0,$ then inequality $\alpha \|x\|_2\leqslant \|x\|_1$ implies that
 $\|x_n-x\|_2\to 0.$
\hfill $\square$

$3^{\circ}.$
{\it If a sequence of elements of a linear space
$\mathcal X$ is fundamental with respect to one of the two equivalent norms, then it is fundamental with respect to both norms.}

\doc\,
Let $\alpha \|x\|_2\leqslant \|x\|_1\leqslant \beta\|x\|_2\;(\alpha>0,\beta >0)$. Assume that sequence $\{x_n\} _{n=1}^{\infty}$
is fundamental with respect to norm $\|\cdot\|_2$, that is,
$\|x_n-x_m\|_2\to 0\;(n,m\to\infty).$ Then it follows from inequality $\|x\|_1\leqslant \beta\|x\|_2$ 
that
$\|x_n-x_m\|_1\to 0\; (n,m\to\infty),$ i.e.\,the sequence is fundamental with respect to norm $\|\cdot\|_1$. If
this sequence is fundamental with respect to norm $\|\cdot\|_1,$, then we prove in the same way, 
using inequality $\alpha \|x\|_2\leqslant \|x\|_1\leqslant$, that it is fundamental with respect to norm $\|\cdot\|_2.$
\hfill $\square$

$4^{\circ}.$
{\it If space $\mathcal X$ is complete with respect to one of the two equivalent norms, then it is
complete with respect to both norms.}

\doc\,
Let $\alpha \|x\|_2\leqslant \|x\|_1\leqslant \beta\|x\|_2\;(\alpha>0,\beta >0).$ Suppose that space $\mathcal X$ is complete with respect to norm $\|\cdot\|_1$, and suppose that sequence $\{x_n\}_{n=1}^{\infty}$ is
fundamental with respect to norm $\|\cdot\|_2$, i.e.\,$\|x_n-x_m\|_2\!\to \!0
(n,m\!\to~\infty).$ According to {\bf{3}},
it is fundamental with respect to norm $\|\cdot\|_1$. Therefore, there exists $x\in \mathcal X$ such that $\|x_n-x\|_1\to 0\;(n\to\infty)$. In view of inequality $\alpha \|x\|_2\leqslant \|x\|_1$, we have $\|x_n-x\|_2\to 0\;
(n\to\infty),$ which means that space $\mathcal X$ is complete with respect to norm
 $\|\cdot\|_2.$ We argue in exactly the same way if space $\mathcal X$ is complete with respect to
norm $\|\cdot\|_2.$
\hfill $\square$

Let $\mathcal X$~ be a finite-dimensional linear space ($\dim\,\mathcal X=n$). Let $\{g_k\}_{k=1}^n$ be
its basis, let $\{e_k\}_{k=1}^n$~ be a  fixed basis in $\mathbb R_2^n.$
Each element of $x\in\mathcal X$ can be expanded uniquely in basis $\{g_k\}_{k=1}^n\!\!:\;$
$x=\sum\limits_{k=1}^nx_kg_k,\; x_k\in\mathbb R$.
Let us associate with this $x$ the element $\hat x\in\mathbb R_2^n$ according to the rule $\hat x=\sum\limits_{k=1}^nx_ke_k$.
This correspondence is an isomorphism.

\begin{lem}
\label{lm1}
{\it There are positive numbers $\alpha$ and $\beta$ such that the following inequality is true:}
\begin{equation}\label{equi}
\alpha\|\hat x\|_{\mathbb P_2^n}\leqslant \|x\|_{\mathcal X}\leqslant \beta\|\hat x\|_{\mathbb P_2^n}.
\end{equation}
\end{lem}
\doc\,
Using the triangle inequality for the norm, and the H\"{o}lder inequality ($p=2$), we obtain
$$
\|x\|_{\mathcal X}=\|\sum\limits_{k=1}^nx_kg_k\|_{\mathcal X}\leqslant \sum\limits_{k=1}^n|x_k|\ |g_k\|_{\mathcal X}\leqslant
\beta\left(\sum\limits_{k=1}^n|x_k|^2\right)^{1/2}=\beta\|\hat x\|_{ \mathbb R^n_2 },
$$
where $\beta=\left(\sum\limits_{k=1}^n\|g_k\|^2\right)^{1/2}$. This proves the right inequality in (\ref{equi}).
The assertion {\bf{1}} and what has already been proved imply the inequalities
\begin{equation}\label{contnorm}
\bigl|\|x'\|_{\mathcal X}-\|x''\|_{\mathcal X}\bigr|\leqslant \|x'-x''\|_{\mathcal X} \leqslant \beta\|\hat x'-\hat x''\|_{\mathbb R_2^n}.
\end{equation}
Consider function $f(x)\doteq \|x\|_{\mathcal X}$ defined on set
$$
S=\{x:\; \|\hat x\|_{\mathbb P_2^n}=1\}\subset\mathcal X;
$$
The set $S$ can be viewed as the unit sphere in space $\mathbb P_2^n.$
It follows from (\ref{contnorm}) that this function is continuous. Since $S$ is a bounded closed set in $\mathbb R_2^n,$ there is
a point on the sphere where
$f$ is equal to its infimum. 
We denote this value by $\alpha$. Since $S$ does not contain zero, we have $\alpha>0$.
Thus, for any $x\in\mathcal X$,
\begin{equation*}
\|\frac{x}{\|\hat x\|_{\mathbb P_2^n}}\|_{\mathcal X}\geqslant\alpha,
\end{equation*}
which yields the left inequality in (\ref{equi}).
\hfill $\square$

\medskip

$5^\circ.$
{\it All norms in $\mathbb R^n$ are equivalent.}

\doc\,
Let us take as $\mathcal X$ in (\ref{equi}) the space $\mathbb R^n$ endowed with any norm.
\hfill $\square$

$6^\circ.$
{\it Let $\mathcal X$~ be a finite-dimensional linear space.
All norms in $\mathcal X$ are equivalent.}

\doc\,
Let $\|\cdot\|_1,\,\|\cdot\|_2$~ be two arbitrary norms in $\mathcal X.$
According to inequality (\ref{equi}), we have
$$
\alpha\|\hat x\|_{\mathbb R_2^n}\leqslant \|x\|_1\leqslant \beta\|\hat x\|_{\mathbb R_2^n},\quad
\gamma\|\hat x\|_{\mathbb R_2^n}\leqslant \|x\|_2\leqslant \delta\|\hat x\|_{\mathbb R_2^n},
$$
which gives us $\frac{\alpha}{\delta}\|x\|_2\leqslant \|x\|_1\leqslant \frac{\beta}{\gamma}\|\|_2.$
\hfill $\square$
\begin{teo}
\label{th1}
{\it All finite-dimensional LNS are complete.}
\end{teo}
\doc\,
Let $\mathcal X$~be a  LNS, $n=\dim\, \mathcal X,$ let $\{g_k\}_{k=1}^n$~ be a  basis in $\mathcal X,
 \{x^{(n)}\}_{n=1}^{\infty}\subset\mathcal X$~ is an arbitrary fundamental sequence. 
Let us construct sequence $\{\hat x^{(n)}\}_{n=1}^{\infty}\subset \mathbb R_2^n$ as above.
It follows from the left inequality in (\ref{equi})  that
$$
\alpha \|\hat x^{(n)}-\hat x^{(m)}\|_{\mathbb R_2^n}\leqslant \|x^{(n)}-x^{(m )}\|_{\mathcal X}.
$$
This means that sequence $\{\hat x^{(n)}\}_{n=1}^{\infty}$ is fundamental in the LNS $\mathbb R_2^n$ and,  since $\mathbb R_2^n$ is complete, has limit in it. Thus, there exists $x\!\in\!\mathcal X$ 
such that $\|\hat x^{(n)}\!-\!\hat x\|_{\mathbb R_2^n}\!\to~0\; (n\to \infty)$.
The right inequality in (\ref{equi}) implies that $\|x^{(n)}-x\|_{\mathcal X}\leqslant\beta\|\hat x^{(n)}-\hat x\|_{\mathbb R_2^n},$
i.e. $\|x^{(n)}-x\|_{\mathcal X}\to 0\quad (n\to\infty).$
This means that $\mathcal X$ is complete.
\hfill $\square$

\medskip

A {\it closed} linear manifold in $\mathcal X$ is called a {\it subspace} of $\mathcal X.$ We draw reader's attention
to the difference between the concept of a subspace in the theory of metric spaces and in the theory of LNS. In the theory of metric spaces, the concept
a subspace is not an embedded closure.

Let us consider an example. Let $P[a,b]$ denote the linear space of algebraic polynomials with the same norm as in ${\bf{C}}[a.b].$
Obviously, $P[a,b]$ is a linear manifold in ${\bf{C}}[a.b]$, but $P[a,b]$ is not closed in ${\bf{C}}[ a.b]$. Indeed,
$x_n(t)=\sum\limits_{k=0}^{n}\frac{t^k}{k!}\rightrightarrows x(t)=e^t$, i.e. $x_n $ converges to $x$ in ${\bf{C}}[a.b];$ \;
$x_n\in P[a,b]\linebreak
(n\in\mathbb N)$, but $x\notin P[a,b].$ Therefore,
 $P[a,b]$ is not a subspace of ${\bf{C}}[a,b].$

\medskip

$7^\circ.$
{\it All finite-dimensional linear manifolds in an infinite-dimensional LNS $\mathcal X$ are closed, and thus they are subspaces.}

\doc\,
Let $\mathcal Y\subset\mathcal X$ be a finite-dimensional linear manifold. By Theorem 1, it is closed. Therefore, by Theorem \ref{mprth1},
it is complete.
\hfill $\square$

\medskip

$8^\circ.$
{\it The set $P_n[a,b]$ of polynomials of degree at most $n$~ is a subspace of the LNS} ${\bf{C}}[a,b].$

\doc\,
Each polynomial 
$p(t)=a_0t^n+a_1t^{n-1}+\dots +a_{n-1}t+a_n$ is defined by an $(n+1)$-dimensional vector
$(a_0,a_1,\ldots ,a_n)$ of its coefficients. Therefore, $P_n[a,b]$~ is a $(n+1)$-dimensional linear manifold in ${\bf{C}}[a,b].$
By assertion {\bf{7}}, $P_n[a,b]$~ is a subspace of ${\bf{C}}[a,b].$
\hfill $\square$



\bigskip

\poonkt{Direct product. Direct sum. Isometric isomorphism} \ 

The {\it direct product} $\mathcal X\times\mathcal Y$ of LNSs $\mathcal X$ and $\mathcal Y$ is a linear space
$$
\mathcal Z=\{z=(x,y):\;x\in\mathcal X,\;y\in\mathcal Y\}\quad\mbox{with norm}\quad \|z\|=\ |x\|_{\mathcal X}+\|y\|_{\mathcal Y}
$$
Let $\mathcal X$ and $\mathcal Y$ be two LNSs such that $\mathcal X\bigcap\mathcal Y=\{0\}$ (these can be two subspaces of the same LNS).
The {\it direct sum} $\mathcal W\doteq \mathcal X\dot+\mathcal Y$ is the linear space
$$
\mathcal W=\{w=x+y:x\in\mathcal X,\;y\in\mathcal Y\}\quad\mbox{with norm}\quad \|w\|=\|x\| _{\mathcal X}+\|y\|_{\mathcal Y}.
$$

$1^\circ.$
{\it $\mathcal X\times\mathcal Y\;\;\bigl(\mathcal X\dot+\mathcal Y\bigr)$ is a Banach space if and only if
spaces $\mathcal X$ and $\mathcal Y$ are Banach.}

\doc\,
Let $\mathcal Z=\mathcal X\times\mathcal Y$ be a Banach space. Let 
$\{x_n\}_{n=1}^{\infty}$ and $\{y_n\}_{n=1}^{\infty}$ be fundamental sequences in $\mathcal X$ and $\mathcal Y$, respectively.
Let $z_n\!=\!(x_n,y_n)$.
Then
$$
\|z_n-z_m\|_{\mathcal Z}\!=\!\|x_n-x_m\|_{\mathcal X}\!+\!\|y_n-y_m\|_{\mathcal Y} \to \!0\,(n,m\to\infty).
$$
Therefore, sequence $\{z_n\}_{n=1}^{\infty}$ is fundamental in $\mathcal Z.$ Since this space is complete,
there exists $z=(x,y)\in \mathcal Z,\;z_n\to z$. Since 
$$
0\gets \|z_n-z\|_{\mathcal Z}=\|x_n-x\|_ {\mathcal X}+\|y_n-y\|_{\mathcal Y},
$$
we have $\|x_n-x\|_{\mathcal X}\to~0,
\|y_n-y\|_{\mathcal Y}\to 0$. This means that spaces $\mathcal X$ and $\mathcal Y$ are complete.

The reverse implication is proved by arguing similarly: completeness of $\mathcal X$ and $\mathcal Y$ implies completeness of $\mathcal Z.$
Similar arguments work for the direct sum.
\hfill $\square$

The LNSs $\mathcal X$ and $\mathcal Y$ are called {\it isometrically isomorphic} (we write $\mathcal X\cong\mathcal Y$)
if there exists a mapping $\varphi: \mathcal X \to \mathcal Y$ (an isometric isomorphism) with the properties:

1)\;$\varphi$~is an isomorphism $\bigl($recall: \, this means that $\varphi$ is bijective and
\begin{equation*}
\varphi(x+y)=\varphi(x)+\varphi(y),\;\varphi(\lambda x)=\lambda\varphi(x)
\end{equation*}
for all $x,y\in \mathcal X$ and $\lambda\in \mathbb P\bigr);$

2)\;$\varphi$~ is an isometry; this means that $\|\varphi(x)\|_{\mathcal Y}=\|x\|_{\mathcal X}.$

Obviously, the direct sum is isometrically isomorphic to the direct product.

\medskip

$2^\circ.$
{\it Isometric isomorphism~ is a continuous mapping.}

\doc\,
Let $\varphi: \mathcal X\to \mathcal Y\;(\mathcal X,\,\mathcal Y$~ are LNSs) be an isometric isomorphism and let $x_n\to x$ in $\mathcal X. $
Since $\varphi (x_n)-\varphi (x)=\varphi(x_n-x),$ we have
$\|\varphi (x_n)-\varphi (x)\|_{\mathcal Y}=\\=
\|\varphi (x_n-x)\|_{\mathcal Y}= \|x_n-x\|_{ \mathcal X}\to 0\;(n\to\infty),$ i.e. $\varphi (x_n)\to \varphi (x)$
in $\mathcal Y.$ This means that $\varphi$ is continuous
\hfill $\square$

$3^\circ.$
{\it The inverse of an isometric isomorphism is also an isometric isomorphism.}

\doc\,
Let $\varphi: \mathcal X\to \mathcal Y\;(\mathcal X,\,\mathcal Y$~ are LNSs) be an isometric isomorphism. Since $\varphi$ is bijective, 
there exists its inverse, which is a bijective mapping
$\varphi ^{-1}: \mathcal Y\to \mathcal X.$

Let us show that $\varphi ^{-1}$~ is an isomorphism. Let $y_i\in \mathcal Y,\,x_i\doteq \varphi ^{-1}(y_i);$ 
which means $y_i=\varphi (x_i)\,(i=1,2).$ Then
\begin{multline*}
\alpha y_1+\beta y_2=\alpha \varphi(x_1)+\beta \varphi (x_2)=\varphi (\alpha x_1+\beta x_2)\Rightarrow \\
\Rightarrow\varphi ^{-1}(\alpha y_1+\beta y_2)=
\alpha x_1+\beta x_2=\alpha \varphi ^{-1}(y_1)+\beta \varphi ^{-1}(y_2).
\end{multline*}
It remains to show that $\varphi ^{-1}$~ is an isometry. Let $y\in \mathcal Y,$
$ x\doteq \varphi ^{-1}(y)$; this means that $y=\varphi (x).$ We have
$
\|\varphi ^{-1}(y)\|_{\mathcal X}=\|x\|_{\mathcal X}=\\
=\|\varphi (x)\|_{\mathcal Y}=\|y\|_{\mathcal Y}.
$
\hfill $\square$

\medskip

$4^\circ.$
{\it If $\mathcal X\cong\mathcal Y,$ then $\mathcal X$ is complete if and only if the $\mathcal Y$ is complete.}

\doc\,
Let $\mathcal X\cong\mathcal Y$, and let $\varphi: \mathcal X\to \mathcal Y$~ be the corresponding isometric isomorphism. 
Since
$\|\varphi(x)-\varphi(y)\|_{\mathcal Y}=\|x-y\|_{\mathcal X},$ i.e.\,the isometric isomorphism preserves distances between elements,
a fundamental sequence in space $\mathcal X\;\bigl(\mathcal Y\bigr)$ corresponds
 to a fundamental
sequence in space $\mathcal Y\;\bigl(\mathcal X\bigr)$; and a convergent sequence in space
$\mathcal X\;\bigl(\mathcal Y\bigr)$ corresponds to a convergent
sequence in space $\mathcal Y\;\bigl(\mathcal X\bigr).$ This yields the required.
\hfill $\square$

Isometric isomorphism  between LNSs is obviously an {\it equivalence relation}. Such spaces may differ
from each other only by the nature of their elements. However, in the LNS theory, isometrically isomorphic spaces are treated as identical.

Let $\mathcal X$~ be an LNS with norm $\|\cdot \|.$ The smallest Banach space $\widehat{\mathcal X}$ that contains $\mathcal X$ as a subspace is called a {\it completion} of the LNS $\mathcal X$. For example, $\mathbb R$~ is a completion of
$\mathbb Q.$

The following completion theorem is meaningful only for incomplete space.

\begin{teo}\label{nhprth2A}
{\it Every LNS has a unique completion up to an isometric isomorphism.
Moreover, the original LNS is dense in its completion: $cl\,\mathcal X=\widehat{\mathcal X}$.}
\end{teo}
\doc\,
To prove this theorem, we only need to add to the proof of Theorem \ref{mprth2} the verification that the classes 
of equivalent fundamental sequences constructed there constitute a
linear space, which is fairly obvious.
\hfill $\square$

As in the case of metric spaces, a completion can often be described by starting directly from the definition (see, 
for example, the remark at the end of paragraph {\bf{1.4}}).
\vskip 2pt

Let $\mathcal X$ be an LNS, let $\mathcal Y$~ be its subspace. The {\it quotient space}
$\mathcal X/\mathcal Y=\{\hat x:\;x,y\in\hat x\Leftrightarrow x-y\in\mathcal Y\}$ is, as is well known, a linear space
(its elements $\hat x$ are called cosets with respect to the subspace $\mathcal Y$, the zero element is $\mathcal Y.$ itself)
Let $\hat x\in\mathcal X/\mathcal Y$. Put $\|\hat x\|_{\mathcal X/\mathcal Y}=\underset{x\in\hat x}{\inf} \|x\|_{\mathcal X}.$

$5^\circ.$
{\it If $\mathcal X$ is a Banach space, then $\mathcal X/\mathcal Y$ is also a Banach space.}

\doc\,
Let $\hat x$~ be the zero element of the quotient space $\mathcal X/\mathcal Y,$ i.e.
$\hat x=\mathcal Y; \|\hat x\|_{\mathcal X/\mathcal Y}\!=\!\underset{x\in \mathcal Y}{\inf}\|x\|_{\mathcal X}=~0 .$
Conversely, if $\|\hat x\|_{\mathcal X/\mathcal Y}=0,$ then there exists a sequence $\{x_n\}\subset\hat x$ such that
$x_n\to~0.$ Since, due to the closedness of $\mathcal Y$, the coset $\hat x=x+\mathcal Y$ is also closed, it contains
the limit of this sequence, that is, $0\in\hat x.$ This means that $\hat x=0.$
The identity axiom is proved.

Since $\lambda \hat x=\{y:y=\lambda x,\,x\in \hat x\},$ we have
$$
\|\lambda \hat x\|_{\mathcal X/\mathcal Y}= \underset{x\in \lambda \hat x}{\inf}\|x\|_{\mathcal X}= \underset {x\in \hat x}{\inf}\|\lambda x\|_{\mathcal X}=
|\lambda| \underset{x\in \hat x}{\inf}\|x\|_{\mathcal X}=|\lambda|\cdot\|\hat x\|_{\mathcal X/\mathcal Y}.
$$
It remains to verify the validity of the triangle inequality. Since $\hat x+\hat y=\widehat{x+y},$ where $x\in\hat x,\,y\in\hat y,$ we have
from the definition of the norm:
$$
\|\hat x+\hat y\|_{\mathcal X/\mathcal Y}=\|\widehat{x+y}\|_{\mathcal X/\mathcal Y}\leqslant \|x+y\ |_{\mathcal X}\leqslant \|x\|_{\mathcal X}+\|y\|_{\mathcal X}.
$$
We now pass on the right side to the infimums. We arrive at the required inequality.

Thus, the axioms of the norm are satisfied.

Let us prove completeness of the quotient space.

Let $\psi :\mathcal X\to \mathcal X/\mathcal Y$ be a natural homomorphism of the space $\mathcal X$
to the factor space $\mathcal X/\mathcal Y,$ i.e., $\psi (x)=\hat x=x+\mathcal Y.$ For any $\hat x\in\mathcal X/\mathcal Y$,
there exists $x\in\mathcal X$ such that $\psi (x)=\hat x,$ and, by the property of the infimum,
\begin{equation}\label{fpr}
\|\hat x\|_{\mathcal X/\mathcal Y}\geqslant \frac{1}{2}\|x\|_{\mathcal X}.
\end{equation}

Let $\mathcal X$~ be a complete LNS, 
$\{\hat x_n\}_{n=1}^{\infty}\subset \mathcal X/\mathcal Y$~ is a fundamental sequence. This means that
$\|\hat x_n-\hat x_m\|_{\mathcal X/\mathcal Y}\to~0\;(m,n\to~\infty).$ Since $\hat x_n-\hat x_m =\widehat{x_n-x_m},$
there exist
$$
x_n\in \hat x_n,\,x_m\in\hat x_m,\,x_n-x_m\in\widehat{x_n-x_m}
$$
such that, due to (\ref{fpr}),
$
0\gets\|\hat x_n-\hat x_m\|_{\mathcal X/\mathcal Y}\geqslant \frac{1}{2}\|x_n-x_m\|_{\mathcal X}.
$
Thus, sequence $\{x_n\}_{n=1}^{\infty}$ is fundamental in complete space $\mathcal X$. Therefore, there exists
$x\!\in \!\mathcal X$ such that $\|x_n-x\|_{\mathcal X}\to~0$  $(n\!\to~\infty)$. We put $\hat x=\psi (x).$
From the definition of the norm in $\mathcal X/\mathcal Y$, we have
$$
\|\hat x_n-\hat x\|_{\mathcal X/\mathcal Y}=\underset{z\in \hat x_n-\hat x}{\inf}\|z\|_{\mathcal X }\leqslant \|x_n-x\|_{\mathcal X}\to~0,
$$
that is, $\;\hat x_n\to\hat x\in \mathcal X/\mathcal Y.$ The quotient space is complete.
\hfill $\square$



\bigskip

\poonkt
{Splitting of B-spaces. Compactness of constituent
subspaces} \ \phantom{01234567890123456789}

Let $X$~ be a B-space, $X_1,X_2$ are its subspaces 
$\bigl(X_1\bigcap X_2=\{0\}\bigr)$, $X=X_1\dotplus X_2, A\subset X$ and
$\;A_i\subset X_i$~ are  ``projections'' of set $A$ onto $X_i\;(i=1,2),\;$ i.e.
\begin{multline*}
A_1=\{x_1\in X_1:\; (\forall\,x\in A)\;\;(\exists\,x_2\in X_2),\;\;x=x_1+x_2\},\\
A_2=\{x_2\in X_2:\; (\forall\,x\in A)\;\;(\exists\,x_1\in X_1),\;\;x=x_1+x_2\}
\end{multline*}

Let us ask ourselves the following question: what is the relationship between properties of closedness
(compactness, relative compactness) of sets $A$ and $A_1,A_2$?
\begin{lem}
\label{locked}
{\it If $A$ is closed, then $A_i\;\;(i=1,2)$ are also closed.}
\end{lem}
\doc\,
Let $x_i^{*}$~ be limit points of sets $A_i$ and let $\{x^{(n)}_i\}_{n=1}^\infty\subset A_i$ be sequences converging to
these limit points $(i=1,2).$ We put $x^{(n)}(t)\doteq \linebreak\doteq 
 x^{(n)}_1(t)+x^{(n)}_2(t) .$ Then $x^{(n)}(t)\to x^{*}_1(t)+x^{*}_2\doteq x^{*}(t).$
Since $A$ is closed, we have $x^{*}\in A,$ is closed, and by the definition of sets $A_i\;\;x^{*}_i\in A_i,$ i.e. $A_i\;\; (i=1,2)$ are closed.
\hfill $\square$

Already the example
\begin{multline*}
X=\mathbb R^2,\;\; A=\{(x_1,x_2):x_1^2+x_2^2\leqslant 1\}
\setminus \{(x_1,x_2):\;x_2=\sqrt{1-x_1^2},\\-1<x_1<0\}
\end{multline*}
shows that the converse is not true: \; the sets $A_1,\,A_2$ are closed, but the set $A$ is not.
\begin{lem}
\label{otnkomp}
{\it If $A_i\subset X_i\;\;(i=1,2)$ are relatively compact, then $A$~ is relatively compact.}
\end{lem}
\doc\,
Let sets $A_i\;(i=1,2)$ be relatively compact, let $\varepsilon$ be arbitrary, let
$E_i\doteq \{y^{k}_i\}_{k=1}^{m_i}\;\;(i=1,2)$~ be the corresponding $\dfrac{\varepsilon}{2 }$-nets\;
for $A_i$ in $X_i\;\;(i=1,2).$ We set $y^{jk}\doteq y^{j}_1+y^{k}_2\;\; (j=1,\ldots,m_1,\linebreak
k=1,\ldots,m_2).\;$ Then $H\doteq \{y^{jk}\}_{j=1,\ldots,m_1,\;k=1,\ldots,m_2}\subset X$~ is an $\varepsilon$-net for $A$ to $X.$

Indeed, let $x\in A.$ Put $x=x_1+x_2\;\;(x_i\in A_i).$ There are
$y^{k_{i}}_i\in E_i$ such that $\|x_i-y^{k_{i}}_i\|_{X_i}<\dfrac{\varepsilon}{2}.$ Since $y^{k_1,k_2}\in H,$ we have
\begin{multline*}
\|x-y^{k_1,k_2}\|_X=\|x_1+x_2-y^{k_1}_1-y^{k_2}_2\|_X=\|(x_1-y^{k_1}_1)+( x_2-y^{k_2}_2)\|_{X}\leqslant \\
\leqslant \|x_1-y^{k_1}_1\|_X+\|x_2-y^{k_2}_2\|_X=\|x_1-y^{k_1}_1\|_{X_1}+\|x_2- y^{k_2}_2\|_{X_2}<\varepsilon.
\end{multline*}
By Theorem \ref{hausd} (Hausdorff), $A$ is relatively compact in $X.$
\hfill $\square$

It can be seen directly from the proofs that the assertions of Lemmas \ref{locked} and \ref{otnkomp} extend to an arbitrary number $m$ of the
constituent subspaces. Let $X$~ be a $B$-space, $X_1,\ldots,X_m$~ are its subspaces such that
$X_i\bigcap X_i=\{0\}$ for $i\ne j$, and
$$X=X_1\dotplus\ldots\dotplus X_m, \quad
A\subset X,$$
and $A_i\subset X_i$ are the ``projections'' $A$ onto the subspace $X_i\;(i=
1,\ldots,m),\;$ i.e.
\begin{multline*}
A_i=\{x_i\in X_i:\; (\forall\,x\in A)\;\;(\exists\,x_j\in X_j),\;j\ne i,\;\;x=\sum\limits_{j=1}^m x_j \},\\
i=1,\ldots,m .
\end{multline*}
Thus, we have the following final result:
\begin{teo}
\label{kompitog}
{\it If $A$~ is closed and $A_i\;\;(i=1,\ldots,m)$~ are relatively compact, then the sets $A,\,A_1,\ldots, A_m\;$ are compact.}
\end{teo}



\bigskip

\poonkt
{Examples of B-spaces of sequences}\phantom{01234567890123456789} \

Consider the following linear normed spaces of sequences of numbers (see \cite[VIII.2.8-11]{dansch}):

${\bf{bv}}\doteq \{x=(x_1,x_2,\ldots):\,\|x\|\doteq |x_1|+\sum\limits_{k=1}^{\infty} |x_{k+1}-x_k|<+\infty\};$

${\bf{bv_0}}\doteq \{x=(x_1,x_2,\ldots):\,\|x\|\doteq \sum\limits_{k=1}^{\infty}|x_{k +1}-x_k|<+\infty,\;\underset{k\to \infty}{\lim}x_k=0\};$

${\bf{bs}}\doteq \{x=(x_1,x_2,\ldots):\,\|x\|\doteq \underset{n}{\sup}\left |\sum\limits_{k =1}^{n}x_k\right |<+\infty\};$

${\bf{cs}}\doteq \{x=(x_1,x_2,\ldots):\, series\; \sum\limits_{k=1}^{\infty}x_k \;converges, 
\|x\|\doteq \underset{n}{\sup}\left |\sum\limits_{k=1}^{n}x_k\right |\};$

\medskip

Below we prove a number of results for these spaces.

\medskip

$1^\circ.$
{\it There following $($set-theoretic$)$ inclusions take place:
$$
{\bf{l_1}}\subset {\bf{bv_0}}\subset {\bf{bv}}\subset {\bf{c}},\qquad {\bf{l_1}}\subset {\bf{ cs}}\subset {\bf{c_0}},\qquad
{\bf{cs}}\subset {\bf{bs}}\subset {\bf{l}}_{\infty};
$$}

\doc\,
Let us first show that for any LNS $\mathcal X$ and any sequence $\{x_n\}_{n=1}^{\infty}\subset\mathcal X$ of its elements,
the convergence of the series $\sum\limits_{n=1}^{\infty}\|x_{n+1}-x_n\|$ implies that this sequence is fundamental.

Indeed, by Cauchy's convergence principle, it follows from the convergence of the series that 
for any $\varepsilon>0$ there exists a natural number
$N$ such that for all $n,m>N,\;m>n$\;\;$$\sum\limits_{k=n}^{m-1}\|x_k-x_{k+1} \|<\varepsilon.$$ Using the triangle inequality, we obtain
$$\|x_n-x_m\|\leqslant \sum\limits_{k=n}^{m-1}\|x_k-x_{k+1}\|<\varepsilon,\;$$ which means that the sequence is fundamental
$\{x_n\}_{n=1}^{\infty}.$

Note that the reverse implication does not take place: the fact that the sequence is fundamental does not, in general, 
imply the convergence of the indicated series.
Let us consider an example. The sequence 
$$x^{(n)}=\bigl(1,\frac{1}{2},\ldots,\frac{1}{n},0,\ldots\bigr)
$$ 
converges in space ${\bf{l_2}}$ to 
$$
x=(1,\frac{1}{2},\ldots,\frac{1}{k},\frac{1}{k+1}\ldots),\;
$$ 
since
$$
\|x^{(n)}-x\|=\left(\sum\limits_{k=n+1}^{\infty}\frac{1}{k^2}\right)^{\frac{1}{2}}\to 0
$$
(as the remainder of a convergent series). So, the sequence
$\{x^{(n)}\}_{n=1}^{\infty}$ is fundamental. However, the series \; $\sum\limits_{n=1}^{\infty}\|x^{(n)}-x^{(n+1)}\|=
\sum\limits_{n=1}^{\infty}\frac{1}{n+1}$ diverges.
 
The inclusion ${\bf{bv}}\subset {\bf{c}}$ now follows from the above implication and the completeness 
of space $\mathbb R$. The remaining inclusions are obvious.
\hfill $\square$

\medskip

$2^\circ.$
{\it There is an isometric isomorphism between spaces}
$$
{\bf{bv_0}}\;\text{and}\;{\bf{l_1}}; \quad {\bf{cs}}\;\text{and}\;{\bf{c}}; \quad {\bf{bs}}\;\text{and}\;{\bf{l}}_{\infty};
$$

\doc\,
We define mapping $$\varphi: {\bf{bv_0}}\to {\bf{l_1}}$$ as follows: for $x=(x_1,x_2,\ldots)\in {\bf{bv_0}}$, we put
$$y\doteq \varphi (x)=(x_2-x_1,x_3-x_2,\ldots ).$$
It follows from the definition of space ${\bf{bv_0}}$ that $\varphi (x)\in {\bf{l_1}}$. Since the inverse mapping is given by
$$x=\varphi ^{-1}(y),\quad x_k=
-\sum\limits_{j=k}^{\infty}y_j,\;k=1,2,\dots,$$
we have $x_k\to 0\;(k\to\infty)$ (as the remainder of the convergent series). So, mapping $\varphi$ is bijective. The linearity
$\varphi$ is obvious. Hence $\varphi$~ is an isomorphism.
Since $$\|\varphi (x)\|_{{\bf{l_1}}}=\sum\limits_{k=1}^{\infty}|x_{k+1}-x_k|=\| x\|_{{\bf{bv_0}}},$$ we obtain that $\varphi$ is an isometry.

We define mapping $$\varphi: {\bf{cs}}\to {\bf{c}}$$ as follows: for $x=(x_1,x_2,\ldots)\in {\bf{cs} }$ 
we set
$y\doteq \varphi (x)$, $y_n=\sum\limits_{k=1}^{n}x_k,\;n\in\mathbb N$.
Due to the convergence of the series
$\sum\limits_{k=1}^{\infty}x_k$ there exists the limit $S=\underset{n\to\infty}{\lim}y_n,$ i.e. $\varphi (x)\in {\bf{c}}$.
The inverse mapping is defined by the equality
$$x=\varphi ^{-1}(y),\;x_1=y_1, \quad x_k=y_k-y_{k-1},\; k=2,3,\ldots (y\in c).$$ The $n$-th partial sum of the series
$\sum\limits_{k=1}^{\infty}x_k$ is equal to $y_n$. Thus,
$\varphi$ is bijective and,
due to its obvious linearity, is an isomorphism. Since
$\|\varphi (x)\|_{{\bf{c}}}=\underset{n\in\mathbb N}{\sup}|y_n|=\|x\|_{{\bf{cs}}},$ we obtain that $\varphi$  is an isometry.

The mapping $$\varphi: {\bf{bs}}\to {\bf{l}}_{\infty}$$ is defined as follows: for $x=(x_1,x_2,\ldots)\in
 {\bf{bs}}$, we put
$y\doteq \varphi (x)$, $y_n=\sum\limits_{k=1}^{n}x_k,$ $n\in\mathbb N$. The inverse
mapping is defined as in the previous case. Thus, as above, we obtain that $\varphi$~ is an isomorphism. Since
$$
\|x\|_{{\bf{bs}}}=\underset{n\in \mathbb N}{\sup}\left |\sum\limits_{k=1}^{n}x_k\right | =\underset{n\in\mathbb N}{\sup}|y_n|=
\|\varphi (x)\|_{{\bf{l}}_\infty},
$$
mapping $\varphi$ is an isometry;
\hfill $\square$

\medskip

$3^\circ.$
${\bf{bv}}=\mathbb R\dot +{\bf{bv_0}};$

\doc\,
Let $x=(x_1,x_2,\ldots)\in {\bf{bv_0}}$. From the proof of ${\bf{1}}$ and the completeness
of space $\mathbb R$, we obtain the existence of the limit
$\alpha\doteq\underset{k\to\infty}{\lim}x_k,$
so one can represent $x=\alpha e+y,$ where
$$
e=(1,1,\ldots)\in {\bf{bv}},\; y=(x_1-\alpha,\,x_2-\alpha,\ldots)\in {\bf{bv_0}}).
$$
The uniqueness of this representation follows from the uniqueness of the limit. The required representation follows.
\hfill $\square$

\medskip

$4^\circ.$
{\it ${\bf{bv}},\,{\bf{bv_0}},\,{\bf{bs}},\,{\bf{cs}}$~ are Banach spaces.}

\doc\,
The completeness of spaces ${\bf{bv_0}},\,{\bf{bs}},\,{\bf{cs}}$ follows from statements {\bf{2}} and {\bf{2.3.4}}, the completeness of ${\bf{bv}}$ follows from statement {\bf{3}}.
\hfill $\square$


\newpage

\begin{center}
\begin{large}
{\bf{Chapter II. \;  Functions of bounded variation}}
\end{large}
\end{center}
\addcontentsline{toc}{section}{Chapter II. \; Functions of bounded variation}

\begin{center}
\begin{large}
\section{Monotone functions}
\end{large}
\end{center}


\poonkt{Definitions}\phantom{01234567890123456789} \

Here we are discuss functions ~$f: J\to\R$ defined at each point of an 
interval
$J\subset\R$ of the real axis, and taking real
values.

Let us recall some definitions from mathematical analysis.

A function~$f(\cdot)$ is called {\it increasing} ({\it decreasing})
on interval~$J$ if, for every $t,s\in J,$ $t<s$,
\begin{equation}
\label{ch11}
f(t) \ls f(s) \qq (f(t)\gs f(s)).
\end{equation}
Thus, we will understand the increase (decrease) 
in the non-strict sense. If in (\ref{ch11})
one has strict inequalities for all $t$ and~$s$ specified there, then~$f(\cdot)$ will be called
{\it strictly increasing} ({\it strictly decreasing}).

A function~$f(\cdot)$ is called {\it monotonic} ({\it strictly monotonic})
on an interval~$J$ if it is increasing on this interval
or decreasing (strictly increasing or strictly decreasing).

Let us also recall some properties of the monotone functions. A monotonie function
has finite one-sided
limits
$$
f(t_0+)\doteq \lim\lt_{t\to t_0+0} f(t), \qq
f(t_0-)\doteq \lim\lt_{t\to t_0-0} f(t),
$$
at each interior point of the interval $J$. 
That is, a monotone function can have only first kind discontinuities.

Denote by $\sg_{t_0}(f)$ ($\sg_{t_0}^+(f),\,\sg_{t_0}^-(f)$) the jump
(right jump, left jump) of function~$f(\cdot)$ at point~$t_0:$
$$
\gathered
\sg_{t_0}(f)\doteq f(t_0+)-f(t_0-), \q \sg_{t_0}^+(f)=f(t_0+)-f(t_0), \\
\sg_{t_0}^-(f)\doteq f(t_0)-f(t_0-).
\endgathered
$$
Clearly,
$$
\sg_{t_0}(f)=\sg_{t_0}^+(f)+\sg_{t_0}^-(f).
$$
If~$\sg_{t_0}^+(f)=0$ ($\sg_{t_0}^-(f)=0$), then~$f(\cdot)$ is continuous
on the right (left) at point~$t_0.$ Since~$f(\cdot)$ is monotone, 
either 
$f(t_0-) \ls f(t_0) \ls f(t_0+)$ or $f(t_0-) \gs f(t_0) \gs f(t_0+),$
so if~$\sg_{t_0}(f)=0,$ then~$f(\cdot)$ is continuous at~$t_0.$
If $J=[a,\,b],$ then at point~$a$~($b$) one can only talk about
about $\sg_{a}^+(f)$ ($\sg_{b}^-(f)$) and right (left) continuity.

If~$f(\cdot)$~ is a decreasing function, then $g(t)=-f(t)$~ is an increasing function.
So, below we consider only increasing functions.


\bigskip

\poonkt{Auxiliary assertions} \
\begin{lem}
\label{ch1lem1}
 Let $J=[a,\,b],$ let $f(\cdot)$ be an increasing function, $t_1,\,
\ldots,\,t_n$~ are distinct interior points of $J.$ Then
\begin{equation}
\label{ch12}
\sg_{a}^+(f)+\sum\lt_{k=1}^n \sg_{t_k}(f)+\sg_{b}^-(f) \ls
f(b)-f(a).
\end{equation}
\end{lem}
\doc\,
Without loss of generality, we can assume that~$t_k$ are ordered:
$t_1<t_2<\ldots<t_n.$ Put $t_0\doteq a,$ $t_{n+1}\doteq b$ and select
arbitrarily~$s_k$ so that $t_k<s_k<t_{k+1},$ $k=0,\,1,\,\ldots,\,n.$

Since~$f(\cdot)$ is increasing, we have inequalities
$$
\aligned
{} & f(a+)-f(a) \ls f(s_0)-f(a); \\
{} & f(t_k+)-f(t_k-) \ls f(s_k)-f(s_{k-1}), \q k=1,\,2,\,\ldots,\,n. \\
{} & f(b)-f(b-) \ls f(b)-f(s_n).
\endaligned
$$
Adding these inequalities, we obtain estimate~(\ref{ch12}).
\hfill $\square$

\begin{lem}
\label{ch1lem2}
{\it Let $J=[a,\,b],$ let $f(\cdot)$ be an increasing function. For every $\eps>0$,
there exist only finitely many points of discontinuity of function~$f(\cdot)$ where its 
jump exceeds~$\eps.$} (In other words, for every $\eps>0$ the set
$T_\eps(f)\doteq\{ t\in J:\sg_{t}(f)>\eps \}$~ is finite.)
\end{lem}
\doc\,
Fix arbitrarily~$n$ points
$t_1,\,\ldots,\,t_n$ from~$T_\eps(f).$
Then,
according to Lemma~\ref{ch1lem1},
$$
n\eps< \sum\lt_{k=1}^n \sg_{t_k}(f) \ls f(b)-f(a).
$$
Therefore,
$
n\ls \frac{f(b)-f(a)}{\eps}.
$
And this is what was needed to be proved.
\hfill $\square$

\begin{teo}
\label{ch1thm1}
{\it Let $J=[a,\,b].$ The set~$T(f)$ of points of discontinuity of an increasing
function~$f(\cdot)$ is at most countable. If $\{t_k\}_{k=1}^\infty\subset$\linebreak
$\subset T(f)$ is the set of interior discontinuity points, then}
\begin{equation}
\label{ch13}
\sg_{a}^+(f)+\sum\lt_{k=1}^\fy \sg_{t_k}(f)+\sg_{b}^-(f) \ls f(b)-f (a).
\end{equation}
\end{teo}
\doc\,
The first assertion follows from Lemma~\ref{ch1lem2} and the representation
$\;
T(f)=\bigcup_{n=1}^\fy T_{\fracs 1n}(f).
$\;
Lemma~\ref{ch1lem1} yields that the sequence of partial  sums of the series 
(having positive terms) from~(\ref{ch13}) is bounded. Therefore, taking $n$ in (\ref{ch12})
to~$+\fy,$ we obtain~(\ref{ch13}).
\hfill $\square$

If $J=(a,\,b)$ (in particular, $a=-\fy$ or (and) $b=+\fy$)
and $T(f)=\{t_k\}_{k=1}^\fy,$ then inequality~(\ref{ch13}) takes form
\begin{equation}
\label{ch14}
\sum\lt_{k=1}^\fy \sg_{t_k}(f) \ls \sup_{t\in J} f(t)-\inf_{t\in J} f(t).
\end{equation}

\bigskip

\poonkt{Example of an increasing function 
discontinuous on $\mathbb Q$} 

\phantom{0123456789012345678901234567890123456789}

The fact that the set of points of discontinuity of a monotone function is at most
countable can give an impression that such functions are rather unsophisticated. However, a countable set can have a structure that is  quite complex.
In this regard, let us consider the following example.
\begin{example}
\label{Qraz}
Let $\Q=\{r_1,\,r_2,\,\ldots\}$ be the set of rational numbers. For each~ $t\in\R$, we put
\begin{equation}
\label{ch15}
f(t)=\sum\lt_{k:\,r_k<t} \frac{1}{2^k}
\end{equation}
(here the summation extends to those~ $k$ for which 
${r_k<t}$). Obviously $\inf_{\R} f(t)=0,$ $\sup_{\R} f(t)=\sum\lt_{k=1}^\fy
\fracs{1}{2^k}=1.$
\end{example}

Let us show that

1)~$f(\cdot)$ is strictly increasing on all~$\R;$

2)~$f(\cdot)$ has discontinuity on the right at every rational point 
$r_n\in\Q,$ moreover~$\sg_{r_n}^+(f)=\fracs{1}{2^n};$

3)~$f(\cdot)$ is left continuous on~$\R;$

4)~$f(\cdot)$ is continuous at irrational points.

Indeed, let $t_2>t_1.$ There is a rational number
$t_1<r_{k_0}<t_2$ such that
$$
\gathered
f(t_2)=\!\!\!\sum\lt_{r_k< t_2} \!\!\frac{1}{2^k}=\!\!\!\sum\lt_{r_k< t_1}
\!\!\frac{1}{2^k}+\!\!\!\!\!\sum\lt_{t_1\ls r_k< t_2}
\!\!\!\!\frac{1}{2^k} >
\!\!\!\sum\lt_{r_k<t_1} \!\!\frac{1}{2^k}+
\frac{1}{2^{k_0}}=f(t_1)+\frac{1}{2^{k_0}},
\endgathered
$$
\begin{equation}
\label{ch16}
f(t_2)-f(t_1)>\frac{1}{2^{k_0}}>0.
\end{equation}
Thus, assertion~1) is proved.

Let~$r_n \in \Q$ be arbitrary. For~$t>r_n$ as above, we have
$f(t)>f(r_n)+\linebreak+
\fracs{1}{2^n}.$ Hence for $t\to r_n+0$ we have
$\sg^{+}_{r_n}(f)=
f(r_n+)-f(r_n) \gs \fracs{1}{2^n}.\;$ This proves that~$f(\cdot)$ is discontinuous from the right at every rational point.

Since the sum of jumps over all discontinuity points~$f(\cdot)$ is larger or
equals to $\sum\lt_{k=1}^\fy \sg_{r_k}(f)$ and $\sum\lt_{k=1}^\fy \sg_{r_k}(f)
\gs \sum\lt_{k=1}^\fy \sg_{r_k}^+(f),$ we have by inequality~(\ref{ch14})
$$
1=\sum\lt_{k=1}^\fy \frac{1}{2^k} \ls \sum\lt_{k=1}^\fy \sg_{r_k}^+(f) \ls
\sum\lt_{k=1}^\fy \sg_{r_k}(f) \ls 1.
$$
Therefore, all series above sum to exactly 1,
$\sg_{r_k}^+(f)= \sg_{r_k}(f)=\fracs{1}{2^k},$ $\sg_{r_k}^-(f)=0$, and
$\sg_t(f)=0$ for every irrational point. So, all statements above for $f(\cdot)$ defined by
equality~(\ref{ch15}) are proved.
\begin{zam}
\label{zampr1}
1. One can carry out a similar to (\ref{ch15}) construction of a strictly increasing function, discontinuous
at all points of a prespecified
countable set ${\{t_1,\,t_2,\,\ldots\}\subset J},$ 
{\it left continuous} at points $t_k,$ having
jumps at these points equal to $\alpha_k>0,$ where
$\sum\lt_{k=1}^\infty \alpha_k<\infty$.

2. If we define
\begin{equation}
\label{ch15+}
g(t)=\sum\lt_{k:\,t_k\leqslant t}\alpha_k,
\end{equation}
then we obtain a function $g(\cdot)$ that is discontinuous at all points of the set ${\{t_1,\,t_2,\,\ldots\}},$
{\it is right continuous} at points $t_k$ with jumps $\alpha_k.$
\end{zam}


\bigskip

\poonkt{Continuous and discrete parts} \phantom{01234567890123456789} \

Let~$J=[a,\,b],$ let $f(\cdot)$~ be an increasing function, $T(f)=\{t_k\}$ is
its set of discontinuities. We define function~$f_d(\cdot)$ in the following way. Let~$f_d(a)=0.$ Put for~$t>a$
\begin{equation}
\label{ch17}
f_d(t)\doteq \sg_a^+(f)+\sum\lt_{t_k<t} \sg_{t_k}(f)+\sg_t^-(f).
\end{equation}

We immediately see from definition~(\ref{ch17})  that~$f_d(\cdot)$ jumps
when we pass through the points of discontinuity of $f(\cdot).$

Let us show that $f_d(\cdot)$ is an increasing function.

Indeed, let~$s>t.$ Then
\begin{multline}
\label{ch18}
f_d(s)-f_d(t)=\sum\lt_{t\ls t_k<s} \sg_{t_k}(f)+\sg_s^-(f)-\sg_t^-(f)=\\
=\sum\lt_{t<t_k<s} \sg_{t_k}(f)+f(t+)-f(t-)+\sg_s^-(f)-f(t)+f(t-) =\\
=\sg_t^+(f)+\sum\lt_{t<t_k<s} \sg_{t_k}(f)+\sg_s^-(f)\gs 0,
\end{multline}
which means that $f_d(\cdot)$ is increasing.

Next, put
\begin{equation}
\label{ch19}
f_c(t)\doteq f(t)-f_d(t).
\end{equation}

Let us show that $f_c(\cdot)$~ is an increasing continuous function.

Let~$t\in J$ be arbitrary, let $s>t.$ We obtain from~(\ref{ch18}) and the inequality~(\ref{ch13}),
applied on closed interval $[t,\,s]$:
\begin{equation}
\label{ch110}
f_d(s)-f_d(t) \ls f(s)-f(t),
\end{equation}
so $f_c(t)=f(t)-f_d(t) \ls f(s)-f_d(s)=f_c(s),$ i.e. $f_c(\cdot)$
is indeed an increasing function.

On the other hand, we also see from~(\ref{ch18}) that
$$
\sg_t^+(f)=f(t+)-f(t)\ls f_d(s)-f_d(t)
$$
(all terms are non-negative, so one term does not exceed the sum).
Letting $s\to t+0$ here and in (\ref{ch110}), we obtain
\begin{equation}
\label{ch111}
f(t+)-f(t)=f_d(t+)-f_d(t),
\end{equation}
so $f_c(t+)=f(t+)-f_d(t+)=f(t)-f_d(t)=f_c(t),$ i.e. $f_c(\cdot)$
is continuous on the right at point~$t$. Taking~$s<t,$ interchainging in the argument above $s$ and~$t$, and then taking $s\to t-0,$ we establish
continuity of $f_c(\cdot)$ at $t$ on the left.

Along the way, we proved $\bigl($see, for example, (\ref{ch111})$\bigr)$ that the jumps
(left jumps, right jumps) of $f_d(\cdot)$ and~$f(\cdot)$ coincide.

Definition~(\ref{ch19}) can be rewritten as follows: 
$
f(t)=f_c(t)+f_d(t).
$
Thus, we have proved the following assertion.

\begin{teo}
\label{ch1thm2}
{\it An increasing function can be represented as the sum of an 
increasing continuous function~$f_c(\cdot)$ {\rm(}continuous part~$f(\cdot)${\rm)}
and an 
increasing jump function~$f_d(\cdot)$ {\rm(}discrete
parts~$f(\cdot)${\rm)}.}
The value~$f(\cdot)$ at  point~$a$ is ``inherited'' from its continuous
part: $f_c(a)=f(a).$
\end{teo}

Note that in the example considered in Section 1.2 $f_d(t)=f(t),\\
f_c(t)=~0.$


\bigskip

\bigskip

\poonkt{Cantor's ladder}\phantom{01234567890123456789} \

Here it will be convenient to define a special function (the so-called Cantor's ``ladder''),
which we will need later in the book.
First, we prove an auxiliary assertion.

{\it Let $f(\cdot)$ be an increasing function that takes all values from  set $A$ that is
everywhere dense in~$[f(a),\,f(b)].$ Then~$f(\cdot)$ is continuous on~$[a,\,b].$}

\doc\,
Assume that~$f(\cdot)$ has a discontinuity point $t_0$ in $[a,\,b]$
Since~$f(\cdot)$ is an increasing function, $t_0$ can only be
a jump discontinuity. 
This means that at least one
of intervals $(f(t_0-),\,f(t_0))$ or $(f(t_0),\,f(t_0+))$ does not contain
values of function~$f(\cdot).$ Hence this interval contains no points of set $A,$
which contradicts to the hypothesis that $A$ is dense in $[f(a),\,f(b)].$ Therefore, $f(\cdot)$ is continuous on~$[a,\,b].$
\hfill $\square$
\begin{example}
\label{Kant1}
{We define function $\theta:[0,\,1]\to[0,\,1]$ as follows (see Section 1.7.7).
We set $\theta(t)=\fracs 12$ for~$t\in G_1,$ $\theta(t)=\fracs 14$ for
$t\in(\fracs 19,\,\fracs 29),$ $\theta(t)=\fracs 34$ for $t\in(\fracs 79,\,
\fracs 89),$ and, in general, $\theta(t)=\frac{2i-1}{2^k}$ on the $i$-th  adjacent interval of $k$-th rank (counting
from left to right), $i=1,\,2,\,\ldots,\,2^{k-1}.$
Function $\theta(\cdot)$ is defined on the Cantor open set~$G_0.$
Let us define it on the Cantor perfect set~${\cal K}.$ We put
$\theta(0)\!=\!0,$ ${\theta(t)=\!\sup\lt_{\tau<t,\,\tau\in G_0}\! \theta(\tau)}$ for
the other points in ${\cal K}.$ Clearly, $\theta(1)\!=\!1.$ Now~$\theta(\cdot)$
is defined everywhere on $[0,\,1].$}
\end{example}
Let us show that:

a$)$~$\theta(\cdot)$~ {\it is an increasing function;

b$)$~$\theta(\cdot)$ is continuous on~$[0,\,1].$ }

\doc\,
By construction,~$\theta(\cdot)$ increases (in a non-strict sense) on the set $G_0.$ The extension of $\theta(\cdot)$ to the entire interval $[0,1]$
preserves its property of being increasing, since the supremum can only become larger as the set over which the supremum is taken expands.

Already on ~$G_0$, function $\theta(\cdot)$ takes all binary rational values. Since the set of binary rational numbers is dense
in $[0,1],$ the previous result implies continuity of ~$\theta(\cdot).$
\hfill $\square$


\bigskip

\poonkt{Sets of measure zero}\phantom{01234567890123456789} \

Everywhere below, for the sake of the uniformity of notation of various sums of jumps, we extend function $f:[a,b],$
originally defined on $[a,b],$ to the left of $t=a$ and to the right of $t=b$ as follows:
\begin{equation}
\label{soglprodolz}
f(t)=f(a)\;\;\;\text{on}\;t<a,\quad f(t)=f(b)\;\;\;\text{on}\; t>b.
\end{equation}
This gives us $\sigma_a^-(f)=0$ and $\sigma_b^+(f)=0.$

Let $g:[a,b]\to\mathbb R$~ {\it be an arbitrary increasing} function.
We say that $A\subset [a,b]$~ {\it is a $g$-measure zero set} if for any $\varepsilon>0$ there is an open set
$G=\bigcup\limits_k(a_k,\,b_k)\supset A$ such that $$\mu_g(G)\doteq\sum\limits_k\bigl(g(b_k-)-g(a_k+)\bigr)< \varepsilon.$$ The fact that
$A\subset [a,b]$~ is a set of $g$-measure zero is also denoted by $\mu_g(A)=0.$
If $g(t)=t,$ then we say that $A$ {\it is a set of Lebesgue measure zero}\; ({\it null-set}) and write 
$$\mes(G)\doteq\sum\limits_k(b_k-a_k),\quad \mes(A)=0.$$
Thus, for an open set
$
G=\bigcup\limits_k(a_k,\,b_k)$, $\mes(G)
$ 
is the sum of the lengths of the constituent intervals (we assume from now on that these intervals do not intersect each other, although they may have common endpoints).

\begin{teo}
\label{konschet1}
{\it Let an increasing function $g:[a,b]\to\mathbb R$ be continuous. Then all finite and countable sets from $[a,b]$
are sets of $g$-measure zero.}
\end{teo}
\doc\,
Let, for example, $A=\{t_1,t_2,\ldots\}\subset [a,b]$, and let $\varepsilon>0$ be fixed arbitrarily.
Due to the uniform continuity of $g(\cdot)$, for every $t_k$ there exists an interval $(a_k,b_k)\ni t_k$
such that $g(b_k-)-g(a_k+)=
g(b_k)-g(a_k)<\frac{\varepsilon}{2^k}.$ Put $G=\bigcup\limits_{k=1}^\infty(a_k,\,b_k)$, which yields the required.
\hfill $\square$

{\it The Cantor perfect set $\pK$ $($see p. {\bf{1.7.7}}$)$ is a set of Lebesgue measure zero, $mes\,\pK=0.$}

\doc\,
For notation, see {\bf{1.7.7}}.
Let $\varepsilon>0$ be arbitrary, let $\Gamma\doteq\bigcup\limits_{k=1}^\infty \gamma_k$
be the set of points of the first kind,
enumerated as in {\bf{1.7.7}}. According to statement c) from {\bf{1.7.7}}, the open set
$$
G\doteq \bigcup\limits_{k=1}^\infty\left(\gamma_k-\frac{\varepsilon}{2^k},\,\gamma_k+\frac{\varepsilon}{2^k}\right)
$$ 
covers
$\pK,$ and $\mes(G)<\varepsilon,$ i.e., $\mes\,\pK=0.$
\hfill $\square$

\begin{teo}
\label{konschet2}
{\it Let $g(\cdot)$~ be an arbitrary increasing function. Finite or countable union of sets
$g$-measure zero is a set of $g$-measure zero.}
\end{teo}

\doc\,
Let $A_n,\;n=1,2,\ldots$~ be sets of $g$-measure zero, let $\varepsilon>0$ be arbitrary. We cover each $A_n$ with the open set
$G_n=\bigcup\limits_{k=1}^\infty \bigl(a_{nk},\,b_{nk}\bigr)$ so that
$g\bigl(b_{nk}-\bigr)-g\bigl(a_{nk}+\bigr)<\dfrac{\varepsilon}{2^{n+k}},\;k,n=1 ,2,\ldots$
Then $A=\bigcup\limits_{n=1}^\infty A_n$
is covered by the open set $G=\bigcup\limits_{n=1}^\infty G_n,\;\mu_g(G)<\varepsilon.$~\hfill $\square$


\newpage

\section{Functions of bounded variation}
\label{ch2}


\poonkt{Total variation  of a function }\phantom{01234567890123456789} \

Let~$J=[a,\,b]$, $f:J\to \R.$ Consider a partition $\tau=\{t_k\}_{k=0}^n$
of the interval $[a,\,b],$ $a=t_0<t_1<\ldots<t_n=b.$ Let us consider the sum of the absolute values of the increments
of function~$f(\cdot)$ corresponding to this partition:
$$
v_\tau\doteq v_\tau(f)\doteq v_\tau(f;\,[a,\,b])\doteq \sum\lt_{k=1}^n
|f(t_k)-f(t_{k-1})| \ (\gs 0).
$$
There are two possible cases: 

1)~The set of sums~$\{v_\tau(f)\}_\tau,$ corresponding to all
partitions, is bounded from above, that is, there is a constant~$C$
such that $v_\tau(f)\ls C$ for all~$\tau;$ in this case we say
that function~$f(\cdot)$ has {\it finite total variation}
(finite change) on $[a,\,b]$, and define {\it the total variation of function~$f(\cdot)$
on interval~$[a,\,b]$} as
$$
\bigvee\lt_a^b(f)\doteq \sup_\tau v_\tau(f)
$$
(the supremum here is taken over all partitions~$\tau$). We also say that
$f(\cdot)$ {\it is a function of bounded variation.}

2)~The set of sums $\{v_\tau(f)\}_\tau$ is not bounded,
that is, for any~$C>0$ there is a partition~$\tau$ such that $v_\tau(f)>C$.
Then we say that~$f(\cdot)$ is not a  function of  bounded variation
or that its total variation~$f(\cdot)$ is infinite, $\bigvee\lt_a^b(f)=+\fy.$

Obviously, the total variation of a function is zero if and only if this function is identically equal to a constant.

Suppose that partitions~$\tau'$ and~$\tau$ differ only by one point, so that
$\tau'\!=\!\tau\cup\{\ol t\}$ and~$\ol t\in(t_{k_0-1},\,t_{k_0}).$ Then
by the properties of the absolute value of a real number,
$$
\gathered
|f(t_{k_0})-f(t_{k_0-1})|=|(f(t_{k_0})-f(\ol t))+(f(\ol t)-f(t_{ k_0-1}))| \ls \\
\ls |f(t_{k_0})-f(\ol t)|+|f(\ol t)-f(t_{k_0-1})|.
\endgathered
$$
The number on the left is included in the sum~$v_\tau,$ both terms on the right are included
in sum $v_{\tau'}.$ The rest of the terms in both sums are the same, so
$v_\tau \ls v_{\tau'}.$ Thus, adding new points to some partition
can only increase the sum of the absolute values of the increments. Since $v_\tau(f)\ls \bigvee\lt_a^b(f)$
for any partition, we may assume that a given partition includes all points that are of interest to us.


\bigskip

\poonkt{Example of a continuous function with infinite total variation}\phantom{01234567890123456789}

\begin{example}
\label{CneBV1}
{\it Let~$\mathfrak{c}(t)=0$ for~$t=0,$ $\mathfrak{c}(t)=t\cdot \cos\frac{\pi}{2t}$ for~$0<t\ls 1.$}
\end{example}
Consider the partition~$\tau_n$ of interval $[0,\,1]$ by points
$$
\gathered
0<\frac{1}{2n}<\frac{1}{2n-1}<\ldots<\frac{1}{3}<\frac{1}{2}<1. \\
v_{\tau_n}(\!f\!)\!=\!\left| \frac{1}{2n}\cos\frac{\pi}{2\cdot\frac{1}{2n}}\!-\!0
\right| \!+\!\left| \frac{1}{2n-1}\cos\frac{\pi}{2\cdot\frac{1}{2n-1}}\!-\!
\frac{1}{2n}\cos\frac{\pi}{2\cdot\frac{1}{2n}} \right|\!+\! \ldots \!+\! \\
+\left| \frac{1}{2}\cos\frac{\pi}{2\cdot\frac{1}{2}}-\frac{1}{3}
\cos\frac{\pi}{2\cdot\frac{1}{3}} \right|+ \left| 1\cdot\cos
\frac{\pi}{2\cdot 1}-\frac{1}{2}\cos\frac{\pi}{2\cdot\frac{1}{2}} \right|.
\endgathered
$$
Since $\cos\pi k=(-1)^k,$ $\cos\frac{(2k-1)\pi}{2}=0,$ we have
$$
v_{\tau_n}=\frac{1}{2n}+\frac{1}{2n}+\frac{1}{2n-2}+\frac{1}{2n-2}+
\ldots+\frac{1}{4}+\frac{1}{4}+\frac{1}{2}+\frac{1}{2}=1+\frac{1}{2}+
\ldots+\frac{1}{n}.
$$
Since the partial sum of the harmonic series tends to~$+\fy$ as ${n\to\fy},$
for every $c>0$ there is~$n$ such that $v_{\tau_n}(\mathfrak{c})>c.$
Thus, $\bigvee\lt_0^1(\mathfrak{c})=+\fy.$

\bigskip

\poonkt{Vector space of functions of bounded variation\;} \ \phantom{012345678900123456789012345678901234567890123456789}

Let us introduce some classes of functions of bounded variation and study their basic properties.
\begin{teo}
\label{ch4thm1}
{\it Monotone functions have finite total variation, and} ${\bigvee\lt_a^b(f)=|f(b)-f(a)|}.$
\end{teo}
\doc\,
Let $f(\cdot)$ be an increasing function on $[a,\,b]$, let $\tau$ be  an arbitrary
partition. Since $f(t_k)\gs f(t_{k-1}),$ we have
$$
v_\tau(f)=\sum\lt_{k=1}^n |f(t_k)-f(t_{k-1})|=\sum\lt_{k=1}^n
(f(t_k)-f(t_{k-1}))=f(b)-f(a).
$$
Here $v_\tau(f)$ does not even depend on the partition, so
$\bigvee\lt_a^b(f)=$\linebreak
$=f(b)-f(a).$
\hfill $\square$

\begin{teo}
\label{ch2thm2}
{\it If~$f(\cdot)$ satisfies the Lipschitz condition, then it has finite total
variation, and $\bigvee\lt_a^b(f)\ls L(b-a),$ where~$L$~ is a Lipschitz constant
 of function~$f(\cdot).$}
\end{teo}
\doc\,
The Lipschitz condition means that for all $t$ and $s$ from $[a,\,b]$
$|f(t)-f(s)|\ls L|t-s|.$ For an arbitrary partition~$\tau$ we have
\begin{multline*}
v_\tau(f)=\sum\lt_{k=1}^n |f(t_k)-f(t_{k-1})| \ls \\
\ls \sum\lt_{k=1}^n L(t_k-t_{k-1})=L \sum\lt_{k=1}^n (t_k-t_{k-1})=L( b-a),
\end{multline*}
so the assertion of theorem follows.
\hfill $\square$

\begin{sle}
\label{corol1}
If function~$f(\cdot)$ is differentiable at each point~$[a,\,b]$
and its derivative is bounded, then~$f(\cdot)$~ is a function of bounded variation.
\end{sle}
\doc\,
Let~$|f'(t)|\ls M_1.$ By the Lagrange finite increment formula, 
$$
|f(t)-f(s)|=|f'(\xi)(t-s)|=|f'(\xi)|\cdot|t-s| \ls M_1|t-s|
$$
($\xi$~ is a
point strictly between~$t$ and~$s$),
i.e. $f(\cdot)$ satisfies the Lipschitz condition with constant~$M_1.$
\hfill $\square$

\begin{teo}
\label{ch2thm3}
{\it A function of bounded variation is bounded.}
\end{teo}
\doc\,
Note first that for any two points~$t$ and~$s$ from~$[a,\,b]$, one has inequality
\begin{equation}
\label{ch21+}
|f(t)-f(s)| \ls\bigvee\lt_a^b(f).
\end{equation}
Indeed, let $a\ls t<s\ls b$ and let partition~$\tau$ consists of
points $a,\,t,\,s,\,b$ (for~$t=a$ or~$s=b$ there are only~2 or~3 points).
$$
v_\tau(f)=|f(t)-f(a)|+|f(s)-f(t)|+|f(b)-f(s)|\ls \bigvee\lt_a^b (f)
$$
so (\ref{ch21+}) follows. Let~$t\in[a,\,b]$~ be arbitrary.
We have from~(\ref{ch21+})
$$
|f(t)|-|f(a)| \ls |f(t)-f(a)|\ls \bigvee\lt_a^b(f),
$$
so $|f(t)|\ls|f(a)|+\bigvee\lt_a^b(f).$
\hfill $\square$
\begin{teo}
\label{ch2thm4}
{\it Let~$f(\cdot)$ and~$g(\cdot)$ be functions of bounded\linebreak 
variation on~$[a,\,b].$
Then the following functions also have finite variation:

1)~$\al f(t)+\bt g(t)\q (\al,\,\bt\in\R),$

2)~$f(t)\cdot g(t),$

3)~$f(t)/g(t)$ under the additional assumption that
${g(t)\!\gs\! c\!>\!0}.$}
\end{teo}
\doc\,
Let~$\tau$~ be an arbitrary partition of interval $[a,\,b]$. We have
$$
\gathered
v_\tau(\al f+\bt g)=\sum\lt_{k=1}^n |(\al f(t_k)\!+\!\bt g(t_k))\!-\!
(\al f(t_{k-1})\!+\!\bt g(t_{k-1}))|= \\
=\sum\lt_{k=1}^n |\al(f(t_k)-f(t_{k-1}))+\bt(g(t_k)-g(t_{k-1})) | \ls \\
\ls |\al|\cdot v_\tau(f)+|\bt|\cdot v_\tau(g) \ls |\al|\cdot
\bigvee\lt_a^b(f)+|\bt|\cdot \bigvee\lt_a^b(g).
\endgathered
$$
This implies that $\al f+\bt g$ has finite variation and
the following useful estimate is valid:
$
\bigvee\lt_a^b(\al f+\bt g)\ls |\al|\cdot \bigvee\lt_a^b(f)+|\bt|\cdot \bigvee\lt_a^b(g).
$

Theorem~\ref{ch2thm3} implies that there are constants~$M_f$ and~$M_g$ such
that $|f(t)|\ls M_f$ and $|g(t)|\ls M_g$ for all $t\in[a,\,b].$ We have
$$
\gathered
v_\tau(f\cdot g)=\sum\lt_{k=1}^n |f(t_k)\cdot g(t_k)-f(t_{k-1})\cdot
g(t_{k-1})|= \\
=\sum\lt_{k=1}^n |(f(t_k)g(t_k)-f(t_k)g(t_{k-1}))
+(f(t_k)g(t_{k-1})-f(t_{k-1})g(t_{k-1}))| \ls \\
\ls \sum\lt_{k=1}^n |f(t_k)||g(t_k)-g(t_{k-1})|+\sum\lt_{k=1}^n
|g(t_{k-1})|\cdot|f(t_k)-f(t_{k-1})| \ls \\
\ls M_f v_\tau(g)+M_g v_\tau(f) \ls M_f\bigvee\lt_a^b(g)+
M_g\bigvee\lt_a^b(f).
\endgathered
$$
Therefore, the product $f(\cdot)\cdot g(\cdot)$ has finite total
variation, and
$$
\bigvee\lt_a^b(f\cdot g) \ls M_f \bigvee\lt_a^b(g)+M_g \bigvee\lt_a^b(f).
$$

Let us prove that $\frac{1}{g(t)}$~ is a function of bounded variation. We have
$$
\gathered
v_\tau\left(\frac 1g\right)=\sum\lt_{k=1}^n \left| \frac{1}{g(t_k)}-\frac{1}{g(t_{k-1})}
\right|=\sum\lt_{k=1}^n \frac{|g(t_k)-g(t_{k-1})|}{g(t_k)\cdot g(t_{k-1} )}
\ls \\ \ls\frac{1}{c^2} v_\tau(g) \ls \frac{1}{c^2}\bigvee\lt_a^b(g).
\endgathered
$$
The finiteness of the total variation of the fraction~$\frac{f(t)}{g(t)}$ follows from
representation $\frac{f(t)}{g(t)}=f(t)\cdot\frac{1}{g(t)}$ and statement~2).
\hfill $\square$

Note that condition $g(t)\gs c>0$ cannot be replaced by condition~$g(t)>0.$
For example, put $g(t)=1-\{t\}$ on~$[0,\,1],$ where~$\{t\}$~ is the fractional part
of number $t.$ Obviously, $g(t)>0$ and (prove this) $\bigvee\lt_0^1(g)=2.$
However, function~$\frac{1}{g(t)}$ is not bounded and, therefore,
by Theorem~\ref{ch2thm3}, it cannot be a function of bounded variation.

Let us denote the set of functions of bounded variation on $[a,\,b]$
by ${\bf{BV}}[a,\,b].$ Then assertion~1) of Theorem~\ref{ch2thm4} means that \linebreak
${\bf{BV}}[a,\,b]$~is a {\it vector space}, and statements~1)
and~2) imply that ${\bf{BV}}[a,\,b]$~ is an {\it algebra} (the axioms of algebra
follow from the fact that the operations on functions are defined pointwise).

%

\bigskip

\poonkt{Additivity property of the total variation}\phantom{01234567890123456789} \

Total variation is an {\it additive function of the interval:}
\begin{teo}
\label{ch2thm5}
{\it Let~$c\in(a,\,b).$ Then
\begin{equation}
\label{ch22+}
{\bf{BV}}[a,\,b]={\bf{BV}}[a,\,c]\cap {\bf{BV}}[c,\,b],
\end{equation}
and if $f\in {\bf{BV}}[a,\,b],$ then
\begin{equation}
\label{ch23+}
\bigvee\lt_a^b(f)=\bigvee\lt_a^c(f)+\bigvee\lt_c^b(f).
\end{equation} }
\end{teo}
\doc\,
Let $f\in {\bf{BV}}[a,\,b],$ $\tau'=\{t_k\}_{k=0}^m$~ is a partition of interval $[a,\,c],$ $\tau''=\{s_k\}_{k=0}^n$~ is a partition of interval $[c,\,b].$
Then $\tau=\tau'\cup\tau''$~ is a partition of interval $[a,\,b]$ containing
point~$c=t_m=
s_0,$ and
$$
\gathered
v_{\tau'}(f)+v_{\tau''}(f)\!=\!\sum\lt_{k=1}^m |f(t_k)-f(t_{k-1} )|\!+
\!\sum\lt_{i=1}^n |f(s_k)-f(s_{k-1})|= \\
=v_\tau(f) \ls \bigvee\lt_a^b(f).
\endgathered
$$
Since the terms on the left hand side are non-negative, we have
$$
v_{\tau'}(f) \ls \bigvee\lt_a^b(f) \q\text{u}\q v_{\tau''}(f)\ls
\bigvee\lt_a^b(f),
$$
i.e. $f\in {\bf{BV}}[a,\,c]\bigcap {\bf{BV}}[c,\,b].$ So,
the left hand side in~(\ref{ch22+}) is included in the right hand side. Moreover, passing in
inequality
$$
v_{\tau'}(f)+v_{\tau''}(f)\ls \bigvee\lt_a^b(f)
$$
to the supremum over all partitions~$\tau'$ of interval $[a,\,c]$
and over all partitions $\tau''$ of interval $[c,\,b],$ we obtain inequality
\begin{equation}
\label{ch24+}
\bigvee\lt_a^c(f)+\bigvee\lt_c^b(f) \ls \bigvee\lt_a^b(f).
\end{equation}

Let now $f\in {\bf{BV}}[a,\,c]\bigcap {\bf{BV}}[c,\,b]$ and let $\tau$ be an
arbitrary partition of the closed interval $[a,\,b].$ As was already noted, we may assume without loss of generality that  $c\in\tau$. Let $\tau=\{t_k\}_{k=0}^n$
and~$c=t_p.$ Put $\tau'=\{t_k\}_{k=0}^p$ and $\tau''=\{t_k\}_{k=p}^n$;
 $\tau'$ and~$\tau''$ are partitions of intervals $[a,\,c]$
and~$[c,\,b]$, respectively. Then
$$
\gathered
v_{\tau}(f)=\sum\lt_{k=1}^n |f(t_k)-f(t_{k-1})|=\sum\lt_{k=1}^p
|f(t_k)-f(t_{k-1})|+ \\
+\!\sum\lt_{k=p+1}^n |f(t_k)-f(t_{k-1})|=v_{\tau'}(f)+
v_{\tau''}(f) \ls \bigvee\lt_a^c(f)+\bigvee\lt_c^b(f).
\endgathered
$$
Inequality $v_\tau(f)\ls \bigvee\lt_a^c(f)+\bigvee\lt_c^b(f)$ means
that $f\in {\bf{BV}}_{[a,\,b]}$ and hence the right hand side  in~(\ref{ch22+})
is contained in the left hand side, and passing in this inequality to the supremum,
we obtain the inequality which is opposite to (\ref{ch24+}).

Therefore, equalities~(\ref{ch22+}) and~(\ref{ch23+}) hold.
\hfill $\square$

\begin{sle}
\label{corol2}
{\it If~$f(\cdot)$ is piecewise monotone on ${[a,\,b]}$ {\rm(}that is, if $[a,\,b]$
can be split by points $c_k,$ $a<c_1<\ldots<c_p<b$ into
intervals~$[c_{k-1},\,c_k],$ $k=1,\,\ldots,\,p+1,$ $c_0=a,$ $c_{p+1}=b, $
so that on each such interval $f(\cdot)$ is monotone\rm{)}, then $f\in {\bf{BV}}[a,\,b].$
Moreover,}
$$
\bigvee\lt_a^b(f)=\sum\lt_{k=1}^{p+1} \bigvee\lt_{c_{k-1}}^{c_k}(f).
$$
\end{sle}
\doc\,
Since~$f(\cdot)$ is monotone on~$\;[c_{k-1},\,c_k],$ we have by Theorem~2.1 
$f\in {\bf{BV}}[c_{k-1},\,c_k].$ By induction, equality~(\ref{ch23+}) can be
extended to any finite number of terms.
\hfill $\square$


\bigskip

\poonkt{Criterion for finiteness of total  variation}\phantom{01234567890123456789}

\begin{teo}
\label{ch2thm6}
{\it A function $f:[a,\,b]\to\R$ has finite\linebreak 
total variation if and only if it can be represented as the\linebreak 
difference of two increasing functions.}
\end{teo}
\doc\,
{\it Sufficiency.\;}
Let $f(t)=f_1(t)-f_2(t),$ %
where~$f_1(\cdot)$ and~$f_2(\cdot)$
are increasing functions. By Theorem~4.1, $f_i\in {\bf{BV}}[a,\,b]$ $(i=1,\,2),$
and, by assertion~1) of Theorem~4.4,  $f_1-f_2 \in {\bf{BV}}[a,\,b].$

{\it Necessity.\;}
We associate to every $f\in {\bf{BV}}[a,\,b]$
functions~$f_{\pi}(\cdot)$ and~$f_\nu(\cdot)$ as follows:
\begin{gather}
\label{ch25+}
f_\pi(a)\doteq 0, \q f_\pi(t)\doteq \bigvee\lt_a^t(f) \q\text{for}\q t>a, \\
\label{ch26+}
f_\nu(t)\doteq f_\pi(t)-f(t).
\end{gather}

Let us show that these functions~ are increasing.

Let $a\ls t<s\ls b.$ We apply equality~(\ref{ch23+}) in Theorem~4.5 to
interval $[a,\,s],$ taking point $t$ as $c.$ Then
\begin{multline}
\label{ch27}
f_\pi(s)-f_\pi(t)=\bigvee\lt_a^s(f)-\bigvee\lt_a^t(f)=\bigvee\lt_a^t(f)+
\bigvee\lt_t^s(f)-\bigvee\lt_a^t(f)= \\
=\bigvee\lt_t^s(f) \gs 0.
\end{multline}
Further, using~(\ref{ch27}) and~(\ref{ch21+}), we see that
$$
\begin{gathered}
f_\nu(s)-f_\nu(t)=(f_\pi(s)-f(s))-(f_\pi(t)-f(t))= \\
=(f_\pi(s)-f_\pi(t))-(f(s)-f(t)) \overset{(\ref{ch27})}{=}
\bigvee\lt_t^s(f)-(f(s)-f(t)) \overset{(\ref{ch21+})}{\gs} 0.
\end{gathered}
$$
Thus, $f_\pi(\cdot)$ and~$f_\nu(\cdot)$ are indeed increasing
functions. From definition\linebreak
(\ref{ch26+}) we have
\begin{equation}
\label{ch28}
f(t)=f_\pi(t)-f_\nu(t),
\end{equation}
as needed
\hfill $\square$

\begin{zam}
\label{zam100}
{\it In the representation of a function of bounded variation
as the difference of two increasing functions, both of these functions can be chosen to be strictly
increasing.}
\end{zam}

Indeed, let $g: [a,\,b]\to\R$~ be an arbitrary strictly 
increasing function. By (\ref{ch28}),
$$
f(t)=f_\pi(t)-f_\nu(t)=\bigl(f_\pi(t)+g(t)\bigr)-\bigl(f_\nu(t)+g(t )\bigr);
$$
where the functions in brackets are strictly increasing.

It is also possible to ensure that {\it in this representation both functions are
non-negative.}

\begin{sle}
\label{corol3}
{\it A function of bounded variation can 
can have at most countable set of points of discontinuity.}
\end{sle}
\doc\,
Let $f(t)=f_1(t)-f_2(t),$ where $f_1$ and $f_2$ are
increasing functions. Then $T(f)\subset T(f_1)\cup T(f_2)$ and Corollary~\ref{corol3} follows from Theorem~\ref{ch1thm1}.
\hfill $\square$
\begin{sle}
\label{corol4}
{\it A function of bounded variation can \linebreak
be represented as the sum of
a continuous function of bounded variation and a jump function, which
is also a function of bounded variation.}
\end{sle}
\doc\,
Let us represent~$f(t)$ by the formula~(\ref{ch28}) and apply Theorem~\ref{ch1thm2}
to functions $f_\pi(\cdot)$ and~$f_\nu(\cdot):$
\begin{equation}
\label{ch29}
f_\pi(t)=f_{\pi c}(t)+f_{\pi d}(t), \q
f_\nu(t)=f_{\nu c}(t)+f_{\nu d}(t),
\end{equation}
where~$f_{\pi c}(\cdot)$ ($f_{\nu c}(\cdot)$)~ is the continuous part,
and $f_{\pi d}(\cdot)$ ($f_{\nu d}(\cdot)$)~ is the discrete part~$f_\pi(\cdot)$
($f_\nu(\cdot)$). Subtracting from the first equality in~(\ref{ch29}) the second equality,
we obtain: \;
$$
f_\pi(t)-f_\nu(t)=(f_{\pi c}(t)-f_{\nu c}(t))+(f_{\pi d}(t)-f_{\ nu d}(t)).
$$
Therefore,
\begin{equation}
\label{predst1}
f(t)=f_c(t)+f_d(t),
\end{equation}
where, by Theorem~\ref{ch2thm6}, $f_c(\cdot)$~ is a
continuous function of bounded variation, and~$f_d(\cdot)$~ is a
jump function, which also has finite total variation.
\hfill $\square$

Note that due to the convention $f_d(a)=0$ the representation (\ref{predst1}) is unique.

Note also that~$f_d(\cdot)$ and~$f_c(\cdot)$ can be found
directly using formulas\linebreak
(\ref{ch17}) and~(\ref{ch19}). We define the {\it unit function} as follows:
$$
\mathfrak h_c(t)=\left\{
\begin{array}{lr}
0\;\;\text{for}\;\;t\leqslant c\;\;(a\leqslant c<b),\\
1 \;\;\text{for}\;\; t>c\;\;(a<c<b), \\
1 \;\;\text{for}\;\;t=b\;\;(c=b).\\
\end{array}
\right.
$$
Function $\mathfrak h_c(\cdot)$ is left continuous at all points of interval $[a,b).$
By means of this function, the representation (\ref{ch17}) can be written as
\begin{equation}
\label{predst1fogrvar}
f_d(a)=0,\quad f_d(t)=\sum\limits_{s\in T(f)}\sigma_{s}(f)\mathfrak h_s(t)+\sigma_{t}^-( f)\quad (t>a),
\end{equation}
\begin{sle}
\label{corol5}
{\it A function of bounded variation is Riemann integrable.}
\end{sle}

To conclude this section, we present one more {\it criterion} of finiteness of total variation.
\begin{teo}
\label{ch2dop1}
{\it $f(\cdot)$ is a function of bounded variation on $[a,\,b]$
if and only if there is an increasing
function~$F(\cdot)$ such that for every pair
$t,\,s\in[a,\,b],$ $s>t$, we have inequality $$|f(s)-f(t)| \ls 
F(s)-F(t).$$}
\end{teo}
\doc\,
{\it Sufficiency.\;} Let $\tau=\{t_k\}_{k=0}^\infty$~ be an arbitrary partition of interval $[a,b].$ We have
$$
v_{\tau}(f)=\sum\limits_{k=1}^\infty |f(t_k)-f(t_{k-1})| \leqslant \sum\limits_{k=1}^\infty \bigl(F(t_k)-F(t_{k-1})\bigr)=F(b)-F(a).
$$
Hence $f\in {\bf{BV}}[a,b].$

{\it Sufficiency.\;} Let $f\in {\bf{BV}}[a,b].$ We set
$F(t)\doteq f_{\pi}(t)$. Let $s>t.$ Then $|f(s)-f(t)| \ls \bigvee\lt_{s}^t(f)=f_{\pi}(t)-f_{\pi}(s)=F(t)-F(s).$
\hfill $\square$


\bigskip

\poonkt{Continuous functions of bounded variation\;}\phantom{01234567890123456789} \

{\bf{1.}}\;Let us begin with the following result.
\begin{teo}
\label{ch2thm7}
{\it Let function $f(\cdot)\in {\bf{BV}}[a,\,b]$ be continuous from the right $($continuous from the left$)$
at point $t_0\in[a,\,b)$ $(t_0\in(a,\,b])$. Then~$f_\pi(\cdot)$ is also
continuous from the right $($from the left$)$ at this point.}
\end{teo}
\doc\,
Let~$f(\cdot)$ be continuous from the right at $t_0\in[a,\,b)$. Let $\eps>0$
is arbitrary. Consider a partition $\tau=\{t_k\}_{k=0}^n$ of interval $[t_0,\,b]$
such that
$$
\gathered
\sum\lt_{k=1}^n |f(t_k)-f(t_{k-1})| > \bigvee\lt_{t_0}^b(f)-\frac{\eps}{2}, \qquad
|f(t_1)-f(t_0)|< \frac{\eps}{2}.
\endgathered
$$
(The first inequality follows from the property of the supremum.
If 
the second inequality is not satisfied, then we add a new point~$t_1$ to satisfy the second inequality, while the first inequality remains valid.) We obtain from these inequalities:
$$
\gathered
\bigvee\lt_{t_0}^b(f)< \sum\lt_{k=1}^n |f(t_k)-f(t_{k-1})|+\frac{\eps}{2} <
\eps+\sum\lt_{k=2}^n |f(t_k)-f(t_{k-1})| \ls \\
\ls \eps+\bigvee\lt_{t_1}^b(f)
\Rightarrow \bigvee\lt_{t_0}^b(f)-\bigvee\lt_{t_1}^b(f)<\eps \Rightarrow
\bigvee\lt_{t_0}^{t_1}(f)<\eps
\Rightarrow f_\pi(t_1)-f_\pi(t_0)<\eps.
\endgathered
$$
This proves continuity of $f_\pi(\cdot)$ from the right at point $t_0\in[a,\,b)$.
In the case of continuity from the left, we consider interval $[a,\,t_0]$ and argue similarly.
\hfill $\square$
\begin{sle}
\label{raznsrogovozr}
{\it A continuous function of bounded variation can be represented in
as the difference of increasing continuous functions. 
It is possible to ensure that in this representation both continuous functions are strictly increasing.}
\end{sle}
\doc\,
Let $f(\cdot)$~ be a continuous function of bounded variation.
By the previous theorem,~$f_\pi(\cdot)$ is also continuous on~$[a,\,b],$ and so
is continuous $f_\nu(t)=f_\pi(t)-f(t).$
\hfill $\square$

Regarding continuous functions of bounded variation, see also {\bf{5.4}}.


{\bf{2.}}
Let $G=\bigcup\limits_{k}(a_k,b_k)$~ be an open set, let $f(\cdot)$~ be a continuous function of bounded variation.
Then $f_{\pi}(\cdot)$ and $f_{\nu}(\cdot)$~ are increasing continuous functions and
\begin{multline*}
\mu_{f_{\pi}}(G)=\sum\limits_{k}\bigl(f_{\pi}(b_k)-f_{\pi}(a_k)\bigr)=\sum\limits_{k }\bigvee\limits_{a_k}^{b_k}(f),\qquad \mu_{f_{\nu}}(G)=\\
=\sum\limits_{k}\bigl(f_{\nu}(b_k)-f_{\nu}(a_k)\bigr)=\sum\limits_{k }\Bigl(\bigl(f_{\pi}(b_k)-f(b_k)\bigr)-\bigl(f_{\pi}(a_k)-f(a_k)\bigr)\Bigr)=\\=
\sum\limits_{k}\bigvee\limits_{a_k}^{b_k}(f)-\sum\limits_{k}\bigl(f(b_k)-f(a_k)\bigr)\leqslant 2\mu_{ f_{\pi}}(G).
\end{multline*}
Therefore, a set of $f_{\pi}$-measure zero is also a set 
$f_{\nu}$-measures zero.

{\bf{3.}}
The definition of total variation can be extended to unbounded intervals
$[a,\,+\fy),$ $(-\fy,\,b],$ $(-\fy,\,+\fy)$ as follows:
$$
\gathered
\bigvee\lt_a^{+\fy}(f)\doteq \lim\lt_{b\to+\fy} \bigvee\lt_a^b(f), \q
\bigvee\lt_{-\fy}^b(f)\doteq \lim\lt_{a\to -\fy} \bigvee\lt_a^b(f), \\
\bigvee\lt_{-\fy}^{+\fy}(f)\doteq \bigvee\lt_{-\fy}^a(f)+
\bigvee\lt_a^{+\fy}(f)=\lim\lt_{a\to -\fy \atop b\to +\fy} \bigvee\lt_a^b(f).
\endgathered
$$


\bigskip

\poonkt{Banach space ${\bf{BV}}[a,b]$}\phantom{01234567890123456789} \

For an element $x$ of vector space ${\bf{BV}}[a,b]$, define
\begin{equation}
\label{normbv1}
\|x\|\doteq |x(a)|+\bigvee\lt_a^b(x).
\end{equation}


$1^{\circ}.\;$
{\it Equality} (\ref{normbv1}) {\it defines a norm on ${\bf{BV}}[a,b],$ i.e. ${\bf{BV}}[a,b]$~ is an LNS.}

\doc\,
Axiom of identity. If $x=0,$ then, clearly, $\|x\|=0.$ Let $\|x\|=0$. Then simultaneously
$|x(a)|=0$ and $\bigvee\lt_a^b(x)=0$. The second equality means that $x(t)\equiv C=\const,$ and the first equality means that $C=0,$ i.e.\,$x=0.$

Positive homogeneity and the triangle inequality follow from the first assertion of Theorem \ref{ch2thm4}.
\hfill $\square$


$2^{\circ}.\;$
{\it ${\bf{BV}}[a,b]$~ is a Banach space.}

\doc\,
Let $\{x_n(\cdot)\}_{n=1}^{\infty}$~ be a fundamental sequence of functions
of bounded variation, that is,
$$
(\forall\varepsilon>0)\ (\exists N\in\mathbb N)\ (\forall n,m>N)
\bigl(||x_n-x_m||<\varepsilon\bigr).
$$
Denote for brevity $y_{nm}(t)=x_n(t)-x_m(t).$ Then
$$
\gathered
\bigl|x_n(t)-x_m(t)\bigr|=\bigl|y_{nm}(t)\bigr|\leqslant \bigl|y_{nm}(t)-y_{nm}(a)\ bigr|+\bigl|y_{nm}(a)\bigr|\leqslant \\
\leqslant \bigl|y_{nm}(a)\bigr|+\bigvee\limits_a^b(y_{nm})= ||y_{nm}||=||x_n-x_m||<\varepsilon.
\endgathered
$$
Thus, for every $t\in[a,\,b]$, the sequence $\{x_n(t)\}_{n=1}^{\infty}$ of real numbers is fundamental, so there exists pointwise limit
$x(t)=\lim\limits_{n\to\infty}x_n(t).$

Now, two facts need to be established: a)\; $x\in {\bf{BV}}[a,\,b]$ and b)\; $x_n\to x$
$(n\to\infty)$ in ${\bf{BV}}[a,\,b].$

Let $\tau=\{t_k\}_{k=0}^p$~ be an arbitrary partition of $[a,\,b]: $\\
$ a=t_0<t_1<\cdots<t_p= b.$
$$
\bigl|y_{nm}(a)\bigr|+\sum\limits_{k=1}^{p}\bigl|y_{nm}(t_k)-y_{nm}(t_{k-1}) \bigr|<||x_n-x_m||<\varepsilon.
$$
Letting $m\to\infty,$ we get
$$
\bigl|x_n(a)-x(a)\bigr|+\sum\limits_{k=1}^{p}\Bigl|\bigl(x_n(t_k)-x(t_{k})\bigr) -\bigl(x_n(t_{k-1})-x(t_{k-1})\bigr)\Bigr|\leqslant \varepsilon.
$$
This implies that $x_n(\cdot)-x(\cdot)$~ is a function of bounded variation and
\begin{equation}\label{bv101}
\bigl|x_n(a)-x(a)\bigr|+\bigvee\limits_a^b(x_m-x)=||x_n-x||\leqslant \varepsilon.
\end{equation}
Since $x(t)=x_n(t)-\bigl(x_n(t)-x(t)\bigr),$ we have  $x(\cdot)\in {\bf{BV}}[a,\,b],$and in view of (\ref{bv101})
$x_n\to x$ in ${\bf{BV}}[a,\,b].$
\hfill $\square$


$3^{\circ}.\;
${\it Space ${\bf{BV}}[a,b]$ is not separable.}

\doc\,
Put
$$
S\doteq \{x_c:\,[a,b]\to \mathbb R:\; x_c(t)=0,\;t\leqslant c,\;x_c(t)=1,\;t>c,\;c\in [a,b]\}.
$$
Set $S\subset {\bf{BV}}[a,\,b]$ is uncountable and $\|x_c-x_{c'}\|=2$ for $c\ne c'.$
By Theorem \ref{mprth6}, ${\bf{BV}}[a,\,b]$ is not separable.
\hfill $\square$

In view of Theorem \ref{ch2thm3}, the norm in space ${\bf{BV}}[a,b]$ can be defined differently.
Namely, define
\begin{equation}
\label{normbv2}
\|x\|^*\doteq \underset{t\in [a,b]}{\sup}|x(t)|+\bigvee\lt_a^b(x).
\end{equation}
Since both terms have the properties of a norm, (\ref{normbv2}) is indeed a norm on the
space of functions of bounded variation.


$4^{\circ}.\;$
{\it Norms $\|\cdot\|$ and $\|\cdot\|^*,$ defined by equalities } (\ref{normbv1}) {\it and} (\ref{normbv2}) {\it respectively, are equivalent.}

\doc\,
Since the value of a function at a point cannot exceed its supremum on the entire interval, we have $\|x\|\leqslant \|x\|^*.$ On the other hand,
$$
|x(t)|-|x(a)|\leqslant |x(t)-x(a)|\leqslant \bigvee\limits_a^b(x)\quad \Longrightarrow |x(t)|\leqslant | x(a)|+\bigvee\limits_a^b(x),
$$
which yields
\begin{multline*}
\underset{t\in [a,b]}{\sup}|x(t)|\leqslant |x(a)|+\bigvee\limits_a^b(x) \Longrightarrow \|x\|^*\leqslant
|x(a)|+2\bigvee\limits_a^b(x)\leqslant 2\|x\|.
\end{multline*}
Thus, $\frac{1}{2}\|x\|^*\leqslant \|x\|\leqslant \|x\|^*.$
\hfill $\square$

\medskip


$5^{\circ}.\;$
Define $${\bf{CBV}}[a,b]\doteq {\bf{C}}[a.b]\bigcap {\bf{BV}}[a,b]$$ with norm (\ref{normbv2}) :\\
$\|x\|_{{\bf{CBV}}}\doteq \underset{t\in [a,b]}{\max}|x(t)|+\bigvee\lt_a^b(x) .$

${\bf{CBV}}[a,b]$~ is a\;{\it Banach space and hence is a subspace of Banach space} ${\bf{BV}}[a,b].$

\doc\,
Let $x*\in {\bf{BV}}[a,b]$~ be a limit point of 
${\bf{CBV}}[a,b]$. Then there is a sequence $\{x_n\}_{n=1}^\infty\subset {\bf{CBV}}[a,b],
x_n\to x^*\;(n\to\infty)$ in ${\bf{BV}}[a,b].$\linebreak
By assertion \textbf{4}, $x_n(t)\rightrightarrows x^*(t),\;$ and hence $x_n(\cdot)$ is continuous, i.e.
$x^*\in {\bf{CBV}}[a,b]$. Thus, ${\bf{CBV}}[a,b]$ is closed.
By Theorem  \ref{mprth1}, ${\bf{CBV}}[a,b]$ is complete and is a subspace of ${\bf{BV}}[a,b].$
\hfill $\square$


\bigskip

\section{Functions of bounded variation-2}
\label{ch3}


\hskip 1mm
\poonkt{Boundary. Module. Superposition}\phantom{01234567890123456789} \

We begin with a few simple results.



{\bf{1.}}\, Let $A\subset [a,b],$ let $\partial A$~ be the boundary of set $A,\;\chi_A(t)=\kappa_A(t)=\linebreak
=1$ for $t\in A,\;\;\chi_A(t)=0,\;\;\kappa_A(t)=-1\;$ for $t\notin A\;(\chi_A(\cdot)$ is the
characteristic function of set $A).$

{\it Functions $\chi_A,\,\kappa_A$ have finite total variation if and only if $\partial A$~ is a finite set.}

\doc\,
Let~$\pt A$ contain~$m$ points, let $\tau$~ be an arbitrary partition of interval $[a,\,b].$ All terms in the sum $v_\tau(\chi_A)=
 \sum\lt_{k=1}^n |\chi_A(t_k)-\chi_A(t_{k-1})|$ are either equal to 0 (when points~$t_{k-1}$
and~$t_k$ either both belong to $A,$ or both belong to $[a,\,b] \setminus A$),
or equal~1 (when one of points~$t_{k-1}$ or~$t_k$ belongs to~$A,$
and the other one belongs to~$[a,\,b] \setminus A$). Therefore, $v_\tau (\chi_A) \ls m,$
which means that the total variation of $\chi_A(\cdot)$ is finite. Analogous argument works for $\kappa_A.$

Now let~$\pt A$~ be an infinite set. Let us choose a partition
$\tau=\{t_k\}_{k=0}^n$ which alternates points belonging to~$A,$
and points not belonging to $A.$ Since any neighborhood of a bundary point
contains points of both types, such a partition exists for any
large~$n.$ For such partition,
\begin{multline*}
v_\tau (\chi_A)=\sum\lt_{k=1}^n |\chi_A(t_k)-\chi_A(t_{k-1})|=n,\\
 v_\tau (\kappa_A)=\sum\lt_{k=1}^n |\kappa_A(t_k)-\kappa_A(t_{k-1})|=2n.
\end{multline*}
Therefore, $\bigvee\lt_a^b (\chi_A)=+\fy,\;\;\bigvee\lt_a^b (\kappa_A)=+\fy$
\hfill $\square$

For example, Dirichlet function $D(t)=\chi_{\mathbb Q}(t)$ has infinite total variation.


{\bf{2.}} {\it If~$f(\cdot)$ is a function of bounded variation, then
$|f(\cdot)|$ is also a function of bounded variation, and
$\bigvee\lt_a^b(|f|) \ls \bigvee\lt_a^b(f).$ The converse is in general not true.}

\doc\,
Let~$f(\cdot)$~ be a function of bounded variation on~$[a,\,b]$.
We have for any partition~$\tau$\;\;
$$\tau(|f|)=\sum\lt_{k=1}^n \Big| |f(t_k)|-|f(t_{k-1})| \big| \ls \sum\lt_{k=1}^n \Big| f(t_k)-f(t_{k-1}) \Big|=v_\tau (f) \ls \bigvee\lt_a^b (f).
$$
Therefore, $|f(\cdot)|$ has finite total variation on~$[a,\,b]$
and $\bigvee\lt_a^b (|f|) \ls \bigvee\lt_a^b (f).$

Let $A\subset [a,b]$, and let $\partial A$~ be an infinite set.
Then 
$\bigvee\lt_a^b (\kappa_A)=+\fy,$\; but $\bigvee\lt_a^b (|\kappa_A|)=0.$
\hfill $\square$

However, {\it if~$f(\cdot)$ is continuous, then the converse is also true,
where} $\bigvee\lt_a^b (|f|)=\bigvee\lt_a^b (f).$

\doc\,
It remains to prove the opposite
to the previously proven inequality.
So, let~$|f(\cdot)|$ have finite total variation on~$[a,\,b],$
and let $\tau$ be an arbitrary partition of this closed interval,
$v_\tau (f)\!=\!\sum\lt_{k=1}^n\! |f(t_k)-f(t_{k-1})|.$  \linebreak
Let us divide all terms in this sum
into two groups. The first group includes the terms for  which
$f(t_k) \cdot f(t_{k-1}) \gs 0.$ For the terms of this group
$$
\big| f(t_k)-f(t_{k-1}) \Big|=\Big| |f(t_k)|-|f(t_{k-1})| \Big|.
$$
The second group includes the terms for which
${f(t_k)\cdot f(t_{k-1})<0}.$
For each of these terms, due to the continuity of $f(\cdot)$
there is at least one point
${t_k' \in (t_{k-1},\,t_k)}$ where~$f(t_k')=0.$ For the terms in this group 
\begin{align*}
\big| f(t_k)-f(t_{k-1}) \Big| & \ls\Big| f(t_k) \Big|+\Big| f(t_{k-1}) \Big| \\
&=
\Big| |f(t_k)|-|f(t_k')| \Big|+\Big| |f(t_k')|-|f(t_{k-1})| \Big|.
\end{align*}\

Let us denote the partition obtained by adding to~$\tau$ the points~$t_k'$,
corresponding to the terms of the second group, by $\tau'.$ Then
$$
v_\tau (f)={\sum}^{'}+{\sum}^{''} \ls v_{\tau'}(|f|) \ls \bigvee\lt_a^b (|f| ).
$$
This implies that~$f(\cdot)$ has finite total variation on~$[a,\,b]$
and $\bigvee\lt_a^b (f) \ls \bigvee\lt_a^b (|f|).$ Together with the opposite inequality proved earlier
this yields the required.
\hfill $\square$


{\bf{3.}} Is the superposition of functions of bounded variation a function of bounded variation?
Let us consider the following example.
\begin{example}
\label{superp}
{\it Both functions $f(t)=\sqrt t$ and $g(t)=t^2\sin\,t^2$ for $t>0$,
$g(0)=0,$  have finite total variation on $[0,1]$} (the first one is strictly increasing, the second one has bounded derivative),
{\it and the function $f(g(t))=t\left|\sin\,\frac{1}{t}\right|$ for $t>0,\,f(g(0))=0 ,$ has infinite total variation on $[0,1].$}
To see this, it suffices to show, in view of assertion {\bf{2}},  that {\it function
$h(t)\doteq t\,\sin\,\frac{1}{t}$ for $t>0,\,h(0))=0$ has infinite total variation.}
For this purpose, we take the following partition $\tau:$
$$
0<\frac{2}{2n+1}<\frac 1n<\frac{2}{2n-1}<\frac{1}{n-1}<\ldots<\frac 23 < 1,
$$
and then argue as in paragraph \textbf{4.2}.
\end{example}

Nevertheless, the following assertion is true.

\begin{teo}
\label{superpoz}
{\it Let~$f(\cdot)$ have finite total variation on~$[a,\,b],$ let $g(\cdot)$ be a
strictly increasing function continuous on~$[\al,\,\bt]$, $g(\al)=a,$
$g(\bt)=b.$ Then $F(t)=f(g(t))$~ is a function of bounded variation on
interval~$[\al,\,\bt],$ with $\bigvee\lt_\al^\bt(F)=\bigvee\lt_a^b(f).$}
\end{teo}

\doc\,
Let $\sg=\{\nu_k\}_{k=0}^n$~ be an arbitrary partition of closed interval \linebreak
$[\al,\,\bt].$ We set $t_k=g(\nu_k)$ $(k=0,\,1,\,\ldots,\,n).$ Since~$g(\cdot)$ is strictly increasing,
$\tau=\{t_k\}_{k=0}^n$~ is a partition of interval $[a,\,b].$ Then
$$
\gathered
v_\sg(F)\!=\!\sum\lt_{k=1}^n \left|F(\nu_k)-F(\nu_{k-1})\right|\!
=\!\sum\lt_{k=1}^n \left|f(g(\nu_k))-f(g(\nu_{k-1}))\right|\!= \\
=\!\sum\lt_{k=1}^n \left|f(t_k)-f(t_{k-1})\right|\!=\!v_\tau (f) \ls \bigvee\lt_a^b(f).
\endgathered
$$
Therefore, $F(\cdot)$ has finite total variation on~$[\al,\,\bt],$
and $\bigvee\lt_\al^\bt (F) \ls\linebreak\ls 
\bigvee\lt_a^b (f).$

Now, let $\tau=\{t_k\}_{k=0}^n$~ be an arbitrary partition of 
interval~$[a,\,b].$ The continuity and strict monotonicity of function~$g(\cdot)$
imply that there exists the inverse function~$g^{-1}(\cdot),$
which is continuous and is strictly increasing on interval~$[a,\,b].$
Put $\nu_k=g^{-1}(t_k).$ Then $\sg=\{\nu_k\}_{k=1}^n$~ is a  partition of ~$[\al,\,\bt]$, and $t_k=g(\nu_k).$ Therefore,
$$
\gathered
v_\tau (f)\!=\!\sum\lt_{k=1}^n |f(t_k)-f(t_{k-1})|\!
=\sum\lt_{k=1}^n
|f(g(\nu_k))-f(g(\nu_{k-1}))|\!= \\
=\!\sum\lt_{k=1}^n |F(\nu_k)-F(\nu_{k-1})|\!=\!v_\sg (F) \ls \bigvee\lt_\al^\bt (F),
\endgathered
$$
i.e. $\bigvee\lt_a^b (f) \ls \bigvee\lt_\al^\bt (F).$ Together with the opposite inequality proved earlier, this implies that
$\bigvee\lt_\al^\bt (F)=\bigvee\lt_a^b (f).$
\hfill $\square$


\bigskip

\poonkt{Representations of function and its total variation}\phantom{01234567890} \

Recall that convention (\ref{soglprodolz}) applies  everywhere below.


\medskip

{\bf{1.}}
Using Corollary \ref{corol4} and representation (\ref{ch17}) $\bigl($see also (\ref{predst1fogrvar})$\bigr)$, we see that
$f_d(\cdot)$ can be represented, up to the values at the interior points of discontinuity, as the sum of series
\begin{equation}
\label{var21}
f_d(t)=\sum\lt_{k=1}^\fy \sigma_{c_k}(f) \mathfrak h_{c_k}(t)
\end{equation}
$\bigl($recall, $T(f)=\{c_1,c_2,\ldots\}$~ is the set of discontinuity points of function $f(\cdot);$ if this set is finite,
then the series in (\ref{var21}) is replaced by a finite sum$\bigr)$. The representation (\ref{var21}) thus takes place for
$t\in [a,b]\setminus T(f)$. According to (\ref{ch17}), $f_d(c_{k_0})$ differs from the sum of the series in (\ref{var21}) by the term
$\sigma_{c_{k_0}}^-(f)$, $|\sigma_{c_{k_0}}^-(f)|\leqslant \bigvee\lt_a^b(f)$.
Let us show that
\begin{equation}
\label{predstvar1}
\bigvee\lt_a^b(f_d)=\sum\lt_k^\infty\;\bigl(|\sigma^-_{c_k}(f)|+|\sigma^+_{c_k}(f)|\bigr)\doteq S.
\end{equation}
(The series in (\ref{predstvar1}) converges because $f_d(\cdot)$~ is a function of bounded variation.)

First, let $T(f)=\{c_1,c_2,\ldots,c_n\}.$ Without loss of generality, we may assume that
$
c_1<c_2<\ldots<c_n.\;\;
$
We can also assume that $c_1=a,\;c_n=b.$ For $n=1,2$ representation (\ref{predstvar1}) is obvious. Let $n\geqslant 3.$
Since function $f_d(\cdot)$ is constant between the points of discontinuity, taking as a partition $\tau'=\{t_k\}_{k=0}^{2(n-1)},$ where
\begin{multline}
\label{razbextr}
t_0=c_1<t_1<c_1=t_2,\ldots,\;t_{2k}=c_{k+1}\;(k=0,1,\ldots,n-1),\\
t_{2(k-1)}<t_{2k-1}<t_{2k}\;\;(k=1,\ldots,n-1),
\end{multline}
and taking into account convention (\ref{soglprodolz}), we have
$$
\bigvee\lt_a^b(f_d)=\underset{\tau}{sup}\;v_{\tau}(f_d)=v_{\tau'}(f_d)=\sum\lt_k^n\bigl(| \sigma^-_{c_k}|+|\sigma^+_{c_k}|\bigr).
$$
i.e.\,in this case representation (\ref{predstvar1}) is true.

For a countable $T(f)$, the series in the representation of $f_d(\cdot)$ is majorized by a convergent series of real numbers
$\bigvee\lt_a^b(f)+\sum\lt_{k=1}^\fy |\sigma_k(f)|$, so it converges absolutely and uniformly.
Define at the points of continuity
\begin{equation}
\label{vtochnepr}
g_n(t)\doteq\sum\lt_{k=1}^n \sigma_k(f) \rho_{c_k}(t),
\end{equation}
and at the discontinuity points $c_{j}$
\begin{equation}
\label{vtochrazr}
g_n(c_j)\doteq\sum\lt_{k=1}^n \sigma_k(f) \rho_{c_k}(c_j)+\sigma_{c_j}^-(f)\qquad (j\in\mathbb N ).
\end{equation}
Then
\begin{equation}
\label{ravnsh}
g_n(t) \rightrightarrows f_d(t)\qquad (n\to\fy).
\end{equation}
Moreover, for any partition $\tau$ (indepdendent of $n$),
\begin{equation}
\label{ravnsh1}
v_{\tau}(g_n)\to v_{\tau}(f_d)\;\;\;(n\to \infty),
\end{equation}
and for any natural $n$
\begin{equation}
\label{ravnsh2}
v_{\tau}(g_n)\leqslant v_{\tau}(f_d).
\end{equation}

For a fixed $n$, we construct a partition $\tau'_n=\{t_k\}_{k=0}^{2(n-1)},$ where $t_k$ are defined in the same way as in (\ref{razbextr}).
Then, according to what was proved above,
\begin{multline}
\label{shpodposl}
\bigvee\lt_a^b(g_n)\!=\!v_{\tau'_n}(g_n)\!=\!\sum\lt_{k=1}^{2(n-1)}\,| g_n(t_k)-g_n(t_{k-1})|\!=\\=
\sum\lt_{k=1}^n\,\bigl(|\sigma^-_{c_k}|+|\sigma^+_{c_k}|\bigr)\to S
\end{multline}
and
\begin{equation}
\label{shpodpos2}
\bigvee\lt_a^b(g_n)\leqslant S
\end{equation}
In view of (\ref{ravnsh2}) --- (\ref{shpodpos2}), we have inequality  $\bigvee\lt_a^b(f_d)\geqslant S.$

Suppose that
 $\bigvee\lt_a^b(f_d)>S.$ Then there exists $\varepsilon >0$ such that $\bigvee\lt_a^b(f_d)>S+\varepsilon.$ \;\;
There also exists partition $\tau_{\varepsilon}$ such that $v_{\tau_{\varepsilon}}(f_d)>S+\dfrac{\varepsilon}{2}.$
It follows from (\ref{ravnsh1}) that there exists $n_0$ such that
$v_{\tau_{\varepsilon}}\bigl(g_{n_0}\bigr)\geqslant S+\dfrac{\varepsilon}{2},$
which contradicts to (\ref{shpodpos2}). So, equality (\ref{predstvar1}) holds.
\hfill $\square$

If a function of bounded variation is such that its values at points of discontinuity are situated
(non-strictly) between its one-sided limits at these points, then the total variation of this function does not depend on these values 
(see exercise \ref{ch2upr6}; see also 
\cite[exercise 2.6]{derr08}). In this case, the left and the right jumps of this function are of the same signs, so
$|\sigma^-_{c_k}|+|\sigma^+_{c_k}|=|\sigma^-_{c_k}(f)+\sigma^+_{c_k}(f)|=| \sigma_k(f)|,$ and the equality (\ref{predstvar1}) takes form
\begin{equation}
\label{predstvar1+}
\bigvee\lt_a^b (f_d)=\sum\lt_{k=1}^\fy |\sigma_k(f)|.
\end{equation}
For example, this is the case for the total variation of $f_d(\cdot)$ if $f(\cdot)$ is continuous on the left or on the right at each point of discontinuity.


{\bf{2.}}\;The following theorem corresponds to representation (\ref{predst1}).
\begin{teo}
\label{predstvar0}
{\it Let~$f(\cdot)$~be a function of bounded variation. Then}
\begin{equation}
\label{var22}
\bigvee\lt_a^b(f)=\bigvee\lt_a^b(f_c)+\bigvee\lt_a^b(f_d).
\end{equation}
\end{teo}

\doc\,
The finiteness of the variations of the components $f_c(\cdot)$ and $f_d(\cdot)$ is shown
in Corollary \ref{corol4}. Let us prove equality (\ref{var22}).

Consider first the case when $f(\cdot)$ has a single point
of discontinuity $c$, say, $c \in (a,b)$. 
Let $h>0$ be sufficiently small. By Theorem \ref{ch2thm5},
\begin{equation}
\label{sumvar}
\bigvee\limits_a^b(f)= \bigvee\limits_a^{c-h}(f)+\bigvee\limits_{c-h}^{c+h}
(f)+\bigvee\limits_{c+h}^b(f).
\end{equation}
Since $f(\cdot)$ is continuous on $[a,c)$ and $(c,b]$, the total variation
$f(\cdot)$ on $[a,c-h]$ and $[c+h,b]$ is the same as the total variation
$f_c(\cdot)$ on these segments. For this reason, the first and third terms
in~(\ref{sumvar}) tend, as $h\to 0+$, to
$\bigvee\limits_{a}^{c} (f_{c})$ and $\bigvee\limits_c^b (f_{c})$, respectively. It is also clear that the continuity of $f$ on the indicated intervals yields equality
$$
\underset{h\to 0+}{\lim}\bigvee\limits_{c-h}^{c+h}(f)=|\sigma_c^+(f)|+|\sigma_c^-(f)|= \bigvee\limits_a^b(f_d).
$$
The required equality is now obtained from (\ref{sumvar}) by taking $h\to 0+$.

Let $T(f)=\{c_1,c_2,\ldots,c_n\}$~ be a finite set. Let us sub-divide interval $[a,b]$
so that each sub-interval contains one point of discontinuity.
Applying Theorem \ref{ch2thm5} and what was proved above, we obtain the sought equality.

Finally, let $T(f)=\{c_k\}_{k=1}^{\infty}$.
Define $g_n(\cdot)$ as in (\ref{vtochnepr}) or in (\ref{vtochrazr}), respectively,
and put
$\widetilde{g}_n(t)=f_c(t)+g_{n}(t)$.
Then
$$
f(t)-\widetilde{g}_n(t)=f_d(t)-g_{n}(t)=\sum\limits_{k=n+1}^{\infty}\sigma_k(f)\rho_{c_k}(t),
$$
$$
\bigvee\limits_a^b\bigl(f-\widetilde{g}_n\bigr)=\sum\limits_{k=n+1}^{\infty}|\sigma_{c_k}^+(f)|+ |\sigma_{c_k}^-(f)|\to 0
$$
as $n\to \infty$ (as the remainder of a convergent series). Therefore,
\begin{equation}
\label{schod1}
\left|\bigvee\limits_a^b(f)-\bigvee\limits_a^b(\widetilde{g}_n)\right|\leqslant \bigvee\limits_a^b(f-\widetilde{g}_n)\to 0\quad (n\to \infty).
\end{equation}
$\bigl($Inequality in (\ref{schod1}) is a consequence of inequality $\bigl||\alpha|-|\beta|\bigr|\leqslant |\alpha-\beta|\bigr)$.
Thus,
\begin{equation}
\label{schod2}
\bigvee\limits_a^b(\widetilde{g}_n)\to \bigvee\limits_a^b(f)\;\;(n\to\infty).
\end{equation}
Moreover (see {\bf{5.2.1}}),
$$
\bigvee\limits_a^b(g_{n})=
\sum\limits_{k=1}^n\bigl(|\sigma_{c_k}^+(f)|+|\sigma_{c_k}^-(f )|\bigr)\to \bigvee\limits_a^b(f_d)\quad (n\to\infty).
$$
Since $g_{n}(\cdot)$ has finitely many points of discontinuity, we have
$$
\bigvee\limits_a^b(\widetilde{g}_{n})=\bigvee\limits_a^b(f_{c})+\bigvee \limits_a^b(g_{n}).
$$
The required equality now follows by taking $n\to\infty.$
\hfill $\square$


{\bf{3.}}\;The following question arises in connection with the proof of Theorem \ref{predstvar0}. 
Suppose that a sequence of functions of bounded variation $\{f_n(\cdot)\}_{n=1}^\infty$  
converges to a function of bounded variation~$f(\cdot).$ Does $\bigvee\limits_a^b(f_n)\;always$
converge to $\bigvee\limits_a^b(f)?$

Consider the following example.
\begin{example}
\label{pila}
\begin{multline*}
f_n(t)=\frac{1}{2^{n-1}}-\Big|t-\frac{1}{2^{n-1}}\Big|\;\;\text{when }\;\;t\in\Big[0,\,\frac{1}{2^{n-2}}\Big],\\
f_n(t)=f_n(t-\frac{1}{2^{n-1}})\;\;\text{when}\;\;t\in \Big(\frac{1}{2 ^{n-2}},\,2\Big],\;n=1,2,\ldots
\end{multline*}
\end{example}
Here
$
\max\lt_{[0,\,2]} f_n(t)=\frac{1}{2^{n-1}} \to 0,  \text{i.e.}\
f_n(t) \rightrightarrows f(t) \equiv 0,
$\;
but $$\bigvee\limits_0^2 (f_n)=2,
n=1,2,\ldots,\quad \bigvee\limits_0^2 (f)=0.$$

However, the following statement is true.
\begin{teo}
\label{shvar}
{\it Let $f_n(t)\to f(t)\;\;(n\to\infty,\;t\in [a,b])$
 and $\bigvee\limits_{a}^{b}(f_n)\leqslant \bigvee\limits_{a}^{b}(f)\;\;(n\in\mathbb N).$ Then}
 $\bigvee\limits_{a}^{b}(f_n)\to \bigvee\limits_{a}^{b}(f)\;\;(n\to\infty).$
\end{teo}
\doc\,
For every $\delta>0$ there is a partition $\tau_{\delta}=\{s_k\}_{k=1}^m$ such that
$v_{\tau}(f)\doteq\sum\limits_{k=1}^m|f(s_k)-f(s_{k-1})|>\bigvee\limits_{a}^{b} (f)-\delta.$
We fix such a partition and consider $v_{\tau}(f_n).$
Due to the pointwise convergence, $v_{\tau}(f_n)\to v_{\tau}(f)$ as $n\to \infty$
(we passed to the limit in a finite sum for each fixed $s_k$).
There is $N_0\in\mathbb N$ such that, for all $n>N_0$, $v_{\tau}(f_n)>\bigvee\limits_{a}^{b}(f)-\delta. $
So, for all
$n>N_0\;\;\;\bigvee\limits_{a}^{b}(f)-\delta<v_{\tau}(f_n)\leqslant \bigvee\limits_{a}^{b}( f_n) \leqslant \bigvee\limits_{a}^{b}(f).\;$
Since $\delta$ is arbitrary, this yields the assertion of the theorem.
\hfill $\square$


\medskip

{\bf{4.}}\;
Let us consider Example \ref{Qraz} in greater detail (see Remark \ref{zampr1}).
Let
$$
T=\{t_1,t_2,\ldots\}\subset [a.b]\;\;\bigl(T=\{t_1,t_2,\ldots, t_n\}\subset [a.b]\bigr)
$$
be an at most countable set, and let numbers $\alpha_k>0$ be such that the series
$\sum\lt_{k=1}^\infty \alpha_k$ converges. Let us put
\begin{equation}
\label{Qraz10}
f_1(a)=0,\quad f_1(t)=\sum\lt_{t_k<t}\alpha_k\;\;\;(t>a),\quad f_2(t)=\sum\lt_{a <t_k\leqslant t}\alpha_k.
\end{equation}
As was mentioned in \textbf{\ref{zampr1}}, an increasing jump function $f_1(\cdot)\;\;\bigl(f_2(\cdot)\bigr)$ is left continuous on
$[a,b)\;\;\bigl($right continuous of $[a,b]\bigr).$ Let us generalize this example.
\begin{example}
\label{ObobQraz}
Let real numbers $\alpha_k$ be such that the series
$\sum\lt_{k=1}^\infty \alpha_k$ converges absolutely. We can represent
$\alpha_k=\alpha_k^+-\alpha_k^-,\;\alpha_k^+>0,\;\alpha_k^->0 $ so that the series
$\sum\lt_{k=1}^\infty \alpha_k^+$ and $\sum\lt_{k=1}^\infty \alpha_k^-$ converge.
Let us define as in (\ref{Qraz10}) the functions $f_{i+}(\cdot)$ and $f_{i-}(\cdot)\;\;(i=1,2)$ 
using $\alpha_k ^+$ and $\alpha_k^-$, respectively, and
let us put $f_i(t)\doteq f_{i+}(t)-f_{i-}(t)\;\;(i=1,2).$ Then the jump function $f_1(\cdot)\;\; \bigl(f_2(\cdot)\bigr)$ is left continuous on
$[a,b)\;\;\bigl($right continuous $[a,b]\bigr).$
\end{example}

Let $f(\cdot)$~ be a function of bounded variation.
We define
\begin{equation}
\label{-30}
f_+(t)\doteq \sum\lt_{s\in T(f),\;a<s\leqslant t}\sigma_s^-(f),\qquad f_-(t)\doteq f(t) -f_+(t).
\end{equation}
Then, clearly, $T(f_+)=T(f_-)=T(f)$. Function $f_+(\cdot)$ is a jump function whose jumps
are equal to the left jumps of function $f(\cdot)$:
$$
\sigma_{c_k}(f_+)=
\sigma_{c_k}^-(f)\quad (k=1,2,\ldots).$$
According to Example \ref{ObobQraz}, function $f_+(\cdot)$ is right continuous.

Using successively representations (\ref{-30}), (\ref{predst1}) and (\ref{ch17}), we obtain
\begin{multline*}
f_-(t)=f(t)-f_+(t)=f_d(t)+f_c(t)-f_+(t)=\\=
\sigma_a(f)+\sum\lt_{s\in T(f),\;s<t}\sigma_s(f)+\sigma_{t}^-(f)+f_c(t)-
\sum\lt_{s\in T(f),\;a<s\leqslant t}\sigma_s^-(f)=\\=
f_c(t)+\sum\lt_{s\in T(f),\;s<t}\sigma_s^+(f).
\end{multline*}
In view of what was said in Example \ref{ObobQraz}, $f_-(\cdot)$ is left continuous.

From definitions (\ref{-30}) we obtain 
\begin{equation}
\label{-10}
f(t)=f_+(t)+f_-(t).
\end{equation}
Thus, we arrive at the following assertion.
\begin{teo}
\label{predstvar01}
{\it A function of bounded variation can be represented as the sum of left-continuous and right-continuous functions so that}
\begin{equation}
\label{-20}
\bigvee\lt_a^b(f)=\bigvee\lt_a^b(f_+)+\bigvee\lt_a^b(f_-).
\end{equation}
\end{teo}

\doc\,
It remains to prove equality (\ref{-20}).

Let us represent $f_-(\cdot)$ as a sum of discrete and continuous parts according to (\ref{predst1}):
$$
f_-(t)=\bigl(f_{-}\bigr)_{d}(t)+\bigl(f_{-}\bigr)_{c}(t).
$$
Then, according to (\ref{-10}) $\bigl($since $\bigl(f_+\bigr)_d(t)=f_+(t)\bigr)$,
$$
f(t)=f_+(t)+\bigl(f_{-}\bigr)_{d}(t)+\bigl(f_{-}\bigr)_{c}(t)=f_d(t )+\bigl(f_{-}\bigr)_{c}(t),
$$
i.e. $f_c(t)=\bigl(f_{-}\bigr)_{c}(t).$ Now, by Theorem \ref{predstvar0}, (\ref{-20}) is true.
\hfill $\square$


\bigskip

\poonkt{Banach space ${\bf{H}}[a,b]$\;}\phantom{0123456789012345678901234567890} \

{\bf{1.}}
Consider the vector space ${\bf{H}}[a,b]$ of functions 
$x\in {\bf{BV}}[a,b]$ admitting representation (see (\ref{predst1fogrvar}))
\begin{equation}
\label{funskach1}
x(a)=0,\quad x(t)=\sum\limits_{s\in T(x)}\sigma_{s}(x)\mathfrak h_s(t)+\sigma_{t}^-( x)\quad (t>a),
\end{equation}
where we define a norm by the formula
\begin{equation}
\label{normH}
\|x\|=\sum\lt_{s\in T(x)}\;\bigl(|\sigma^-_{s}(x)|+|\sigma^+_{s}(x) |\bigr)\;\left(=\bigvee\lt_a^b(x)\right).
\end{equation}

The normed space ${\bf{H}}[a,b]\;$ {\it is isometrically isomorphic to the space of sequences} $\mathbf{l}_1.$

Indeed, the isometric isomorphism $\Phi:{\bf{H}}[a,b]\to \mathbf{l}$ is defined by \;\;$\bigl(T(x)=\{c_1,c_2,\ldots\}\bigr)$
\begin{equation}
\label{izometrizomorph}
\Phi(x)=\bigl(\sigma_{c_1}^-(x),\,\sigma_{c_1}^+(x),\,\sigma_{c_2}^-(x),\,\sigma_ {c_2}^+(x),\ldots\bigr)\;\;
\bigl(\in \mathbf{l}\bigr).
\end{equation}
The bijectivity of the mapping $\Phi$ is not violated by the fact that the same sequence from $l$ corresponds to an infinite
set of jump functions having the same jumps possibly at different points of discontinuity since the ``identifier'' $T(x)$
always singles out the function we need. The existence of an isometric isomorphism ${\bf{H}}[a,b]\simeq \mathbf{l}$ yields completeness of space ${\bf{H}}[a,b]$.
Thus, the Banach space ${\bf{H}}[a,b]$ consists of all possible jump functions
and is a subspace of ${\bf{BV}}[a,b],$ and
\begin{equation}
\label{razlBVvprsum}
{\bf{BV}}={\bf{CBV}}[a,b]\dotplus {\bf{H}}[a,b] .
\end{equation}


{\bf{2.}}
It was shown in subsection {\bf{5.2.2}} that the convergence of total variations (\ref{schod2}) is ensured by the convergence
\begin{equation}
\label{schod3}
\bigvee\limits_a^b(f-\widetilde{g}_n)\to 0\quad (n\to\infty)\quad \bigl(\text{see}\;\;(\ref{schod1}) \bigr).
\end{equation}
Example \ref{pila} shows that {\it the uniform convergence $\widetilde{g}_n(t)\rightrightarrows f(t)\;\;(n\to\infty)$ is not sufficient to ensure convergence
$(\ref{schod2})$.}

{\it The uniform convergence is not sufficient even for finiteness of the total variation of the limit function.}

Let us construct the corresponding example by departing from Example \ref{CneBV1}.
\begin{example}
\label{{CneBV2}}
{\it Let
$\mathfrak{c}(t)=t\,\cos\,\frac{\pi}{2t}$ for $0<|t|\leqslant 1,\;\;\mathfrak{c}(0)=0 ,$ \\
$g_n(t)=\left\{
\begin{array}{lcr}
0\;\text{for}\;\;|t|<\frac{1}{2n-1},\\
\mathfrak{c}(t)\;\text{for}\;\;\frac{1}{2n-1}\leqslant |t|\leqslant 1,\\
\end{array}
\;\;n\in\mathbb N.
\right.
$}
\end{example}
Then
$g_n(\cdot) \;(n\in\mathbb N)$\; have finite total variations since these are piecewise monotone functions. According to
{\bf{4.2}} (see example~\ref{CneBV1})\;\;\;$\bigvee\limits_{-1}^1(\mathfrak{c})=\linebreak=
+\infty,$ and for all $t\in [-1,\,1] \;\;|g_n(t)-\mathfrak{c}(t)|\leqslant \frac{1}{2n-1}\to 0
(n\to\infty),$ i.e. $g_n(t)\rightrightarrows \mathfrak{c}(t)\;\;(n\to\infty).$


\bigskip

\poonkt{Space ${\bf{CH}}[a,b]\;$} \phantom{01234567890123456789} \

Put
\begin{equation}
\label{predstgot0}
{\bf{CH}}[a,b]\doteq {\bf{C}}[a,b]\dot + {\bf{H}}[a,b].
\end{equation}
The functions $x \in {\bf{CH}}[a,b]$ are thus uniquely representable in the form (see \ref{predst1fogrvar})
$\bigl($see also (\ref{predst1})$\bigr)$:
\begin{equation}
\label{predstgot}
x(t)=x_{\gd}(t)+x_{\gc}(t),
\end{equation}
where $x_{\gd}\in {\bf{H}}[a,b],$ and continuous function $x_{\gc}(\cdot)$ can now have {\it infinite total variation}.
Let us also put
$$
\|x\|_{{\bf{CH}}}\doteq \underset{t\in [a,b]}{\sup}|x(t)|+\sum\limits_{t\in T( x)}|\sigma_t(x)|;
$$
${\bf{CH}}$~ is a Banach space with respect to this norm (see  {\bf{2.3}}).
Thus, {\it we have strict inclusions}
\begin{equation}
\label{vkl2_}
\underline{{\bf{H}}[a,b]\subset{\bf{BV}}[a,b]}\subset{\bf{CH}}[a,b]
\end{equation}
where {\it the underlined inclusion is an embedding of Banach spaces.}

We obtain from  Theorem \ref{predstvar01} and the definition of space ${\bf{CH}}$ the following result.
\begin{sle}
\label{predstvar01+}
{\it A function $f\in{\bf{CH}}$ can be represented as the sum of left-continuous and right-continuous functions so that}
\begin{equation}
\label{-20+}
\bigvee\lt_a^b(f_{\gd})=\bigvee\lt_a^b\bigl(f_{\gd}\bigr)_{+}+\bigvee\lt_a^b\bigl(f_{\gd} \bigr)_{-}.
\end{equation}
\end{sle}

In the representations (\ref{predstgot0}),\,(\ref{predstgot}), the discrete part of a function 
is written explicitly $\bigl($see\,(\ref{funskach1}), and let us recall that the series in (\ref{funskach1}) converges absolutely and uniformly$\bigr).$
According to the discussion in subsection {\bf{5.3.1}},
\begin{equation}
\label{hcl}
{\bf{CH}}[a,b]\simeq {\bf{C}}[a,b]\dot+{\bf{l}}.
\end{equation}
Therefore, one can define the discrete part of a function by specifying the set of points of discontinuity $T(x)=\{c_1,\,c_2,\ldots\}$ and the sequence of jumps $\bigl($see (\ref{izometrizomorph})$\bigr):$
\begin{equation}
\label{hcls}
\sigma(x)\doteq\bigl(\sigma_{c_1}^-(x),\,\sigma_{c_1}^+(x),\,\sigma_{c_2}^-(x),\,\sigma_{c_2}^+(x),\ldots\bigr)\in {\bf{l}}.
\end{equation}
At the same time, one should keep in mind that by specifying a countable set $T(x)$ one already fixed the method 
of enumeration of its elements (we will say:
{\it particularly numbered set}).
If we interchange two elements of a 
particularly numbered set, we obtain {\it another} 
particularly numbered set.
In what follows, we will assume that the notation $T(x)$ always means 
particularly numbered set.

Thus, each function from ${\bf{CH}}[a,b]$ defines a triple
\begin{equation}
\label{hcltr1}
\bigl(x_{\gc}(\cdot),\,T(x),\,\sigma(x)\bigr).
\end{equation}
Conversely, each triple $\bigl(y(\cdot),\,T,\,a\bigr),$ where $y\in {\bf{C}}[a,b],\;\;T= \{c_1,c_2,\ldots\}$ is at most countable
specifically enumerated set
$a=(a_1,a_2,\ldots)\in {\bf{l}}$ defines a function $x\in{\bf{CH}}[a,b]$ for which
\begin{equation}
\label{hcltr2}
x_{\gc}(t)=y(t),\;\; T(x)=T,\;\; \sigma_{c_j}^-(x)=a_{2j-1},\;\sigma_{c_1}^+(x)=a_{2j},\;\;\;j=1,2,\ldots
\end{equation}


\bigskip

\poonkt{Preservation of convergence\;} \phantom{01234567890123456789} \

{\bf{1.}}\,
{\bf{Let}} $g(\cdot),\,g_n(\cdot)\in {\bf{CH}}[a,b]$ $(n=1,2,\ldots)$. It will be convenient to assume that all these functions have common set
of points of discontinuity $T\;\Bigl($one can always put $T\doteq T(g)\bigcup\left(\bigcup\limits_{k=1}^\infty T(g_n)\right)\Bigr).$
Assume that convergence
\begin{equation}
\label{tochsx}
g_n(t)\to g(t),\quad\text{with}\;\; n\to\infty\qquad (t\in [a,b])
\end{equation}
takes place.
Consider the representation (\ref{predstgot})
\begin{equation}
\label{predst00}
g(t)=g_{\gd}(t)+g_{\gc}(t)
\end{equation}
of a function from {\bf{CH}}[a,b] \;(recall that functions of bounded variation have the same representation, see (\ref{predst1})).

Note that the pointwise convergence (\ref{tochsx}), generally speaking, does not entail pointwise convergence 
for the components $g_{\gd}$+$g_{\gc}(t)$.
In this regard, let us consider the following example.
\begin{example}
\label{nesh}
Let
$p_n(t),\,n\in\mathbb N$~be continuous functions of bounded variation,
$p_n(a)=0;\;p\in {\bf{H}},\;p(a)=0,\;\;p_n(t)\to p(t)$
$(n\to~\infty);\;\;q_n(t),\linebreak
n\in\mathbb N$~ are
``step'' functions,\;$q_n(a)=0,
q_n(t)\to q(t)\;\;
(n\to~ \infty),\;q$~ is a continuous function of bounded variation, $q(a)=0$. Further,
$$
g_n(t)\doteq p_n(t)+
q_n(t),\;\;g(t)\doteq p(t)+q(t).$$ Then $g_n(t)\to g(t), $
however, for $n\to\infty$
$$
(g_{n})_{\gc}(t)=p_n(t)\to p(t)=g_{\gd}(t),\;\;\;(g_{n})_{\gd}(t)=q_{n}(t)\to q(t)=g_{\gc}(t).
$$
\end{example}

Let us indicate a specific example from the previous ``family'' of sequences.

Let $[a,b]=[0,1],$ and let sequence $\{c_k\}_{k=0}^{\infty}$ be such that
$c_0=0,\linebreak
_{k+1}>c_k,\;c_k\to 1\;(k\to \infty)$. Put \\
$p_n(t)=\left\{
\begin{array}{lcr}
c_1\left(\frac{t}{c_1}\right)^n\;\text{for}\;t\in [0,c_1], \\
c_{k-1}+(c_k-c_{k-1})\left(\frac{t-c_{k-1}}{c_k-c_{k-1}}\right)^n\;for\;t\in (c_{k-1},\,c_k], k=2,\ldots,n,\\
c_{n}+(1-c_n)\left(\frac{t-c_n}{1-c_n}\right)^n\;for\;t\in (c_{n},\, 1]
\end{array}
\right.
$\\
($p_n$~ are continuous strictly increasing functions). Then \\
$p(t)=\left\{
\begin{array}{lcr}
c_{k-1}\;for\;t\in [c_{k-1},\,c_k), k=1,2,\ldots,\\
1\;for\;t=1;
\end{array}
\right.
\;\;p
$~ is an increasing (in non-strict sense) function from ${\bf{H}},\;T(p)=\{c_1,c_2,\ldots\}$.
Let further\;
$q(t)=t(1-t),\;q_n(t)=q\left(\frac{k}{n}\right),\;t\in\left[\frac{k-1 }{n},\,\frac{k}{n}\right),\;
k=1,2,\ldots,n,\;q_n(1)=0;\linebreak
n=1,2\ldots$ Then we even have $q_n(t)\rightrightarrows q(t)\;(n\to\infty).$

We say that in the representation (\ref{predstgot}) {\it convergence is preserved} if convergence (\ref{tochsx})
implies convergence
\begin{equation}
\label{sohrsch}
\bigl(g_{n}\bigr)_{\gd}(t)\to g_{\gd}(t),\qquad\bigl(g_{n}\bigr)_{\gc}(t)\to g_{\gc}(t)\quad (n\to\infty).
\end{equation}
When needed, we will specify the type of preserved convergence.
Note also that this term refers to specific sequences
$$
\bigl\{g_n(\cdot)\bigr\}_{n=1}^\infty,\qquad \bigl\{\bigl(g_{n}\bigr)_{\gc}(\cdot)\bigr \}_{n=1}^\infty,\qquad
\bigl\{\bigl(g_n\bigr)_{\gd}(\cdot)\bigr\}_{n=1}^\infty.
$$


\newpage

{\bf{2.\;Uniform convergence of sequence $\{g_n(\cdot)\}_{n=1}^\infty$}}

\begin{lem}
\label{unisochrSM}
{\it If convergence in $(\ref{tochsx})$ is uniform on $[a,b],$ then convergence of right jumps \;\;$\bigl($left jumps,\;
jumps$\bigr)$ is also uniform on $T\;\;([a,b]):$ } 
$$
\sigma_t^+(g_n)\;\rightrightarrows\;\sigma_t^+(g),\;\;\;\sigma_t^-(g_n)\;\rightrightarrows\;\sigma_t^-(g),
\;\;\;\sigma_t(g_n)\;\rightrightarrows\;\sigma_t(g)\;\;\;(n\to\infty).
$$
\end{lem}
\doc\,
Let $g_n(t)\rightrightarrows g(t)$ on $[a,b].$ This means:
$$
\bigl(\forall\,\varepsilon>0\bigr)\quad\bigl(\exists\,n_0\in\mathbb N\bigr)\quad\bigl(\forall\,n\geqslant n_0\bigr)\quad
\bigl(\forall\,t\in [a,b]\bigr)\quad \bigl|g_n(t)-g(t)\bigr|<\dfrac{\varepsilon}{2}.
$$
Let us show that the sequence of right jumps also converges uniformly on $T\;\;([a,b]).$
\begin{multline*}
\bigl|\sigma_t^+(g_n)-\sigma_t^+(g)\bigl|=\bigl|(g_n(t+)-g_n(t)-g(t+)+g(t)\bigl|=\\=
\bigl|(g_n(t+)-g(t+))-(g_n(t)-g(t))\bigl|\leqslant\bigl|g_n(t+)-g(t+)\bigl|+\bigl| g_n(t)-g(t)\bigl|
\end{multline*}
The smallness of the second term is evident from the inequality in the previous line.
Since $g_n(t+)=\underset{s\to t+}{lim}\;g_n(s),$, we have by Theorem \ref{SM}
$$
\underset{n\to\infty}{\lim}\;g_n(t+)=
\underset{n\to\infty}{\lim}\;\underset{s\to t+}{\lim}\;g_n(s)=
\underset{s\to t+}{\lim}\;\underset{n\to\infty}{\lim}\;g_n(s)=\underset{s\to t+}{\lim}\;g(s)= g(t+).
$$
In detail, there is number $n_1$ such that for $n>n_1$ the first term will also be less than $\dfrac{\varepsilon}{2}.$
Let $n>max\,\{n_0,\,n_1\}.$ Then $\;\bigl|\sigma_t^+(g_n)-\sigma_t^+(g)\bigl|<\varepsilon,\;\;t\in T$
(it can also be assumed that $\;t\in [a,b]).$
Therefore,
\begin{equation}
\label{unisk0}
\sigma_t^+(g_n)\rightrightarrows \sigma_t^+(g)\;\;\;\text{on}\;\; T\;\;(\text{to}\;\; [a,b]).
\end{equation}
The uniform convergence of left jumps (jumps) is proved in exactly the same way.
\hfill $\square$

However, the uniform convergence of the sequence of jumps is not sufficient for the preservation of convergence.
Consider the following example (see subsection {\bf{1.3.3}} for $p=1$).
\begin{example}
\label{ex1.3.3}
Let $[a,b]=[0,1],\;c_k=1-\dfrac{1}{k},\;k=1,2,\ldots$
$$
\alpha_{kn}=\dfrac{1}{n}\;\;\text{for}\;\;k\leqslant n,\;\;\alpha_{kn}=0\;\;\text{ for}\;\;k>n,
$$
$$
\bigl(g_n\bigr)_{\gd}(t)=g_n(t)=\sum\limits_{c_k<t}\alpha_{kn},\qquad \bigl(g_n\bigr)_{\gc} (t)=0\;\;\;\;(t\in [0,1]).
$$
\end{example}
Here
$$
\sigma^{-}_{c_k}(g_n)=0,\;k=1,2,\ldots\,,0\leqslant\sigma^{+}_{c_k}(g_n)=\sigma_{c_k }(g_n)\leqslant\dfrac{1}{n}\to 0\;\;(n\to\infty)
$$
(the sequence of jumps tends to zero uniformly with respect to $k$), $g_{\gd}(t)=0\;(t\in [a,b]).$ However, $(g_n)_{\gd}( b)=1,$
i.e. $g_{\gd}(t)\nrightarrow g_{\gd}(t).$

\begin{teo}
\label{unisochrSMt}
{\it Let convergence in\;$(\ref{tochsx})$ be uniform on $[a,b],$ and let the series
$\sum\limits_{k=1}^\infty \bigl|\sigma_{c_k}(g_n)\bigr|$
converge uniformly in $n\in\mathbb N.$ Then we have the preservation of $($uniform$)$ convergence
in the representation} \;$(\ref{predstgot}).$
\end{teo}
\doc\,
It suffices to show the uniform convergence of sequence $\{(g_n)_{\gd}(\cdot)\}_{n=1}^\infty.$
By Lemma \ref{unisochrSM}, the sequence of jumps converges\linebreak 
uniformly in $k\in\mathbb N$
(right jumps, left jumps). By our assumptions, the series 
$\sum\limits_{k=1}^\infty \sigma_{c_k}(g_n)\mathfrak h_{c_k}(t)$
converges absolutely and uniformly in
$n\in\mathbb N$ and $t\in [a,b]$.
Therefore, one can pass to the limit as $n\to\infty$ under the series sign. Thus,
$\;(g_n)_{\gd}(t)\rightrightarrows g_{\gd}(t)$ on $\;[a,b].$
\hfill $\square$


\bigskip

{\bf{3.\;Quasi-uniform convergence}}

We begin with some auxiliary statements.
Let $\!\{c_{kn}\}_{k,n=1}^\infty$~be a sequence of  real numbers in two indices, $c:\mathbb N\times\mathbb N\!\to\!\mathbb R.$
\begin{lem}
\label{Quarz}
{\it Let $c_{kn}\to c_{k\,\bullet}$ for $n\to\infty\;\;(k=1,2,\ldots),$ assume that the series
\begin{equation}
\label{sochrsch5}
\sum\limits_{k=1}^\infty c_{kn}\;\;\;(n=1,2,\ldots)\qquad\text{and}\qquad\sum\limits_{k=1} ^\infty c_{k\,\bullet}
\end{equation}
converge. Then
\begin{equation}
\label{sochrsch6}
\underset{n\to\infty}{lim}\;\sum\limits_{k=1}^\infty c_{kn}=\sum\limits_{k=1}^\infty c_{k\,\bullet }
\end{equation}
if and only if the first series in $($\ref{sochrsch5}$)$ converges quasi-uniformly with respect to $n\in\mathbb N.$}
\end{lem}
\doc\,
Let us put, for $k=1,2,\ldots$,
\begin{multline*}
\sigma_k(0)\doteq c_{k\,\bullet},\;\;\sigma_k\left(\frac{1}{n}\right)\doteq c_{kn}\;\;(n=1 ,2,\ldots),\\
\sigma_k(t)\doteq c_{kn}+\theta_n(t)\bigl(c_{k,n+1}-c_{kn}\bigr)\;\;\text{for}\;\;t \in\left(\frac{1}{n+1},\,\frac{1}{n}\right),
\end{multline*}
where
$$
\theta_n(t)\doteq \dfrac{t-\frac{1}{n}}{\frac{1}{n+1}-\frac{1}{n}}=(n+1)(1 -nt)\quad (0\leqslant\theta_n(t)\leqslant 1).
$$
Clearly, $\sigma_k(\cdot)$ are continuous on $(0,1].$ Since
$$
\underset{t\to 0+}{\lim}\;\sigma_k(t)=\underset{n\to \infty}{lim}\;\Bigl(c_{kn}+\theta_n(t)\bigl( c_{k,n+1}-c_{kn}\bigr)\Bigr)=
c_{k\,\bullet}=\sigma_k(0),
$$
functions $\sigma_k(\cdot)$ are continuous on $[0,1].$

Let $S_{\bullet\,n}$ and $S$~ be the sums of the first and second  series 
in (\ref{sochrsch5}), respectively,
$$
\sigma(0)\doteq S,\;\;\sigma(t)\doteq\sum\limits_{k=1}^\infty \sigma_{k}(t),\;\;t\in(0 ,1].
$$
Since
$$
\sigma(t)=S_{\bullet\,n}+\theta_n(t)\bigl(S_{\bullet\,n+1}-S_{\bullet\,n}\bigr),\;\; \;t\in \left(\frac{1}{n_+1},\,\frac{1}{n}\right),
\;\;n=1,2,\ldots,
$$
function $\sigma(\cdot)$ is piecewise linear and continuous in $(0,1].$ In view of equalities
\begin{equation}
\label{sochrsch7}
\underset{t\to 0+}{\lim}\;\sigma(t)=\underset{n\to \infty}{\lim}\;S_{\bullet\,n}=S=\sigma(0)
\end{equation}
the validity of the limit relation (\ref{sochrsch6}) is equivalent to the continuity of $\sigma(\cdot)$ for $t=0,$ and hence continuity of this function on $[0,1].$ By the Arzela Theorem \cite[sect. 432]{ficht2}\;(see also {\bf{1.10.6.2$^{\circ}$}})
the validity of the limit relation (\ref{sochrsch6}) is equivalent to the quasi-uniform convergence of the series
$\sum\limits_{k=1}^\infty\sigma_k(t)$ to $[0,1].$

{\it The latter is equivalent to the quasi-uniform convergence of the series
\begin{equation}
\label{sochrsch8}
\sum\limits_{k=1}^\infty\sigma_k\left(\frac{1}{n}\right)=\sum\limits_{k=1}^\infty c_{kn}
\end{equation}
with respect to $n\in\mathbb N.$}

Indeed, since $\left\{\frac{1}{n}\right\}_{n=1}^\infty\subset [0,1],$ in one direction this assertion is obvious.

Let series (\ref{sochrsch8}) converge quasi-uniformly. This means (see subsection {\bf{1.10.6.2$^{\circ}$}},
here  we use the  formulation of quasi-uniform convergence in terms of series),
$$
\bigl(\forall\,\varepsilon>0\bigr)\;\;\bigl(\forall\,k_0\in\mathbb N\bigr)\;\;
\bigl(\exists\,p\in\mathbb N\bigr)\;\;\bigl(\exists\,k_1,\,\ldots,\,k_p\,(\geqslant k_0)\bigr)\;\ ;\bigl(\forall\,n\in\mathbb N\bigr)\;\;
$$
$$
\underset{i=1,\,\ldots,\,p}{min}\;\left|\sum\limits_{k=k_i}^\infty c_{kn}\right|<\varepsilon.
$$
Since $\sigma_k(\cdot)$ are piecewise linear,  we have for $t\in [0,1]$
$$
\left|\sum\limits_{k=k_i}^\infty\sigma_k(t)\right|\leqslant \left|\sum\limits_{k=k_i}^\infty c_{kn}\right|,\quad i=1,\,\ldots,\,p.
$$
That is, series $\sum\limits_{k=1}^\infty\sigma_k(t)$ also converges quasi-uniformly.
\hfill $\square$

It will be convenient for what follows to reformulate Lemma \ref{Quarz} in the language of sequences.
Let $\bigl\{C_{kn}\bigr\}_{k,n=1}^\infty$~ be a sequence of real numbers in two indices, $C:\mathbb N\times\mathbb N\to\mathbb R.$
\begin{lem}
\label{Quarzperef}
{\it Assume that for each $n\in\mathbb N$
\begin{equation}
\label{lembis}
C_{kn}\to C_{\bullet\,n}\quad (k\to\infty)
\end{equation}
and that for every $k\in\mathbb N$ there exists the limit
$C_{k\,\bullet}=\underset{n\to\infty}{\lim}\;C_{kn}.$ Then
the limits $\underset{n\to\infty}{\lim}\;C_{\bullet\,n}$ and $\underset{k\to\infty}{\lim}\;C_{k\,\bullet}$ exist and are equal, that is,
\begin{equation}
\label{sochrsch62}
\underset{n\to\infty}{\lim}\;\underset{k\to\infty}{\lim}\; C_{kn}=\underset{k\to\infty}{\lim}\;\underset{n\to\infty}{\lim}\;C_{kn},
\end{equation}
if and only if convergence in $($\ref{lembis}$)$ is quasi-uniform with respect to $n\in\mathbb N.$}
\end{lem}

Let $g(\cdot),\;g_n(\cdot)\in{\bf{CH}}[a,b]\;(n=1,2,\ldots).$
\begin{lem}
\label{Quarz1}
{\it We have
\begin{multline}
\label{sohrpm}
\underset{n\to\infty}{\lim}\;\underset{s\to t\pm}{\lim}\;g_n(s)=\underset{s\to t\pm}{\lim}\;\underset{n\to\infty}{\lim}\;g_n(s)\\
\Bigl(=\underset{s\to t\pm}{\lim}\;g(s)\;\;\bigl(t\in T(g)\bigr)\Bigr).
\end{multline}
if and only if convergence in $($\ref{tochsx}$)$ is quasi-uniform on\;$[a,b].$}
\end{lem}
\doc\,
Let $t\in [a,b)\;\;\bigl(t\in (a,b]\bigr),\;\;s_k\to t+\;\;\bigl(s_k\to t-\bigr)$ for $k\to\infty.$
In view of Lemma \ref{lembis}, it suffices to put
$C_{kn}\doteq g_n(s_k).$
\hfill $\square$
\begin{sle}
\label{sohrsk}
{\it
We have
\begin{equation}
\label{sochrsk1}
\underset{n\to\infty}{lim}\;\sigma_t^{\pm}(g_n)=\sigma_t^{\pm}(g),\quad
\underset{n\to\infty}{lim}\;\sigma_t(g_n)=\sigma_t(g)\;\;\;(t\in T(g),
\end{equation}
if and only if the convergence in $($\ref{tochsx}$)$ is quasi-
uniform on $[a,b]$.
In which case the convergence in $(\ref{sochrsk1})$ is also quasi-
uniform on\;$[a,b].$}
\end{sle}
\doc\,
We only need to prove the last statement.
Let $g_n(t)\to g(t)$ quasi-uniformly on $[a,b].$ This means that $g_n(t)\to g(t)$ for all $\;t\in [a, b]$ and
\begin{equation}
\label{shsost0}
\bigl(\forall\,\varepsilon>0\bigr)\;\bigl(\forall\,n_0\in\mathbb N\bigr)\;
\bigl(\exists\,p\in\mathbb N\bigr)\;\bigl(\exists\,n_1,\,\ldots,\,n_p\,(\geqslant n_0)\bigr)\;\bigl( \forall\,t\in [a,b]\bigr)
\end{equation}
$$
\underset{i=1,\,\ldots,\,p}{min}\;|g_{n_i}(t)-g(t)|<\dfrac{\varepsilon}{2}.
$$
Index $i_{0}$, where the minimum is attained, depends on $t,\;\;i_{0}=i_{0}(t).$ Let $t\in T$ and, for example, $s\to t+.$
For each such $s$, the minimum $a_{n_i}(s)\doteq |g_{n_i}(s)-g(s)|$ is attained at index $1\leqslant i_0(s)\leqslant p$;
and although this index may change when $s$ changes, the inequality 
$a_{n_i}(s)<\dfrac{\varepsilon}{2}$ continues to hold for all $s$ sufficiently close to $t.$ Therefore,
\begin{equation}
\label{shsost}
\underset{i=1,\,\ldots,\,p}{min}\;|g_{n_i}(t+)-g(t+)|\leqslant \dfrac{\varepsilon}{2},
\end{equation}
i.e. if (\ref{shsost0}) holds,
\begin{multline*}
\bigl|\sigma_t^+(g_n)-\sigma_t^+(g)\bigl|=\bigl|(g_n(t+)-g_n(t)-g(t+)+g(t)\bigl|=\\=
\bigl|(g_n(t+)-g(t+))-(g_n(t)-g(t))\bigl|\leqslant\bigl|g_n(t+)-g(t+)\bigl|+\bigl| g_n(t)-g(t)\bigl|<\varepsilon.
\end{multline*}
Consequently, convergence $\sigma_t^+(g_n)\to\sigma_t^+(g)$ is quasi-
uniform on $[a,b].$ The quasi-uniform convergence of
left jumps (jumps) is proved similarly.
\hfill $\square$

Let us now state the main result of this section.

\begin{teo}
\label{osntQU}
{\it Assume that the following conditions are satisfied:

a$)$\; $g_n(t)\to g(t)\;\;(n\to\infty)$ quasi-uniformly
on $[a,b];$

b$)$\;the series\;\;$\sum\limits_{k=1}^\infty\bigl|\sigma_{c_k}(g_n)\bigr|$\: converges
uniformly with respect to\;$n\in\mathbb N.$

Then the preservation of $($quasi-uniform$)$ convergence takes place
in the 
representation} \;$(\ref{predstgot}).$
\end{teo}
\doc\,
Condition a)\; of the theorem is equivalent to the assertions of Lemmas \ref{Quarz},\;\ref{Quarzperef},\;\ref{Quarz1} and Corollary \ref{sohrsk}.
Condition\;b)\; of the theorem implies absolute and uniform with respect to $n\in\mathbb N$ and $t\in [a,b]$ convergence of the series
$$
\sum\limits_{k=1}^\infty \sigma_{c_k}(g_n)\mathfrak h_{c_k}(t),
$$
that is, the possibility of passing to the limit under the sign of this series as $n\to\infty.$ Therefore,
$(g_n)_{\gd}(t)\to g_{\gd}(t))\;$ quasi-uniformly on $[a,b]$ and, hence,
$\;(g_n)_{\gc}(t)\to g_{\gc}(t)\;\;(t\in [a,b]).$ By Theorem \ref{Arzela} (Arzela) 
the last convergence is
quasi-uniform on $[a,b].$
\hfill $\square$


\bigskip

\poonkt{Continuous functions of bounded variation-2} \phantom{01234} \

{\bf{1.}}\;{\it If function~$f(\cdot)$ has finite total variation
on~$[a,\,b]$ and is continuous at point~$b\;\;($at point~$a),$ then
$\bigvee\lt_a^b(f)=\underset{t\to b-0}{\lim}\bigvee\lt_a^t(f)$ \\
$\left(\bigvee\lt_a^b(f)=\underset{t\to a+0}{\lim}\bigvee\lt_t^b(f)\right).$
These limits are infinite if and only if the
variation of $f(\cdot)$ is infinite.}

\doc\,
This assertion follows immediately  from Theorem \ref{ch2thm7} and the properties of function $f_{\pi}(\cdot).$
\hfill $\square$


\medskip

{\bf{2.}}
Let $\om_{\tau k}(f)\doteq\underset{t\in [t_{k-1},t_k]}{\sup}\,f(t)-\underset{t\in [t_ {k-1},t_k]}{\inf}\,f(t)$denote
the oscillation of function $f(\cdot)$ on the interval
$[t_{k-1},t_k]$ of a partition $\tau=\{t_k\}_{k=0}^n.$

{\it If $f(\cdot)$~ is a continuous function of bounded variation, then} 
$$
\bigvee\lt_a^b(f)=\underset{d(\tau)\to 0}{lim} v_\tau(f).
$$
\doc\,
After adding a point $\overline t$ to partition $\tau=\{t_k\}_{k=0}^n$
the sum~$v_\tau (f)$ can only increase. Since 
$$
\big| f(t_k)-f(t_{k-1}) \big| \ls\big| f(t_k)-f(\ol t) \big|+
\big| f(\ol t)-f(t_{k-1}) \big| \ls 2\om_{\tau k}(f),
$$
adding a point to the partition can increase $v_\tau(f)$
no more than by twice the oscillation of the function on the interval that contains this new point.

Let $\eps>0$ be arbitrary. Let us fix a partition $\tau^*=\{t_k^*\}_{k=0}^m$
such that the sum $v^* \doteq v_{\tau^*}(f)$ satisfies inequality
$$
v^* > \bigvee\lt_a^b (f)-\eps \doteq A.
$$
Since~$f(\cdot)$ 
is also uniformly continuous on~$[a,\,b]$, there exists $\dl>0$ such that
if $|t'-t''|<\dl,$ then $\big| f(t')-f(t'') \big|<\frac{v^*-A}{4m}.$
Let~$\tau$~ be an arbitrary partition whose diameter is $d(\tau)<\dl$
and $v \doteq v_\tau (f).$ Then $\om_{\tau k}<\frac{v^*-A}{4m}.$

Consider partition $\tau_0 \doteq \tau \cup \tau^*,$ $v_0=v_{\tau_0} (f).$
Since we have added at most~$m$ points to~$\tau$, we have
$v_0-v < 2\cdot \frac{v^*-A}{4m} \cdot m= \linebreak =\frac{v^*-A}{2}.$
Therefore, $$v>v_0-\frac{v^*-A}{2} \gs v^*-\frac{v^*-A}{2}=\frac{v^*+A}{2}> A
\;\;(=\bigvee\lt_a^b (f)-\eps).$$ Thus, for $d(\tau)<\dl$\;\;\;\;
$\bigvee\lt_a^b (f)-\eps<v \ls \bigvee\lt_a^b (f),$
as claimed.
\hfill $\square$


{\bf{3.}}
{\it Let $f\in {\bf{BV}}[a,b].$ Then $\bigvee\lt_a^b(f)=\underset{\tau}{\sup}\sum\lt_{k =1}^n\om_{\tau k}(f).$}

\doc\,
Given an arbitrary partition $\tau=\{t_k\}_{k=0}^n$ of interval $[a,\,b]$
let us denote $\Om_\tau(f)=\sum\lt_{k=1}^n \om_{\tau k}(f)$. Clearly,
$$
v_\tau (f)=\!\sum\lt_{k=1}^n\! |f(t_k)--f(t_{k-1})| \ls \sum\lt_{k=1}^n\om_{\tau k}(f)=\Om_\tau (f).
$$ 
Hence $\bigvee\lt_a^b (f) \ls \sup\lt_\tau\Om_\tau (f).$
It follows from the properties of the infimum and the supremum that for any $\eps>0$
and any~$k$ ${(1\ls k\ls n)}$ there are points $\xi_k$ and $\eta_k,$
belonging to interval $(\Dl_k)\doteq [t_{k-1},\,t_k],$ such that
$f(\xi_k)>\sup\lt_{(\Dl_k)} f(t)-\fracs{\eps}{2n},$
${f(\eta_k) < \inf\lt_{(\Dl_k)} f(t)+\fracs{\eps}{2n}}$.
Thus,
$$
\gathered
\Om_\tau (f)=\sum\lt_{k=1}^n \Big(\sup\lt_{(\Dl_k)} f(t)-
\inf\lt_{(\Dl_k)} f(t)\Big) \ls \sum\lt_{k=1}^n \Big( f(\xi_k)-f(\eta_k)
+\frac{\eps}{n} \Big)= \\
=\sum\lt_{k=1}^n \Big| f(\xi_k)-f(\eta_k) \Big|+\eps \ls v_{\tau'}(f)+\eps
\ls\bigvee\lt_a^b (f)+\eps,
\endgathered
$$
where the partition~$\tau'$ is obtained by adding to partition~$\tau$ all
points~$\xi_k$ and~$\eta_k$ ($k=1,\,\ldots,\,n$). Since $\eps$ was arbitrary,
we have $\Om_\tau (f) \ls \bigvee\lt_a^b (f).$ 
Together the opposite inequality proved earlier, we obtain the required.
\hfill $\square$

\medskip

{\it If $f(\cdot)$~ is a continuous function of bounded variation, then} 
$$
\bigvee\lt_a^b(f)=\underset{d(\tau)\to 0}{lim} \sum \limits_{k=1}^\infty\om_{\tau k}(f).
$$

\doc\,
It was shown above that for an arbitrary~$\eps$
$$
v_\tau (f) \ls \Om_\tau (f) \ls v_{\tau'}(f)+\eps,
$$
where partition~$\tau'$ was obtained by adding to~$\tau$ at most
$2n$ points. The assertion now follows from {\bf{2.}}
\hfill $\square$


\medskip

{\bf{4.}}
{\it Let~$f(\cdot)$ be continuous on~$[a,\,b],$ assume that the derivative~$f'(\cdot)$
exists everywhere except perhaps finitely many points of $[a,b]$, $f'(\cdot)$ is Riemann integrable in the proper sense.
Then~$f\in {\bf{BV}}[a,b],$ and} $\bigvee\lt_a^b(f)=\int\lt_a^b |f'(t)|dt.$

\doc\,
Consider partition $\tau=\{t_k\}_{k=1}^n$ containing all points where the derivative does not exist, so that~$f(\cdot)$ is differentiable on
all intervals~$(t_{k-1},\,t_k).$ This allows us to apply on each of these
intervals the Lagrange finite increment formula, so
$$
v_\tau (f)=\sum\lt_{k=1}^n \big| f(t_k)-f(t_{k-1}) \big|=\sum\lt_{k=1}^n
\big| f'(\xi_k) \big| (t_k-t_{k-1}),
$$
where~$\xi_k \in (t_{k-1},\,t_k)$ $(k=1,\,\ldots,\,n).$ The sast sum
is the corresponding to the partition~$\tau$ Riemann sum for function~$|f'(t)|$ on interval~$[a,\,b]$.
We obtain from the result in {\bf{2.}}
$$
\gathered
\bigvee\lt_a^b (f)=\underset{d(\tau)\to 0}{\lim}v_\tau (f)=
\underset{d(\tau)\to 0}{\lim}\sum\lt_{k=1}^n \big| f'(\xi_k) \big| (t_k-t_{k-1})=\int\lt_a^b \big| f'(t) \big|\,dt.\;\;\hfill \square
\endgathered
$$


\bigskip

\poonkt{Rectifiability of a Jordan curve}\phantom{01234567890123456789} \

The concept of total variation of a function was introduced in connection with the need to resolve the issue of 
existence of a finite length for a Jordan curve, i.e. a curve $x=\vfi(t),$ $y=\psi(t),$ $z=\chi(t),$ $\al\ls t\ls \bt,$ 
where functions $\vfi,\, \psi,\,\chi$ are assumed to be continuous on $[\al,\,\bt].$ 

\begin{teo}[Jordan's theorem]
\label{Jord}
{\it A Jordan curve is rectifiable if and only if the corresponding functions $\vfi(\cdot)$, $\psi(\cdot)$ and
$\chi(\cdot)$ have finite total variation on}~$[\al,\,\bt].$
\end{teo}
\doc\,
{\it{Necessity.}}
Let $\tau=\{t_k\}_{k=0}^n$~ be an arbitrary partition of interval $[\al,\,\bt],$ let $P_\tau$~ be the length of the broken line
inscribed in the curve corresponding to the partition~$\tau.$ Since the curve
is rectifiable, we have
\begin{align*}
P_\tau & =\!\sum\lt_{k=1}^n\!\! \sqrt{\big(\vfi(\nu_k)\!-\!
\vfi(\nu_{k-1})\big)^2\!+\!\big(\psi(\nu_k)\!-\!\psi(\nu_{k-1})\big)^ 2\!+\!
\big(\chi(\nu_k)\!-\!\chi(\nu_{k-1})\big)^2} \\
& \ls \ell,
\end{align*}
 where $\ell$~ is the length
of the curve. Clearly, 
\begin{align*}
&\big| \vfi(\nu_k)-\vfi(\nu_{k-1})\big| \\
&\ls \sqrt{\big(\vfi(\nu_k)\!-\!
\vfi(\nu_{k-1})\big)^2\!+\!\big(\psi(\nu_k)\!-\!\psi(\nu_{k-1})\big)^ 2\!+\!
\big(\chi(\nu_k)\!-\!\chi(\nu_{k-1})\big)^2},
\end{align*}
where
$$
v_\tau (\vfi)=\sum\lt_{k=1}^n \big| \vfi(\nu_k)-\vfi(\nu_{k-1}) \big| \ls
P_\tau\ls\ell,
$$
i.e. $\vfi(\cdot)$ has finite total variation on~$[\al,\,\bt]$ and
$\bigvee\lt_\al^\bt (\vfi) \ls \ell.$ We prove similarly the finiteness of the total variation of $\psi(\cdot)$ and $\chi(\cdot).$

{\it{Sufficiency.}} We first prove that for any $\{t_k\}_{k=1}^\infty$
$$
\sqrt{\sum\lt_{k=1}^m t_k^2} \ls \sum\lt_{k=1}^m |t_k|.
\eqno{(*)}
$$
Since, clearly,
$$
\Big( \sum\lt_{k=1}^n |t_k| \Big)^2=\sum\lt_{k=1}^m \sum\lt_{j=1}^m
|t_k| \cdot |t_j| \gs \sum\lt_{k=1}^m |t_k|^2,
$$
we have inequality~$(*)$.

Let $\vfi(\cdot),\,\psi(\cdot),\,\chi(\cdot)$~ be functions of bounded
variations. Given an arbitrary partition~$\tau$, we have, due to~$(*)$,
\begin{align*}
& P_\tau \\
&=\sum\lt_{k=1}^n\sqrt{\big(\vfi(\nu_k)\!-\!\vfi(\nu_{k-1})\big)^2\ !+\!
\big(\psi(\nu_k)\!-\!\psi(\nu_{k-1})\big)^2\!+\!\big(\chi(\nu_k)\!-\!
\chi(\nu_{k-1})\big)^2}  \\
& \ls \sum\lt_{k=1}^n \Big( \big|\vfi(\nu_k)\!-\!\vfi(\nu_{k-1})\big|\!+\!
\big|\psi(\nu_k)\!-\!\psi(\nu_{k-1})\big|\!+\!\big|\chi(\nu_k)\!-\!\chi( \nu_{k-1})\big|
\big) \\
& =v_\tau(\vfi)+v_\tau(\psi)+v_\tau(\chi) \ls \bigvee\lt_\al^\bt (\vfi)+
\bigvee\lt_\al^\bt (\psi)+\bigvee\lt_\al^\bt (\chi).
\end{align*}
Since~$P_\tau$ increases as we decrease the diameter of the partition,
the resulting inequality means that~$P_\tau$ has finite limit
$\ell \ls \bigvee\lt_\al^\bt (\vfi)+\bigvee\lt_\al^\bt (\psi)
+\bigvee\lt_\al^\bt (\chi).$
\hfill $\square$


\vskip 1cm
\begin{center}
{\Large\bf Exercises}
\end{center}

\up
\label{ch2upr6}
For what values of $c$ the total variation of function $f$ defined on interval~$[-1,\,1]$ by $f(t)=t-\sign t$
for~$0<|t|\ls 1,$ $f(0)=c$  is minimal?

\up
\label{ch2upr9}
Express the total variation of  function $F(t)=k\cdot f(t)+m$ in terms of the total
function variation of $f(\cdot).$

\up
\label{ch2upr11}
Let~$f(\cdot)$ and~$g(\cdot)$~ be functions of bounded variation. Prove
that $$h(t)\!=\!\max\{f(t),\,g(t)\}, \quad t\!\in\![a,\,b]$$ and 
$${m(t)=\min\{f(t),\,g(t)\}}, \quad t\in[a,\,b]$$ are also functions of bounded variation.

\up
\label{ch2upr30}
Prove that if~$f(\cdot)$ is a function of bounded variation on $[a,\,b],$
then function~$F(t)=\frac{1}{t-a}\int\lt_a^t f(s)ds$ for~$t>a,$ $F(a)=0$
has finite total variation.

\up
\label{ch2upr33}
Let $g(0)=0,$ $g\Big(\frac{1}{2n-1}\Big)=0,$ $g\Big(\frac{1}{2n}\Big)=1 ,$
$g(\cdot)$ is linear on each interval $\Big[\frac{1}{2n+1},\,
\frac{1}{2n}\Big]$ and $\Big[\frac{1}{2n},\,\frac{1}{2n-1}\Big],$
$n=1,\,2,\,\ldots$ Does $g(\cdot)$ have finite total variation on~$[0,\,1]$?

\up
\label{ch2upr34}
Let $g(0)=0,$ $g\Big(\frac{1}{2n-1}\Big)=0,$ $g\Big(\frac{1}{2n}\Big)=
\frac{1}{2n},$ $g(\cdot)$ is linear on each interval $\Big[
\frac{1}{2n+1},\,\frac{1}{2n}\Big]$
and $\Big[\frac{1}{2n},\,\frac{1}{2n-1}\Big],$
$n=1,\,2,\,\ldots$ Does $g(\cdot)$ have finite total variation on~$[0,\,1]$?

\up
\label{ch2upr35}
Let $g(0)=0,$ $g\Big(\frac{1}{2^n}\Big)=0,$ $g(t)=\frac{1}{2^n}$ in the middle of 
interval $\Big(\frac{1}{2^{n+1}},\,\frac{1}{2^n}\Big),$ $g(\cdot)$
is linear on intervals $\Big[\frac{1}{2^{n+1}},\,\frac{3}{2^{n+2}}\Big]$
and $\Big[\frac{3}{2^{n+2}},\,\frac{1}{2^n}\Big].$  Does $g(\cdot)$ have finite total variation on~$[0,\,1]$?

\up
\label{ch2upr36}
Let $\{t_n\}_{n=1}^\fy$ be a sequence of real numbers such that
${t_1=1},$ ${t_{n+1}<t_n},$ ${n=1,\,2,\,\ldots},$ ${t_n\to +0};$
$${g(0)=0}, \quad {g(t_n)=0}, \quad {g\Big(\! \frac{t_n+t_{n+1}}{2}\! \Big)=t_n},$$
$g(\cdot)$ is linear on each interval
$\Big[t_{n+1},\,\frac{t_n+t_{n+1}}{2}\Big]$
and $\Big[\frac{t_n+t_{n+1}}{2},\,t_{n+1}\Big].$ What are the conditions on
sequence~$\{t_n\}_{n=1}^\fy$ so that the total variation of ~$g(\cdot)$
is finite?

\up
\label{ch2upr38}
Show that  function $g(t)=(n+1)\Big(t-\frac{1}{n+1}\Big)$
on  ${\frac{1}{n+1}<t \ls \frac{1}{n}},$ $n=1,\,2,\,\ldots,$
$\ g(0)=0,$ has infinite total variation on~$[0,\,1].$

\up
\label{ch2upr39}
Let $\{t_n\}_{n=1}^\fy$ be the same sequence of real numbers as in problem \ref{ch2upr36}. Put
$g(t)=\frac{t_n}{t_n-t_{n+1}}(t-t_{n+1})$ for $t_{n+1}<t\ls t_n,$ \linebreak
$n=1,\,2,\,\ldots,$ $g(0)=0.$ What condition need to be imposed on $\{t_n\}_{n=1}^\fy$ so that $g(\cdot)$
had finite variation on~$[0,\,1]$?

\up
\label{ch2upr49}
Prove that function $f(t)=t^\al\sin\frac{1}{t^\bt}$ for $0\!<\!t\! \ls\! 1,$ $f(0)=0,$
$\bt\!>\!0$ has finite total variation if and only if
$\al>\bt.$

\up
\label{ch2upr46}
Let~$\vfi(\cdot)$ be proper Riemann integrable,
$f(t)=\linebreak=
\int\lt_a^t \vfi(s)ds,$ $t\in[a,\,b].$ Prove that
 $\bigvee\lt_a^b(f)=\int\lt_a^b |\vfi(s)|ds.$ Prove that
the same is true if~$\vfi(\cdot)$ is absolutely integrable on~$[a,\,b]$
in the improper sense.

\up
\label{ch2upr47}
Let~$f(\cdot)$ be continuously differentiable on~$[a,\,b].$ Prove that
$$
\frac{d}{dt}\bigvee\lt_a^t(f)=|f'(t)|.
$$

\up
\label{ch2upr42}
Let~$f(\cdot)$ be a function of bounded variation on~$[0,\,1].$ By the rectangles formula,
$\int\lt_0^1 f(t)dt \approx \frac 1n \sum\lt_{k=1}^n
f\Big( \frac kn \Big).$ Estimate the error in this approximate equality: show that
$$
\left| \int\lt_0^1 f(t)dt -\frac 1n\sum\lt_{k=1}^n f\Big(\frac kn\Big)\right|
\ls \frac 1n \bigvee\lt_0^1(f).
$$

\newpage

\begin{center}
\begin{Large}
{\bf{Chapter III. \; Riemann---Stieltjes integral}}
\end{Large}
\end{center}
\addcontentsline{toc}{section}{Chapter III. \; Riemann---Stieltjes integral}

\begin{center}
\begin{large}
\section{Definition of the Riemann--Stieltjes integral}
\end{large}
\end{center}

\hskip 1mm
\poonkt{Definitions}  \

Let us start with the following definition.

Let $f,\,g:[a,\,b]\to\R.$ Consider a partition $\tau=\{t_k\}_{k=0}^n$
of the interval $[a,\,b],$ $a=t_0<t_1<\ldots<t_n=b$. Let $d(\tau)=
\max\lt_{1\ls k\ls n}(t_k-t_{k-1})$~ be the diameter of the partition. In each of
the intervals $[t_{k-1},\,t_k]$ select a point $\xi_k$ and consider the sum
\begin{equation}
\label{ch31}
\gS_\tau \doteq \gS_\tau(f,\,g) \doteq \sum\lt_{k=1}^n f(\xi_k)\bigl(g(t_k)-g(t_{k-1} )\bigr).
\end{equation}
This sum is called {\it Stieltjes integral sum}. It 
resembles Riemann integral sum: instead of $\Dl t_k \doteq t_k-t_{k-1}$
we have increments $\Dl g_k\doteq g(t_k)-g(t_{k-1})$.
If $g(t)\equiv t$, then the Stieltjes sum becomes the Riemann sum.

We say that the sum $\gS_\tau(f,\,g)$ has limit as $D(\tau)\to 0$
if there exists a real number $J$ such that for every $\eps>0$ there is
$\dl>0$ such that for any partition $\tau$ whose diameter
$d(\tau)<\dl$ we have $|\gS_\tau(f,\,g)-J|<\eps,$
regardless of the choice of points $\xi_k.$ The notation for this limit is the usual one:
\begin{equation}
\label{ch32}
J=\lim\lt_{d(\tau)\to 0} \gS_\tau(f,\,g).
\end{equation}

It is easy to see that the limit~(\ref{ch32}) has all the properties
of the limit of a sequence of real number. In particular, the Bolzano-Cauchy principle of convergence
is formulated in the following form.

{\it For the limit~$($\ref{ch32}$)$ to exist, it is necessary and sufficient that for
any $\eps>0$ there exists $\dl>0$ such that for any partitions~$\tau',\,
\tau''$ inequalities $d(\tau')<\dl,$ $d(\tau'')<\dl$ imply inequality}
$|\gS_{\tau'}-\gS_{\tau''}|<\eps.$

If the indicated limit exists, then we say that  the function~$f(\cdot)$ is
integrable with respect to the function~$g(\cdot)$ in the sense of Stieltjes. This limit
is called the {\it Stieltjes integral} (more precisely, the Riemann~--- Stieltjes integral) of function $f(\cdot)$ with respect to function~$g(\cdot)$ over the
segment~$[a,\,b].$ Thus,
\begin{multline*}
(S)\!\int\lt_a^b f(t)\,dg(t)=(RS)\!\!\!\int\lt_a^b f(t)\,dg(t)=
\int\lt_a^b f(t)\,dg(t)\doteq \lim\lt_{d(\tau)\to 0} \gS_\tau(f,\,g).
\end{multline*}
In this case~$f(\cdot)$ is called the {\it integrated} function, and~$g(\cdot)$ is called
the {\it integrating} function.

Obviously, if $g(t)\equiv t$, then $(S)\!\int\lt_a^b f(t)\,dg(t)=
({R})\!\int\lt_a^b f(t)\,dt$. Therefore, the Stieltjes integral is a generalization
of the Riemann integral.

Note also that for $g(t)\equiv c= \const$ $(S)\!\int\lt_a^b f(t)\,dg(t)=0$
$\bigl($all sums (\ref{ch31}) are equal to~0$\bigr)$ for any partition, and for $f(t)={\cal K}=
\const$ 
$(S)\!\int\lt_a^b f(t)\,dg(t)={\cal K}(g(b)-g(a)),$
since for any partition~$\tau\quad{}$
\begin{multline*}
\gS_\tau(f,\,g)=\sum\lt_{k=1}^n {\cal K}(g(t_k)-
g(t_{k-1}))={\cal K}\sum\lt_{k=1}^n (g(t_k)-g(t_{k-1}))=\\
={\cal K}(g(t_n)-g(t_0))={\cal K}(g(b)-g(a)).
\end{multline*}

We warn the reader that the symbol  $dg(t)$ is not a differential:
function $g(t)$ is, generally speaking, not assumed to be differentiable. 
Expression $f(t)\,dg(t)$ only follows the form of the term in the
sum~(\ref{ch31}): $f(\xi_k)\,\Dl g_k.$

In what follows, we omit the sign~$(S)$ in front of the integral.

It follows directly from the definition~(\ref{ch32}), the linearity of sums~(\ref{ch31})
and the linearity of limit
(\ref{ch32}) that the Stieltjes integral is linear with respect to the
integrated and integrating functions:
$$
\gathered
\int\lt_a^b (\al f(t)+\bt F(t))\,d(\gam g(t)+\dl G(t))=\al\gam\int\lt_a^b
f(t)\,dg(t)+ \\
+\al\dl\int\lt_a^b f(t)\,dG(t)+\bt\gam\int\lt_a^b F(t)\,dg(t)+
\bt\dl\int\lt_a^b F(t)\,dG(t)
\endgathered
$$
($\al,\,\bt,\,\gam,\,\dl$~ are real numbers).


\newpage

\poonkt{Basic properties of the integral\;}\phantom{012345678901234567890} \

Let us establish here some properties of the integral introduced in the previous subsection.
The {\it additivity of the Stieltjes integral} as a function of an interval takes has the following form.

\begin{teo}
\label{ch3thm1}
{\it If the integral $\int\lt_a^b f(t)\,dg(t)$ exists, then
for every $c\in(a,\,b)$
the integrals $\int\lt_a^c f(t)\,dg(t)$ and $\int\lt_c^b f(t)\,dg(t)$ exist as well, and 
the following equality holds:}
\begin{equation}
\label{ch34}
\int\lt_a^b f(t)\,dg(t)\!=\!\int\lt_a^c f(t)\,dg(t)\!+\!\int\lt_c^b f(t)\ ,dg(t)
\end{equation}
\end{teo}
\doc\,
Assume that the integral $\int\lt_a^b f(t)\,dg(t)$ exists. Let $\eps>0$ be 
arbitrary, and let $\dl>0$ be from the principle of convergence. Further, let 
$\tau_1',\,\tau_1''$~ be some partitions of interval $[a,\,c]$
and $\tau_2$ a partition of interval $[c,\,b]$ such that $d(\tau_1')<\dl,$
$d(\tau_1'')<\dl,$ $d(\tau_2)<\dl.$ Then $\tau'\doteq \tau_1'\cup\tau_2,$
$\tau''=\tau_1''\cup\tau_2$~ are partitions of interval $[a,\,b]$ such that
$d(\tau')<\dl,$ $d(\tau'')<\dl.$ Therefore, by the principle of convergence,
$|\gS_{\tau'}-\gS_{\tau''}|<\eps.$ If in the terms in the sums~$\gS_{\tau'}$
and~$\gS_{\tau''}$ corresponding to the interval~$[c,\,b]$ we choose identical~$\xi_k,$ then $\gS_{\tau'}-\gS_{\tau''}=\gS_{\tau_1'}-
\gS_{\tau_1''}.$ Therefore,
$$
|\gS_{\tau_1'}-\gS_{\tau_1''}|=|\gS_{\tau'}-\gS_{\tau''}|<\eps,
$$
which, by the principle of convergence, yields the existence of the integral
$\int\lt_a^c f(t)\,dg(t).$
The existence of the integral $\int\lt_c^b f(t)\,dg(t)$  is proved in the same way.

To prove equality~(\ref{ch34}) it suffices to include point $c$
in all partitions, split sum~(\ref{ch31}) into two sums which correspond to the intervals $[a,\,c]$ and~$[c,\,b]$, respectively, and pass to the limit
as $d(\tau)\to 0$. The existence of the two limits in the right-hand side was
proved above.
\hfill $\square$

\begin{zam}
{\it The converse to Theorem~\ref{ch3thm1} is in general not true: the existence
of the integrals $\int\lt_a^c f(t)\,dg(t)$ and
$\int\lt_c^b f(t)\,dg(t)$ does not imply
the existence of the integral $\int\lt_a^b f(t)\,dg(t).$} 
\end{zam}

In this regard, let us consider the following example.
\begin{example}
\label{notexist}
Let functions $f,\,g:[-1,\,1]\to\R$ be defined as follows:
$$
f(t)=\begin{cases} 0 & \text{when}\q -1\ls t\ls 0, \\ 1 & \text{when}\q 0<t\ls 1;
\end{cases}\
g(t)=\begin{cases} 0 & \text{when}\q -1\ls t<0, \\ 1 & \text{when}\q 0\ls t\ls 1.
\end{cases}
$$
\end{example}
The integrals $\int\lt_{-1}^0 f(t)\,dg(t)$ and $\int\lt_0^1 f(t)\,dg(t)$
exist and are equal to $0$ since  $\gS_{\tau'}=\gS_{\tau''}=0$
($\tau',\,\tau''$~ are arbitrary partitions of intervals
$[-1,\,0]$ and~$[0,\,1]$, respectively). Let~$\tau$~ be an arbitrary partition of
interval $[-1,\,1]$ not containing~0. Let $0\in(t_{k_0-1},\,t_{k_0}).$
Then
$$
\gS_\tau(f,\,g)=f(\xi_{k_0})(g(t_{k_0})-g(t_{k_0-1}))=f(\xi_{k_0}).
$$
Hence we see that $f(\xi_{k_0})=0$ if $\xi_{k_0}<0$ and $f(\xi_{k_0})=1$
if $\xi_{k_0}>0.$ This means that the limit~$\gS_\tau(f,\,g)$
$\Bigl($that is, the integral $\int\lt_a^b f(t)\,dg(t)\Bigr)$ \linebreak
does not exist.

The non-existence of the integral in the previous example is related to the fact that 
functions~$f(\cdot)$ and~$g(\cdot)$ have a common point of discontinuity.


\begin{teo}
\label{ch3thm2}
{\it The existence of one of the integrals $\int\lt_a^b f(t)\,dg(t)$ or
$\int\lt_a^b g(t)\,df(t)$ implies the existence of the other, and yields equality}
\begin{equation}
\label{ch35}
\int\lt_a^b \!\! f(t)\,dg(t)\!+\!\!\int\lt_a^b\!\! g(t)\,df(t)\!=\!f(b)g(b)\!-\!f(a)g(a)
\!\doteq\! f(t)g(t)\Big|_a^b
\end{equation}
{\rm (}the integration by parts formula{\rm )}.
\end{teo}
\doc\,
Assume that the integral $\int\lt_a^b g(t)\,df(t)$ exists. Consider
an arbitrary partition~$\tau$ of interval $[a,\,b]$, and fix points
$\xi_k\in [t_{k-1},\,t_k]$ arbitrarily. Thus, we have
$$
a=t_0\ls\xi_1\ls t_1\ls\xi_2\ls\ldots\ls\xi_n\ls t_n=b.
$$
Let $\xi_0=a,$ $\xi_{n+1}=b.$
\begin{multline*}
\gS_\tau(f,\,g)=\sum\lt_{k=1}^n f(\xi_k)(g(t_k)-g(t_{k-1}))=\\
=f(\xi_1)g(t_1)+f(\xi_2)g(t_2)+\ldots+f(\xi_n)g(t_n)-\\
-(f(\xi_1)g(t_0)+f(\xi_2)g(t_1)+\ldots+f(\xi_n)g(t_{n-1})).
\end{multline*}
Let us add to the right side
$$
f(\xi_0)g(t_0)-f(\xi_{n+1})g(t_n)+f(b)g(b)-f(a)g(a)=0.
$$
Then
\begin{multline*}
\gS_\tau(f,\,g)=-\Big(g(t_0)(f(\xi_1)-f(\xi_0))+g(t_1)(f(\xi_2)-f(\xi_1) )+\ldots \\
\ldots+g(t_n)(f(\xi_{n+1})-f(\xi_n))\Big)+f(t)g(t)\Big|_a^b=\\
=-\gS_{\tau'}(g,\,f)+f(t)g(t)\Big|_a^b,
\end{multline*}
where $\tau'=\{\xi_k\}_{k=0}^{n+1}$~ is a partition of the interval $[a,\,b].$ The fact that some of the ``neighbors'' $\xi_k$ can be equal to
each other does not play any role since the coresponding terms
in the above sum will be equal to 0. So, we get the equality
$$
\gS_\tau(f,\,g)=-\gS_{\tau'}(g,\,f)+f(t)g(t)\Big|_a^b.
$$
Since the limit of sums~$\gS_{\tau}(g,\,f)$ exists by our assumption,
the limit\linebreak
$\gS_{\tau'}(f,\,g)$ also exists ($d(\tau)\to 0$
is equivalent to $d(\tau')\to 0$). Hence we get
equality~(\ref{ch35}).
\hfill $\square$


\poonkt{Existence criterion \;}\;\;\phantom{012345678901234567890} \

Let us present a necessary and sufficient condition for the 
existence of the Stieltjes integral, similar to the
criterion for the existence of the Riemann integral.

We will assume that~$g(\cdot)$~ is an {\it increasing function}.
Let $\tau=\{t_k\}_{k=0}^n$\linebreak
be a partition of interval $[a,\,b],$
$\Dl g_k=g(t_k)-g(t_{k-1})\gs 0.$ Let us define for the Stieltjes integral
the lower~$(s_\tau)$ and the upper~$(S_\tau)$ {\it Darboux sums}:
$$
s_\tau=s_\tau(f,\,g)\doteq \sum\lt_{k=1}^n m_k\Dl g_k, \qq
S_\tau=S_\tau(f,\,g)\doteq \sum\lt_{k=1}^n M_k\Dl g_k,
$$
where
$$
m_k=\inf\lt_{[t_{k-1},\,t_k]} f(t), \qq
M_k=\sup\lt_{[t_{k-1},\,t_k]} f(t).
$$
Obviously, 1)~$s_\tau\ls \gS_\tau \ls S_\tau;$
2)~$s_\tau=\inf\lt_{\xi_k}\gS_\tau,$ $S_\tau=\sup\lt_{\xi_k} \gS_\tau,$ where
the infimum/supremum are taken over all possible choices of 
points~$\xi_k$ in a fixed partition. Moreover, 3)~when new points are added to the partition
the lower sum can only increase, and the upper sum can only
decrease; 4)~a lower sum does not exceed an upper sum,
even corresponding to another partition: if~$\tau'$ and~$\tau''$~ are
two arbitrary partitions, then $s_{\tau'}\ls S_{\tau''}.$

Let us prove property 3) for the lower sums; for the upper sums the proof is similar.
It suffices to inspect what happens to the lower sum after we add a new point to the partition. Let~$\tau$~ be an arbitrary partition,
$\ol t\not\in\tau,$ $\ol t\in(t_{k_0-1},\,t_{k_0}),$ $\tau'=\tau\cup\{\ol t\}.$
Then the term $m_{k_0}\Dl g_{k_0}$ in the sum~$s_\tau$ corresponds to
two terms:
$$
m_{k_0}'(g(\ol t)-g(t_{k_0-1}))+m_{k_0}''(g(t_{k_0})-g(\ol t)),
$$
where $m_{k_0}'$ and $m_{k_0}''$~ are infimums of $f(\cdot)$
on the intervals
$[t_{k_0-1},\,\ol t]$ and $[\ol t,\,t_{k_0}]$, respectively.

Since, clearly, $m_{k_0}'\gs m_{k_0},$ $m_{k_0}''\gs m_{k_0}$
(when we make the set smaller the infimum can only increase), we have
$$
m_{k_0}'(g(\ol t)-g(t_{k_0-1}))+m_{k_0}''(g(t_{k_0})-g(\ol t)) \gs
m_{k_0} (g(t_{k_0})-g(t_{k_0-1})).
$$
Since the other terms do not change as we transition from partition~$\tau$ to
partition~$\tau'$, assertion~3) is proved.

Let~$\tau'$ and~$\tau''$~ be two arbitrary partitions. Denote 
$\tau=\tau'\cup \tau''.$ Applying successively properties~3), 1)
and again~3), given that~$\tau$ is smaller than~$\tau'$ and~$\tau'',$
we obtain
$$
s_{\tau'} \ls s_\tau \ls S_\tau \ls S_{\tau''},
$$
which proves property~4).

Since the set~$\{s_\tau\}_\tau$ of lower sums is bounded from above
(for example, by any upper sum),
there exists $J_*=\sup\lt_\tau s_\tau$. We have $s_\tau\ls J_* \ls S_{\tau'}$
for all $\tau$ and~$\tau'.$ The set~$\{S_\tau\}_\tau$ of upper sums
is bounded from below by $J_*,$ so there exists $J^*=\inf\lt_\tau S_\tau,$
and
\begin{equation}
\label{ch36}
s_\tau \ls J_* \ls J^* \ls S_{\tau'}
\end{equation}
for all $\tau$ and~$\tau'.$ The numbers~$J_*$ and~$J^*$ are called
the {\it lower} and the {\it upper Stieltjes integrals}, respectively.


\begin{teo}
\label{ch3thm3}
{\it Let~$g(\cdot)$~ be an increasing function. Integral $\int\lt_a^b f(t)\,dg(t)$ exists if and only if
$$
\lim\lt_{d(\tau)\to 0} (S_\tau-s_\tau)=0.
$$}
\end{teo}
\doc\,
{\it Sufficiency.}
Let $\lim\lt_{d(\tau)\to 0} (S_\tau-s_\tau)=0$
and let $\eps>0$ be arbitrarily. Then there is~$\dl>0$ such that for~$D(\tau)<\dl
\quad{}$ $S_\tau-s_\tau<\eps,$ and it follows from~\ref{ch36} that $J^*-J_* \ls
S_\tau-s_\tau<\eps.$ Since~$\eps$ is arbitrary, this means that
$J^*-J_*=0.$ Denote $J\doteq J^*=J_*.$ Then $s_\tau \ls J\ls S_\tau.$
From this and property~1) of Darboux sums, we obtain that $|\gS_\tau-J|<\eps,$
i.e. $J$ is the integral $\int\lt_a^b f(t)\,dg(t).$

{\it Necessity.} Suppose that integral $\int\lt_a^b f(t)\,dg(t)$ exists,
let $\eps>0$ be arbitrary. By the definition of the integral, there exists~$\dl>0$ such that
that if~$d(\tau)<\dl,$ then $|\gS_\tau-J|<\fracs{\eps}{3}$ or, put differently, $J-\fracs{\eps}{3}<\gS_\tau < J+\fracs{\eps}{3}.$ The latter, 
property 2) of Darboux sums, yield $J-\fracs{\eps}{3}\ls s_\tau \ls
S_\tau \ls J+\fracs{\eps}{3},$ i.e. $S_\tau-s_\tau \ls\fracs 23 \eps<\eps.$
Hence $\lim\lt_{d(\tau)\to 0}(S_\tau -s_\tau)=0.$
\hfill $\square$

Theorem~\ref{ch3thm3} admits the following two equivalent reformulations.

Denote by $\om_{\tau k}(f)=M_k-m_k$ the oscillation of function~$f(\cdot)$ on the interval~$[t_{k-1},\,t_k].$ 
Then $S_\tau-s_\tau=\sum\lt_{k=1}^n \om_{\tau k}(f)\,\Dl g_k.$


\begin{teo}[6.3']
{\it Let~$g(\cdot)$~ be an increasing function. Integral
$\int\lt_a^b f(t)\,dg(t)$ exists if and only if
$$
\lim\lt_{d(\tau)\to 0} \sum\lt_{k=1}^n \om_{\tau k}(f)\,\Dl g_k=0.
$$}
\end{teo}


\begin{teo}[6.3'']
{\it Let~$g(\cdot)$~ be an increasing function. Stieltjes integral
$\int\lt_a^b f(t)\,dg(t)$ exists if and only if the lower~$J_*$ and the upper~$J^*$ Stieltjes integrals coincide. }
\end{teo}


\poonkt{Sufficient conditions for the existence\;}\;\phantom{012345678901234567890} \

\label{ch3poonkt4}
We provide several sufficient conditions for the existence of the Stieltjes integral.


\begin{teo}
\label{ch3thm4}
{\it If~$f(\cdot)$ is continuous on~$[a,\,b],$ and $g(\cdot)$~ is of bounded
variation on~$[a,\,b],$ then integral $\int\lt_a^b f(t)\,dg(t)$ exists.}
\end{teo}
\doc\,
First, let $g(\cdot)$ be an increasing function. Then we may assume that
$g(b)\!>\!g(a),$ because if $g(b)\!=\!g(a),$ then~$g(t)\!\equiv \! g(a),$ and the integral
exists (and is equal to~0, see~subsection {\bf{6.1}}).
Fix an arbitrary~$\eps>0.$ 
By Cantor's theorem, $f(\cdot)$ is uniformly continuous on~$[a,\,b].$
This means that there exists~$\dl>0$ such that if~$d(\tau)<\dl$, then
$\om_{\tau k}(f)<
\eps/(g(b)-g(a)).$ Therefore,
$$
\sum\lt_{k=1}^n \om_{\tau k}(f)\,\Dl g_k <\frac{\eps}{(g(b)-g(a))}
\sum\lt_{k=1}^n \Dl g_k=\eps.
$$
By Theorem~$3.3'$, integral $\int\lt_a^b f(t)\,dg(t)$ exists.

In the general case, by Theorem~\ref{ch2thm6}, $g(t)=g_1(t)-g_2(t),$
where~$g_1(\cdot)$ and~$g_2(\cdot)$ are
increasing functions. By what has been proved above, integrals
$\int\lt_a^b f(t)\,dg_1(t)$ and $\int\lt_a^b f(t)\,dg_2(t)$ exist, so integral $\int\lt_a^b f(t)\,dg(t)$ exists and is equal to the difference of the previous two integrals.
\hfill $\square$

\medskip

Let~$\mathcal X$ and~$\mathcal Y$~ be two classes of functions $f:[a,\,b]\to\R.$
We say that
$(\mathcal X,\,\mathcal Y)$~ is a {\it pair of existence classes for the Stieltjes integral} (we also say: an {\it RS-pair}) if
any $f(\cdot)\in\mathcal X$ is integrable with respect to any~$g(\cdot)\in \mathcal Y.$

\medskip

In these terms, Theorem~3.4 states that $({\bf{C}}[a,\,b],\,{\bf{BV}}[a,\,b])$~ is a
pair of existence classes of the Stieltjes integral (an RS-pair), and taking into account Theorem~\ref{ch3thm2},
the same can be said about the pair $({\bf{BV}}[a,\,b],\,{\bf{C}}[a,\,b])$.

\medskip

A pair $(\mathcal X,\,\mathcal Y)$ of existence classes for the Stieltjes integral is called {\it exact},
if none of these classes can be extended without narrowing the other one. In other words, this pair is exact
if the fact that $f$ is integrable with respect to any $g\in\mathcal Y$ implies that $f\in\mathcal X,$ and if the integral
$\int\lt_a^b f(t)\,dg(t)$ exists for any~$f\in\mathcal X,$ then $g\in\mathcal Y.$
Below (see Theorem \ref{ch3thm9}) we will show that the pair $\bigl({\bf{C}}[a,\,b],\,{\bf{BV}}[a,\,b ]\bigr)$
is exact, and hence so is the pair  $\bigl({\bf{BV}}[a,\,b],\,{\bf{C}}[a,\,b]\bigr)$.


\begin{teo}
\label{ch3thm5}
{\it If~$f(\cdot)$ is Riemann integrable in the proper sense, and~$g(\cdot)$
satisfies the Lipschitz condition, then integral $\int\lt_a^b f(t)\,dg(t)$
exist.}
\end{teo}
\doc\,
Assume that for all 
$t,\,s\in[a,\,b]$ $(L>0)$ we have
\begin{equation}
\label{ch37}
|g(t)-g(s)| \ls L|t-s|.
\end{equation}

First, let $g(\cdot)$~ be an increasing function, and let $\tau$~ an
arbitrary partition~$[a,\,b].$ Then $\Dl g_k=g(t_k)-g(t_{k-1})\ls
L\Dl t_k$ and
$$
\sum\lt_{k=1}^n \om_{\tau k}(f)\Dl g_k \ls L \sum\lt_{k=1}^n \om_{\tau k}(f)
\Dl t_k \to 0 \quad (d(\tau)\to 0)
$$
according to the Riemann's integrability criterion applied to $f(\cdot)$. 
So $\sum\lt_{k=1}^n \om_{\tau k}(f) \Delta g_k \to 0$ as $d(\tau)\to 0$.
By Theorem~$3.3'$, integral $\int\lt_a^b f(t)\,dg(t)$ exists.

In the general case, put
$$
g(t)=Lt-(Lt-g(t))=g_1(t)-g_2(t).
$$
Function $g_1(t)=Lt$ is increasing (this is a linear function with positive angular
coefficient) and satisfies Lipschitz' condition   with constant~$L$
$\bigl($with the equality in~(\ref{ch37})$\bigr)$.

Let us show that $g_2(t)=Lt-g(t)$ also increases and satisfies 
Lipschitz' condition with  constant~$2L.$

Let~$s>t.$ By \ref{ch37}, $g_2(s)-g_2(t)=Ls-g(s)-Lt+g(t)=L(s-t)- \linebreak -(g(s)-g(t ))\gs 0$. Therefore,
$$
|g_2(s)-g_2(t)|\ls L|s-t|+|g(s)-g(t)| \ls 2L|s-t|.
$$

By what was proved above, there exist integrals $\int\lt_a^b f(t)\,dg_1(t)$ and \\
$\int\lt_a^b f(t)\,dg_2(t).$ Therefore, 
integral 
$\int\lt_a^b f(t)\,dg(t)$ exists and is equal to the difference of the previous two integrals.
\hfill $\square$

Denote by ${\bf{Rm}}[a,\,b]$ the class of functions that are Riemann integrable in the proper sense, and by ${\bf{Lip}}[a,\,b]$ the class of Lipschitz continuous functions. Then Theorems~\ref{ch3thm5} and~\ref{ch3thm2} yield that
$\bigl({\bf{Rm}}[a,\,b],\,{\bf{Lip}}[a,\,b]\bigr)$ and $\bigl({\bf{Lip} }[a,\,b],\,{\bf{Rm}}[a,\,b]\bigr)$ are RS-pairs.

As the following theorem shows, these pairs are not exact.


\begin{teo}
\label{ch3thm6}
{\it If~$f(\cdot)$ is Riemann integrable in the proper sense, and~$g(\cdot)$
can be represnted in the form
\begin{equation}
\label{ch38}
g(t)=c+\int\lt_a^t \vfi(s)ds,
\end{equation}
where~$\vfi(\cdot)$ is absolutely integrable on~$[a,\,b]$ at least in
the improper
sense, then integral $\int\lt_a^b f(t)\,dg(t)$ exists and
we have equality}
\begin{equation}
\label{ch39}
(S)\!\int\lt_a^b f(t)\,dg(t)=(R)\!\int\lt_a^b f(t)\vfi(t)\,dt.
\end{equation}
\end{teo}
\doc\,
We carry out the proof in several steps.

1. First, let~$\vfi(\cdot)$ be absolutely integrable in the proper sense.
Then~$\vfi(\cdot)$ is bounded on~$[a,\,b],$ that is, there exists
constant~$\cal K$ such that $|\vfi(s)|\ls {\cal K}$ for all
$s\in[a,\,b].$ For any~$t,\,s\in[a,\,b]$
$$
|g(s)-g(t)|=\left|\int\lt_t^s \vfi(v)\,dv\right| \ls \int\lt_t^s |\vfi(v)|\,dv\ls
{\cal K} |s-t|,
$$
hence, in this case~$g(\cdot)$ satisfies Lipschitz condition and
the existence of integral $\int\lt_a^b f(t)\,dg(t)$ follows from the previous
theorems.

2. Let~$\vfi(t)\gs 0$ be integrable only in the improper sense.
We restrict our consideration to the case when $\vfi(\cdot)$ has only one singularity. Let~$\vfi(\cdot)$
be integrable in the proper sense on any interval~$[a,\,h],$ $h<b$
and let $\vfi(t)\to +\fy$ for~$t\to b-0,$ $\int\lt_a^b \vfi(t)\,dt<+\fy.$

Fix an arbitrary~$\eps>0$ and let us find $\dl>0$ such that
\begin{equation}
\label{ch310}
\int\lt_{b-\dl}^b \vfi(s)\,ds <\frac{\eps}{2\Om},
\end{equation}
where~$\Om=\sup\lt_{[a,\,b]} f(t)-\inf\lt_{[a,\,b]} f(t)$~ is the oscillation of
$f(\cdot)$ on $[a,\,b];$ $|\Om|<+\fy,$ since $f(\cdot)$ is bounded as a Riemann integrable in the proper sense function. Consider
a partition~$\tau$ of interval $[a,\,b]$ such that $d(\tau)<\fracs{\dl}{2}$. 
For this partition, we represent 
$$
\sum\lt_{k=1}^n \om_{\tau k}(f) \Dl g_k=\sum\nolimits'+\sum\nolimits''
$$
where sum~$\sum'$ corresponds to the intervals $[t_{k-1},\,t_k]$
that lie entirely in~$[a,\,b-\fracs{\dl}{2}],$ the $\sum''$~  corresponds the other intervals
(i.e.\,those that lie in $[b-\dl,\,b]$). For the latter
$$
\om_{\tau k}(f)\Dl g_k \ls \Om \Dl g_k=\Om \int\lt_{t_{k-1}}^{t_k}
\vfi(\nu)\,d\nu,
$$
so due to (\ref{ch310})
\begin{equation}
\label{ch311}
\sum\nolimits'' < \Om \int\lt_{b-\dl}^b \vfi(\nu)\,d\nu <\frac{\eps}{2}.
\end{equation}

Function $\vfi(\cdot)$ is integrable in
the proper sense on interval~$[a,\,b-\fracs{\dl}{2}]$, so there exists $\dl_1,$ $0<\dl_1 \ls \fracs{\dl}{2},$ such that
if $d(\tau)<\dl_1,$ then $\sum'<\fracs{\eps}{2}.$ This, together with~(\ref{ch311}),
implies that
$
\sum\lt_{k=1}^n \om_{\tau k}(f) \Dl g_k<\eps.
$
By Theorem~$3.3'$, integral $\int\lt_a^b f(t)\,dg(t)$ exists.

3. In the general case, we put
$
\vfi^+(t)\doteq \frac{|\vfi(t)|+\vfi(t)}{2} \gs 0, \q
\vfi^-(t)\doteq \\\doteq \frac{|\vfi(t)|-\vfi(t)}{2} \gs 0, \quad
g_1(t)=c+\int\lt_a^t \vfi^+(s)\,ds, \q
g_2(t)=\int\lt_a^t \vfi^-(s)\,ds.
$
Then $\vfi(t)=\vfi^+(t)-\vfi^-(t),$ $g(t)=g_1(t)-g_2(t).$
By what was proved above, there exist integrals $\int\lt_a^b f(t)\,dg_1(t)$
and $\int\lt_a^b f(t)\,dg_2(t),$ and hence integral $\int\lt_a^b f(t)\,dg(t)$
exists and is equal to the difference of the previous two.

4. Let us now prove equality~(\ref{ch39}). Note that the integral in the right-hand side of (\ref{ch39}) exists under the assumptions of the theorem. Let us denote it by~$J.$
First, let $\vfi(t)\gs 0$, which means that $g(\cdot)$ is increasing.
As was already proved, the integral in the left-hand
side of (\ref{ch39}) exists, so by Theorem~$3.3'$, given an arbitrary $\eps>0$,
there exists~$\dl>0$ such that for any partitions $\tau$ such that
$d(\tau)<\dl$ we have $\sum\lt_{k=1}^n \om_{\tau k}(f) \Dl g_k <\eps.$

Let~$\tau$~ be such a partition. We successively obtain that
$$
\aligned
{} & \big|\gS_\tau(f,\,g)-J\big|=\Big|\sum\lt_{k=1}^n f(\xi_k)(g(t_k)-g(t_ {k-1}))
-\int\lt_a^b f(s)\vfi(s)\,ds\Big|= \\
{} & =\Big|\sum\lt_{k=1}^n f(\xi_k) \int\lt_{t_{k-1}}^{t_k} \vfi(s)\,ds-
\sum\lt_{k=1}^n \int\lt_{t_{k-1}}^{t_k} f(s)\vfi(s)\,ds\Big|= \\
{} & =\!\Big|\sum\lt_{k=1}^n \!\int\lt_{t_{k-1}}^{t_k}\! (f(\xi_k)-
f(s))\vfi(s)\,ds\Big|\!\ls\! \sum\lt_{k=1}^n\! \int\lt_{t_{k-1}}^{t_k}\!
|f(\xi_k)-f(s)| \vfi(s)\,ds \!\ls \\
{} & \ls \sum\lt_{k=1}^n \om_{\tau k}(f) \int\lt_{t_{k-1}}^{t_k}
\vfi(s)\,ds=\sum\lt_{k=1}^n \om_{\tau k}(f) \Dl g_k <\eps.
\endaligned
$$
(We put the constant on  the interval~$[t_{k-1},\,t_k]$ value~$f(\xi_k)$
under the integral sign, write down the difference of the integrals as the integral of the difference
and use inequality $|f(\xi_k)-f(s)|\ls \om_{\tau k}(f)$ for
$s\in[t_{k-1},\,t_k].$) Thus,
$\lim\lt_{d(\tau)\to 0} \gS_\tau(f,\,g)=J.$ Since, by definition,
$\lim\lt_{d(\tau)\to 0} \gS_\tau(f,\,g)=\linebreak=
\int\lt_a^b f(t)\,dg(t),$
this proves equality~\ref{ch39} for~$\vfi(t)\gs 0$.

In the general case, we introduce again~$\vfi^+(\cdot),$ $\vfi^-(\cdot),$ $g_1(\cdot)$
and~$g_2(\cdot)$ (as above in part 3 of the proof). Having obtained equalities~(\ref{ch39})
for~$\vfi^+(\cdot),$ $g_1(\cdot)$ and~$\vfi^-(\cdot),$ $g_2(\cdot),$
we subtract the second equality from the first one. As a result, we arrive at~(\ref{ch39}).
\hfill $\square$

\begin{sle}
\label{Slet6}
{\it Let~$f(\cdot)$ be Riemann integrable in the proper sense, let $g(\cdot)$
be continuous and such that its derivative~$g'(\cdot)$ exists at all points, except,
possibly, on some finite set.
Further, assume that $g'(\cdot)$ is absolutely integrable $($at least in the improper sense$).$
Then integral $\int\lt_a^b f(t)\,dg(t)$ exists and}
\begin{equation}
\label{ch312}
(S)\!\int\lt_a^b f(t)\,dg(t)=(R)\!\int\lt_a^b f(t)g'(t)\,dt.
\end{equation}
\end{sle}

This corollary follows from the representation
$$
g(t)=g(a)+\int\lt_a^t g'(s)\,ds.
$$


\poonkt{Reduction to a Riemann integral and a finite sum\;}\;\phantom{01234567890} \

We define the identity function~$\mathfrak h_c(\cdot)$ as in Section {\bf{4.5}}.
Note that $\bigvee_a^b(\mathfrak h_c)=1.$

Let $f(\cdot)$ be continuous on~$[a,b]$. Integral 
$\int\lt_a^b\!\! f(t)d\mathfrak h_c(t)$ exists since $\mathfrak h_c(\cdot)$ is of bounded variation.

1. Let $a<c<b.$ Since we already proved that the integral exists,
partitions~$\tau$ can be chosen so that they do not contain point~$c,$ i.e.
let $c\in(t_{k_0-1},\,t_{k_0}).$ Then, by the continuity of $f(\cdot)$,
$$
\gS_\tau(f,\,\mathfrak h_c)=f(\xi_{k_0})(1-0) \to f(c) \quad (d(\tau)\to 0).
$$

2. If~$c=a,$ then
$$
\gS_\tau(f,\,\mathfrak h_a)=f(\xi_1)(1-\mathfrak h_a(a))=f(\xi_1)
\to f(a) \quad (d(\tau)\to 0).
$$

3. Similarly, if $c=b$, then
$$
\gS_\tau(f,\,\mathfrak h_b)=f(\xi_n)(\mathfrak h_b(b)-0)=f(\xi_n)
\underset{d(\tau)\to 0}{\to} f(b).
$$
Thus, in all three cases,
\begin{equation}
\label{ch313}
\int\lt_a^b f(t)\,d\mathfrak h_c(t)=f(c).
\end{equation}
As we have seen above, the integral does not depend on the values of the integrating function
at its points of discontinuity that lie in the interior of the interval of integration.
Therefore, up to these values, we can write the jump function~$g_d(\cdot)$ using~$\mathfrak h_c(\cdot)$
as follows $\bigl($according to convention (\ref{soglprodolz})$\bigr)$:
\begin{equation}
\label{ch314}
g_d(t)=\sum\lt_{t_k\in T(g)} \sg_{t_k}(g) \mathfrak h_{t_k}(t)
\end{equation}
(recall that $T(g)=\{t_k\}_k$ denotes the set of points of discontinuity of $g(\cdot)$).
Combined with the corollary of Theorem~\ref{ch3thm6}, representations~(\ref{ch17})
and~(\ref{ch314}) allows us to state the following result, which can be used to calculate the Stieltjes integral.

{\it Let~$f(\cdot)$ be continuous on~$[a,\,b],$ let $g(\cdot)$ have at most
finitely many points of discontinuity, $T(g)=\{t_1,\,\ldots,\,t_n\}$. Assume that everywhere on $[a,\,b],$ except, perhaps, some finite set,
there exists derivative~$g'(\cdot)$ which is absolutely
Riemann integrable, at least in the improper sense. Then}
\begin{equation}
\label{ch316}
(S)\!\int\lt_a^b f(t)\,dg(t)=({\cal R})\!\int\lt_a^b f(t)g'(t)\,dt+
\sum\lt_{k=1}^n f(t_k)\cdot \sg_{t_k}(g).
\end{equation}

\doc\,
We define $g_d(t)$ by formula~(\ref{ch314}), that is, we set
$$g_d(t)=\sum\lt_{k=1}^n \sg_{t_k}(g)\mathfrak h_{t_k}(t), \quad g_c(t)=g(t) -g_d(t).$$
Clearly, $g_c(\cdot)$ is continuous. Also, $g_c'(\cdot)$ exists
at all points where $g'(\cdot)$ exists, and at these points
$g_c'(t)=g'(t)$
(because $g_d'(t)=0$ at these points). Therefore,
$g_c(\cdot)$ satisfies the conditions of the corollary of Theorem~\ref{ch3thm6}, and
$
\int\lt_a^b f(t)g_c'(t)\,dt=\int\lt_a^b f(t)g'(t)\,dt.
$
Consequently, using representation $g(t)=g_c(t)+g_d(t),$ equality~(\ref{ch313})
and the linearity of the integral, we obtain
\begin{multline*}
\int\lt_a^b f(t)\,dg(t)=\int\lt_a^b f(t)\,d(g_c(t)+g_d(t))=\int\lt_a^b f(t)\ ,dg_c(t)+\int\lt_a^b f(t)\,dg_d(t)=\\
+\int\lt_a^b f(t)g_c'(t)\,dt+\int\lt_a^b f(t)\,d\sum\lt_{k=1}^n \sg_{t_k}(g) \mathfrak h_{t_k}(t)=\int\lt_a^b f(t)g'(t)\,dt+\\
+\sum\lt_{k=1}^n \sg_{t_k}(g) \int\lt_a^b f(t)\,d\mathfrak h_{t_k}(t)=\int\lt_a^b f(t )g'(t)\,dt+\sum\lt_{k=1}^n \sg_{t_k}(g) f(t_k).
\end{multline*}
\hfill $\square$


\poonkt{Integral over infinite intervals \;}\;\phantom{0123456789012345678901234567890} \

Stieltjes integrals over infinite intervals are defined in the usual way:
$$
\gathered
\int\lt_a^{+\fy} f(t)\,dg(t) \doteq \lim\lt_{b\to +\fy} \int\lt_a^b
f(t)\,dg(t); \\
\int\lt_{-\fy}^b f(t)\,dg(t) \doteq \lim\lt_{a\to -\fy} \int\lt_a^b
f(t)\,dg(t); \\
\int_{-\fy}^{+\fy} f(t)\,dg(t) \doteq \lim\lt_{a\to -\fy \atop b\to +\fy}
\int\lt_a^b f(t)\,dg(t).
\endgathered
$$


\begin{center}
\begin{large}
\section{Further properties of Riemann --- Stieltjes integral}
\end{large}
\end{center}

\poonkt{First exact RS-pair \;}\;\phantom{0123456789012345678901234567890} \

Here we will supplement Theorem \ref{ch3thm1} and prove the exactness of pair $\bigl({\bf{C}}[a,b],\,{\bf{BV}}[a,b]\bigr)$ 
(see above item {\bf{6.4}}).

\medskip

\begin{teo}
\label{ch3thm7}
{\it If at point~$c\in(a,\,b)$ one of the functions~$f(\cdot)$ or~$g(\cdot)$
is continuous and the other one is bounded in a neighborhood of $c$, then the
existence of integrals $\int\lt_a^c f(t)\,dg(t)$ and~$\int\lt_c^b
f(t)\,dg(t)$ implies the existence of integral} $\int\lt_a^b f(t)\,dg(t).$
\end{teo}
\doc\,
Let $\tau=\{t_k\}_{k=0}^n$~ be an arbitrary partition~$[a,\,b].$ If
$c\in \tau$ $(c=t_m),$ then, denoting $\tau'=\{t_k\}_{k=0}^m,$
$\tau''=\{t_k\}_{k=m}^n,$ we obtain $\gS_\tau(f,\,g)=
\gS_{\tau'}(f,\,g )+\gS_{\tau''}(f,\,g).$ This yields the existene of the limit
as $d(\tau)\to 0$,
\begin{equation}
\label{ch3.7.d1}
\lim\lt_{d(\tau)\to 0} \gS_{\tau}(f,\,g)=\int\lt_a^b f(t)\,dg(t),
\end{equation}
and equality~(\ref{ch313}).

Now, let $c\not\in\tau$. Denote $\tau^{\circ}=\tau\cup\{c\}.$
According to what was said above, $\gS_{\tau^{\circ}}(f,\,g) \to \int\lt_a^b f\,dg\;\;\;$
($d(\tau) \to 0$). Let $c\in(t_{k-1},\,t_k).$ Instead of the term
$f(\xi_k)(g(t_k)-g(t_{k-1}))$ from~$\gS_{\tau}$, $\gS_{\tau^{\circ}}$
includes the sum $f(\xi_k')(g(t_k)-g(t))+f(\xi_k'')(g(c)-g(t_{k-1}))$
($t_{k-1}\ls \xi_k'' \ls c \ls \xi_k' \ls t_k$).
Since the existence of the limit of $\gS_{\tau^{\circ}}$
as $d(\tau)\to 0$ has been already established, we can choose $\xi_k'=\xi_k''=c,$
so
$$
\big| \gS_{\tau^{\circ}}-\gS_\tau \big|=\big| f(\xi_k)-f(c) \big|
\cdot\big| g(t_k)-g(t_{k-1})\big|.
$$
Here one of the factors can be made arbitrarily small due to
the continuity, while the other one remains bounded. Therefore, the difference
$\big| \gS_{\tau^{\circ}}-\gS_\tau \big|$ can be made as small as needed
provided that $d(\tau)$ is sufficiently small. This implies the existence
of limit~(\ref{ch3.7.d1}) in this case as well.
\hfill $\square$

\begin{teo}
\label{ch3thm8}
{\it If integrated $(f)$ function and an integrating $(g)$ function have a common point of discontinuity, 
then integral 
$\int\lt_a^b f(t)\,dg(t)$ does not exist.}
\end{teo}
\doc\,
Let $c$ be a common point of discontinuity for $f$ and $g$, $a<c<b$, $g(c+)\ne g(c-).$ We will consider only those
partitions~$\tau$ that do not contain point~$c.$ For such partitions, we
consider two sums~$\gS_{\tau'}$ and~$\gS_{\tau''}$: in one of these sums we include $\xi_k=c$ ($t_{k-1}<c<t_k$) and in the other one $\xi_k>c$ or $\xi_k<c$
depending on whether~$f(\cdot)$ is discontinuous at point~$c$ on the right or on the left.
As in the proof of the previous theorem, we obtain
\begin{equation}
\label{ch3.7.d2}
\big| \gS_{\tau'}-\gS_{\tau''} \big|=\big| f(\xi_k)-f(c) \big|
\cdot\big| g(t_k)-g(t_{k-1})\big|.
\end{equation}
Since as $d(\tau)\to 0$ the right side can be made as close to numbers
$$
\big| f(c+)-f(c) \big|\cdot \big| g(c+)-g(c-) \big| \,\text{or}\,
\big| f(c-)-f(c) \big|\cdot \big| g(c+)-g(c-)\big|
$$
as ond needs, and
since at least one of these numbers is strictly positive, integral
$\int\lt_a^b f\,dg$ cannot exist.

If $g(c-)=g(c+)\ne g(c),$ or $c=a$ and $g(a+)\ne g(a),$ or $c=b$
and $g(b-)\ne g(b),$ then, on the contrary, we include point~$c$ in all partitions
($c=t_k$), and again choose~$\xi_k$ differently in sums~$\gS_{\tau'}$
and~$\gS_{\tau''}.$ As a result, we again obtain equality~(\ref{ch3.7.d2}), which again implies that integral $\int\lt_a^b f\,dg$ cannot exist. 
\hfill $\square$

\begin{teo}
\label{ch3thm9}
{\it The pair $\bigl({\bf{C}}[a,\,b],\,{\bf{BV}}[a,\,b]\bigr)$ of Stieltjes integral existence classes is exact.}
\end{teo}
\doc\,
The class of integrated functions cannot be extended without narrowing the class
integrating functions. If~$f(\cdot)$ has a discontinuity at point~$c$, then it is
obviously not integrable with respect to the function of bounded variation~$\mathfrak h_c(\cdot)$
(see Theorem~\ref{ch3thm8}).

The class of integrating functions cannot be extended without narrowing the class
integrated functions. Let ${\bigvee\lt_a^b (g)=+\fy}.$
There exists point $c\in[a,\,b]$ such that in any of its neighborhoods the total variation of function $g(\cdot)$
is infinite. Without loss of generality, we may assume that $c=b.$
Then we can construct such sequence $\{t_n\}_{n=0}^\fy,$
$t=t_0<t_1<\ldots<t_n<b,$ $t_n\to b,$ that $\sum\lt_{k=1}^\fy \big|
g(t_k)-g(t_{k-1})\big|=+\fy.$ There also exists a
sequence $\{c_k\}_{k=1}^\fy$ such that $c_k\to 0,$ $c_k>0$
$(k=0,\,1,\,2,\,\ldots),$ but
$$
\sum\lt_{k=1}^\fy c_k \big| g(t_k)-
g(t_{k-1}) \big|=+\fy
$$
(one can put $c_n=1\Big/ \sum\lt_{k=1}^n \big|
g(t_k)-g(t_{k-1})\big|,$ since %
it is known
(\cite[sect. 375]{ficht2}) that both series $\sum\lt_{n=1}^\fy a_n/s_n$ ($s_n=\sum\lt_{k=1}^n a_k$), $\sum\lt_{n=1}^\fy a_n$ diverge).
We define a continuous function~$f(\cdot)$ as follows: 
\begin{multline*}
f(t_k)=c_k \cdot \sign \bigl(g(t_k)-g(t_{k-1})\bigr)\;\;\;k=0,\,1,\,2,\,\ldots,\;\;f(b)=0,\\
f(t)=f(t_{k-1})+\frac{f(t_k)-f(t_{k-1})}{t_k-t_{k-1}}(t-t_{k -1})\;\;(t_{k-1}<t<t_k), k=1,\,2,\,\ldots 
\end{multline*}
Then, for $\tau_n=\{t_k\}_{k=0}^n$,
$$
\gS_{\tau_n} (f,\,g)=\sum\lt_{k=1}^n f(t_k) (g(t_k)-g(t_{k-1}))=\sum\lt_{k=1}^n c_k \big| g(t_k)-g(t_{k-1}) \big| \to +\fy
$$
as $n\to\fy,$ so integral $\int\lt_a^b f\,dg$ does not exist.
\hfill $\square$


\poonkt{Reduction to the Riemann integral\;}\phantom{0123456789012345678901234567890} \

Here we discuss Stieltjes integral with respect to an increasing integrating function.
\begin{teo}[On the change of variable]
{\it Let $f(\cdot)$ be continuous on~$[a,\,b],$ $g(\cdot)$ increasing
on~$[a,\,b],$ \;$\vfi(\cdot)$ continuous strictly increasing on~$[\al,\,\bt]$,
$\vfi(\al)=a,$ $\vfi(\bt)=b.$ Then}
$$
\int\lt_a^b f(t)\,dg(t)=\int\lt_\al^\bt f(\vfi(s))\,dg(\vfi(s)).
$$
\end{teo}
\doc\,
Let $F(s)=f(\vfi(s)),$ $G(s)=g(\vfi(s)).$ $F(\cdot)$ is continuous
on~$[\al,\,\bt]$ as a superposition of continuous functions, and~$G(\cdot)$~ is a
function of bounded variation due to Theorem \ref{superpoz}. In these notations
the sought equality takes form
$$
\int\lt_a^b f(t)\,dg(t)=\int\lt_\al^\bt F(s)\,dG(s).
\eqno{(*)}
$$

Let $\tau=\{t_k\}_{k=1}^n$~ be a partition of interval $[a,\,b],$
$\xi_k\in [t_{k-1},\,t_k],$ $s_k=\vfi^{-1}(t_k),$ $\mu_k=\vfi(\xi_k)$.
Then $\tau'=\{s_k\}_{k=0}^n$~ is a partition of interval $[\al,\,\bt]$
(see the proof of Theorem \ref{superpoz}). Obviously, $d(\tau)\to 0$ if and only if
$d(\tau')\to 0.$ We have
$$
\gathered
\gS_\tau (f,\,g)=\sum\lt_{k=1}^n f(\xi_k) \big( g(t_k)-g(t_{k-1}) \big)= \\
=\sum\lt_{k=1}^n f(\vfi(\mu_k)) \big( g(\vfi(\nu_k))-g(\vfi(\nu_{k-1})) \big) = \\
=\sum\lt_{k=1}^n F(\mu_k) \big( G(s_k)-G(s_{k-1}) \big)=\gS_{\tau'}(F,\, G).
\endgathered
$$
Hence, taking $d(\tau)\to 0$, we obtain equality~$(*).$
\hfill $\square$

\begin{teo}[On the reduction to  Riemann integral]
\label{ch3svedkrim}
{\it Let\linebreak
$f(\cdot)$ be continuous and let $g(\cdot)$ be strictly increasing.
Then
$$
(S)\!\int\lt_a^b f(t)\,dg(t)=({\cal R})\!\!\int\lt_{g(a)}^{g(b)}
f(g^{-1}(s))\,ds,
$$
where~$g^{-1}(\cdot)$~ is the inverse function of~$g(\cdot)$.}
\end{teo}
\doc\,
Since~$g(\cdot)$ is strictly increasing, the inverse
function $g^{-1}$, acting from $[g(a),\,g(b)]$ to $[a,\,b]$, exists, it is increasing
(in a~non-strict sense, if~$g(\cdot)$ is not continuous). If~$g(\cdot)$
is continuous at point~$t',$ then~$g^{-1}(\cdot)$ is continuous at point $s'=g(t')$.
If~$g(\cdot)$ has a discontinuity at point~$t''$, that is, the interval
$[g(t''-),\,g(t''+)]$ does not degenerate to a point, then $g^{-1}(s)=t''$
for all $s\in[g(t''-),\,g(t''+)]$. Thus $g^{-1}(\cdot)$
is continuous on~$[g(a),\,g(b)].$ Therefore, function $F(s)=f(g^{-1}(s))$
is also continuous on this interval. This ensures the existence of the integral
on the right-hand side of the equality being proved. (The existence of the integral in the left-hand side  was established earlier.)

If~$\tau_t$~ is a partition of interval~$[a,\,b],$ $\tau_s$~ is the corresponding
partition of interval $[g(a),\,g(b)],$ then, generally speaking, from $d(\tau_t)\to 0$ 
it does not follow that $d(\tau_s)\to 0$. This situation
differs from what was considered in Theorem \ref{ch3thm9}. So,  let us consider an arbitrary
partition $\tau_s=\{s_k\}_{k=0}^n$ of interval $[g(a),\,g(b)]$ and put
$t_k=g^{-1}(s_k)$. Some $t_k$ can coincide,
so we will remove 
the ``extra'' elements from $\{g^{-1}(s_k)\}_{k=0}^n,$
so that the remaining points form a partition of interval $[a,\,b],$ which we
denote by $\tau_t$. 
Now, we have implication $d(\tau_s)\to 0 \q\Rightarrow\q d(\tau_t)\to 0.$
The a priori established existence of the integrals in both parts
of the equality being proved allows us to choose 
points~$\xi_k$ in the most convenient way.
\begin{multline*}
\gS_{\tau_t}(f,\,g)=\sum\lt_{k=1}^n f(t_k) \big( g(t_k)-g(t_{k-1})\big)= \\
=\sum\lt_{k=1}^n f(g^{-1}(s_k))(s_k-s_{k-1})=\gS_{\tau_s}(F).
\end{multline*}
(We kept the ``old'' notation for sum~$\gS_{\tau_t}$ because for
those $t_k$ that coincide the corresponding terms are equal to zero; $\gS_{\tau_s}(F)$~ is the
Riemann sum for~$F(\cdot)$). Letting $d(\tau_s)\to 0$ in $\gS_{\tau_k}(f,\,g)=\gS_{\tau_s}(F),$ we obtain the required
equality.
\hfill $\square$


\newpage

\poonkt{Analogue of Leibniz' formula\;}\phantom{012345678901234567890} \

The following two statements can be called a {\it generalization of the Leibniz formula for the derivative of a product.}
\begin{teo}
\label{ch3thm10}
{\it Let $f,g\in {\bf{CBV}}[a,b],$ let $h\in {\bf{C}}[a,b]$ or $h\in { \bf{BV}}[a,b]$. Then}
\begin{equation}
\label{diff1}
\int\limits_a^b h(t)\,d\bigl(g(t)f(t)\bigr)=\int\limits_a^b h(t)g(t)\,df(t)+\int\limits_a^b h(t)f(t)\,dg(t).
\end{equation}
\end{teo}
\doc\,
Since $h$ (respectively $hg$ and $hf$) is either continuous or has finite total variation, and $gf$
(respectively $f$ and $g$)~ is a continuous function of bounded variation, then by Theorem \ref{ch3thm4} all three integrals in (\ref{diff1})
exist. So, we can choose partition $\tau=\{t_k\}_{k=0}^n$ of interval $[a,b]$ and point $\xi_k$ in the way that is the most convenient for us.

Due to the uniform continuity of function $f$, for every
$\varepsilon >0$ there exists $\delta >0$ such that, if the diameter of partition
$d(\tau)<\delta$, then the oscillation $\omega_{\tau k}(f) <\frac{\varepsilon}{M_h\bigvee\limits_a^b(g)}$, 
where $M_h=\underset{t\in [a,b]}{\sup}|h(t)|$.
We put $\xi_k=t_k$ in all three sums below. Let $d(\tau)<\delta$. Then
\begin{multline*}
\triangle_{\tau}\doteq\bigl|\mathfrak S_{\tau}(h,gf)-\mathfrak S_{\tau}(hg,f)-\mathfrak S_{\tau}(hf,g)\ bigr|=\\
=\left|\sum\limits_{k=1}^n h(t_k)\bigl(g(t_k)f(t_k)-g(t_{k-1})f(t_{k-1})\bigr )\right.-\\
-\sum\limits_{k=1}^n h(t_k)g(t_k)\bigl(f(t_k)- f(t_{k-1})\bigr)-\left.\sum\limits_{k= 1}^n h(t_k)f(t_k)\bigl(g(t_k)-g(t_{k-1})\bigr)\right|=\\
=\left|\sum\limits_{k=1}^n h(t_k)\bigl(f(t_k)-f(t_{k-1})\bigr)\bigl(g(t_k)-g(t_{ k-1})\bigr)\right|\leqslant \\
\leqslant M_h\sum\limits_{k=1}^n \omega_{\tau k}(f)\bigl|g(t_k)-g(t_{k-1})\bigr|\leqslant 
M_h\cdot 
\frac{\varepsilon}{M_h\bigvee\limits_a^b(g)}\bigvee\limits_a^b(g)=\varepsilon
\end{multline*}
Therefore,
$$
\underset{d(\tau)\!\to\! 0}{\lim}\triangle_{\tau}\!=\!\left|\int
\limits_a^b h(t)\,dg(t)f(t)\!-\!\int\limits_a^b h(t)g(t)\,df(t)\!-\!
\int\limits_a^b h(t)f(t)\,dg(t)\right|\!=\!0,
$$
that is, equality (\ref{diff1}) is valid.
\hfill $\square$
\begin{teo}
\label{ch3thm11}
{\it Let any two among three functions $f,g,h$~ be continuous functions of bounded variation, and let the third one be continuous or have finite total
variation. Then equality} (\ref{diff1}) holds.
\end{teo}
\doc\,
The case when $h$ is continuous or has finite total variation was considered
in the previous theorem. Now, let $g$ be
a continuous function or a function of bounded variation, and let
$f,h\in CBV[a,b]$. In this case, in each pairs of functions
$(h,gf),(hg,f),(hf,g)$ one is continuous and the other one is a function
of bounded variation, therefore, by Theorem~\ref{ch3thm4},
all three integrals in (\ref{diff1}) exist. In the notations used in the proof of the previous
theorem we have again
\begin{equation}\label{diff2}
\triangle_{\tau}\leqslant \sum\limits_{k=1}^n |h(t_k)||f(t_k)-f(t_{k-1})||g(t_k)-g(t_ {k-1})|<\varepsilon,
\end{equation}
which is equivalent to equality (\ref{diff1}). Finally,
let $f$ be a continuous function or a function of bounded variation, while $g,h\in {\bf{CBV}}$. If $f\in {\bf{C}}[a,b]$, 
then our reasoning remains the same.
If $f\in {\bf{BV}}[a,b]$, then we use the uniform continuity of 
function $g$ and, choosing a partition of interval $[a,b]$ with a sufficiently small
diameter, we ensure that the oscillation $\omega_{\tau k}(g)$ of $g$
satisfies
$\omega_{\tau k}(g) <\frac{\varepsilon}{M_h\bigvee\limits_a^b(f)}$.
Thus, inequality~(\ref{diff2}) holds,
which is equivalent to equality~(\ref{diff1}).
\hfill $\square$


\begin{center}
\begin{large}
\section{Passing to the limit in Riemann --- Stieltjes integral\;}
\end{large}
\end{center}

\hskip 1mm
\poonkt{Mean value theorem\;}\phantom{0123456789012345678901234567890} \

First, let us show that {\it the Mean value theorem and the Estimation theorem}
are valid for the Stieltjes integral.
\begin{teo}[Mean value theorem]
\label{ch4thm1sr}
{\it Let~$f(\cdot)$~ be a bounded function, let $m\ls f(t)\ls M$, $g(\cdot)$~ be
an increasing function, and assume that integral $\int\lt_a^b f(t)\,dg(t)$ exists.
Then there exists~$\mu\in[m,\,M]$ such that
\begin{equation}
\label{ch41}
\int\lt_a^b f(t)\,dg(t)=\mu(g(b)-g(a)).
\end{equation}
If, moreover,~$f(\cdot)$ is continuous, then there exists~$\xi\in[a,\,b]$ such that
\begin{equation}
\label{ch42}
\int\lt_a^b f(t)\,dg(t)=f(\xi)(g(b)-g(a)).
\end{equation}
If in addition to the above conditions $g(\cdot)$~ is a
strictly increasing function, then in equality} $(\ref{ch42})$ one has $a<\xi<b.$
\end{teo}
\doc\,
If~$g(b)=g(a)$, then~$g(t)\equiv \const$, the integral vanishes, and in~\eqref{ch41}
one can take any~$\mu$. Let $g(b)>g(a)$ and let $\tau$~ be an arbitrary
partition of interval $[a,\,b]$. It is easy to see that
$$
m(g(b)-g(a))\ls \gS_\tau(f,\,g) \ls M(g(b)-g(a)).
$$
Since the integral exists by our assumption, passing here to the limit as ${D(\tau)\to 0}$,
we obtain
$$
m(g(b)-g(a))\!\ls\!\int\lt_a^b f(t)\,dg(t)\!\ls\!M(g(b)-g(a) )\;\;\text{and}\;\;m\!\ls\!\frac{\int\lt_a^b f(t)\,dg(t)}{g(b)-g(a) }\!\ls\!M.
$$
Putting here $\mu\doteq \frac{\int\lt_a^b f(t)\,dg(t)}{g(b)-g(a)}$,
we arrive at \eqref{ch41}.

Let~$f(\cdot)$ be continuous, $m\doteq\inf\lt_{[a,\,b]} f(t)$,
$M\doteq\sup\lt_{[a,\,b]} f(t)$. By the Intermediate value theorem
for continuous functions, there is ${\xi\in[a,\,b]}$ such that~$f(\xi)=\mu$,
which implies equality \eqref{ch42}.

If $f(t) \equiv\const,$ then the last assertion of the theorem is obvious.
Let $f(t)\not\equiv \const$, let $g(\cdot)$ be strictly increasing on $[a,\,b],$
and let $\al$ and~$\bt$ ($a\!<\!\al\!<\!\bt\!<\!b$) be such that
$m' \doteq \inf\lt_{[\al,\,\bt]}
f(t)\!>\!m,$ $M'\doteq \sup\lt_{[\al, \bt]} f(t)\!<\!M.$ Then
$$
\begin{cases}
m\big(g(\bt)-g(\al)\big)<m' \big(g(\bt)-g(\al)\big) \ls \int\lt_\al^\bt
f(t)\,dg(t) \ls \\
\q \ls M' \big(g(\bt)-g(\al)\big)< M\big(g(\bt)-g(\al)\big), \\
m\big(g(\al)-g(a)\big) \ls \int\lt_a^\al f(t)\,dg(t) \ls
M\big(g(\al)-g(a)\big), \\
m\big(g(b)-g(\bt)\big) \ls \int\lt_\bt^b f(t)\,dg(t) \ls
M\big(g(b)-g(\bt)\big).
\end{cases}
$$
Adding up these inequalities, we obtain
$$
m\big(g(b)-g(a)\big) < \int\lt_a^b f(t)\,dg(t) < M\big(g(b)-g(a)\big),
$$
i.e. $\mu=\frac{\int\lt_a^b f(t)\,dg(t)}{g(b)-g(a)} \in (m,\,M),$
hence $a<\xi<b.$
\hfill $\square$

\begin{teo}[Estimation theorem]
\label{ch4thm2}
Let $f\in {\bf{C}}{[a,\,b]}$,
$g\in {\bf{BV}}{[a,\,b]}$. Then
\begin{equation}
\label{ch43}
\left| \int_a^b f(t)\,dg(t) \right| \ls \int_a^b |f(t)|\,dg_\pi (t)\ls
M_f \cdot \bigvee_a^b (g),
\end{equation}
where $g_\pi(t)=\bigvee_a^t(g)$ 
(see~\S~4), $M_f=\max\lt_{[a,\,b]}|f(t)|$ .
\end{teo}
\doc\,
By Theorem \ref{ch3thm4}, both integrals in (\ref{ch43}) exist. Given an arbitrary
partition~$\tau$, we have
\begin{multline*}
|\gS_\tau(f,\,g)|=\left| \sum\lt_{k=1}^n f(\xi_k) \Dl g_k \right| \ls
\sum\lt_{k=1}^n |f(\xi_k)|\cdot|g(t_k)-g(t_{k-1})| \ls \\
\ls \sum\lt_{k=1}^n |f(\xi_k)|\cdot \bigvee_{t_{k-1}}^{t_k}(g)=
\sum\lt_{k=1}^n |f(\xi_k)| (g_\pi(t_k)-g_\pi(t_{k-1})) \ls \\
\ls M_f \sum\lt_{k=1}^n (g_\pi(t_k)-g_\pi(t_{k-1}))=M_f\cdot \bigvee_a^b(g).
\end{multline*}
Passing to the limit $d(\tau)\to 0$ in inequality
$$
\left|\gS_\tau(f,\,g)| \ls |\gS_\tau(|f|,\,g_\pi)\right| \ls M_f \bigvee_a^b(g),
$$
we obtain (\ref{ch43}).
\hfill $\square$

\begin{sle}
\label{slthoc}
{\it Assume that the conditions of Theorem~\ref{ch4thm2} are
satisfied and let partition~$\tau$ be such that
$\om_{\tau k}(f)<\gam$ for all $k=1,\,\ldots,\,n$. Then}
\begin{equation}
\label{ch44}
|\gS_\tau(f,\,g)-\int\lt_a^b f(t)\,dg(t)|<\gam\cdot \bigvee_a^b(g).
\end{equation}
\end{sle}
\doc\,
\begin{multline*}
|\gS_\tau(f,\,g)-\int\lt_a^b\!\! f(t)\,dg(t)|\!=\!\left|
\!\sum\lt_{k=1}^n \! f(\xi_k) \Dl g_k\!-\!\sum\lt_{k=1}^n\!\!
\int\lt_{t_{k-1}}^{t_k}\! f(t)\,dg(t) \right|= \\
=\left| \!\sum\lt_{k=1}^n\!\! \int\lt_{t_{k-1}}^{t_k}\!
(f(\xi_k)-f(t))\,dg(t) \right|
\!\ls\! \sum\lt_{k=1}^n \!\!\int\lt_{t_{k-1}}^{t_k}\!
|f(\xi_k)-f(t)|\,dg_\pi(t) \!\ls \\
\ls\! \sum\lt_{k=1}^n\! \om_{\tau k}(f) \!\int\lt_{t_{k-1}}^{t_k}\! dg_\pi(t)<
\gam \!\sum\lt_{k=1}^n \!(g_\pi(t_k)-g_\pi(t_{k-1}))=\gam \bigvee_a^b(g).\;
\end{multline*}
\hfill $\square$

\poonkt{Limit theorems\;}\phantom{01234567890123456789012345678901234567890} \

Let us start with two {\it classical limit theorems.} Later we will prove some other theorems of this type.

\begin{teo}
\label{ch4thm3}
{\it Let $f_n\in {\bf{C}}{[a,\,b]}$, $n=1,\,2,\,\ldots$,
$f_n(t)\to f(t)\quad{}$ $(n\to\fy)$
uniformly on~$[a,\,b]$, $g\in {\bf{BV}}{[a,\,b]}$. Then
\begin{equation}
\label{ch45}
\lim\lt_{n\to\fy} \int\lt_a^b f_n(t)\,dg(t)=\int\lt_a^b f(t)\,dg(t).
\end{equation} }
\end{teo}
\doc\,
By the well-known theorem in analysis, $f\in {\bf{C}}{[a,\,b]}$. Therefore,
the integral in the right-hand side of (\ref{ch45}) exists. The uniform
convergence of
sequence $\{f_n(\cdot)\}_{n=1}^{\fy}$ means that for
an arbitrary~$\eps>0$ there exists $N\in\N$ such that for all~$n>N$
and for all~$t\in[a,\,b]\;\;$
$
|f_n(t)-f(t)|<\linebreak<
\frac{\eps}{\bigvee_a^b(g)}.
$
Therefore, by Theorem~\ref{ch4thm2}, for~$n>N$
\begin{multline*}
\left| \int\lt_a^b f_n(t)\,dg(t)-\int\lt_a^b f(t)\,dg(t) \right|=\left|
\int\lt_a^b (f_n(t)-f(t))\,dg(t) \right| \ls \\
\ls \frac{\eps}{\bigvee\lt_a^b(g)}
\cdot\bigvee\lt_a^b(g)=\eps,
\end{multline*}
which is equivalent to (\ref{ch45}).
\hfill $\square$
\begin{teo}[Helly]
\label{ch4thm4}
{\it Let $f\in {\bf{C}}[a,\,b],$
$g_n\in {\bf{BV}}[a,\,b],$
$n \geq 1$, $\bigvee_a^b(g_n)\ls {\cal K},$ $n=1,\,2,\,\ldots,$
and let $g_n(t)\to g(t)$ for every $t\in[a,\,b].$ Then}
\begin{equation}
\label{ch46}
\lim\lt_{n\to\fy} \int\lt_a^b f(t)\,dg_n(t)=\int\lt_a^b f(t)\,dg(t).
\end{equation}
\end{teo}
\doc\,
Let us first prove that the limit function~$g(\cdot)$ also has finite total
variation, and that $\bigvee_a^b(g)\ls {\cal K}.$

Let~$\tau$~ be an arbitrary partition of $[a,\,b]$, $\tau=\{t_k\}_{k=0}^m$.
For every $n=1,\,2,\,\ldots$, we have
$$
v_\tau(g_n)=\sum\lt_{k=1}^m |g_n(t_k)-g_n(t_{k-1})| \ls \bigvee_a^b(g_n)
\ls {\cal K}.
$$
Thus, given an arbitrary partition~$\tau$, we have estimate
$$
\sum\lt_{k=1}^m |g_n(t_k)-g_n(t_{k-1})| \ls {\cal K},
$$
therefore, due to the pointwise convergence of sequence
$\{g_n(\cdot)\}_{n=1}^\fy$ as $n\to\fy$, we obtain
$
\sum\lt_{k=1}^m |g(t_k)-g(t_{k-1})| \ls {\cal K}.
$
So, $g \in {\bf{BV}}[a,\,b]$\linebreak and $\bigvee_a^b(g)\ls {\cal K}$.

Let $\eps\!>\!0$ be arbitrary. The uniform continuity of $f(\cdot)$ on~$[a,\,b]$
implies that there exists $\dl>0$ such that if  
$D(\tau)<\dl$,
then $\om_{\tau k}(f)<\fracs{\eps}{3{\cal K}}$. Let us fix such a 
partition~$\tau$. According to the corollary of Theorem~\ref{ch4thm2},
\begin{multline*}$$
\left|\sg_\tau(f,\,g)-\int\lt_a^b f(t)\,dg(t)\right|<\frac{\eps}{3{\cal K}}\cdot \bigvee_a^b(g)\ls \frac{\eps}{3} \\
\left|\sg_\tau(f,\,g_n)-\int\lt_a^b f(t)\,dg_n(t)\right|<\frac{\eps}{3{\cal K}}\cdot \bigvee_a^b(g_n)\ls \frac{\eps}{3}.
\end{multline*}
By the pointwise convergence of sequence~$\{g_n(\cdot)\}_{n=1}^\fy$,
$$
\lim\lt_{n\to\fy} \sg_\tau(f,\,g_n)=\sg_\tau(f,\,g),
$$
so there is~$N\in\N$ such that for~$n>N$
$$
|\sg_\tau(f,\,g)-\sg_\tau(f,\,g_n)|<\fracs{\eps}{3}.
$$
Hence, for~$n>N$,
$$
\gathered
\left| \int\lt_a^b f(t)\,dg_n(t)-\int\lt_a^b f(t)\,dg(t) \right|\ls
\left| \int\lt_a^b f(t)\,dg_n(t)-\sg_\tau(f,\,g_n) \right|+ \\
+|\sg_\tau(f,\,g_n)-\sg_\tau(f,\,g)|+\left| \sg_\tau(f,\,g)-\int\lt_a^b f(t)\,dg(t) \right|
< \frac{\eps}{3}+\frac{\eps}{3}+\frac{\eps}{3}=\eps.
\endgathered
$$
Thus, (\ref{ch46}) is true.
\hfill $\square$


\poonkt{Reduction to a Riemann integral and a series\;}\phantom{0123456789012334567890} \

Helly's theorem allows us to obtain a more general than (\ref{ch316})
formula for calculating Stieltjes integral in the case when $T(g)=
\{t_1,\,t_2,\,\ldots\}$, i.e.\,the set of points of discontinuity of the integrating function~$g(\cdot)$, is countable.

{\it Assume that $f\in {\bf{C}}[a,\,b],$ $g'(\cdot)$
exists everywhere except, perhaps, some countable
set, $g'(\cdot)$ is Riemann integrable. Then}
\begin{equation}
\label{ch47}
(S)\!\int\lt_a^b f(t)\,dg(t)=({R})\!\!\int\lt_a^b f(t)g'(t)\,dt+
\sum\lt_{k=1}^\fy f(t_k) \sg_{t_k}(g).
\end{equation}
\doc\,
It follows from our assumptions that $g\in {\bf{BV}}[a,\,b].$
In view of~(\ref{ch316}), it suffices to prove that
$\;
\int\lt_a^b f(t)\,dg_d(t)=\sum\lt_{k=1}^\fy f(t_k) \sg_{t_k}(g).
$

Let $g_{dn}(t)=\sum\lt_{k=1}^n \sg_{t_k}(g)\mathfrak h_{t_k}(t)$.
Then $\bigl($see (\ref{predstvar1})$\bigr)$
$$
\gathered
\bigvee_a^b(g_{dn})=\sum\lt_{k=1}^n \bigl(|\sigma^-_{c_k}|+|\sigma^+_{c_k}|\bigr) 
\ls\sum\lt_{k=1}^\fy \bigl(|\sigma^-_{c_k}|+|\sigma^+_{c_k}|\bigr)\leqslant\bigvee_a^b(g)(\doteq {\cal K}),\\
|g_{dn}(t)-g_d(t)| \ls\sum\lt_{k=n+1}^\fy |\sg_{t_k}(g)|
\underset{n\to\fy}{\to} 0,
\endgathered
$$
that is, we even have the uniform convergence of sequence\linebreak
$\{g_{dn}(\cdot)\}_{n=1}^\fy$ to~$g_d(\cdot)$.
By Helly's theorem, taking into account (\ref{predstvar1}),
$$
\gathered
\int\lt_a^b f(t)\,dg_d(t)=\lim\lt_{n\to\fy} \int\lt_a^b f(t)\,dg_{dn}(t)
=\lim\lt_{n\to\fy} \sum\lt_{k=1}^n \sg_{t_k}(g) f(t_k)= \\
=\sum\lt_{k=1}^\fy \sg_{t_k}(g) f(t_k).
\endgathered
$$


\poonkt
{A generalization of limit theorem \ref{ch4thm3}}
\begin{teo}
\label{doparcel}
{\it Let $f_n,\,f\in{\bf{C}}[a,b]$, $$f_n(t)\to f(t) \quad \text { as }n\to\infty, \quad t\in [a,b],$$ let $g\in {\bf{BV}}[a,b]$ and assume that there exists $M>0$ such that $|f_n(t)|\leqslant M .$
Then  equality $(\ref{ch45})$ holds:}
$\;\;
\lim\lt_{n\to\fy} \int\lt_a^b f_n(t)\,dg(t)=\int\lt_a^b f(t)\,dg(t).
$\end{teo}
\doc\,
Let first $g(\cdot)$ be strictly increasing. By Theorem \ref{ch3svedkrim}, in view of Arzela's limit theorem from \cite[p. 745]{{ficht2}},
$$
(S)\!\int\lt_a^b f_n(t)\,dg(t)=({R})\!\!\int\lt_{g(a)}^{g(b)} f_n(g^{-1}(s))\,ds\quad (n\in\mathbb N).
$$
Since functions $ f_n\bigl(g^{-1}(s)\bigr)\;(n\in\mathbb N)$ are, clearly, $R$-integrable, by the aforementioned Arzela limit theorem
$$
({R})\!\!\int\lt_{g(a)}^{g(b)} f_n(g^{-1}(s))\,ds\to ({R })\!\!\int\lt_{g(a)}^{g(b)} f(g^{-1}(s))\,ds=(S)\!\int\lt_a^b f (t)\,dg(t)\;\;(n\to\infty).
$$
In the general case, we represent $g(\cdot)$ as the difference of strictly increasing functions (see Remark \ref{zam100}). For each function the result is already proved, and so the required equality follows. 
\hfill $\square$


\poonkt{Variable limit of integration\;}\phantom{01234567890123456789012345678901234567890} \

We note here some properties of the Stieltjes integral with variable upper limit of integration (the indefinite Stieltjes integral).
\begin{teo}
\label{ch4thm6}
{\it Let~$f(\cdot)$ be continuous on~$[a,\,b]$ and let~$g(\cdot)$~ be a function of bounded
variation. Then function
$
G(t)=\int\lt_a^t f(s)\,dg(s)
$

a$)$ has finite total variation;

b$)$ continuous on the right $($left$)$ at the points of right $($left$)$ continuity of function $g(\cdot);$ furthermore,}
\begin{equation}
\label{sigma}
\sigma_t^{+}(G)=f(t)\sigma_t^{+}(g),\quad \sigma_t^{-}(G)=f(t)\sigma_t^{-}(g), \quad\sigma_t(G)=f(t)\sigma_t(g).
\end{equation}
\end{teo}
\doc\,
a)\, Given an arbitrary partition $\tau=\{t_k\}_{k=0}^n$ 
of interval $[a,\,b]$, we have
\begin{multline*}
v_\tau(G)=\sum\lt_{k=1}^n \big|G(t_k)-G(t_{k-1})\big|=\sum\lt_{k=1}^n
\left| \int\lt_{t_{k-1}}^{t_k} f(s)\,dg(s) \right| \ls
\sum\lt_{k=1}^n M_f \bigvee\lt_{t_{k-1}}^{t_k} (g)=\\=M_f \bigvee\lt_a^b (g), \q\text {where}\q M_f=\max\lt_{[a,\,b]} |f(t)|.
\end{multline*}
This implies that~$G(\cdot)$ has finite total variation, and that
$\bigvee\lt_a^b (G) \ls M_f \bigvee\lt_a^b (g).$

b)\,Let first $g(\cdot)$ be an increasing function. Let $h>0.$
By Theorem~4.1,
$
G(t+h)-G(t)=\int\lt_a^{t+h} f(s)\,dg(s)-\int\lt_a^t f(s)\,dg(s)=\int\lt_t^{t+h} f(s)\,dg(s)=\\
=f(\xi) \big( g(t+h)-g(t) \big).
$
Hence, as $h\to 0+$, due to the continuity of $f(\cdot)$,
$$
G(t+)-G(t)=f(t)\big( g(t+)-g(t) \big).
\eqno{(A)}
$$
Similarly, for $h<0$ and $h\to 0-$,
$$
G(t)-G(t-)=f(t)\big( g(t)-g(t-) \big).
\eqno{(B)}
$$
Let us add up $(A)$ and~$(B)$. We obtain
$$
G(t+)-G(t-)=f(t)\big(g(t+)-g(t-)\big).
\eqno{(C)}
$$

In the general case, we represent $g(t)=g_\pi(t)-g_\nu(t),$ and, having obtained
equalities~$(A),\,(B),\,(C)$ for~$g_\pi(\cdot)$ and~$g_\nu(\cdot),$ let us subtract
from the second from the first one. As a result, we obtain equalities~$(A),\,(B),\,(C)$ in the general case.
From this we see that equalities (\ref{sigma}) hold and that $G(\cdot)$ is right (left) continuous.
\hfill $\square$

\begin{zam}
\label{zamrsvar}
{\it We will show in Theorem \ref{teoravvar0} below that, in fact,}
$$
\bigvee\lt_a^b (G)=\int\limits_a^b |f(t)|\,dg_{\pi}(t).
$$
\end{zam}


\poonkt
{Substitution\;}
\begin{teo}
\label{upr314}
{\it Assume that one of the following conditions is satisfied:

1) \;$f\in {\bf{BV}}[a,b],\;g,h\in {\bf{C}}[a,b]$ or

2)\; $f\in {\bf{C}}[a,b],\;$ $g\in {\bf{CBV}}[a,b],h\in {\bf{BV}}[a, b]$. Then}
\begin{equation}\label{subst}
\int\limits_a^b h(t)\,d \,\left(\int\limits_a^t g(s)\,df(s)\right)=\int\limits_a^b h(t)g(t)\,df(t).
\end{equation}
\end{teo}
\doc\,
Let $f\in {\bf{BV}}[a,b],\;g,h\in {\bf{C}}[a,b].$ Denote $F(t)\doteq \int\limits_a^t
g(s)\,df(s)$. According to Theorem \ref{ch4thm6},
$F\in {\bf{BV}}[a,b]$, and by Theorem \ref{ch3thm4} the integrals
in (\ref{subst}) exist. Let $\varepsilon >0$ be arbitrarily small. Then
by Theorem~$6.3'$ there exists $\delta_1>0$ such that for all partitions
$\tau$ whose diameter is $d(\tau)<\delta_1$ we have inequality
\begin{equation}\label{subst2}
\sum_{k=1}^{n} \omega_{\tau k}(g)(f_{\pi}(t_k)-f_{\pi}(t_{k-1}))<\frac{\varepsilon}{3M_h},
\end{equation}
where $f_{\pi}(t)=\bigvee_a^t(f)$ is an increasing function;
regarding notation $M_h$, see the proof of Theorem \ref{ch3thm10}.

For such partitions,
$$
|\mathfrak S_{\tau}(h,F)-\mathfrak S_{\tau}(hg,f)|=\left|\sum_{k=1}^n h(\xi_k)\bigl(F(t_k )-F(t_{k-1})\bigr)-\right.
$$
$$
\left.-\sum_{k=1}^n h(\xi_k)g(\xi_k)\bigl(f(t_k)-f(t_{k-1})\bigr)\right|\leqslant
\sum_{k=1}^n\int_{t_{k-1}}^{t_k}\big|h(\xi_k)g(t)-
$$
$$
-h(\xi_k)g(\xi_k)\big|\,df_{\pi}(t) \leqslant M_h\sum_{k=1}^n \omega_{\tau k}(g)\bigl(f_ {\pi}(t_k)-f_{\pi}(t_{k-1})\bigr)<\frac{\varepsilon}{3}
$$

Further, it follows from the definition of the RS-integral that there exists $\delta_2$ such that
if $D(\tau)<\delta_2$ then
\begin{equation}\label{subst3}
\left |\int_a^b h\,dF-\mathfrak S_{\tau}(h,F)\right |<\frac{\varepsilon}{3},\quad
\left |\int_a^b hg\,df-\mathfrak S_{\tau}(hg,f)\right |<\frac{\varepsilon}{3}.
\end{equation}

Let $D(\tau)<min\{\delta_1,\delta_2\}.$ Then we have inequalities
(\ref{subst3}) and inequality
$$
|\mathfrak S_{\tau}(h,F)-\mathfrak S_{\tau}(hg,f)|<\frac{\varepsilon}{3}.
$$
Consequently,
$$
\left |\int_a^b h(t)\,d\left(\int_a^t g(s)\,df(s)\right)-
\int_a^b h(t)g(t)\,df(t)\right |\leqslant
$$
$$
\leqslant \left |\int_a^b h(t)\,dF(t)-\mathfrak S_{\tau}(h,F)\right |+
\big|\mathfrak S_{\tau}(h,F)-\mathfrak S_{\tau}(hg,f)\big|+
$$
$$
+\left |\mathfrak S_{\tau}(hg,f)-\int_a^b h(t)g(t)\,df(t)\right |<\frac{\varepsilon}{3}+\frac {\varepsilon}{3}+\frac{\varepsilon}{3}=\varepsilon.
$$

Since $\varepsilon$ is arbitrarily small, this argument proves equality~(\ref{subst})

Now let condition 2) be satisfied. All integrals in (\ref{subst}) exist in this case as well.
Let us apply Theorem \ref{ch3thm2} to the inner integral on the left-hand side of (\ref{subst}). Then
$$
\int_a^b h(t)\,d\left(\int_a^t g(s)\,df(s)\right)=\int_a^b\,h(t)\,d\left(gf|_{ a}^{t}-\int_a^t f(s)\,dg(s)\right)=
$$
$$
=\int_a^b h(t)\,dg(t)f(t)-\int_a^b\,h(t)d\left(\int_a^t f(s)\,dg(s)\right).
$$
Applying Theorem \ref{ch3thm10}  to the first term on the right-hand side of the last equality,
and applying to the second term what was already proved, we obtain equality (\ref{subst}) in this case as well.
\hfill $\square$

\bigskip


\poonkt
{Change of the order of integration\;}
\begin{teo}
\label{upr413}
{\it Assume that one of the following conditions is satisfied:
$$
1)f\in {\bf{BV}}[a,b],g\in {\bf{BV}}[c,d],h(\cdot,s)\in {\bf{C}} [a,b],h(t,\cdot)\in {\bf{C}}[c,d]
$$
$(s\in [c,d],t\in [a,b])$ or
$$
2)f\in {\bf{C}}[a,b],g\in {\bf{C}}[c,d],h(\cdot,s)\in {\bf{BV}} [a,b],h(t,\cdot)\in {\bf{BV}}[c,d]
$$
$(s\in [c,d],t\in [a,b])$,  
and there exists $M>0$ such that $|h(t,s)|\leqslant M$
$(t\in [a,b],s\in [c,d])$. Then}
\begin{equation}\label{chRSBVC}
\int\limits_a^b \left(\int\limits_c^d h(t,s)\,dg(s)\right)\,df(t)=\int\limits_c^d \left(\int\limits_a^ b h(t,s)\,df(t)\right)\,dg(s).
\end{equation}
\end{teo}
\doc\,
Let condition 1) be satisfied. Denote
$$
H_g(t)=\int_c^d h(t,s)\,dg(s),\quad H_f(s)=\int_a^b h(t,s)\,df(t).
$$
By Theorem 7.5, $H_g\in BV[a,b], H_f\in BV[c,d]$, hence all integrals in (\ref{chRSBVC})
exist.

Let $\tau_n=\{t_{kn}\}_{k=0}^{m_n}$ be a sequence
of partitions of interval $[a,b]$ such that $d(\tau_n)\to 0$ as $n\to \infty$.
Consider Stieltjes sums for the left-hand side integral (\ref{chRSBVC}) corresponding to these partitions:
$$
\mathfrak S_{\tau_n}(H_g,f)=\sum\limits_{k=1}^{m_n}H_g(\xi_{kn})
\bigl(f(t_{kn})-f(t_{(k-1)n})\bigr)=
$$
$$
\!=\!\sum\limits_{k=1}^{m_n}\int_c^d h(\xi_{kn},s)\,
dg(s)\bigl(f(t_{kn})\!-\!f(t_{(k-1)n})\bigr)\!=\!\int_c^d \mathfrak
S_{\tau_n}(h(\cdot,s),f)\,dg(s).
$$
Since for $n\to\infty\quad \mathfrak S_{\tau_n}(h(\cdot,s),f)\to \int_a^b h(t,s)\,df(t)$ (for all $s\in [c,d]$) and
$|\mathfrak S_{\tau_n}(h(\cdot,s),f)|\leqslant M\bigvee\limits_a^b(f)$,
then, according to Exercise \ref{ch4upr12},
$$
\mathfrak S_{\tau_n}(H_g,f)\to \int\limits_c^d H_f(s)\,dg(s)
$$
for $n\to \infty$. This proves equality (\ref{chRSBVC}).

With the help of Theorem 6.2, condition~2) is reduced to what was already proved.
\hfill $\square$


\vskip 1cm
\begin{center}
{\Large\bf Exercises}
\end{center}

Regarding solutions to exercices \ref{ch4upr4}--\ref{ch4upr12} see \cite{derr08}, regarding solution to \ref{ch3upr12} see \cite{Kond37}.

\medskip

\up
\label{ch4upr4}
Prove the second Mean value theorem for the Riemann integral: if~$f(\cdot)$
is Riemann integrable in the proper sense, and~$g(\cdot)$ is monotone
on~$[a,\,b]$, then there exists $\xi\in[a,\,b]$ such that
$$
\int\lt_a^b f(t)g(t)\,dt=g(a)\int\lt_a^\xi f(t)\,dt+g(b)\int\lt_\xi^b
f(t)\,dt.
$$
If~$g(\cdot)$ is strictly monotone, then~$\xi\in(a,\,b)$.

\up
\label{ch4upr10}
Prove the following strengthening of Theorem~4.3 (weakening of the assumptions of a theorem
while maintaining the statement means its strengthening). Assume that a sequence
of continuous functions~$\{f_n(\cdot)\}_{n=1}^\fy$ 
converges for every $t\in[a,\,b]$ to a continuous function~$f(\cdot)$ and there is a constant~$M$
such that~$|f_n(t)|\ls M$ for all~$t\in[a,\,b]$. Let $g(\cdot)$~ be a function
of bounded variation. Then
$$
\lim\lt_{n\to\fy} \int\lt_a^b f_n(t)\,dg(t)=\int\lt_a^b f(t)\,dg(t).
$$

\up
\label{ch4upr12}
Using Theorem 4.3, Theorem \ref{ch4thm4} and Exercise~\ref{ch4upr10},
formulate a theorem on the passage to the limit of the form
$$
\lim\lt_{n\to\fy} \int\lt_a^b f_n(t)\,dg_n(t)=\int\lt_a^b f(t)\,dg(t).
$$
Prove this theorem.

\up
\label{ch3upr12}
Prove V.\,Kondurariy's theorem: if~$f(\cdot)$ and~$g(\cdot)$ satisfy conditions
$$
|f(t)-f(s)|\ls L|t-s|^\al,\quad|g(t)-g(s)|\ls K|t-s|^\bt,
$$
$\al+\bt>1$ ($t,\,s\in[a,\,b]$), then integral
$\int\lt_a^b f(t)\,dg(t)$ exists.


\newpage

\begin{center}
\begin{Large}
{\bf{Chapter IV. \;Classes of functions}}
\end{Large}
\end{center}
\addcontentsline{toc}{section}{Chapter IV. \;Classes of functions}

\begin{center}
\begin{large}
\section{Absolutely continuous functions}
\end{large}
\end{center}

\poonkt
{Definition. Basic properties\;}\phantom{12345678901234567890} \

A function $f:[a,\,b]\to\R$ is said to be absolutely continuous on interval $[a,\,b]$
if for any $\varepsilon>0$ there exists $\delta>0$ such that for any finite
systems of pairwise non-overlapping intervals
$
\bigcup\lt_{k=1}^n(a_k,\,b_k)\subset[a,\,b]
$
one has
\begin{equation}
\label{par21}
\sum\lt_{k=1}^n(b_k-a_k)<\delta\Longrightarrow\sum\lt_{k=1}^n|f(b_k)-f(a_k)|<\varepsilon
\end{equation}
Let us show that {\it in this definition a finite system of intervals can be replaced by a countable one.}

Let $f(\cdot)$ be absolutely continuous in the sense of the above definition
and let $\delta$ correspond to $\frac{\varepsilon}{2}$ in (\ref{par21}). Next,
let
$\{(a_k,\,b_k)\}_{k=1}^{\infty}$ be a countable system of
non-overlapping intervals such that
$
\sum\limits_{k=1}^{\infty} (b_k-a_k)<\delta.
$
Then, for any $n\in\mathbb N$,
$$
\sum\limits_{k=1}^{n} (b_k-a_k)<\delta, \quad \mbox{and therefore}\quad
\sum\limits_{k=1}^{n} \mid f(b_k)-f(a_k)\mid <\frac{\varepsilon}{2}.
$$
Since $n$ is arbitrary, this means that $\sum\limits_{k=1}^{\infty}\mid f(b_k)-f(a_k)\mid
\ls\frac{\varepsilon}{2}<\varepsilon.$

Suppose that for every $\epsilon >0$ there is
$\delta >0$ such that, for any counting system of non-overlapping
intervals, we have the implication
\begin{equation}\label{ilsur-count}
\sum\limits_{k=1}^{\infty} (b_k-a_k)<\delta \Rightarrow \sum\limits_{k=1}^{\infty} \mid f(b_k)-f(a_k)\mid<\varepsilon.
\end{equation}
Let $\{(a_k,\,b_k)\}_{k=1}^{\infty}$ be a countable system
non-overlapping intervals such that for any $n\in\mathbb N$ we have
inequality $\sum\limits_{k=1}^{n} (b_k-a_k)<\frac{\delta}{2}.$
Then inequality
$\sum\limits_{k=1}^{\infty} (b_k-a_k)\ls \frac{\delta}{2}<\delta$ holds. In view of implication
 (\ref{ilsur-count}), inequality
$$
\sum\limits_{k=1}^{n} \mid f(b_k)-f(a_k)\mid <\sum\limits_{k=1}^{\infty}
\mid f(b_k)-f(a_k)\mid <\varepsilon
$$
holds as well. This means that $f(\cdot)$ is absolutely continuous in the sense of the above definition.
\hfill $\square$

If we take ${n=1}$ in the definition of absolute continuity, then we get the
definition of uniform continuity. Therefore, {\it an absolutely continuous function is uniformly continuous.}


\begin{teo}
\label{par2thm1}
Let $f(\cdot)$ and $g(\cdot)$ be absolutely continuous on $[a,\,b],$ $\al,\ \bt\in\R.$
Then
\begin{enumerate}
\item[a)] the function $F(t)=\al f(t)+\bt g(t)$ is absolutely continuous;
\item[b)] the function $G(t)=f(t)\cdot g(t)$ is absolutely continuous;
\item[c)] the function $H(t)=f(t)/g(t)$ is absolutely continuous if $g(t)\ne 0$ for all $t\in[a,\,b].$
\end{enumerate}
\end{teo}
\doc\,
Assertions a) and b) are proved in the same way as the analogous results for the functions of bounded variation
(see the proof of Theorem \ref{ch2thm4}).
Let us prove assertion c).
Since $g(\cdot)$ does not vanish on $[a,\,b]$, we may assume that
$g(t)>0$ on $[a,\,b].$ The continuity of $g(\cdot)$ implies the existence of 
$\sigma>0$ such that $g(t)\gs\sigma$ for all $t\in[a,\,b].$ (If $\sigma=0,$ then by
the Weierstrass theorem we can find a point $t_0$ such that $g(t_0)=\sigma=0,$ which would contradict
to inequality $g(t_0)>0.$)

Let $\varepsilon>0$ be arbitrary and let $\delta>0$ be such that implication
 (\ref{par21}) holds with $\varepsilon$ replaced by $\varepsilon\cdot\sigma^2.$ Then
$$
\gathered
\sum\lt_{k=1}^n\left|\dfr1{g(b_k)}-\dfr1{g(a_k)}\right|=\sum\lt_{k=1}^n
\dfr{|g(b_k)-g(a_k)|}{g(b_k)\cdot g(a_k)}\ls\\
\ls\dfr1{\sigma^2}\sum\lt_{k=1}^n|g(b_k)-g(a_k)|<
\dfr{\varepsilon\cdot\sigma^2}{\sigma^2}=\varepsilon.
\endgathered
$$

Therefore, function $\dfr1{g(\cdot)}$ is absolutely continuous. Now, it follows from assertion b) that
function $H(t)=f(t)\cdot\dfr1{g(t)}$ is also absolutely continuous.
\hfill $\square$

Denote by $AC[a,\,b]$ the set of all absolutely continuous functions on $[a,\,b]$. The previous theorem implies
that $AC[a,\,b]$~ is a {\it vector space.} Since the addition and the multiplication by scalar functions are defined pointwise,
that is, it boils down to addition, multiplication and division of real numbers, the axioms of vector space follow 
from the properties of addition and multiplication of real numbers.


\hskip 1mm
\poonkt
{Sufficient conditions. Necessary conditions \;} \phantom{12345678901234567890} \

The following two theorems indicate the ``place'' of absolutely continuous functions among 
other classes of functions.
\begin{teo}
\label{par2thm2}
{\it If function $f(\cdot)$ satisfies the Lipschitz condition on~$[a,\,b]$, then it is absolutely
continuous on $[a,\,b].$ }
\end{teo}
\doc\,
Let $f(\cdot)$ satisfy the Lipschitz condition on $[a,\,b]$ with constant
$L$, let $\varepsilon >0$ be arbitrary,
$\delta =\frac{\varepsilon}{L}$ and let $\sum\limits_{k=1}^{n} (b_k-a_k)<\delta.$
Then
$$
\sum\limits_{k=1}^{n} \mid f(b_k)-f(a_k)\mid \ls L\sum\limits_{k=1}^{n} (b_k-a_k)<
L\frac{\varepsilon}{L}=\varepsilon.
$$
\hfill $\square$
\begin{sle}
\label{par2corol1}
{\it If $f(\cdot)$ has bounded derivative, then it is absolutely continuous.}
\end{sle}
\doc\,
Let $\mid f'(t)\mid \ls M.$ Using the Lagrange finite increment formula
($t,s\in [a,\,b]$):
$$
\mid f(t)-f(s)\mid =\mid f'(\xi)(t-s)\mid \ls M\mid t-s\mid (\xi\in (t,s)\;
\mbox{or}\;\xi\in (s,t)),
$$
we obtain that the Lipschitz condition is satisfied by $f$ with Lipschitz constant $L=M.$
\hfill $\square$
\begin{teo}
\label{par2thm3}
{\it A function that is absolutely continuous on~$[a,\,b]$ has finite total variation
on}~$[a,\,b]\;\;\bigl(AC[a,\,b]\subset {\bf{CBV}}[a,\,b]$; at the moment, we consider this as an inclusion in the set-theoretic sense$\bigr).$
\end{teo}
\doc\,
Let $f\in AC[a,\,b].$ Let us fix $\delta>0$ from the definition
of absolute continuity corresponding to $\varepsilon=1.$ Then, for any
finite system of non-overlapping intervals $\bigcup\lt_{k=1}^n(a_k,\,b_k)$, we have
\begin{equation}
\label{par22}
\sum\lt_{k=1}^n(b_k-a_k)<\delta\Longrightarrow\sum\lt_{k=1}^n|f(b_k)-f(a_k)|<1
\end{equation}

By (\ref{par22}), for any partition $\tau_k=\{t_{kj}\}_{j=0}^{p_k}$ of interval $[t_{k-1},\,t_k]$, we have
$v_{\tau_k}(f)=\sum\lt_{j=1}^{p_k}|f(t_{kj})-f(t_{k,j-1})|<1\;\; $
(regarding notation $v_{\tau_k}(f)$, see section {\bf{4.1}}). Therefore, for an arbitrary partition $\tau'$ of interval $[a,\,b]$ (without loss of generality, we can assume that
$\tau_k\subset\tau',$
see~{\bf{4.1}}), we have
$v_{\tau'}(f)=
\sum\lt_{k=1}^n v_{\tau_k}(f)\ls m.\;$
Therefore, $f\in {\bf{BV}}[a,\,b]$ and $\bigvee\lt_a^b(f)\ls m.$
\hfill $\square$

Let ${\bf{AC}}={\bf{AC}}[a,\,b]$~ be the normed linear space of absolutely continuous functions with the norm induced from ${\bf{BV}}[a,\,b]$
(or, which is the same, from ${\bf{CBV}}[a,\,b]).\;$ It can be proved (we will not do it at the moment since we do not have use for it)
that ${\bf{AC}}$~ is a Banach space, and hence, by virtue of Theorem \ref{mprth1}, is a subspace of ${\bf{BV}}\;\bigl($
and also of ${\bf{CBV}}\bigr).$


\bigskip

\poonkt
{Representation in the form of a difference\;} \phantom{12345678901234567890} \

A function has finite total variation if and only if it can be represented as the difference of two increasing functions (Theorem \ref{ch2thm6}). 
It is possible to establish a similar criterion for the absolutely continuity of a function.
\begin{teo}
\label{par2thm4}
{\it An absolutely continuous function can be\linebreak
represented as a difference of two increasing
absolutely continuous\linebreak
functions.}
\end{teo}
\doc\
Let $f(\cdot)$ be absolutely continuous on $[a,\,b].$ Then\linebreak
$f\in {\bf{CBV}}[a,\,b]\;$
and according to (\ref{ch28})
$f(t)\!=\!f_{\pi}(t)-f_{\nu}(t),$ where
$f_{\pi}(\cdot)$ and $f_{\nu}(\cdot)$ are increasing continuous functions,
$f_{\pi}(a)=0,\;\; f_{\pi}(t)=\bigvee\lt_a^t(f)\;$ for $\;t>a.$ Therefore, it remains to prove that
function $f_{\pi}(\cdot)$ is absolutely continuous on $[a,\,b].$

Let $\varepsilon>0$ be arbitrary and let $\delta>0$ be such that implication~(\ref{par21}),
with $\varepsilon$ replaced by $\varepsilon/2$ holds. Then, since $f_{\pi}(\cdot)$ is increasing,
$$\sum\lt_{k=1}^n\left(f_{\pi}(b_k)-f_{\pi}(a_k)\right)=\sum\lt_{k=1}^n\bigvee \lt_{a_k}^{b_k}(f)=
\sum\lt_{k=1}^n\sup\lt_{\tau_k}v_{\tau_k}(f),$$
where $\tau_k=\{t_{kj}\}_{j=0}^{m_k}$
is a partition of interval $[a_k,\,b_k],$ $k=1,2,\cdots,n.$

Since
$$
\sum\lt_{k=1}^n\sum\lt_{j=1}^{m_k}(t_{kj}-t_{kj-1})=
\sum\lt_{k=1}^n(b_k-a_k)<\delta,
$$
we have $\;\sum\lt_{k=1}^nv_{\tau_k}(f)<\dfr{\varepsilon}{2}\;$
and, therefore,
$\;\;\sum\lt_{k=1}^n\sup\lt_{\tau_k}v_{\tau_k}(f)\ls\dfr{\varepsilon}2<\varepsilon.\;$
Consequently,
$\sum\lt_{k=1}^n(f_{\pi}(b_k)-f_{\pi}(a_k))<\varepsilon,$
i.e. $f_{\pi}\in AC[a,\,b].$

Since $f_{\nu}(t)=f_{\pi}(t)-f(t),$ we also have $f_{\nu}\in AC_{[a,\,b]}. $
\hfill $\square$


\begin{center}
\begin{large}
\section{Proper functions}
\end{large}
\end{center}


\bigskip

\poonkt{Definition and basic properties\;} \ \phantom{01234567890123456789012345678901234567890123456789} 

{\bf{1.}}
Let ${\bf{R}}[a,b]$ denote the set of functions $f:[a,b]\to \mathbb R$
having finite one-sided limits $f(a+),f(b-),f(t+), f(t-)$ at each point $t\in (a,b)$. We will call the functions
from ${\bf R}[a,b]$ {\it proper}. 
The original sources (see \cite{dje64},\,\cite{Schw}; see also \cite{Tvr93}) use term ``regulated functions'',
which does not seem to us very fortunate.
The Russian translation of book \cite{Schw}, article \cite{Tolst84}; see also 
\cite{derr00},\, \cite{derr08}\,\cite{derr18},\,\cite{derr19},\, \cite{DK05} and \cite{BarRod22}) use 
Russian term ``pravilniye funktsii'', which  is translated to English as ``proper functions''.
Continuous, piecewise continuous functions,
functions of bounded variation, functions from ${\bf{CH}}$ are  proper functions.

Let us consider an example. 

\begin{example}
\label{RminusBV}
Let $\{t_n\}_{n=1}^\fy$  be a sequence of real numbers such that
${t_1=1},$ ${t_{n+1}<t_n},$ ${n=1,\,2,\,\ldots},$ ${t_n\to +0}.$
{\it We assume
$${f(0)=0},\quad f(t)=\dfrac{(t-t_{n+1})t_n}{t_n-t_{n+1}}\;$$ for 
$t_{n+1}<t\leqslant t_n,\;n=1,2,\ldots$}
If the series $\sum\limits_{n=1}^\infty t_n$ converges, then $f(\cdot)$~ is a function of bounded variation, and
$\bigvee\limits_0^1(f)=2\sum\limits_{n=1}^\infty t_n$.
If this series diverges, then the total variation of $f(\cdot)$ is infinite, $f(\cdot)$ is discontinuous
(right continuous) and $f\in{\bf{R}}[a,b]\setminus {\bf{BV}}[a,b].$
\end{example}

\begin{teo}
\label{ch8thm1}
{\it  Proper functions are bounded} $\bigl($i.e. \linebreak
 ${\bf{R}}[a,b]\subset {\bf{M}}[a,b]\bigr)$
\end{teo}
\doc\,
Let $f\in {\bf{R}}[a,b]$. Suppose that
$f\notin {\bf{M}}[a,b].$ Then there is a sequence $
\{t_n\}_{n=1}^{\infty}\subset [a,b]$ such that
$f(t_n)\geqslant n$, and its subsequence $\{t_{n_k}\}_{k=1}^{\infty}$ such that 
$t_{n_k}\to t^*+$ or $t_{n_k}\to t^*-\;\bigl(t^*\in [a,b]\bigr).$ It follows that $f (t^*+)=\infty$
or $f(t^*-)=\infty$. This contradicts to the definition of proper functions, which 
proves boundedness of $f(\cdot).$
\hfill $\square$

\begin{lem}
\label{ch8lem1}
{\it Let $f\in {\bf R}[a,b],\;c\in (a,b]\quad \bigl(c\in[a,b)\bigr)$. Then}
$$
\underset{t\to c+0}{\lim}f(t\pm)=f(c+)\qquad \left(\underset{t\to c-0}{\lim}f(t\pm) =f(c-)\right).
$$
\end{lem}
\doc\,
Let us prove, for example, that $\underset{t\to c-}{\lim}f(t+)=f(c-)$. Function $\tilde f(t)=f(t+)$ is right 
continuous and is different from $f(\cdot)$
only at the points of discontinuity. If $c$ is a point of continuity of $f(\cdot)$ (and hence of $\tilde f(\cdot)$), then $\tilde f(c-)=\tilde f(c)= f (c)=f(c-)$,
that is, the statement is true. Let $c$ be a discontinuity point of $f(\cdot)$ and $\tilde f(\cdot)$. There is a sequence
$\{t_k\}_{k=1}^{\infty}$ of points of continuity of $f(\cdot)$ such that $t_k\to c-0$:
$$
(\forall\varepsilon >0)\quad (\exists N)\quad (\forall k>N)\Longrightarrow (|\tilde f(c-)-\tilde f(t_k)|<\varepsilon).
$$
However, $\tilde f(t_k)=f(t_k)$, so $|\tilde f(c-)- f(t_k)|<\varepsilon$ and,
hence
$$\tilde f(c-)=\underset{k\to \infty}{\lim}f(t_k)=\linebreak=f(c-).$$
The other assertions are proved similarly.
\hfill $\square$

\medskip

\begin{teo}
\label{ch8thm2}
{\it The set $T(f)$ of discontinuity points of a  proper function $f(\cdot)$ is at most countable.}
\end{teo}
\doc\,
Put
$$
\mathfrak{s}_t(f)\doteq \max \{|\sigma_t(f)|,\,|\sigma_t^{+}(f)|,\,|\sigma_t^{-}(f)| \}.
$$
Let $T_{\varepsilon}(f)$ denote 
the set of points of discontinuity of $f(\cdot)$ on $[a,b]$ where $\mathfrak{s}_t(x)\geqslant \varepsilon.$
Suppose that $T_{\frac{1}{n}}(f)$ is an infinite set and that $t_{*}$ is its limit point. Say, $t_{*}>a$.
Without loss of generality, we may assume that there is an increasing sequence $\{t_k\}_{k=1}^{\infty}$ of points
from $T_{\frac{1}{n}}(f),\;t_k\to t_{*}-0$ as $k\to \infty$. Passing to the limit  $k\to \infty$ in inequality $\mathfrak{s}_t(f)\geqslant\frac{1}{n}$
 using Lemma \ref{ch8lem1}, we obtain $0\geqslant \frac{1}{n}$, which is impossible.
Hence $T_{\frac{1}{n}}(f)$ is a finite set.
The assertion of the theorem follows from the representation $T(f)=\bigcup\limits_{n=1}^{\infty}T_{\frac{1}{n}}(f).$
\hfill $\square$

\begin{sle}
\label{ch8corol1}
{\it  Proper functions are $R$-integrable.}
\end{sle}

\begin{teo}
\label{shodnsk1}
{\it Let $f\in {\bf{R}}[a,b]$, let $t^*$~ be a limit point of set $T(f)$ and let sequence $\{t_k\} _{k=1}^\infty\subset T(f)$ be such that
$t_k\to t^*\;\;(k\to\infty).$ Then}
\begin{equation}
\label{skachki0}
\mathfrak{s}_{t_k}(f)\to 0\;\;\;\;(k\to\infty).
\end{equation}
\end{teo}
\doc\,
According to the proof of Theorem \ref{ch8thm2},
there is only a finite number of points of discontinuity  of $f$ where $\mathfrak{s}_{t_k}(f)$ is greater than a given value.
Let $\varepsilon_n\to 0\;$ for $n\to\infty$, let $\varepsilon>0$ be arbitrary. Let us find a number $\mathcal N$ such that for all
$n>\mathcal N\;\;\;\varepsilon_n<\varepsilon.$  For such $n$, for all terms of sequence $\{t_k\}_{k=1}^\infty$, except, perhaps,
finitely many of them, inequality $\mathfrak{s}_{t_k}(f)<\varepsilon$  is satisfied. 
Hence $\mathfrak{s}_{t_k}(f)\to 0\;\;(k\to\infty).$
\hfill $\square$


{\bf{2.}}\;
Denote
$$
\hat f(t)=\max\{f(t+),f(t-)\},\quad \check f(t)=\min\{f(t+),f(t-)\}.
$$
\begin{teo}
\label{ch8thm3}
{\it Let $f\in {\bf R}[a,b]$. There are points 
$t^{*},\linebreak
t_{*}\in [a,b]$ such that}
$\;\;
\underset{t\in [a,b]}{\sup}\hat f(t)=\hat f(t^{*}),\quad \underset{t\in [a,b]}{\inf }\check f(t)=\check f(t_{*}).
$
\end{teo}
\doc\,
The proof basically follows the proof of the Weierstrass theorem.
Suppose, for instance, that 
$\underset{s\!\in\! [a,b]}{\sup}\hat f(s)>\hat f(t)$ for all
$t\!\in\! [a,b]$. Consider function \;\;
$g(t)\!=\dfrac{1}{\underset{s\in [a,b]}{\sup}\hat f(s)\!-\!\hat f(t)}.\;$ According to our
assumption, the denominator does not vanish,
$g(t)>0,\;g\in \mathbf R[a,b]$. By Theorem 1, there exists 
constant $K>0$ such that for all $t\in [a,b]\;g(t)\leqslant K$. This implies that for all $t\in [a,b]$ inequality
$\hat f(t)\leqslant\underset{s\in [a,b]}{\sup}\hat f(s)-\frac{1}{K}$  is satisfied. 
But this inequality contradicts to the definition of the supremum. So, there is point
$t^*\in [a,b]$ satisfying the conditions of the theorem. The existence of point $t_*$ is proved in exactly the same way.
\hfill $\square$


\bigskip

\poonkt{Criterion for the properness of a function \;} \ \phantom{01234567890123456789012}

Recall that function $f:[a,b]\to\mathbb R$ is called a {\it step function} if $f(t)=c_k$ for $t\in (t_{k-1},t_k)\quad (k=1,2,...,n\;$ for some $n\in\mathbb N ,$\;
$a=t_0<t_1<...<t_n=b).$ Obviously, any step function is proper.

\begin{teo}
\label{ch8thm4}
{\it For $f:[a,b]\to\mathbb R$ to be  proper,
it is necessary and sufficient that there exists a sequence of
step functions that converges to $f(\cdot)$ uniformly.}
\end{teo}
{\it{Necessity.}}
Let $f:[a,b]\to\mathbb R$ be a  proper function. For any $n\in\mathbb N$ there is only a finite set
$T_{\frac{1}{n}}(f)$ of points of discontinuity where $\mathfrak{s}_t(f)$ is greater than or equal to $\frac{1}{n}$
(see the proof of Theorem \ref{ch8thm2}). We cover each point of set $T_{\frac{1}{n}}(f)$ with an open interval centered at this point and
with radius so small that each such interval contains
only one point from set $T_{\frac{1}{n}}(f)\bigr)\;\;\Bigl($here and in what follows, if
$a\in T_{\frac{1}{n}}(f)\;\;\left(b\in T_{\frac{1}{n}}(f)\right),$ then instead of open intervals we will consider the corresponding half-intervals$\Bigr).$
Denote the union of thus obtained intervals by $G$.

We cover each point of set $[a,b]\setminus T_{\frac{1}{n}}(f)$ with an open interval in such a way that
\begin{equation}\label{Dje}
|f(t')-f(t'')|<\frac{1}{n}
\end{equation}
for any two points $t',t''$ of this interval. Note that none of these intervals contains points from $T_{\frac{1}{n}}(f).$
Combined with the intervals from the union $G$, we obtain an open cover of closed interval $[a,b]$. By the Borel --- Lebesgue  lemma, 
we can extract from this cover a finite subcover that contains all intervals from the union $G.$ After removing the latter, there remains
a finite cover $I_1,\ldots,I_p$ of set $[a,b]\setminus T_{\frac{1}{n}}(f)$ with property (\ref{Dje}).
We may also assume that the intervals $I_k=(t_{k,1},\,t_{k,2})$ of this finite cover are numbered in the ascending order.

For the endpoints of the intervals $I_k$ that coincide with points\linebreak
$t\in T_{\frac{1}{n}}(f)$ we set $f_n(t)=f(t).$
Also, if $t_{k',1}\in T_{\frac{1}{n}}(f)\;\;\Bigl(t_{k'',2}\in T_{\frac{ 1}{n}}(f)\Bigr),$ then we set
\begin{multline*}
f_n(t)=f\left(\dfrac{t_{k',1}+t_{k',2}}{2}\right)\;\;\text{for}\;\;t\in \bigl(t_{k',1},\,t_{k',2}\bigr)\\
\Bigl(f_n(t)=f\left(\dfrac{t_{k'',1}+t_{k'',2}}{2}\right)\;\;\text{for}\; \;t\in \bigl(t_{k'',1},\,t_{k'',2}\bigr)\Bigr),
\end{multline*}
Further, put for some $\;k_0$
$$
\Gamma_{k_0}\doteq (t_{k_0,1},\,t_{k_0,2})\bigcap (t_{k_0+1,1},\,t_{k_0+1,2})\ne\varnothing\;\;\text{and}\;\; c_{k_0}\in \Gamma_{k_0}.
$$
Note that $\Gamma_{k_0}\bigcap T_{\frac{1}{n}}(f)=\varnothing.$
Set
$$
f_n(t)=f\bigl(c_{k_0}\bigr)\;\;\text{for}\;\;t\in \bigl(t_{k_0,1},\,t_{k_0+1, 2}\bigr).
$$
As a result, we obtain that $f_n(\cdot)\;(n\in\mathbb N)$ are step functions, and, by  (\ref{Dje}), we have
$\;\;|f_n(t)-f(t)|<\frac{1}{n}\to 0\;$ as $\;n\to \infty,\;$ i.e., $\; f_n(t)\rightrightarrows f(t).$

{\it{Sufficiency.}}
Let $f_n(\cdot),\;n=1,2,...$ be step functions and let 
$f_n(t)\rightrightarrows f(t)$ as $n\to \infty$. Given an arbitrary $\varepsilon >0$,
there is a number $n$ such that inequality $|f_n(t)-f(t)|<\frac{\varepsilon}{3}$ holds for all $t\in [a,b]$.
 Further, let $t_0\in (a,b].$ There exists $c<t_0$ such that, for any points $t',t''\in (c,t_0)$,
$$
|f(t')-f(t'')|\leqslant |f(t')-f_n(t')|+|f_n(t')-f_n(t'')|+|f_n(t' )-f(t'')|<\varepsilon.
$$
By virtue of the Cauchy convergence principle, $f(t_0-)$ exists. The existence of the right limit for any
$t_1\in [a,b)$ is proved in the same way. Thus, $f(\cdot)$ is a   proper function.
\hfill $\square$


\bigskip

\poonkt
{Banach space ${\bf{R}}[a,b]$\;} \ \phantom{12345678901234567890}

The set ${\bf{R}}[a,b]$ is, obviously, a vector space over field $\mathbb R$ with respect to the operations of addition and
multiplications by scalars from $\mathbb R.$
Let us define the following norm in ${\bf{R}}[a,b]$:
\begin{gather}
\label{normprav}
\|x\|\doteq\underset{t\in [a,b]}{\sup}|x(t)|.
\end{gather}
The axioms of the norm are easily verified. It is also obvious that the convergence with repsect to this norm is the uniform convergence on $[a,b]$.

\begin{teo}
\label{ch8thm5}
{\it The linear normed space ${\bf{R}}[a,b]$ is a Banach space with respect to norm $(\ref{normprav}).$}
\end{teo}
\doc\,
Let $\{x_n\}_{n=1}^{\infty}\subset {\bf R}[a,b]$~ be a fundamental sequence. Arguing in the same way as in the proof of completeness of space 
$\mathbf{M}[a,b]$ (see section \textbf{1.3}) 
we prove that for each $t\in [a,b]$ there exists the limit $x(t)=\underset{n\to \infty}{\lim}x_n(t), $ 
and that $x_n(t)\rightrightarrows x(t)$ on $[a,b].$
That is,
$$
(\forall \varepsilon >0)\;(\exists N\in \mathbb N\; (\forall n>N)\;(\forall t\in [a,b])\;\left(|x_n( t)-x(t)|<\frac{\varepsilon}{4}\right).
$$
By our assumption, for every $n\in\mathbb N$, the limit $\underset{s\to t+}{\lim}x_n(s)$ exists for every $t\in [a,b).$ Hence for $t \in[a,b)$
$$
(\exists \delta >0)\;(\forall s>t,\;s-t<\delta) \;\left(|x_{n_0}(s)-x_{n_0}(t+)|<\frac{ \varepsilon}{4}\right),
$$
where $n_0>N$ is fixed. Let $s',s''>t,\;s'-t<\delta,\;s''-t<\delta$ (so also $|s'-s''|<\delta$) .
\begin{multline*}
|x(s')-x(s'')|\leqslant |x(s')-x_{n_0}(s')|+|x_{n_0}(s')-x_{n_0}(t+) |+\\+|x_{n_0}(t+)-x_{n_0}(s'')|+|x_{n_0}(s'')-x(s'')|<
\frac{\varepsilon}{4}+\frac{\varepsilon}{4}+\frac{\varepsilon}{4}+\frac{\varepsilon}{4}=\varepsilon.
\end{multline*}
By  the Cauchy convergence principle, there exists $x(t+)$. We proved similarly the existence of the limit on the left at every
$t\in (a,b].$ Therefore,
$x$~ is a proper function. Since the convergence in ${\bf R}[a,b]$ is the uniform convergence, this
proves completeness of space ${\bf R}[a,b]$.
\hfill $\square$

Denote by $\mathcal H[a,b]$ the vector space of step functions, equipped with  norm (\ref{normprav}) (not to be confused with the Banach space {\bf{H}}[a,b], see {\bf{5.3}}).
From Theorems \ref{mprth1},\;\ref{ch8thm1},\;\ref{ch8thm4} and \ref{ch8thm5}, we immediately obtain:

\begin{teo}
\label{ch8podpr}
{\it The vector space $\mathcal H[a,b]$ is dense in 
$B$-space ${\bf{R}}[a,b]$ $\bigl(cl\mathcal{H}[a,b]={\bf{R}}[a,b]\bigr)$.
${\bf{R}}[a,b]$~ is a subspace of $B$-space ${\bf{M}}[a,b].$}
\end{teo}

We also mention the following important properties.
\begin{teo}
\label{ch8thmdop}
{\it
Let $\{x_n\}_{n=1}^\infty\subset{\bf{R}}[a,b],\;x_n(t)\rightrightarrows x(t)$. Then
$\underset{n\to\infty}{\lim}\;\underset{s\to t\pm}{\lim}x_n(s)=x(t\pm)$ and
hence}\;\;$\bigl(t\in (a,b)\bigr),$
$$
\underset{n\to\infty}{\lim}\sigma_t(x_n)=\sigma_t(x)\;\;
\underset{n\to\infty}{\lim}\sigma_a^{+}(x_n)=\sigma_a^{+}(x),\;\;
\underset{n\to\infty}{\lim}\sigma_b^{-}(x_n)=\sigma_b^{-}(x).
$$
\end{teo}
\doc\,
Theorem \ref{ch8thm5} implies that $x\in{\bf{R}}[a,b],$ and the uniform\linebreak 
convergence means that $x_n\to x\;(n\to\infty)$ in ${\bf{R}}[a,b].$
By Theorem \ref{SM} (the Shatunovskii --- Moore lemma \cite[I.6.6]{dansch}) the limits can be interchanged:
$\;\;
\underset{n\to\infty}{\lim}\;\underset{s\to t\pm}{\lim}x_n(s)=
\underset{s\to t\pm}{\lim}\;\underset{n\to\infty}{\lim}x_n(s)=
x(t\pm).$
\hfill $\square$

Since for proper functions all arguments in the proof of Lemma \ref{unisochrSM} remain valid
(they rely only on Theorem \ref{SM}), we obtain the following result.
\begin{sle}
\label{uniskachR}
{\it Under conditions of Theorem \ref{ch8thmdop},
$$
\sigma_t(x_n)\rightrightarrows \sigma_t(x),\qquad \sigma_t^{\pm}(x_n)\rightrightarrows \sigma_t^{\pm}(x)
$$
at} $n\to \infty.$
\end{sle}


\begin{teo}
\label{ch8nsep}
{\it ${\bf R}[a,b]$~ is not separable.}
\end{teo}

\doc\,
Consider again the set from Section 1.6,
$S=\bigcup\limits_{c\in[a,b]}\{x_c(\cdot)\},$ where $x_c(t)=0$ for $t\ne c,
x_c(c)=1.$ Obviously, $S\subset {\bf{R}}[a,b],$ and for
any distinct $x$ and $y\quad \|x-y\|=1.$ By Theorem \ref{mprth6}, ${\bf{R}}[a,b]$ is not separable.
\hfill $\square$


\bigskip

\poonkt{On the $rl$-representation of a  proper function\;} \ \phantom{01234567890123456789}

{\bf{1.}}\;By analogy with space ${\bf{CH}}\;\bigl(\subset{\bf{R}}\bigr)$
(see Section ${\bf{5.4}},$ definitions (\ref{predstgot0}), (\ref{predstgot}) and the assertion formulated in Section {\bf{5.4}}), 
we define the following normed linear space $\bigl(\subset {\bf{R}}\bigr).$

Let $T\!=\!\{c_1,c_2,\ldots\}$~be at most a countable, 
particularly numbered set;
the triple $\;(y,T,a)\;$ is
such that $\;y(\cdot)\;$ runs through space \;${\bf{C}}$, and $a=
(a_1,a_2,\ldots)$ runs through the Banach space
$$
{\bf{cs}}\doteq \left\{x=(x_1,x_2,\ldots):\;\text{series}\;
\sum\limits_{k=1}^{\infty}x_k \;\text{converges},\;\|x\|\doteq \underset{n}{\sup}\left|\sum\limits_{k =1}^{n}x_k\right|\right\}
$$
\cite[c. 260]{dansch}.
Let us denote the vector space of functions of the form
\begin{equation}
\label{xyz}
x(t)\doteq y(t)+z(t),
\end{equation}
where 
\begin{multline}
\label{zazt}
z(a)=0,\;z(t)\doteq \sum\limits_{c_j<t} s_j+p_t,\quad p_j\doteq a_{2j-1},\;q_j\doteq a_{2j}, \;s_j\doteq p_j+q_j,\\
j=1,2,\ldots,\;\;\text{(when}\;\; t\notin T\;\;\;p_t=q_t=s_t=0),
\end{multline}
by $\mathcal{C}s=\mathcal{C}s[a,b]$.
Then
$
x_{\mathfrak c}(t)=y(t),\;\; T(x)=T,\;\; \sigma_{c_j}^-(x)=a_{2j-1},\linebreak
\sigma_{c_j}^+(x)=a_{2j},\;\;j=1,2,\ldots,\;\;\sigma(x)=a.\;\;
$
Definition (\ref{xyz}) corresponds to representation
(\ref{predstgot}), definition (\ref{zazt})~--- to representation
(\ref{funskach1}). Now the series in (\ref{zazt}) may also converge non-absolutely.

By our construction, we have (set-theoretic) inclusion $\mathcal{C}s\subset {\bf{R}}.$
Vector space $\mathcal{C}s$ endowed with norm
$\;
\|x\|_{{\bf{Cs}}}\doteq\|y\|_{{\bf{C}}}+\|a\|_{{\bf{cs}}}
$\;
is a Banach space ${\bf{Cs}}={\bf{Cs}}[a,b].$ Obviously, $x\in{\bf{Cs}}$ satisfies inequality
$\|x\|_{{\bf{R}}}\leqslant\|x\|_{{\bf{Cs}}}.$


{\bf{2.}}\;We say that function $x:[a,b]\to \mathbb R$ is $rl$-representable (has an $rl$-representation) if
(see (\ref{-10}))
\begin{equation}
\label{rlpredst}
x(t)=x_+(t)+x_-(t),
\end{equation}
where $x_+(\cdot)\;\;\bigl(x_-(\cdot)\bigr)$ is right-continuous\;\;(left-continuous).

The fact that a function has an $rl$-representation does not guarantee its properness:
\begin{example}
\label{rlnor}
{\it Let}
$$
T(x)\doteq\{0,\ldots t_{n+1},t_n,\ldots t_2,t_1\;(=1)\},\;t_{n+1}<t_n,\;\; t_n\to 0\;\;(n\to\infty),
$$
$
x(t)=\left\{
\begin{array}{lcr}
\sqrt{1-t^2}\;\;\text{for}\;\;-1\leqslant t\leqslant 0,\\
(-1)^n\;\;\text{for}\;\;t_{n+1}<t\leqslant t_n,\;\;n=1,2,\ldots \\
\end{array}
\right.
$
\end{example}
Obviously, $x(\cdot)$ is left-continuous, which means that it has \linebreak
$rl$-representation (the second term is zero), but $x\notin {\bf{R}}$ (because $x(0+)$ does not exist).
Moreover, a function that is left (right) continuous (and therefore $rl$-representable) need not be bounded.

Denote by $\mathcal{RL}\!=\!\mathcal{RL}[a,b]$ the vector space of
functions $x:[a,b]\to \mathbb R\;$ representable in the form (\ref{rlpredst}).
Theorem \ref{predstvar01} and Corollary \ref{predstvar01+}
imply inclusions ${\bf{BV}}\subset \mathcal{RL},\;{\bf{CH}}\subset \mathcal{RL}.$
Since for $x\in {\mathcal{C}s}$ Definition (\ref{-30}) and all subsequent
calculations remain valid, we have inclusion
$\;\mathcal{C}s\subset \mathcal{RL}.$
Of course, we do not have equalities (\ref{-20}) or (\ref{-20+}) here.

The listed subsets of $\mathcal{RL}$ do not exhaust this space. Thus, the function considered in Example \ref{RminusBV}, provided that the series there diverges, belongs to set $\mathcal{RL}\setminus \mathcal{C}s$ (it is left continuous itself).
Let us point out some other examples of functions from $\mathcal{RL}\setminus \mathcal{C}s.$ We can modify Example \ref{RminusBV}
by changing the values of the function at its points of discontinuity.
Let $\alpha$ be a real parameter, let a sequence of real numbers
$\{t_n\}_{n=1}^\infty$ be such that
${t_1<1},\;{t_{n+1}<t_n},$ ${n=1,\,2,\,\ldots},$ ${t_n\to +0}$ and such that the series $\sum\limits_{n=1}^\infty t_n$ diverges.
\begin{example}
\label{RminusBVdop}
{\it We assume
${x_{\alpha}(0)=0},\;x_{\alpha}(t)=\dfrac{(t-t_{n+1})t_n}{t_n-t_{n+1}} \;$ for \linebreak
$t_{n+1}<t<t_n,\;x_{\alpha}(t_n)=\alpha t_n,\;n=1,2,\ldots$}
\end{example}
For $\alpha=1$, we obtain the former function which is continuous on the left. When $\alpha=0$, this function becomes continuous on the right.
For $0<\alpha<1$, its value at a point of discontinuity is strictly between the left and the right limits. For $\alpha<0$ and $\alpha>1$,
its value at a point of discontinuity lies outside of the interval $[x_{\alpha}(t-),\,x_{\alpha}(t+)].$ It is important that for any
$\alpha\;\;\;x_{\alpha}(t_n)\to 0\;(n\to\infty)$. If this condition is violated, then the function is no longer proper. We also note
that the function will remain proper if its values are changed at a finite number of points.


{\bf{3.}}\;Let us give an example of a proper function $x(\cdot)$ that is discontinuous at every rational points, i.e.,
$T(x)=\mathbb Q \bigcap (0,1)$, but has an $rl$-representation despite the divergence of the series of jumps 
$\;\bigl($series of type (\ref{zazt})$\bigr )$.
Note that here $T(x)$ is everywhere dense on $[0,1]$).

The following table contains \textit{all} rational numbers of interval $(0,1)$ {\it once}.
The line with number $m$ contains proper fractions with numerator $m$, provided 
that the numbers that have already been listed in the previous lines
are omitted  (for example, in the second line all fractions with an even denominator are omitted, in the third line we omitted all fractions
with denominator that is a multiple of three, etc.)
$$
\begin{array}{cccccc}
\frac{1}{2}&\frac{1}{3}&\frac{1}{4}&\frac{1}{5}&\frac{1}{6}&\ldots \\
\frac{2}{3}&\frac{2}{5}&\frac{2}{7}&\frac{2}{9}&\frac{2}{11}&\ldots \\
\frac{3}{4}&\frac{3}{5}&\frac{3}{7}&\frac{3}{8}&\frac{3}{10}&\ldots \\
\frac{4}{5}&\frac{4}{7}&\frac{4}{9}&\frac{4}{11}&\frac{4}{13}&\ldots \\
\frac{5}{6}&\frac{5}{7}&\frac{5}{8}&\frac{5}{9}&\frac{5}{11}&\ldots \\
\ldots &\ldots &\ldots &\ldots &\ldots &\ldots \\
\end{array}
$$
Let us denote by $t_{mn}$ the number in the $m$th row and the $n$th column.
Note that the series $\sum\limits_{n=1}^\infty t_{mn},\;m=1,2,\ldots$ {\it diverge}. Let us construct functions
$x_m(\alpha_m,\,\cdot)\;\;(\alpha_m\in\mathbb R,
m=1,2,\ldots)$ as in Example \ref{RminusBVdop}.
By our construction, function $x_m(\alpha_m,\,\cdot)$ is discontinuous at points $t_{mn}\;(n=1.2.\ldots)$ and continuous at the other points
of interval $[0,1]\;(m=1,2,\ldots).$ Note that at each rational point of interval $(0,1)$ only one of
functions $x_m(\alpha_m,\,\cdot)$  is discontinuous. Moreover, $x_m(\alpha_m,\,\cdot)\in \mathcal{RL}\;(m=1,2,\ldots).$
\begin{example}
\label{RminusBVdop+}
{\it Let us put}
$$
x(\alpha,\,t)\doteq\sum\limits_{m=1}^\infty\dfrac{x_m(\alpha_m,\,t)}{2^m}\qquad
\bigl(\alpha=(\alpha_1,\alpha_2,\ldots)\bigr).
$$
\end{example}
The series converges uniformly, so $x(\alpha,\cdot)$ is continuous at all
irrational points (see Section {\bf{1.10.3)}}) and discontinuous
at all rational points (from interval (0,1)). Moreover, $x(\alpha,\cdot)\in\mathcal{RL}.$

As was already noted,
$
x_m(\alpha_m,\,t)=\bigl(x_m\bigr)_{+}(\alpha_m,\,t)+\bigl(x_m\bigr)_{-}(\alpha_m,\,t),
$\;\linebreak
Therefore, it suffices to put \;\;
$$
x_{\pm}(\alpha,\,t)\doteq\sum\limits_{m=1}^\infty\dfrac{\bigl(x_m\bigr)_{\pm}(\alpha_m,\,t) }{2^m}.
$$
Since 
the series converge uniformly, we obtain the desired representation.


\medskip

{\bf{4.}}\;Below we prove the main assertion of this section.
\begin{teo}
\label{rlpr}
${\bf{R}}\subset \mathcal{RL}.$
\end{teo}
\doc\;
Let $x\in {\bf{R}}.$ If $T(x)$ is a finite set, then we may assume that 
$x\in{\bf{CH}},$ and the required inclusion is proved.
Therefore, we may assume in what follows that $T(x)=\{t_1,\,t_2,\,\ldots\}$ 
is an arbitrary countable 
particularly numbered set.
Function $y(t)=x(t+)\;\;\bigl(y(t)=x(t-)\bigr)$ is obviously right-
continuous $\bigl($left-continuous$\bigr).$ Consider
function (mapping $\mathcal{W}$)
$$
(\mathcal{W}x)(t)\doteq x(t)-\dfrac{1}{2}\bigl(x(t+)+x(t-)\bigr).
$$
It is easy to see that at the points of continuity of
$x(\cdot)$ we have $(\mathcal{W}x)(t)=0,$ and for $t\in T(x)$
$$(\mathcal{W}x)(t)=x(t)-\dfrac{1}{2}\bigl(x(t+)+x(t-)\bigr)=\dfrac{1}{2 }\Bigl(\sigma_t^{-}(x)-\sigma_t^{+}(x)\Bigr)\doteq \delta_t(x).$$

If we denote by $\widetilde H=\widetilde H[a,b]$ the vector space of functions that differ from the function identically equal to zero only on at most
countable set, then
$\mathcal{W}:{\bf{R}}\to\widetilde H,$ and since the right side of $\delta_t(x)$ is also equal to zero at the points of continuity 
of  function $x(\cdot),$ we have
\begin{equation}
\label{Wrl}
(\mathcal{W}x)(t)=\delta_t(x)\quad (t\in [a,b]).
\end{equation}
Since
$\;
x(t)=\dfrac{1}{2}x(t+)+\dfrac{1}{2}x(t-)+(\mathcal{W}x)(t),
$\;
then to prove inclusion $x\in \mathcal{RL}$ it is enough to show that
\begin{equation}
\label{Wrl2+}
\mathcal{W}x\in \mathcal{RL}.
\end{equation}
Let us prove
(\ref{Wrl2+}). Assume that for an arbitrary $n$
$$
\Tt_n\doteq\{t_1,\ldots,t_{n}\}\subset \{t_1,\ldots,t_{n},t_{n+1}\}=\Tt_{n+1}.
$$
Note that $t_{n+1}$ can fall into any of the intervals determined by points of set $\Tt_n.$ Let us select
$t_{jn}^{'}>t_j\;(j=1,\ldots,n)$ so that the interval $(t_j,\,t_{jn}^{'})$ contains no points of set $ \Tt_n.$
When passing from $n$ to $n+1$, only one of numbers $t_{jn}^{'}$ can change (if $t_{n+1}<\underset{1\leqslant j\leqslant n}{\min}\{t_j\}$
or $t_{n+1}>\underset{1\leqslant j\leqslant n}{\max}\{t_j\}$) or two (otherwise).
For the other points $t_{j,n+1}^{'}=t_{jn}^{'}.$
(When selecting points $t_{jn}^{'}$ we take into account convention (\ref{soglprodolz}).) We put
$$
r_{jn}(t)=\left\{
\begin{array}{lcr}
0\;\;\text{for}\;\; t<t_j\;\;\text{and}\;\; t\geqslant t_{jn}^{'},\\
\delta_{t_j}(x)\dfrac{t-t_{jn}^{'}}{t_j-t_{jn}^{'}}\;\;\text{for}\;\; t_j\leqslant t<t_{jn}^{'},\\
\end{array}
\right.
$$
$$
l_{jn}(t)=\left\{
\begin{array}{lcr}
0\;\;\text{for}\;\; t\leqslant t_j\;\;\text{and}\;\;t\geqslant t_{jn}^{'},\\
-\delta_{t_j}(x)\dfrac{t-t_{jn}^{'}}{t_j-t_{jn}^{'}}\;\;\text{for}\;\; t_j<t<t_j^{'},\\
\end{array}
\right.
\;\;\;j=1,\ldots,n.
$$
By our construction, $r_{jn}(\cdot)$ $\bigl(l_{jn}(\cdot)\bigr)$ $(j=1,2,\ldots,n)$ are right-continuous (left-continuous) functions.
Since $x(\cdot)$ is a proper function, we have by Theorem \ref{shodnsk1} $\;\delta_{t_n}(x)\to 0\;\;(n\to\infty) ,$ and so (since $|r_{nn}(t)|, |l_{nn}(t)|$ are less than $|\delta_{t_n}(x)|$),
\begin{equation}
\label{control+}
r_{nn}(t),\,l_{nn}(t)\rightrightarrows 0\;\;\;\;\text{at}\;\;n\to\infty.
\end{equation}
Moreover, in a neighborhood of $t_n\;\;\bigl($see (\ref{Wrl})$\bigr)$,
$$
r_{nn}(t)+l_{nn}(t)=\left\{
\begin{array}{lcr}
\delta_{t_n}(x)\;\;\text{for}\;\;t=t_n,\\
0\;\;\text{for}\;\;t\ne t_n\\
\end{array}
\right.
\quad =\bigl(\mathcal W\bigr)(t).
$$
We put further
$$
y_{+}(t)=\sum\limits_{n=1}^{\infty}r_{nn}(t),\qquad y_{-}(t)=\sum\limits_{n=1}^ {\infty}l_{nn}(t).
$$
In each of the partial sums of these series, for every $t\in [a,b]$ there can be only one term that is not identically equal to zero, therefore,
due to (\ref{control+}), the series converge uniformly. Moreover, $y_{+}(\cdot)\;\;\bigl(y_{-}(\cdot)\bigr)$ is right- (left- continuous),
and since, obviously, $\;y_{+}(t)+y_{-}(t)=
(\mathcal{W}x)(t),$ relation (\ref{Wrl2+}) and the theorem is proved.
\hfill $\square$

\small

\medskip

\poonkt{Other properties of  proper functions\;} \ \phantom{01234567890123456789}

{\bf{1.}}\,
{\it Let $f(\cdot)$~ be a  proper function, $F(t)=\int\limits_a^t f(s)\,ds.$ Then}
$$
F'_{\pm}(t)\doteq \underset{h\to \pm 0}{\lim}\frac{1}{h}\Bigl(F(t+h)-F(t)\Bigr )=f(t\pm).
$$

\doc\,
Let $h>0$. By the Mean value theorem for the Riemann integral,
$$
\Delta F(t,h)\doteq F(t+h)-F(t)=\int\limits_t^{t+h} f(s)\,ds=\mu(t)\cdot h,
$$
where $\underset{[t,t+h]}{\inf}f(s)\leqslant \mu(t)\leqslant
\underset{[t,t+h]}{\sup}f(s)$.
Therefore, $$\underset{[t,t+h]}{\inf}f(s)\leqslant
\frac{1}{h}\Delta F(t,h)
\leqslant \underset{[t,t+h]}
{\sup}f(s),$$ and taking $h\to +0$ we obtain the required.
The reasoning for $h<0$ is similar.
\hfill $\square$

\medskip


{\bf{2.}}\,
{\it Let $f(\cdot)$~ be a   proper function. Then}
$$
\int\limits_a^b (f(t+h)-f(t))\,dt=(f(b\pm)-f(a\pm))h+o(h)\quad \mbox{when }\quad h\to\pm 0.
$$

\doc\,
Let $h>0$. We set $f(t)=f(b)$ for $t>b$. Due to the previous result
$\int\limits_a^b \bigl(f(t+h)-f(t)\bigr)\,dt=\int\limits_{a+h}^{b+h} f(t)\,dt-\int\limits_a^b f(t)\,dt=\linebreak-
\int\limits_b^{b+h}f(t)\,dt-\int\limits_a^{a+h} f(t)\,dt=$
$(f(t+)-f(a+))h+o(h)$ for $h\to+0.$
For $h<0$, arguing similarly, we obtain the assertion with ``lower'' signs.
\hfill $\square$

\medskip

{\bf{3.}}\,
{\it Let $f(\cdot)$ be a  proper function that is RS-integrable with respect to
$g\in BV[a,b],\;$
$c\in [a,b]$. Then
$\;
\int\limits_c^{c+h} f(t)\,dg(t)=f(c\pm)\bigl(g(c\pm)-g(c)\bigr)+\gamma(h) ,
$\;\linebreak
where $\gamma(h)\to 0$ for $h\to\pm 0.$ }

\doc\,
Let $h>0$. By Theorem \ref{ch4thm1},
\begin{equation}\label{regRS}
\int\limits_c^{c+h} g(t)\,df(t)=\mu(f(c+h)-f(c)),
\end{equation}
where $\underset{[c,c+h]}{\inf}g(t)\leqslant \mu\leqslant
\underset{[c,c+h]}{\sup}g(t)$.
This implies that as $h\to +0$
$\mu \to g(c+)$ and, therefore,
$\mu =g(c+)+\alpha (h)$, where $\alpha (h)\to 0$ as $h\to +0$. 
Similarly, $f(c+h)=f(c+)+\beta (h)$, where $\beta (h)\to 0$ as $h\to +0$. 
Substituting these equalities into (\ref{regRS}), we obtain the required with 
$$\gamma (h)=\alpha (h)(f(c+)-f(c))+\beta (h)g(c+)+\alpha (h)\beta (h).$$

For $h<0$, similar arguments give us the statement with lower signs.
\hfill $\square$


\medskip

{\bf{4.}}\,
{\it The composition of two  proper functions may not be a  proper function.}

\doc\,
Let $f(t)=tcos\,\frac{\pi}{2t}$ for $t\ne 0,\;f(0)=0$ (see the example from Section 4.2), $ g (t)=sign\, t$. Function $f(\cdot)$
is continuous on $[-1,1]$, and therefore is
 proper on this interval. Function $g(\cdot)$ is also  proper on this interval. However, function $F(t)=g(f(t))={\rm sign}\,f(t)$
is not  proper on this segment since the limits $F(0+),\;F(0-)$ do not exist.
\hfill $\square$

\medskip

{\bf{5.}}\,
Let $h:[a,b]\times [c,d]\to\mathbb R$~ be a function of two variables. Functions 
$h(\cdot,s)\;(s\in [c,d])$ and $h(t,\cdot)$ are called sections of $h(\cdot,\cdot).$
The following assertion will be useful in what follows.

{\it Let sections} $h(\cdot,s)\;(s\in [c,d]),$ $h(t,\cdot)\;(t\in [a,b])$~ be
{\it  proper functions. Then there is 
constant $M$ independent of $t$ and $s$ such that} $|h(t,s)|\leqslant M$ {\it for all} $(t,s)\in [a,b]\ times[c,d].$

\doc\,
Suppose that the opposite is true: for any natural $n$ there are $t_n$ and $s_n$ such that
$|h(t_n,s_n)|\geqslant n$. Without loss of generality, we may assume that

a)\;$t_n\to t_{*},\;\;s_n\to s_{*}\;\;\bigl(n\to \infty,\;(t_{*},s_{*}) \in [a,b]\times [c,d]\bigr);$

b)\;sequence $\{|h(t_n,s_n)|\}_{n=1}^{\infty}$ is increasing and hence \linebreak
$|h(t_n,s_n)|\to+\infty\;\;(n\to \infty);$

c)\; $t_n\to t_{*}-$ or $t_n\to t_{*}+,\;s_n\to s_{*}-$ or $s_n\to s_{*}+\;(n\to \infty ),$\\
otherwise it is possible to pass in a)-c) to the corresponding subsequence.
Let, for example, $t_n\to t_{*}-,
s_n\to s_{*}-$. Since in any neighborhood of point
$(t_{*},s_{*})$ function $|h(t_n,s_n)|$ takes arbitrarily large values, there are no finite limits
$|h(t-,s_{*})|$ and $|h(t_{*},s-)|,$ which contradicts to the  properness of sections
$h(\cdot,s),\;h(t,\cdot).$
\hfill $\square$
\begin{zam}\label{zam1}
{\it In the absence of  properness, the boundedness of both sections, generally speaking, does not imply boundedness of } $h(\cdot,\cdot).$
\end{zam}

In this regard, consider the following example. Let $[a,b]=[c,d]=[0,1];\\ h(t,s)=\left\{
\begin{array}{lcr}
tg\,\dfrac{2\pi\,ts}{(t+1)(s+1)},\;\;\bigl((t,s)\in[0,1)^2\bigr) ,\\
0,\quad t=1,s=1;
\quad
\end{array}
\right.
$
then
$$
\left|h(t,s)\right|\leqslant R(s)=\left\{
\begin{array}{lcr}
tg\,\dfrac{\pi\,s}{(s+1)}\;\bigl(s\in [0,1)\bigr);\\
0,\;s=1;
\end{array}
\right.
$$
$$
\left|h(t,s)\right|\leqslant Q(t)=\left\{
\begin{array}{lcr}
tg\,\dfrac{\pi\,t}{(t+1)}\;\bigl(t\in [0,1)\bigr),\\
0,\;t=1.
\end{array}
\right.
$$
However, $h(t,s)\to +\infty$ for $t\to 1,\;s\to1.$

\normalsize

\bigskip


\begin{center}
\begin{large}
\section{$\sigma$-continuous functions}
\end{large}
\end{center}


\medskip

\poonkt
{Spaces ${\bf N}[a,b]$ and ${\bf{Rm}}[a,b]$\;} \  \phantom{012345678901234567890}

Denote by ${\bf{N}}[a,b]$ the set of bounded functions 
$x:[a,b]\to\mathbb R$ having at most countable set $T(x)$ of points of discontinuity.
Obviously, ${\bf{N}}[a,b]$~ is a vector space. Let us endow it with norm (\ref{normprav}). Functions $x\in {\bf{N}}[a,b]$
will be called $\sigma$-continuous.
Denote by ${\bf{Rm}}[a,b]$ the vector space of Riemann-integrable (and hence bounded)
functions with norm (\ref{normprav}).
Theorems \ref{ch8thm1},\, \ref{ch8thm2}, Corollary \ref{ch8corol1} and the Lebesgue theorem (a criterion for
$R$-integrability of a function) imply inclusions of normed linear spaces
\begin{equation}
\label{vkl1}
{\bf{R}}[a,b]\subset{\bf{N}}[a,b]\subset{\bf{Rm}}[a,b]\subset{\bf{M}}[ a,b].
\end{equation}

{\it Inclusions} (\ref{vkl1}) {\it are strict and are inclusions of Banach spaces.}

\doc\,
The strictness of the inclusions is obvious. We only need to prove completeness of spaces ${\bf{Rm}}[a,b]$ and ${\bf{N}}[a,b]$.

Let us show that ${\bf{Rm}}[a,b]\;\Bigl({\bf{N}}[a,b]\Bigr)$ is closed in ${\bf{M}}[ a,b]\;\Bigl({\bf{Rm}}[a,b]\Bigr).$
Let $\{x_n\}_{n=1}^\infty\subset{\bf{M}}[a,b],\;x_n(t)\rightrightarrows x(t)\;(n\to\infty)$.
Due to completeness of ${\bf{M}}[a,b]$, the limit $x\in{\bf{M}}[a,b],$ and due to the uniform convergence we have the inclusion
\begin{gather} \label{PN}
T(x)\subset\bigcup\limits_{n=1}^\infty T(x_n).
\end{gather}
Indeed, let $C\doteq [a,b]\setminus\bigcup\limits_{n=1}^\infty T(x_n)$. The points of this set are the points of continuity of all elements of this sequence. Due to the uniform convergence of the sequence, these are  points of continuity of function $x$, which means that
the points of set $C$ do not contain points from $T(x)$. This proves inclusion (\ref{PN}).

From this inclusion, Theorem \ref{konschet2}  and the aforementioned Lebesgue theorem (the fact that a countable union of countable
sets is a countable set), 
it follows that if $\{x_n\}_{n=1}^\infty\subset{\bf{Rm}}[a,b]\;\;\Bigl(\subset{\bf{N}}[ a,b]\Bigr),$
then $x\in{\bf{Rm}}[a,b]\\
\Bigl(\in{\bf{N}}[a,b]\Bigr)$. Hence ${\bf {Rm}}[a,b]\;\;
\Bigl({\bf{N}}[a,b]\Bigr)$ is closed in
${\bf{M}}[a,b]\;\;\Bigl({\bf{Rm}}[a,b]\Bigr).$ And this means that
${\bf{Rm}}[a,b]\;\;\Bigl({\bf{N}}[a,b]\Bigr)$~ is a Banach space (note: we are talking about norm ( \ref{normprav})).
\hfill $\square$

Inclusions (\ref{vkl1}) imply that {\it space ${\bf N}[a,b]$ is not separable}.

\small


\medskip

\poonkt{Space ${\bf{N_1}}[a,b]\;$} \ \phantom{01234567890123456789}\phantom{12345678901234567890}

To construct a set dense in ${\bf{N}}[a,b]$ (see Corollary \ref{plotFvN} below),
consider a space that is wider than ${\bf{N}}[a,b]$. To this end, we recall 
some elements of the classification of
functions by R. Baire (see \cite[XV.1]{nat97}). According to this classification, the zeroth class includes continuous functions,
the first class includes the functions that can be represented as the pointwise limits of continuous functions, the second class includes functions
representable as the pointwise limits of the functions of the first class, etc. The vector space of {\it bounded} functions that are not higher than the first
class endowed with  norm (\ref{normprav}) is denoted by ${\bf{N_1}}[a,b].$ Let us list some properties of R. Baire classes.

1.\;{\it The limit of a uniformly convergent sequence of functions of class at most some $k$ is a function of class at most } $k$ (Theorem 4 in \cite[XV.1]{nat97}).

2.\;{\it For a function $x:[a,b]\to \mathbb R$ to be a function not higher than the first class it is necessary and sufficient 
that} for any real $\;y\;$ the Lebesgue sets
$$
E_y^{-}(x)\doteq\{t:x(t)<y\},\quad E_y^{+}(x)\doteq\{t:x(t)<y\}
$$
{\it were the sets of type} $F_{\sigma}$ (Lebesgue's theorem \cite[XV.1]{nat97}).

3.\;{\it Characteristic function of a closed set $F\subset [a,b]$
$(F\ne [a,b])$ is a function of the first class } \cite[XV.3]{nat97}.

4.\;{\it A function that has a finite $(\geqslant 1)$ or a countable set of points of discontinuity is a function 
of the first class} \cite[XV.3]{nat97}.

Let us consider some examples.
\begin{example}
\label{Bklass1}
Let $\mathcal K$~ be Cantor's perfect set.
The set of points of the second kind (that is, the Cantor set without the endpoints of the adjacent intervals) will be denoted by $\mathcal S.$ Let
$$
w_{1}(t)\doteq \chi_{\mathcal K}(t),\;\;w_{2}(t)\doteq \chi_{\mathcal S}(t).
$$
Then \;$T\bigl(w_{1}\bigr)=T\bigl(w_{2}\bigr)=\mathcal K,\;$
therefore $w_{1},\;w_{2}\notin {\bf{N}}[0,1]$ (property 4),
and $w_{1}\in {\bf{N_1}}[a,b]\setminus {\bf{N}}[a,b]$ (property 3).
\end{example}
\begin{example}
\label{Bklass2}
Let $\mathcal K_{\alpha}$~ be a perfect set of Lebesgue measure
$\alpha>0$ (see \cite[c. 101]{derr08}). We set $w_{3}(t)\doteq \chi_{\mathcal K_{\alpha}}(t).$ Then $T(w_{3})=\mathcal K_{\alpha}.$
Consequently, $w_{3}\in {\bf{N_1}}[a,b]\setminus {\bf{N}}[a,b]$ (property 3), and by Lebesgue's R-integrability criterion
 $w_{3}\notin \setminus {\bf{Rm}}[a,b].$
\end{example}

Property 4 implies inclusion ${\bf{N}}[a,b]\subset {\bf{N_1}}[a,b]$. Function $w_1(\cdot)$ from Example \ref{Bklass1}
shows that this inclusion is strict. Functions
$$
w_2\in {\bf{Rm}}[a,b]\setminus {\bf{N_1}}[a,b],\;\;w_3\in {\bf{N_1}}[a,b] \setminus {\bf{Rm}}[a,b]
$$
show the relationship between classes ${\bf{N_1}}[a,b]$ and ${\bf{Rm}}[a,b].$

Let us show that ${\bf{N_1}}[a,b]$ is a Banach space.

Let $\{x_n\}_{n=1}^\infty$~ be a fundamental sequence from this space. This sequence is also fundamental in the Banach space
${\bf{M}}[a,b]$ of bounded functions. Therefore, there exists function $x\in {\bf{M}}[a,b]$ such that $x_n(t)\rightrightarrows x(t).$
By property 1, $x\in {\bf{N_1}}[a,b].$
\hfill $\square$


\medskip

\poonkt{Dense sets in ${\bf{N_1}}[a,b]\;$} \ \phantom{01234567890123456789}

Let us call a {\it $F_{\sigma}$-partition} of interval $[a,b]$ a collection
$$
J_m\doteq \{B_1,B_2,\ldots,B_m\}
$$
of sets of type $F_{\sigma}$ (see Sections {\bf{1.9.2}} and {\bf{1.9.3}}) provided that
$$
[a,b]=\bigcup\limits_{k=1}^m B_k,\qquad B_j\bigcap B_k=\varnothing\;\;\;\;\;\text{if}\;\;\;\;\; j\ne k.
$$
We call function $x:[a,b]\to\mathbb R\;$ a {\it $F_{\sigma}$-step} function if it
is constant on every set $B_k$ from some $F_{\sigma}$-partition $J_m.$
The set of $\;F_{\sigma}$-step functions with norm \; (\ref{normprav})\; will be denoted by
$\mathcal F_1[a,b] \bigl(\supset \mathcal H[a,b]\bigr).$

Each $F_{\sigma}$-partition $J_m=\{B_1,B_2,\ldots,B_m\}$
corresponds to $m$-
dimensional vector space $\mathcal J_m$ of $F_{\sigma}$-step functions.
A function from $\mathcal J_m$ has representation $x(t)=\sum\limits_{k=1}^m b_k\chi_{B_k}(t),$ where 
$b_k$ is  the value of $x(\cdot)$ on set
$B_k\;(k=1,2,\ldots,m).$

What does a sequence of functions from $\mathcal F_1[a,b]$ look like? 
There is a sequence of $F_{\sigma}$-partitions
$$
\{J_{m_n,n}\}_{n=1}^{\infty}=\{B_{1,n},\,B_{2,n},\ldots,B_{m_n,n}\},
$$
some of which (or even all) may be the same. The required sequence is
$$
x_n(t)=\sum\limits_{k=1}^{m_n} b_{k,n}\chi_{B_{k,n}}(t).
$$
When increasing the number of the term of a $F_{\sigma}$-sequence, the partition may change or may stay the same.
\begin{teo}
\label{plotF1vN1}
The set $\mathcal F_1[a,b]$ is {\it dense in } ${\bf{N_1}}[a,b].$
\end{teo}

\doc\,
The proof of this theorem (in some other terms) is essentially scattered in book \cite[XV.3]{nat97}.
We present here the main arguments in order to make it easier for the reader to follow the proof.

Let $x\in {\bf{N_1}}[a,b].$ According to property 2 of the previous subsection,
for any real $y$ the Lebesgue sets
$$
E_y^{-}(x)=\{t:x(t)<y\}\quad \text{and}\quad E_y^{+}(x)=\{t:x(t)>y\}
$$
are sets of type $F_{\sigma}.$
Divide interval $\bigl[-\|x\|,\,\|x\|\bigr]$ into $n$ equal parts by fixing points
 $$
 y_0=-\|x\|<y_1<y_2<\ldots<y_n=\|x\|,\quad \left(y_{k+1}-y_k=\dfrac{2\|x\|}{n }\right)
 $$
 and put
 $$
 A_0=E_{y_1}^{-},\;\;A_k=E_{y_{k-1}}^{+}\bigcap E_{y_{k+1}}^{-}\;(k= 1,2,\ldots,n-1)\;\;A_n=E_{y_{n-1}}^{+}.
 $$
 As was already noted, these are sets of type $F_{\sigma}$, and $[a,b]=\bigcup\limits_{k=0}^n A_k.$ 
 According to statement 5 from subsection {\bf{1.9.3}}
 there exists an $F_{\sigma}$-partition \linebreak
 $\{B_0,B_1,\ldots,B_n\}$ such that $B_k\subset A_k\;\;(k=0,1,\ldots,n).$
 We put $x_n(t)=y_k$ for $t\in B_k\;\;(k=0,1,\ldots,n).$ Then $x_n\in \mathcal F_1[a,b],$ and for $t\in [a,b],\;t\in B_k\subset A_k$
 $x_n(t)=y_k,\;\;y_{k-1}<x(t)<y_{k+1}.$ It follows that
 $$
 |x_n(t)-x(t)|<\dfrac{2\|x\|}{n},
 $$
 that is, the sequence of $F_{\sigma}$-step functions $x_n(\cdot)$ converges uniformly to $x(\cdot)$ as $n\to\infty$
 \hfill $\square$

If we denote $\mathcal F[a,b]\doteq \mathcal F_1[a,b]\bigcap {\bf{N}}[a,b],$ then we arrive at the following statement.
\begin{sle}
\label{plotFvN}
Set $\mathcal F[a,b]$ is {\it dense in} space  ${\bf{N}}[a,b].$
\end{sle}

\normalsize

In conclusion, we note that {\it there are strict inclusions}
\begin{multline}
\label{vkl2}
\underline{{\bf{H}}[a,b]\subset{\bf{BV}}[a,b]}\subset{\bf{CH}}[a,b]\subset\\
\subset\underline{{\bf{R}}[a,b]\subset{\bf{N}}[a,b]\subset{\bf{Rm}}[a,b]\subset{\bf {M}}[a,b],}
\end{multline}
where {\it underlined inclusions are embeddings of Banach spaces.}


\newpage

\begin{center}
\begin{Large}
{\bf{Chapter V. \; $RS$-integrability}}
\end{Large}
\end{center}
\addcontentsline{toc}{section}{Chapter V. \; $RS$-integrability.\; Exact $RS$-pairs}

\begin{center}
\begin{large}
\section{Criterion for $RS$-integrability}
\end{large}
\end{center}


\bigskip

\poonkt{Lebesgue's theorem}

\begin{teo}[Analogue of Lebesgue's theorem]
\label{ch10thm1+}
{\it Let $f$~ be a bounded function on $[a,b]$ $(f\in M[a,b]),$ let $g$ be an increasing function.
The integral
\begin{equation}
\label{rs}
(RS)\!\int\limits_a^b f(t)\,dg(t),
\end{equation}
exists if and only if the set $T(f)$ of points of discontinuity of function~$f$ is a set of $g$-measure zero.}
\end{teo}
\doc\,
{\bf{Sufficiency.}}
Let $\Omega (f)$ be the oscillation of $f(\cdot)$ on the whole interval $[a,\,b]$ and let $\varepsilon >0$ be arbitrary. Let us cover
$T(f)$ by such open set $G=\bigcup (a_k,\,b_k)$ that
$\mu_g(G)<\frac{\varepsilon}{2\Omega(f)}$ (see Section 3.6).
The set $F=[a,\,b]\setminus G$ is closed and $f(\cdot)$ is continuous, and thus uniformly continuous
on $F.$ Therefore, there exists $\delta >0$ such that for every $t_0\in F$ the oscillation
$\omega (f)$ on set $[t_0-\delta,\,t_0+\delta]\bigcap F$ satisfies
inequality $\omega (f)<\frac{\varepsilon}{2(g(b)-g(a))}.$

Consider an arbitrary partition $\tau =\{t_k\}_{k=0}^n$ of interval $[a,\,b]$
having diameter $d(\tau)<\delta$. Let $\mathcal J'_1,...,\mathcal J'_p$ be
the intervals $[t_{k-1},\,t_k]$ of this partition that contain at least one
point of $F,$ and let $\mathcal J''_1,\dots, \mathcal J''_q\;(p+q=n)$ be
the remaining intervals of partition $\tau$. The intervals of the second type are contained in $G.$

Since the set of points of discontinuity of $g(\cdot)$ is at most countable, we may assume without loss of generality that all $t_k$~ are
points of continuity $g(\cdot).$ Let us split the sum
$S\doteq\sum\limits_{k=1}^n \omega _k(f)\triangle g_k=\sum\limits_{k=1}^n \omega _k(f)\mu _{g}(\mathcal J_k)$
into two parts:
$$
S=\sum\limits_{i=1}^p \omega'_i(f)\mu_{g}(\mathcal J'_i)+ \sum\limits_{j=1}^q \omega _j''( f)\mu_{g}(\mathcal J''_j).
$$
($\omega'_i(f)\;(\omega''_j(f))$ are the oscillations on the intervals of the first (the second) type). By our construction,
$$
S<\frac{\varepsilon}{2(g(b)-g(a))}\sum\limits_{i=1}^p \mu_{g}(\mathcal J'_i)+
\Omega(f)\sum\limits_{j=1}^q \mu_{g}(\mathcal J''_j)<\frac{\varepsilon}{2}+\frac{\varepsilon}{2}= \varepsilon.
$$
By Theorem 6.3', integral (\ref{rs}) exists.

{\bf{Necessity.}}
Assume that integral\; (\ref{rs}) exist. By Theorem 6.8,
$T(g)\bigcap T(f)=\varnothing.$

Let $j\in \mathbb N$ be arbitrary, $T_j(f)=\{t:\,|\sigma_t(f)|>\frac{1}{j}\}.$
By Theorem 6.3', given an arbitrary $\varepsilon >0$, there exists $\delta >0$
such that for any partition $\tau$ whose diameter $d(\tau)$ is less than $\delta$ we have inequality
\begin{equation}\label{crit}
\sum\limits_{k=1}^n \omega _k(f)\triangle g_k<\frac{\varepsilon}{j}.
\end{equation}
Let $\tau$ be one of such partitions
and let $\mathcal J_1,...,\mathcal J_p$ be the intervals of this partition containing
the points of  set $T_j(f)$ as their interior points.
It may turn out that some of the points of partition $\tau$ also belong to $T_j(f).$
However, $\mu_{g}\{t_{k_1},\ldots,t_{k_l}\}=0$ since $t_{k_i}$ are points of continuity of $g(\cdot).$
We obtain from (\ref{crit})
$$
\frac{\varepsilon}{j}>\sum\limits_{k=1}^n \omega _k(f)\triangle
g_k >\frac{1}{j}\sum\limits_{i=1}^p\mu_{g}(\mathcal J_i)=\frac{1}{j}\sum\limits_{i=1}^ p\mu_{g}(\int\,\mathcal J_i)=
$$
$$
=\frac{1}{j}\mu_{g}\left(\bigcup\limits_{i=1}^p \int\,\mathcal J_i\right)\;\;\;\Rightarrow\;\;\mu_ {g}\left(\bigcup\limits_{i=1}^p \int\,\mathcal J_i\right)<\varepsilon.
$$

Thus, set $A_j=T_j(f)\setminus \{t_{k_1},...,t_{k_l}\}$ is covered
by open set $$G_j=\bigcup\limits_{i=1}^p int\,\mathcal J_i,\;
\mu_{g}(G_j)<\varepsilon.$$ Consequently,
$\mu_{g}(A_j)=\mu_{g}\bigl(T_j(f)\bigr)=0$
$(j\in \mathbb N).$ \; Since
$T(f)=\bigcup\limits_{j=1}^{\infty} T_j(f),$ we have $\mu_{g}T(f)=0.$
\hfill $\square$

\begin{sle}
\label{ch8corolthm21}
{\it Let $f(\cdot)$~ be a bounded function on $[a,b]$ \linebreak
$(f\in {\bf{M}}[a,b]),$ and let $g(\cdot)$ be a function of bounded variation $(g\in {\bf{BV}}[a,b]). $
For the existence of the integral $($\ref{rs}$)$
it is necessary and sufficient that the set $T(f)$ of points of discontinuity of function~$f(\cdot)$ be a set of $g_{\pi}$-measure zero.}
\end{sle}

\doc\,
According to Theorem \ref{ch10thm1+} and the result of section  4.6.2, condition $\mu_{g_{\pi}}\bigl(T(f)\bigr)=0$ is equivalent 
to the existence of integrals \\
$(RS)\!\int\limits_a^b f(t)\,dg_{\pi}(t)\;$ and $(RS)\!\int\limits_a^b f(t)\,dg_{\nu }(t),$ and thus of integral (\ref{rs}).
\hfill $\square$

\begin{sle}
\label{ch8corolthm22}
{\it Let $f\in{\bf{N}},\;g\in {\bf{CBV}}[a,b]$. Then integral $($\ref{rs}$)$ exists.}
\end{sle}
(In other words: {\it $\sigma$-continuous functions are $RS$-integrable with respect to continuous functions of bounded variation,} or
{\it $({\bf{N}},\,{\bf{CBV}})$ and $({\bf{CBV}},\;{\bf{N}})$~ are \textit{pairs of existence classes} for the Riemann ---  Stiltjes integral.})

\doc\,
The set of points of discontinuity $T(f)$ is at most countable and 
hence, by Theorem \ref{konschet1}, is a set of $g_{\pi}$-measure zero.
\hfill $\square$


\bigskip

\poonkt{An extension of the Lebesgue theorem\;}\, \

The following result is an extension of Theorem \ref{ch10thm1+} which allows to treat unbounded 
integrated functions.
\begin{teo}
\label{ch10thm2+}
{\it Let $g(\cdot)$ be an increasing function. For the
existence of integral $($\ref{rs}$)$
it is necessary and sufficient that the following conditions hold:

1)\;$f(\cdot)$ is bounded on the complements to a finite number of intervals where
which $g(\cdot)$ is constant; if $g(\cdot)$ is strictly increasing, then
$f(\cdot)$ is bounded on interval} $[a,\,b];$

2)\; $\mu_{g}\bigl(T(f)\bigr)=0.$
\end{teo}
\doc\,
{\bf{Sufficiency.}}
First, let $f(\cdot)$ be bounded on $[a,\,b],$ and let $\Omega (f)$ denote the oscillation of
$f(\cdot)$ on $[a,\,b].$ Let $\varepsilon >0\;$ be arbitrary. Let us cover
$T(f)$  by open set $G=\bigcup (a_k,\,b_k)$ such that
$\mu_{g}(G)<\frac{\varepsilon}{2\Omega(f)}.$ The complement $F=[a,\,b]
\setminus G$ is closed, $f(\cdot)$ is continuous on $[a,b]$, so it is uniformly continuous
on $F.$ Therefore, there exists $\delta >0$ such that for any $t_0\in F$ the oscillation
$\omega (f)$ on set $[t_0-\delta,\,t_0+\delta]\bigcap F$ satisfies
inequality $\omega (f)<\frac{\varepsilon}{2(g(b)-g(a))}.$

Consider an arbitrary partition $\tau =\{t_k\}_{k=o}^n$ of interval $[a,\,b]$
having diameter $d(\tau)<\delta.$ Let $\mathcal J'_1,...,\mathcal J'_p$ be the
intervals $[t_{k-1},\,t_k]$ of this partition that contain at least one
point of $F,$ and let $\mathcal J''_1,..., \mathcal J''_q\;(p+q=n)$ be
the remaining intervals of partition $\tau.$ The intervals of the second type are contained in $G.$

Since the set of points of discontinuity of $g(\cdot)$ is at most countable,
we can, without loss of generality, assume that all $t_k$~ are points of continuity of
$g(\cdot).$ We split the sum $S=\sum\limits_{k=1}^n \omega _k(f)\triangle g_k=
\sum\limits_{k=1}^n \omega _k(f)\mu _{g}(\mathcal J'_k)$ in two:
$$
S=\sum\limits_{i=1}^p \omega'(f)\mu_{g}(\mathcal J'_i)+ \sum\limits_{j=1}^q \omega _j''(f )\mu_{g}(\mathcal J''_j).
$$
($\omega'_i(f)\;(\omega''_j(f))$ are the oscillations of $f$ on the intervals of the first (the second) type). By our construction,
$$
S<\frac{\varepsilon}{2(g(b)-g(a))}\sum\limits_{i=1}^p \mu_{g}(\mathcal J'_i)+
\Omega(f)\sum\limits_{j=1}^q \mu_{g}(\mathcal J''_j)<\varepsilon.
$$
By Theorem 6.3, integral (\ref{rs}) exists.

Let $f(\cdot)$ be unbounded on $[a,\,b]$ but bounded on intervals
$\mathcal J_1,...,\mathcal J_m,\\
\mathcal J_k=[a_k,b_k],\;k=1,..,m,$
which are in the complement of the finite union of intervals where
$g(\cdot)$ is constant. As was proved above, integrals
$(RS)\!\int\limits_{a_k}^{b_k} f(t)\,dg(t)$ exist, and integrals over
intervals where $g(\cdot)$ is constant are equal to zero regardless of
the values of function $f(\cdot)$ on these intervals. Therefore, in this case too,
integral (\ref{rs}) exists and is equal to
$\sum\limits_{k=1}^m (RS)\!\int\limits_{a_k}^{b_k} f(t)\,dg(t).$

{\bf{Necessity.}}
Assume that integral (\ref{rs}) exist. Then, by Theorem 3.3, for
every $\varepsilon>0$ there exists $\delta>0$ such that, for any partition
$\tau=\{t_k\}_{k=1}^{n}$ of interval $[a,\,b]$ having diameter $d(\tau)<\delta,$
one has inequality
\begin{equation}
\label{crit2}
\sum\limits_{k=1}^n \omega _k(f)\triangle g_k<\varepsilon.
\end{equation}

Let $\tau$ be any of these partitions. Suppose that $f(\cdot)$ is not
bounded in a neighborhood of point $t'\in [t_{k_0-1},\,t_{k_0}].$
Then there is a sequence $t'_k\to t'$ such that $f(t'_k)\to \infty.$
Therefore, $ \omega _{k_0}(f)=+\infty$ and hence $\triangle g_{k_0}=0,$
i.e., since the partition is arbitrary, $g(\cdot)$ is constant on some
interval. Thus, if $f(\cdot)$ is unbounded on $[a,\,b],$ then
there exists a finite system of intervals $\mathcal J_1,...,\mathcal J_m$ where $g(\cdot)$ is constant, and on the complement to these
intervals function $f(\cdot)$ is bounded. Consequently, condition 1) of the theorem is satisfied.

Since integral (\ref{rs}) exists, by Theorem 3.8 $T(g)\bigcap T(f)=\varnothing.$

Let $j\in \mathbb N$ be arbitrary, $T_j(f)=\{t;\,|\sigma_t(f)|>\frac{1}{j}\}.$
Given an arbitrary $\varepsilon >0$ we find $\delta >0$
such that for any partition $\tau$ whose diameter is $d(\tau)<\delta$
an inequality of the form (\ref{crit2}) is satisfied, where the right side $\varepsilon$
is replaced by $\frac{\varepsilon}{j}.$ Let $\tau$ be one of these partitions
and let $\mathcal J_1,...,\mathcal J_p$ be the intervals of this partition containing
points of set $T_j(f)$ as their interior points.
It can happen that some of the partition points $\tau,$ for example $t_1,...,t_l,$
also belong to $T_j(f).$
Then, as noted above, $\mu_{g}\{t_1,...,t_l\}=0,$ since
$t_{k_i}$ are points of continuity of $g(\cdot).$
From (\ref{crit2}) we get
$$
\frac{\varepsilon}{j}>\sum\limits_{k=1}^n \omega _k(f)\triangle g_k >\frac{1}{j}\sum\limits_{i=1}^p \mu_{g}(\mathcal J_i)=
\frac{1}{j}\sum\limits_{i=1}^p\mu_{g}(\int\,\mathcal J_i) =\frac{1}{j}\mu_{g}\left( \bigcup\limits_{i=1}^p\int\,\mathcal J_i\right).
$$
Thus, set $A_j=T_j(f)\setminus \{t_{k_1},...,t_{k_l}\}$ is covered by
open set $$G_j=
\bigcup\limits_{i=1}^p \int\,\mathcal J_i,\;\mu_{g}(G_j)<\varepsilon.$$
Therefore, $\mu_{g}(A_j)=0.$ Since $T(f)=\bigcup\limits_{j=1}^{\infty} T_j(f),$ then $\mu_{g} T(f)=0.$
\hfill $\square$

Note that the above proofs of Theorems \ref{ch10thm1+} and \ref{ch10thm2+} basically go back to \cite[13.16]{hild63}.


\hskip 1mm
\poonkt
{Limit theorems for RS-pair $({\bf{N}},{\bf{CBV}})\;$}\phantom{0123456789} \

It is clear from the proofs of Theorems \ref{ch4thm1sr} and \ref{ch4thm2} that for pair \linebreak
({\bf{N}},{\bf{CBV}}) assertion (\ref{ch41}) of Theorem \ref{ch4thm1sr} and estimate (\ref{ch43}) hold.

\begin{teo}
\label{Fst3}
{\it Let $f_n\in{\bf{N}}[a,b],\,n\in\mathbb N,\;f_n(t)\rightrightarrows f(t) 
(n\to \infty),\linebreak 
g\in{\bf{CBV}}[a,b]$. Then}
\begin{gather} \label{Fst31}
\underset{n\to\infty}{\lim}\int_a^b f_n(t)\,dg(t)=\int_a^b f(t)\,dg(t).
\end{gather}
\end{teo}

The proof follows the proof of Theorem
\ref{ch4thm3}
\hfill $\square$
\begin{teo}
\label{RS2}
{\it Let
$g\in{\bf{CBV}}[a,b],\;\;f_n\in{\bf{N}}[a,b],\;\;n\in\mathbb N,\linebreak
f_n(t)\to f(t),\;
t\in [a,b]\;(n\to \infty),\;f\in{\bf{N}}[a,b],\;
\|f_n\|\leqslant M\;(n\in\mathbb N),$
for some $M$. Then relation} (\ref{Fst31}) holds.
\end{teo}
\doc\;
For a strictly increasing function $g(\cdot)$, we argue as in the proof of Theorem \ref{doparcel}. Since in this case functions
$f_n\bigl(g^{-1}(s)\bigr)\;(n\in\mathbb N)$ are obviously $R$-integrable, all arguments of the proof of the cited theorem go through.
In the general case, we represent $g(\cdot)$ as the difference of two strictly increasing functions (see Corollary \ref{raznsrogovozr}).
\hfill $\square$
\begin{teo}\label{RS3}
{\it Let $g_n\in{\bf{N}}[a,b],\;n\in\mathbb N,\;g_n(t)\rightrightarrows g(t)\; 
(n\to \infty),\linebreak
f\in{\bf{CBV}}[a,b];$ then}
\begin{equation}
\label{RS3b}
\underset{n\to\infty}{\lim}\int_a^b f(t)\,dg_n(t)=\int_a^b f(t)\,dg(t).
\end{equation}
\end{teo}
\doc\, According to Theorems
\ref{ch3thm2} and \ref{Fst3}:\;\;$\int_a^b f(t)\,dg_n(t)=\\
=f(t)g_n(t)\left|_a^b-\int_a^b g_n(t)\,df(t)\to f(t)g(t)\right|_a^b-\int_a^ b g(t)\,df(t)=\int_a^b f(t)\,dg(t)\linebreak
(n\to\infty).
$
\hfill $\square$

We represent function $g\in{\bf{CBV}}[a,b]$ in the form (\ref{ch28}). Then
\begin{equation}
\label{lsgpi1}
\mu_{g_{\pi}}(G)=\sum\limits_k\bigvee\limits_{a_k}^{b_k}(g)\qquad \bigl(G=\bigcup\limits_{k} (a_k,\, b_k),
\;\;\text{intervals do not intersect}\bigr).
\end{equation}

\begin{teo}[{Helly's theorem for pair $({\bf{N}},\,{\bf{CBV}})$}]
\label{Chelli2}
{\it Let \;
$$
f\in{\bf{N}}[a,b],\;g_n\in{\bf{CBV}}[a,b]\;(n\in\mathbb N),\;g_n(t) \to g(t)\;\;(t\in [a,b],\;n\to\infty),
$$
$$
g\in {\bf{C}}[a,b];\;\;\bigvee\limits_a^b(g_n)\leqslant V\;\;(V>0).
$$
Then $(\ref{RS3b})$ holds.}
\end{teo}
\doc\,
In the same way as in the proof of Theorem
\ref{ch4thm4}, we show that $g\in{\bf{BV}}[a,b],$ and hence $g\in{\bf{CBV}}[a,b]$ and
$\bigvee\limits_a^b(g_n)\leqslant V.$

Define function $\widehat g:\Delta\to [0,V],\;\;\Delta\doteq\{(t,s):0\ls s\ls t\ls V\},$
$$
\widehat g(s,\,s)=0,\;\;\;\;\text{and for}\;\;s<t\qquad
\widehat g(t,\,s)\doteq\underset{n\in\mathbb N}{\sup}\left\{\bigvee\limits_s^t(g),\,\bigvee\limits_s^t(g_n )\right\}.
$$
Obviously, $\widehat g$ is continuous on $\Delta,$ and hence uniformly continuous on $\Delta$.
For $\widehat g$ the following statements are true:
 
a)\;{\it for fixed $s,\;0\ls s\ls V$ $\widehat g$ is increasing in $t;$ for fixed $t,\;0\ls t\ls V$ $\widehat g$ is decreasing in} $s;$

b)\;{\it for any $\varepsilon>0$ and $c\in (a,b)\;\;\bigl(c=a,\;\;\;c=b\bigr)$ there exists interval $(\alpha,\beta)\ni c$
$\bigl($half-interval $[a,\beta),\;(\alpha,b]\bigr)$
such that} $\widehat g(\beta,\,\alpha)<\varepsilon\;\;\bigl(\widehat g(\beta,\,a)<\varepsilon,\linebreak
\widehat g(b,\,\alpha)<\varepsilon\bigr).$

Let $\Omega (f)$ be the oscillation of $f$ on the whole interval $[a,\,b]$, let $\gamma >0$ be arbitrary and 
$T(f)=\{c_1,c_2,\ldots\}.$ According to a) and b)
and equality $\underset{(t',s')\to (s,s)}{\lim}\widehat g(t',s')=\widehat g(s,s)=0$ (see also \ref{lsgpi1})
there is interval $(a_k,\,b_k)\ni c_k,$ such that $\widehat g(b_k,\,a_k)<\frac{\gamma}{2^k\cdot\Omega(f)},\ ;\;k=1,2,\ldots$.
Thus, $T(f)$ is covered by open set $G=\bigcup\limits_{k=1}^\infty (a_k,\,b_k),$ so that
\begin{gather} \label{ch21}
\mu_{g_\pi}(G)=\sum\limits_{k=1}^\infty\bigvee\limits_{a_k}^{b_k}(g)\ls\sum\limits_{k=1}^\infty\widehat g(b_k,\,a_k)<\frac{\gamma}{\Omega(f)},
\end{gather}
as well as
\begin{gather} \label{ch22}
\mu_{g_{n,\pi}}(G)=\sum\limits_{k=1}^\infty\bigvee\limits_{a_k}^{b_k}(g_n)\ls\sum\limits_{k= 1}^\infty\widehat g(b_k,\,a_k)<\frac{\gamma}{\Omega(f)}.
\end{gather}
$\Bigl($If $c_{k_0}=a\;\bigl(c_{k_1}=b\bigr),$ then we have half-interval $[a,\,b_{k_0})\;\bigl( (a_{k_1},\,b]\bigr);
$
$G$ is relatively open in this case.$\Bigr)$
Set $F=[a,\,b]\setminus G$ is closed and $f(\cdot)$ is continuous, and therefore uniformly continuous
on $F.$ So, there exists $\delta >0$ such that for any $t_0\in F$ the oscillation $\omega (f)$ 
of function $f(\cdot)$ on  set $\left[t_0-\frac{\delta} {2},\,t_0+\frac{\delta}{2}\right]\bigcap F$
satisfies inequality $\omega (f)<\gamma.$

We fix an arbitrary partition $\tau =\{t_k\}_{k=0}^m$ of interval $[a,\,b]$
having diameter $d(\tau)<\delta.$ Let $\mathcal J'_1,...,\mathcal J'_p$ be the
intervals $[t_{k-1},\,t_k]$ of this partition that contain at least one
point of $F,$ and let $\mathcal J''_1,..., \mathcal J''_q\;(p+q=m)$ are
be the remaining intervals of partition $\tau$. The intervals of the second type are contained in $G.$

Then
\begin{multline*}
\left|\gS_\tau(f,\,g)-\int_a^b f(t)\,dg(t)\right|=\left|\sum\limits_{k=1}^m f(\xi_k) \Delta g_k-\sum\limits_{k=1}^m\int_{t_{k-1}}^{t_k} f(t)\,dg(t) \right|=\\=
\left| \sum\limits_{k=1}^m \int_{t_{k-1}}^{t_k}(f(\xi_k)-f(t))\,dg(t) \right|\ls \sum \limits_{k=1}^m \int_{t_{k-1}}^{t_k}|f(\xi_k)-f(t)|\,dg_\pi(t) \ls \\
\ls \sum\limits_{k=1}^m \omega_{\tau k}(f) \int_{t_{k-1}}^{t_k}\, dg_\pi(t)=\sum\limits_ {k=1}^m \omega_{\tau k}(f) \Delta g_{\pi k}\doteq S.
\end{multline*}
Let us split $S$ into two sums:
$$
S=\sum\limits_{i=1}^p \omega'_i(f)\mu_{g_{\pi i}}(\mathcal J'_i)+ \sum\limits_{j=1}^q \omega _j''(f)\mu_{g_{\pi j}}(\mathcal J''_j).
$$
($\omega'_i(f)\;\bigl(\omega''_j(f)\bigr)$ are the oscillations of function $f(\cdot)$ on the intervals of the first (the second) type). 
By our construction, taking into account (\ref{ch21}),
$$
S<\gamma\sum\limits_{i=1}^p \mu_{g_{\pi}}(\mathcal J'_i)+\Omega(f)\sum\limits_{j=1}^q \mu_ {g_{\pi}}(\mathcal J''_j)
<\left(\gamma\bigvee\limits_a^b(g)+\gamma\right)=\gamma (V+1).
$$
Thus, partition $\tau$ satisfies the estimate
\begin{gather} \label{ch23}
\left|\gS_\tau(f,\,g)-\int_a^b f(t)\,dg(t)\right|<\gamma\cdot(V+1)
\end{gather}
Similarly, given (\ref{ch22}), we obtain that
\begin{gather} \label{ch24}
\left|\gS_\tau(f,\,g_n)-\int_a^b f(t)\,dg_n(t)\right|<\gamma\cdot (V+1)\quad (n\in\mathbb N ).
\end{gather}
Let $\varepsilon>0$ be arbitrary. In all previous arguments, we set $\gamma=\frac{\varepsilon}{3(V+1)}.$ Due to the pointwise convergence
$g_n(t)\to g(t)$, we have $\gS_\tau(f,\,g_n)\to \gS_\tau(f,\,g)\linebreak
(n\to\infty),$ so there is a positive integer $N$ such that for $n>N$ we have inequality
\begin{gather} \label{ch25}
|\gS_\tau(f,\,g_n)-\gS_\tau(f,\,g)|<\frac{\varepsilon}{3}.
\end{gather}
According to estimates (\ref{ch23})--(\ref{ch25}), for such $n$
$$
\gathered
\left| \int_a^b f(t)\,dg_n(t)-\int_a^b f(t)\,dg(t) \right|\ls
\left| \int_a^b f(t)\,dg_n(t)-\gS_\tau(f,\,g_n) \right|+ \\
+|\gS_\tau(f,\,g_n)-\gS_\tau(f,\,g)|+\left| \gS_\tau(f,\,g)-\int_a^b f(t)\,dg(t) \right|
< \frac{\varepsilon}{3}+\frac{\varepsilon}{3}+\frac{\varepsilon}{3}=\varepsilon.
\endgathered
$$
\hfill $\square$


\begin{center}
\begin{large}
\section{Exact $RS$-pairs}
\end{large}
\end{center}

\poonkt
{Exactness of pair $\bigl({\bf{Rm}},\,AC\bigr)$ } 
\begin{teo}
\label{RSabs1}
{\it If integral (\ref{rs})
exists for any 
Riemann integrable function $f(\cdot),$ then function $g(\cdot)$ is absolutely continuous on $[a,\,b].$}
\end{teo}
\doc\,
Arguing as in the proofs of Theorems \ref{ch3thm8} and \ref{ch3thm9},
we show that $g(\cdot)$ is a continuous function of bounded variation.

First, assume that $g(\cdot)$ is increasing and integral (\ref{rs}) exists
for any Riemann
integrable function $f(\cdot).$
Assume that $g(\cdot)$ is not absolutely continuous on
$[a,\,b].$ Then there exists $\varepsilon_0>0$ such that for any
$\delta >0$ there exists a finite system of non-intersecting intervals
$\{(a_k,\,b_k)\}_{k=1}^n$ such that that $\sum\limits_{k=1}^n (b_k-a_k)<\delta,$ but
$\sum\limits_{k=1}^n (g(b_k)-g(a_k))\ge \varepsilon_0.$
Consider function $f(t)=\sum\limits_{k=1}^n \sigma_k\,\mathfrak h_{c_k}(t),$ where
$\mathfrak h_c(\cdot)$~ is the identity function (see ~{\bf{4.5}}),\; $c_k=\frac{a_k+b_k}{2}$ and
$\sigma_k>1.$ Then $f(\cdot)$ is Riemann integrable (as a piecewise continuous function), and by Theorem \ref{ch3thm3}
for this $\varepsilon_0$ there exists $\delta_0 >0$ such that if
$\sum\limits_{k=1}^n (b_k-a_k)<\delta_0,$ then for any partition
$\tau=\{t_j\}_{j=1}^m,$ containing points $a_k,\ b_k\;(k=1,..,n)$
whose diameter is less than $\delta_0$ inequalities
$$
\varepsilon_0>\sum\limits_{j=1}^m \omega _j(f)\triangle g_j>\sum\limits_{k=1}^n (g(b_k)-g(a_k))\ge \varepsilon_0 
$$
are satisfied.
The resulting contradiction implies that $g(\cdot)$ is absolutely continuous.

The general case is reduced, with the help of Theorems \ref{ch2thm6} and \ref{par2thm4}, to what has already been proved.
\hfill $\square$

\begin{teo}
\label{RSabs2}
{\it If integral $(\ref{rs})$ exists for any absolutely
continuous function $g(\cdot),$ then $f(\cdot)$ is Riemann integrable.}
\end{teo}
\doc\,
Since $f(\cdot)$ is integrable with respect to any absolutely continuous function, it follows that
it is also integrable with respect to absolutely continuous function $g(t)\equiv t,$
that is, $f(\cdot)$ is Riemann integrable.
\hfill $\square$

\begin{teo}
\label{RSabs3}
{\it $RS$-pair $\bigl({\bf{Rm}}[a,\,b],\,AC[a,\,b]\bigr)$~ is an
exact pair of classes for the existence of integral } (\ref{rs}).
\end{teo}

\doc\,
The fact that $\bigl({\bf{Rm}}_{[a,\,b]},\,AC_{[a,\,b]}\bigr)$ is an $RS$-pair follows from Theorems \ref{ch8thm1}
and \ref{ch8thm2}. The exactness of this pair follows from Theorems \ref{RSabs1} and \ref{RSabs2}.


\bigskip

\poonkt
{Exactness of pair $\bigl({\bf{N}},{\bf{CBV}}\bigr)$\;}\phantom{012345678901234567890} \

Here we will have to use the definitions and notations from Section {\bf{3.6}}. In addition, let us extend  notation
$\mes\,A$ to finite or countable unions of arbitrary intervals from $[a,b]\;\;(b>a)$ (not only of open intervals,
as in {\bf{3.6}} ). Moreover, given an arbitrary set $A\subset [a,b]$ and a continuous increasing function $g(\cdot)$
we will say $\mu_g(A)>0$ if $\mu_g\bigl([a,b]\setminus A\bigr)=0.$ Obviously, if $A\subset B$ and $\mu_g(A) >0,$ then $\mu_g(B)>0.$

Recall that, according to Corollary \ref{ch8corolthm22}, the pair of function classes $\bigl({\bf{N}},{\bf{CBV}}\bigr)$ is an RS-pair.
To prove its exactness, we need a number of auxiliary results.

\begin{lem}
\label{4.l1}
{\it The set  $T(f)$ of points of discontinuity of function $f(\cdot),$ defined on a closed set, is a set of type} $F_{\sigma}.$
\end{lem}
\doc\,
Let $T_{n}=\{t:|\sigma_t(f)|\geqslant \frac{1}{n},\;n\in\mathbb N\}$~ be the set of points of discontinuity of 
function $f (\cdot)$ where the absolute value
of the jump is greater than or equal to $\frac{1}{n}.$ Let $t_*$~ be a limit point of this set. For every $n\in\mathbb N$ there exists
a sequence $\{t_{mn}\}_{m=1}^{\infty}\to t_*,\;$
$|\sigma_{t_{mn}}(f)|\geqslant \frac{1}{n}.$ Passing 
to a subsequence, if needed, we obtain that $|\sigma_{t_*}(f)|\geqslant \frac{1}{n}.$ Therefore, set $T_{n}$ is closed.
The assertion of the lemma now follows from representation $T(f)=\bigcup\limits_{n=1}^{\infty} T_{n}.$
\hfill $\square$

The proof of the following two general results can be found in \cite{nat97}.
\begin{utv}
\label{4.l2}(Cantor --- Bendixson theorem, see \;\cite[II.6]{nat97}.)
{\it An uncountable closed set $F$ 
can be represented in the form $F=P\bigcup D,$ where $P$~ is a perfect set, and $D$~ is only a countable set.}
\end{utv}

\begin{utv}
\label{4.l3}
{\it A non-empty perfect set has the cardinality of continuum.}\;(\;\cite[II.6]{nat97}.)
\end{utv}
\begin{lem}
\label{4.l4}
{\it A non-empty closed null-set $F$ is nowhere dense.}
\end{lem}
\doc\,
Let $(\alpha,\,\beta)$~ be an arbitrary interval. Since the length of this interval is $\beta-\alpha>0,$ and $\mes\,F=0,$ this interval contains points
not belonging to $F.$ Let $t_0$~ be such a point. The complement of set $F$~ is an open set, so there is a neighborhood of point
$t_0$ containing no points of set $F$. The intersection of this neighborhood with interval $(\alpha,\,\beta)$ 
is an interval that is free from points of set $F.$ Therefore, set $F$ is nowhere dense.
\hfill $\square$

We obtain from statements \ref{4.l1}--\ref{4.l4} the following lemma.

\begin{lem}
\label{4.l5}
{\it If set $T(f)$ of points of discontinuity of  
function $f(\cdot)$ is uncountable and $\mes\,T(f)=0,$ then there exists a perfect
nowhere dense set $F\subset T(f).$}
\end{lem}

\begin{lem}
\label{4.l6}
{\it Let $F$~ be a perfect nowhere dense
null-set.\;There is a continuous increasing function $g(\cdot)$
having the following properties:

a$)\; g^{'}(t)=1$ for points $t\in F;$

b$)\; g^{'}(t)=0$ for points $t\in [a,b]\setminus F;$}
\end{lem}
(Compare with \cite[VIII.2]{nat97}.)

\doc\,
For every $n\in \mathbb N$ we cover $F$ with  an open set $G_n\supset F$ such that
$\mes G_n<\frac{1}{2^n}$, and consider function
$s_n(t)\doteq \mes\big(G_n\bigcap[a,\,t]\big).$
This function has the following properties:

1) $s_n(t)\gs 0\ \ (t\in [a,\,b]);$ \; 2) $s_n(\cdot)$ is increasing; \;
3) $s_n(t)<\frac{1}{2^n};$ \; 4) $s_n(\cdot)$ is continuous on $[a,\,b].$

Properties 1) --- 2) hold due to non-negativity and monotonicity of the Lebesgue measure,
property 3) is true by the construction.

Let us prove property 4). Let $h>0.$
\begin{multline*}
s_n(t+h)-s_n(t)=\mes\big(G_n\bigcap [a,\,t+h]\big)-\mes\big(G_n\bigcap [a,\,t]\big )=\\
=\mes G_n\bigcap(t,\,t+h]\ls h,
\end{multline*}
which implies continuity on the right. Continuity on the left is proved similarly.

Let $t_0\in F.$ For a fixed $\;n\in \mathbb N\;$ and sufficiently small
$h>0$ we have $[t_0,\,t_0+h] \subset G_n.$ Hence
\begin{equation*}
s_n(t_0+h)=s_n(t_0)+h,\quad \dfrac{s_n(t_0+h)-s_n(t_0)}{h}=1.
\end{equation*}
We obtain an analogous equality for $h<0.$ Thus, $s'_n(t)=1$ for all
$t\in F.$

If $t_1\notin F$ for a fixed $n\in \mathbb N,$ then $t_1$ falls into one of the \textit{adjacent} intervals
\; (this interval is a \textit{component} of open set $G\doteq [a,b]\setminus F$),
and for sufficiently small $h>0\;$ interval $\;[t_1,\,t_1+h]$ lies entirely in this interval, 
so $s_n(t_1+h)-s_n(t_1)=\mes G_n\bigcap(t_1,\,t_1+h]=0,$ i.e. $s_n^{'}(t_1)=0 .$

Let
\begin{equation}
\label{rjad1}
s(t)\doteq \sum\limits_{n=1}^{\infty} s_n(t).
\end{equation}
Due to property 3), this series converges uniformly. Hence $s(\cdot)$ is a non-
negative continuous increasing function. Next, we set
\begin{equation}
\label{rjad2}
g_n(t)=\dfrac{1}{2^n}s_n(t)\;\;(n\in\mathbb N),\qquad g(t)\doteq \sum\limits_{n=1}^ {\infty} g_n(t).
\end{equation}
Then $g_n^{'}(t)=\frac{1}{2^n}$ for points $t\in F$, and $g_n^{'}(t)=0$ for points $t\ in [a,b]\setminus F\linebreak
(n\in\mathbb N).$
The series in (\ref{rjad2}) can be differentiated term-wise, since the series obtained by differentiation is majorized by a convergent 
series
$\sum\limits_{n=1}^\infty\frac{1}{2^n}\doteq S$ or real numbers, and hence converges uniformly \cite[c. 440]{ficht2}. Thus,
$g^{'}(t)=S=1$ for points $t\in F,\;\; g^{'}(t)=0$ for points $t\in [a,b]\setminus F, $ i.e. 
 function $g(\cdot)$ has properties a) and b).
\hfill $\square$

\begin{lem}
\label{4.l7}
{\it Let set $F$ and function $g(\cdot)$ be as in Lemma \ref{4.l6}.
Then} $\mu_g(F)>0.$
\end{lem}
\doc\,
We put $G=[a,b]\setminus F$. Then $G$~ is an open set, $G=\bigcup\limits_{n=1}^\infty (a_n,\,b_n).$
Due to property b),
$\;
g(b_n-)-g(a_n+)=g(b_n)-g(a_n)=g^{'}(c_n)\bigl(b_n-a_n\bigr)=0,
$
so $\mu_g(G)=0,$ i.e., according to the definition $\mu_g(F)>0.$
\hfill $\square$

\begin{teo}
\label{4.t1}
{\it RS-pairs $\bigl({\bf{N}},\,{\bf{CBV}}\bigr)$ and $\bigl({\bf{CBV}},\,{\bf{ N}}\bigr)$ are exact.}
\end{teo}
\doc\,
If the integrating function $g(\cdot)$ is discontinuous at a point $c_0\in[a,b],$ then integral (\ref{rs}) does not exist for
$f(t)=\mathfrak h_{c_0}(t),$ where function $\mathfrak h_c(\cdot)$ was defined in section {\bf{4.5}},
$\;f\in{\bf {N}},\;$ since the integrated and the
integrating functions in this case have a common point of discontinuity
(see Theorem \ref{ch3thm8}).

Assume that $g$ has infinite total variation, ${\bigvee\lt_a^b (g)=+\fy}.$ There is a
point $c\in[a,\,b]$ such that in any of its neighborhoods the total variation of $g(\cdot)$
is infinite. Without loss of generality, $c=b.$
Then we can construct such a sequence $\{t_n\}_{n=0}^\fy,$
$t=t_0<t_1<\ldots<t_n<b,$ $t_n\to b,$ that $\sum\lt_{k=1}^\fy \big|
g(t_k)-g(t_{k-1})\big|=+\fy.$ There is also a sequence of real numbers
$\{c_k\}_{k=1}^\fy$ such that $c_k\to 0,$ $c_k>0$
$(k=0,\,1,\,2,\,\ldots),$ but at the same time
$$
\sum\lt_{k=1}^\fy c_k \big| g(t_k)-
g(t_{k-1}) \big|=+\fy
$$
(one can take $c_n=1\Big/\sum\lt_{k=1}^n \big|
g(t_k)-g(t_{k-1})\big|,$ since 
the series $\sum\lt_{n=1}^\fy a_n/s_n$ (here $s_n=\sum\lt_{k=1}^n a_k$) diverges together with the series
$\sum\lt_{n=1}^\fy a_n$).

We define continuous function~$f(\cdot)$ in the following way: 
\begin{multline*}
f(t_k)=c_k \cdot \sign (g(t_k)-g(t_{k-1}))\;\;k=0,\,1,\,2,\,\ldots,\;f(b)=0,\\
f(t)=f(t_{k-1})+\frac{f(t_k)-f(t_{k-1})}{t_k-t_{k-1}}(t-t_{k -1})\;
(t_{k-1}<t<t_k),\; k=1,\,2,\,\ldots 
\end{multline*}
Then for $\tau_n=\{t_k\}_{k=0}^n\quad $
$
\gS_{\tau_n} (f,\,g)=\sum\lt_{k=1}^n f(t_k) (g(t_k)-g(t_{k-1}))= \\
=\sum\lt_{k=1}^n c_k \big| g(t_k)-g(t_{k-1}) \big| \to +\fy
$
as $n\to\fy,$ so integral (\ref{rs}) does not exist.
Thus, the class of integrating functions cannot be enlarged without narrowing the class of integrated functions.

Let $f\notin {\bf {N}}$. If $f\notin {\bf {Rm}},$ then $f(\cdot)$ is not integrable with respect to  absolutely 
continuous (that is, from {\bf{CBV}})
functions $g(t)\equiv t.$

Let $f\in {\bf {Rm}}\setminus {\bf {N}}$. Then $T(f)$ is an uncountable null-set, so by Lemma \ref{4.l5} there exists a perfect nowhere dense null-set $F\subset T(f),$ and by Lemma \ref{4.l6}
there exists an increasing continuous (that is, from $\mathbf{CBV}$) function $g(\cdot)$ such that, according to Lemma \ref{4.l7},
$\mu_g(F)>0,$ and hence $\mu_g\bigl(T(f)\bigr)>0.$
By virtue of Theorem \ref{ch10thm1+}, integral (\ref{rs}) does not exist.

So, the class of integrable functions also cannot be enlarged without narrowing the class of integrating functions.
\hfill $\square$

\begin{sle}
\label{4.s1}
{\it RS-pair $\bigl({\bf {R}},\,{\bf{CBV}}\bigr)$ is not exact.}
\end{sle}

to this concept
brought ideas for work

\newpage

\begin{center}
\begin{Large}
{\bf{Chapter VI. \; (*)-Riemann--Stieltjes integral}}
\end{Large}
\end{center}
\addcontentsline{toc}{section}{Chapter VI. \; (*)-Riemann--Stieltjes integral}

\begin{center}
\begin{large}
\section{Definition of (*)-integral}
\end{large}
\end{center}

\hskip 1mm
\poonkt{Definition and basic properties}\phantom{01234567890123456789} \

{\bf{1.}}\; (*)-integral was first introduced and studied in \cite{derr19}
(see also \cite{derr20}); some ideas from earlier works \cite{derr97},\, \cite{derr02},\, \cite{DK06} led to it.

Let $f\in{\bf{N}}{[a,b]},\;$ $g\in{\bf{BV}}{[a,b]},\;a<b.$ Put
\begin{multline}
\label{defali}
(*)\!\int\limits_a^b f(t)\,dg(t)\doteq\int_a^bf(t)\,dg_{\gc}(t)+\\
+\Bigl(f(a)\sigma_a^{+}(g)+\sum\limits_{t\in T(g)\bigcap(a,b)}f(t)\sigma_t(g)+f( b)\sigma_b^-(g)\Bigr).
\end{multline}
The first term in the right-hand side of (\ref{defali}) is the Riemann-Stieltjes integral over interval $[a,b].$
The existence of this integral is ensured by Corollary \ref{ch8corolthm22}. Since the total variation of  function 
$g$ is finite, the series $\sum\limits_{t\in T(g)\bigcap(a,b)}\sigma_t(g)$ converges absolutely.
The boundedness of function $f\;\;\Bigl($and, consequently, the boundedness of sequence
$\{f(t)\}_{t\in T(g)\bigcap(a,b)}\Bigr)\;$ yield the absolute  convergence of the series in (\ref{defali}).

If the integrable function is continuous ($f\in{\bf{C}}{[a,b]}$), then, by the proof of (\ref{ch47}),
\begin{gather}
\label{dop1}
(*)\!\int\limits_a^b f(t)\,dg(t)=(RS)\!\int\limits_a^b f(t)\,dg(t).
\end{gather}
Also, we set by definition
\begin{gather}
\label{dop2}
(*)\!\int\limits_a^a f(t)\,dg(t)=0.
\end{gather}


\newpage

{\bf{2.} The next properties follow directly from definition (\ref{defali}):}

a)\; The (*)-{\it integral is linear with respect to the integrated function $(f)$. 
It is also linear with respect to the integrating function $(g)$}.

b) Let $[c,d]\subset [a,b]\;(c<d)$. The integral $(*)\!\int\limits_c^d f(t)\,dg(t)$ corresponds to interval $[c,d]$ 
and is defined by formula (\ref{defali}), where $a\;(b)$
stands for $c\;(d)$. Of course, the continuous part of  function $g$ on interval $[c,d]$ differs from the continuous part
of this function on interval $[a,b]$ by a constant. This constant is equal to 
$$
m\doteq \sum\limits_{s<c,\,s\in T(g)}\sigma_s(g)+\sigma^- _c(g),
$$
however, this does not affect the value of the Riemann-Stieltjes integral
$\int_c^df(t)\,dg_{\gc}(t)$.

c)\; (*)-{\it integral is additive in the following sense:\; if\; $a<c<b,$ then}
\begin{gather} \label{add}
(*)\!\int\limits_a^b f(t)\,dg(t)=(*)\!\int\limits_a^c f(t)\,dg(t)+(*)\!\int\limits_c^b f(t)\,dg(t)
\end{gather}
(this is the additivity of the integral as a function of an interval). This fact requires a more detailed justification.
As was already noted in part c), we can assume that in the Riemann--Stieltjes integrals
$$
\int_a^b f(t)\,dg_{\gc}(t),\quad \int_a^c f(t)\,dg_{\gc}(t),\quad \int_c^b f(t)\,dg_ {\gc}(t)
$$
we have the same integrating function (actually, in the last integral the integrating function differs by a constant $m,$ but this does not affect the value of the integral). So,
$$
\int_a^b f(t)\,dg_{\gc}(t)=\int_a^c f(t)\,dg_{\gc}(t)+ \int_c^b f(t)\,dg_{\gc} (t);
$$
and since $\sigma_c(g)=\sigma_c^+(g)+\sigma_c^-(g),$
we have
$$
\sum\limits_{t\in T(g)}f(t)\sigma_t(g)=\sum\limits_{t\in T(g),\,t\leqslant c}f(t)\sigma_t( g)+\sum\limits_{t\in T(g),\,t\geqslant c}f(t)\sigma_t(g).
$$
The last two equalities are equivalent to statement (\ref{add}).


{\bf{3.}}
{\it The following estimates hold:} \;
\begin{equation}
\label{tonkgrub}
\left|(*)\!\int\limits_a^b f(t)\,dg(t)\right|\leqslant (*)\!\int\limits_a^b |f(t)|\,dg_{\ pi}(t)\leqslant \|f\|_{\bf{N}}\bigvee_a^b(g)
\end{equation}
$\Bigl(g_{\pi}(t)=\bigvee\limits_a^t(g),$ see (\ref{ch25+})$\Bigr).$

\doc\,
Due to estimate (\ref{ch43}) $\bigl($as was noted in Section {\bf{12.3}}, this estimate is also valid for the pair ({\bf{N}},{\bf{CBV}})$\bigr)$,
the representation (\ref{var22}) and the equality (\ref{var21}), we have
\begin{multline*}
\left|(*)\!\int\limits_a^b f(t)\,dg(t)\right|\leqslant\left|\int_a^bf(t)\,dg_{\gc}(t)\right |+\\
+\left|f(a)\sigma_a^{+}(g)+\sum\limits_{t\in T(g)\bigcap(a,b)}f(t)\sigma_t(g)+f( b)\sigma_b^-(g)\right|\leqslant \\
\int_a^b |f(t)|\,d\bigl(g_{\gc}\bigr)_{\pi}(t)+
\left(|f(a)||\sigma_a^{+}(g)|+\sum\limits_{t\in T(g)\bigcap(a,b)}|f(t)||\sigma_t (g)|+\right.\\
\left.\phantom{\sum\limits_{t\in T(g)\bigcap(a,b)}}+|f(b)||\sigma_b^-(g)|\right)=
(*)\!\int\limits_a^b |f(t)|\,dg_{\pi}(t)\leqslant \|f\|_{\bf{N}}\bigvee_a^b(g).
\end{multline*}
(We also use the obvious relations
$$
|\sigma_t(g)|=\sigma_t(g_{\pi}),\;\;|\sigma_a^{\pm}(g)|=\sigma_t(g_{\pi}^{\pm}), \;\;
\bigl(g_{\gc}\bigr)_{\pi}(t)=\bigl(g_{\pi}\bigr)_{\gc}(t).)\qquad \hfill \square
$$

 
\bigskip

\poonkt
{Integration by parts formula\;}\phantom{0123456780123456789} \

Let $f\in{\bf{CH}}[a,b],\;$ $g\in{\bf{BV}}[a,b]$. Then the expression in the right-hand side of equality 
\begin{multline}
\label{defali2}
(*)\!\int\limits_a^b g(t)\,df(t)\doteq\int_a^bg(t)\,df_{\gc}(t)+\\
+\Bigl(g(a)\sigma_a^{+}(f)+\sum\limits_{t\in T(f)\bigcap(a,b)}g(t)\sigma_t(f)+g( b)\sigma_b^-(f)\Bigr)
\end{multline}
is well defined, which justifies its designation as a *-integral.
\begin{teo}
\label{bypart}
\textit{Let $f\in{\bf{CH}}[a,b],\;$ $g\in{\bf{BV}}[a,b].$ Then the following identity holds $($the integration by parts formula$)$:
\begin{multline}
\label{bp1}
(*)\!\int\limits_a^b f(t)\,dg(t)+(*)\!\int\limits_a^b g(t)\,df(t)=\\
=f(t)g(t)\Big|_a^b-\sum\limits_{t\in T}\Bigl(\sigma^{+}_t(f)\sigma^{+}_t(g)-\sigma^{-}_t(f)\sigma^{-}_t(g)\Bigr),
\end{multline}
where } $T\doteq T(f)\bigcap T(g).$
\end{teo}
\doc\,
Throughout the proof, given any function, we extend it 
to $t<a$ by its value at point $a,$ and to $t>b$ by its value
at point $b$ (see convention (\ref{soglprodolz})). Then
$$
\sigma_a^{+}(x)=\sigma_a(x),\quad \sigma_b^-(x)=\sigma_b(x).
$$
Taking this into account, the out-of-integral terms in (\ref{defali})\, and \, (\ref{defali2}) can be written as
\begin{gather} \label{exd4}
\sum\limits_{t\in T}f(t)\sigma_t(g),\qquad \sum\limits_{t\in T}g(t)\sigma_t(f)
\end{gather}
respectively 
(the reason why these simplifications in the notation were not made  
right away in the definitions (\ref{defali}), (\ref{defali2})
is the need to consider integrals with variable limits, see, for example, Theorem \ref{brw}).


From definitions (\ref{defali}),\;(\ref{defali2}) $\bigl($ taking into account notation (\ref{exd4})$\bigr)$ 
and representations (\ref{ch17}) and~ (\ref{ch19}),
we obtain
\begin{gather} \label{exd3}
(*)\!\int\limits_a^b f(t)\,dg(t)+(*)\!\int\limits_a^b g(t)\,df(t)=F+G+H,
\end{gather}
where (see Theorem 6.2):
$$
F\doteq \int_a^bf_{\gc}(t)\,dg_{\gc}(t)+\int_a^bg_{\gc}(t)\,df_{\gc}(t)\;\Bigl (=g_{\gc}(t)f_{\gc}(t)\Big|_a^b\Bigr),
$$
$$
G\doteq \int_a^bf_{\gd}(t)\,dg_{\gc}(t)+\int_a^bg_{\gd}(t)\,df_{\gc}(t),\quad
H\doteq \sum\limits_{t\in T}f(t)\sigma_t(g)+ \sum\limits_{t\in T}g(t)\sigma_t(f).
$$
By  Theorem 6.2 and the formulas for computing Riemann-Stieltjes integral (see \S\,6),
$$
\int_a^bf_{\gd}(t)\,dg_{\gc}(t)=f_{\gd}(t)g_{\gc}(t)\Big|_a^b-\int_a^bg_{ \gc}(t)\,df_{\gd}(t)
=f_{\gd}(b)g_{\gc}(b)-\sum\limits_{t\in T}g_{\gc}(t)\sigma_t(f).
$$
Similarly,
$$
\int_a^bg_{\gd}(t)\,df_{\gc}(t)\!=\!g_{\gd}(t)f_{\gc}(t)\Big|_a^b\!-\!\int_a^bf_{\gc}(t)\,dg_{\gd}(t)\!=\!g_{\gd}(b)f_{\gc}(b)!-\!
\sum\limits_{t\in T}f_{\gc}(t)\sigma_t(g).
$$
Adding up the last two equalities, we obtain
\begin{multline*}
F+G=g_{\gc}(t)f_{\gc}(t)\Big|_a^b+f_{\gd}(t)g_{\gc}(t)\Big|_a^b+ g_{\gd}(t)f_{\gc}(t)\Big|_a^b-\int_a^bg_{\gc}(t)\,df_{\gd}(t)-\\
-\int_a^bf_{\gc}(t)\,dg_{\gd}(t)\!=\!g(t)f(t)\Big|_a^b-g_{\gd}(b)f_{\gd }(b)\!-\!
\sum\limits_{t\in T}g_{\gc}(t)\sigma_(f)\!-\!\sum\limits_{t\in T}f_{\gc}(t)\sigma_t( g)\!=\\
=g(t)f(t)\Big|_a^b-H+S,\;\text{where}\;S\doteq\sum\limits_{t\in T}\bigl(g_{\gd}(t)\sigma_t(f)+f_{\gd}(t)\sigma_t(g)\bigr)-g_{\gd}(b)f_{\gd}(b)
\end{multline*}
According to formula (\ref{ch17}), taking into account our convention
$$
f_{\gd}(b)=\sum\limits_{t\in T}\sigma_t(f),\quad g_{\gd}(b)=\sum\limits_{t\in T}\sigma_t(g );
$$
$$
f_{\gd}(t)=\sum\limits_{s\in T,\,s<t}\sigma_s(f)+\sigma^{-}_t(f),\quad g_{\gd}(t)=\sum\limits_{s\in T,\,s<t}\sigma_s(g)+\sigma^{-}_t(g).
$$
Substituting these expressions in $S,$ we obtain
\begin{multline*}
S=\sum\limits_{t\in T}\sigma_t(f)\left(\sum\limits_{s\in T,\,s<t}\sigma_s(g)+\sigma^{-}_t( g)\right)+\\+
\sum\limits_{t\in T}\sigma_t(g)\left(\sum\limits_{s\in T,\,s<t}\sigma_s(f)+\sigma^{-}_t(f) \right)-\\-
\sum\limits_{t\in T}\sigma_t(g)\left(\sum\limits_{s\in T,\,s<t}\sigma_s(f)+\sigma_t(f)+\sum\limits_ {s\in T,\,s>t}\sigma_s(f)\right)=\\=
\sum\limits_{t\in T}\sigma_t(f)\sum\limits_{s\in T,\,s<t}\sigma_s(g)+\sum\limits_{t\in T}\sigma_t( g)\sum\limits_{s\in T,\,s<t}\sigma_s(f)+\\
+\sum\limits_{t\in T}\sigma_t(f)\sigma^{-}_t(g)+\sum\limits_{t\in T}\sigma_t(g)\sigma^{-}_t( f)-\\-
\sum\limits_{t\in T}\sigma_t(g)\sum\limits_{s\in T,\,s<t}\sigma_s(f)-\sum\limits_{t\in T}\sigma_t( g)\sigma_t(f)-
\sum\limits_{t\in T}\sigma_t(g)\sum\limits_{s\in T,\,s>t}\sigma_s(f)\\
\end{multline*}
The second and the fifth terms on the right side of the last equality cancel out each other. So, what remains is
\begin{multline}
\label{zn2}
S=\sum\limits_{t\in T}\sigma_t(f)\sum\limits_{s\in T,\,s<t}\sigma_s(g)+
\sum\limits_{t\in T}\sigma_t(f)\sigma^{-}_t(g)+\sum\limits_{t\in T}\sigma_t(g)\sigma^{-}_t(f )-\\-
\sum\limits_{t\in T}\sigma_t(g)\sigma_t(f)-\sum\limits_{t\in T}\sigma_t(g)\sum\limits_{s\in T,\,s> t}\sigma_s(f).
\end{multline}
Since the double series $\sum\limits_{t,s\in T}\sigma_t(f)\sigma_s(g)$ converges absolutely (the partial sums of the double series of absolute values
are bounded by number $\|f\|\cdot \|g\|$, where one of the norms is in  space ${\bf{BV}}[a,b]$ and 
the other one is in space ${\bf{CH}}[a,b]$),
 one can change the order of summation in the last term in (\ref{zn2}). After doing this and re-ordering the index of summation, we arrive at the equality
$$
\sum\limits_{t\in T}\sigma_t(g)\sum\limits_{s\in T,\,s>t}\sigma_s(f)\!\!=\!\!
\sum\limits_{s\in T}\sigma_s(f)\sum\limits_{t\in T,\,t<s}\sigma_t(g)\;\;=\!\!
\sum\limits_{t\in T}\sigma_t(f)\sum\limits_{s\in T,\,s<t}\sigma_s(g).
$$
Thus, in (\ref{zn2}), the first term and the last term cancel out each other. The remaining terms:
\begin{multline*}
S=\sum\limits_{t\in T}\sigma_t(f)\sigma^{-}_t(g)+\sum\limits_{t\in T}\sigma_t(g)\sigma^{-}_t (f)-\sum\limits_{t\in T}\sigma_t(g)\sigma_t(f)=\\=
\sum\limits_{t\in T}\sigma_t(f)\sigma^{-}_t(g)-\sum\limits_{t\in T}\sigma_t(g)\sigma^{+}_t(f )
\end{multline*}
Let us transform the resulting expression:
\begin{multline*}
S=\sum\limits_{t\in T}\bigl(\sigma^{+}_t(f)+\sigma^{-}_t(f)\bigr)\sigma^{-}_t(g)-
\sum\limits_{t\in T}\bigl(\sigma^{+}_t(g)+\sigma^{-}_t(g)\bigr)\sigma^{+}_t(f)=\\ =
\sum\limits_{t\in T}\Bigl(\sigma^{+}_t(f)\sigma^{-}_t(g)+\sigma^{-}_t(f)\sigma^{-} _t(g)-
\sigma^{+}_t(g)\sigma^{+}_t(f)-\sigma^{-}_t(g)\sigma^{+}_t(f)\Bigr).
\end{multline*}
The first and fourth terms cancel out each other. We obtain
\begin{multline*}
S=\sum\limits_{t\in T}\Bigl(\sigma^{-}_t(f)\sigma^{-}_t(g)-\sigma^{+}_t(f)\sigma^{ +}_t(g)\Bigr)=\\=
-\sum\limits_{t\in T}\Bigl(\sigma^{+}_t(f)\sigma^{+}_t(g)-\sigma^{-}_t(f)\sigma^{- }_t(g)\Bigr).
\end{multline*}
Thus, finally,
$$
F+G+H=f(t)g(t)\Big|_a^b-\sum\limits_{t\in T}\Bigl(\sigma^{+}_t(f)\sigma^{+} _t(g)-\sigma^{-}_t(f)\sigma^{-}_t(g).\qquad \qquad \hfill \square
$$
\begin{zam}
\label{zam(2)}
The formula (\ref{bp1}) can be put in exactly the same form as the integration by parts formula for the PS-integral in\linebreak  
\cite[I.4.33]{Stv79}.

Let $f\in{\bf{CH}}[a,b],\;g\in{\bf{BV}}[a,b];\;[c,d] \subset [a,b]$. Then 
\begin{multline}
\label{bp2}
(*)\!\int\limits_c^d f(t)\,dg(t)+(*)\!\int\limits_c^d g(t)\,df(t)=f(t)g(t) \Big|_c^d-\\-
\sum\limits_{t\in T,\,t\in[c,\,d)}\bigl(\sigma^{+}_t(f)\sigma^{+}_t(g)\bigr)+
\sum\limits_{t\in T,\,t\in(c,\,d]}\bigl(\sigma^{-}_t(f)\sigma^{-}_t(g)\bigr).
\end{multline}
\end{zam}


\bigskip

\poonkt
{Indefinite\;*-integral\;}\phantom{012345678901234567890} \

\medskip

{\bf{1.}}\;First, let us take definition (\ref{defali}) and define
\begin{equation}
\label{def1}
\Phi(t)\doteq (*)\!\int\limits_{a}^t x(s)\,dg(s).
\end{equation}
\begin{teo}
\label{brw}
{\it Let $x\in{\bf{N}}[a,b],\;g\in{\bf{BV}}[a,b]$. Then
$\Phi(\cdot)$~ is a function of bounded variation, moreover,}
$$\bigvee\limits_{a}^{b}(\Phi)\leqslant \|x\|_{{\bf{N }}}\bigvee\limits_a^b(g),$$
\begin{equation}
\label{skachkifi}
\sigma^+_t(\Phi)=x(t)\sigma^+_t(g),\;\;\sigma^-(\Phi)=x(t)\sigma^-_t(g),\;\;\sigma_t(\Phi)=x(t)\sigma_t(g).\;\;
\end{equation}
\end{teo}
\doc\,For an arbitrary partition $\{t_k\}_{k=0}^n$, due to estimate \ref{tonkgrub},
\begin{align*}
\sum\limits_{k=1}^n\bigl|\Phi(t_k)-\Phi(t_{k-1})\bigr|& =\sum\limits_{k=1}^n\left|(* )\!\int\limits_{t_{k-1}}^{t_k} x(t)\,dg(t)\right| \\
& \leqslant\|x\|_{{\bf{N}}}\sum\limits_{k=1}^n\bigvee\limits_{t_{k-1}}^{t_k}(g)
=\|x\|_{{\bf{N}}}\bigvee\limits_a^b(g).
\end{align*}
Let $h>0$. Due to (\ref{add}) and (\ref{defali}),
\begin{align*}
& \Phi(t+h)-\Phi(t)=(*)\!\int\limits_t^{t+h} x(t)\,dg(t) =\int_t^{t+h} x( t)\,dg_{\mathfrak c}(t)+x(t)\bigl(g(t+)-g(t)\bigr) \\
&+
\sum\limits_{s\in T(g)\bigcap(t,t+h)}x(s)\bigl(g(s+)-g(s-)\bigr)
+x(t+h)\bigl(g(t+h)-g(t+h-)\bigl).
\end{align*}
The integral term in the right-hand side tends to zero as $h\to 0$. This follows from the estimate above
and the continuity of function $\bigl(g_{\mathfrak c}\bigr)_{\pi}(\cdot)$.
The term $\sum$ in the second line tends to $0$ as $h\to 0+$. This is obvious if $t$ is an isolated point of the set $T(g)$. In the case
$t$ is a limit point of this set then,  by Lemma \ref{ch8lem1}, $\sum$ tends to
$x(t)\cdot\bigl(g(t+)-g(t+)\bigr)=0$ (the first multiple may not have a limit, but it is bounded).
For the same reason, the term 
in the last line also tends to zero. Thus, we get: $\;\sigma^+_t(\Phi)=x(t)\cdot\bigl(g(t+)-g(t)\bigr)$. Considering the difference
$\Phi(t)-\Phi(t-h),$ we see that $\sigma^-(\Phi)=x(t)\cdot\bigl(g(t)-g(t-)\bigr). $
\hfill $\square$

Below we will specify the exact value of $\bigvee\limits_a^b(\Phi)$ (see (\ref{teoravvar})).


 {\bf{2.}}\;Let us proceed to conditions (\ref{defali2}).

 Let $x\in{\bf{BV}}[a,b]$, $
 g\in{\bf{CH}}[a,b]$, let $
 \Phi(\cdot)$ be defined by equality (\ref{def1}).
 According to (\ref{bp1}),
\begin{equation}
\label{bp1+}
\Phi(t)=-A(t)+B(t)-\Gamma(t),
\end{equation}
where
$$
A(t)=(*)\!\int\limits_a^t g(s)\,dx(s),\quad B(t)=x(s)g(s)\Big|_a^t,
$$
$$
G(t)=\sum\limits_{s\in T\bigcap (a,t)}\Bigl(\sigma^{+}_s(x)\sigma^{+}_s(g)-\sigma^{ -}_s(x)\sigma^{-}_s(g)\Bigr).
$$
By Theorem \ref{brw}, $\quad A\in {\bf{BV}},\;B\in {\bf{CH}},\;{\Gamma\in \bf{H}}.$
So,
now we can only guarantee that $\Phi\in{\bf{CH}}\;\Bigl(\subset{\bf{R}}\Bigr).$

According to Theorem \ref{brw},
$$
\sigma^+_t(A)=g(t)\sigma^+_t(x),
\sigma^-_t(A)=g(t)\sigma^-_t(x),
\sigma_t(A)=-g(t)\sigma_t(x).$$
It follows directly from the definitions that
$$\sigma^+_t(B)=x(t+)g(t+)-x(t)g(t),\;\;\sigma^-_t(B)=
x(t)g(t)-x(t-)g(t-),$$
$$\sigma_t(B)=
x(t+)g(t+)-x(t-)g(t-),$$
function $\Gamma(\cdot)$ is continuous at the points of continuity of at least one of functions $x(\cdot)$ or $g(\cdot)$.
At the points $t$ where both functions are discontinuous,
$$
\sigma^+_t(\Gamma)=\linebreak=
-\sigma^+_t(x)\sigma^+_t(g),\quad \sigma^-_t(\Gamma)=\sigma^-_t(x)\sigma^-_t(g), 
$$
$$
 \sigma_t(\Gamma)=\sigma^-_t(x)\sigma^-_t(g)-
 \sigma^+_t(x)\sigma^+_t(g)
$$
(see also Example \ref{ObobQraz}).
According to (\ref{bp1+}),
$$\sigma^+_t(\Phi)=
x(t)\sigma^+_t(g),\;\;\sigma^-(\Phi)=x(t)\sigma^-_t(g),\;\;\sigma_t(\Phi)= x(t)\sigma_t(g).$$

Thus, assertion (\ref{skachkifi}) of Theorem \ref{brw} is valid under conditions (\ref{defali2}) as well.


\bigskip

\poonkt
{Change of the order of integration}

\begin{teo}
\label{chali}
\textit{Assume that one of the following conditions is satisfied:}
$$
I.\;\;f\in{\bf{BV}}[a,b],\;g\in{\bf{BV}}[c,d],\;h( \cdot,s)\in{\bf{R}}[a,b],\;h(t,\cdot)\in{\bf{R}}[c,d]\ ;
$$
or
$$
II.\;\;f\in{\bf{BV}}[a,b],\;g\in{\bf{BV}}[a,b],\;h( \cdot,s)\in{\bf{N}}[a,b],\;h(t,\cdot)\in{\bf{N}}[c,d]\ ;
$$
{\it and there exists $M>0,$ such that $|h(t,s)|\leqslant M \;\; (t\in [a,b],\,s\in [c,d]).$}

{\it Then}
\begin{gather}\label{chang0}
(*)\!\int\limits_a^b\left((*)\!\int\limits_c^d h(t,s)\,dg(s)\right)\,df(t)=
(*)\!\int\limits_c^d\left((*)\!\int\limits_a^bh(t,s)\,df(t)\right)\,dg(s).
\end{gather}
\end{teo}
\doc\,
It suffices to carry out the proof in the case when condition $I$ is satisfied.
If $f$ and $g$~ are continuous functions of bounded variation, then
there is equality
\begin{gather}\label{chRS}
\int_a^b\left(\int_c^d h(t,s)\,dg(s)\right)\,df(t)=\int_c^d\left(\int_a^bh(t,s)\, df(t)\right)\,dg(s)
\end{gather}
(in this case both integrals are treated in the sense of Riemann--Stieltjes).

By the result in subsection {\bf{10.5.5}}, 
conditions of the theorem imply that there exists $M>0$ 
such that $|h(t,s)|\leqslant M \;\; (t\in [a,b],\,s\in [c,d].$
Let condition \;I be satisfied. Then

a) the double series converges absolutely
\begin{gather*}
S\doteq \sum\limits_{s\in T(g),\,t\in T(f)}h(t,s)\sigma_s(g)\sigma_t(f)\;\;\;
\text {and}\;\;\;|S|\leqslant M\bigvee\limits_a^b(f)\bigvee\limits_c^d(g),
\end{gather*}
since the partial sums of the following double series with non-negative terms \\ 
$\sum\limits_{s\in T(g),\,t\in T(f)}|h(t,s)|\cdot |\sigma_s(g )|\cdot |\sigma_t(f)|$
are bounded by $M\bigvee\limits_a^b(f)\bigvee\limits_c^d(g)$;

b) by the Pringsheim theorem (see, for example, \cite[2.51]{WW1}) the corresponding repeated series converge absolutely, and
\begin{gather*}
\sum\limits_{t\in T(f)}\sum\limits_{s\in T(g)}h(t,s)\sigma_t(f)\sigma_s(g)=
\sum\limits_{s\ in T(g)}\sum\limits_{t\in T(f)}h(t,s)\sigma_s(g)\sigma_t(f).
\end{gather*}

c) the series $\sum\limits_{s\in T(g)}h(t,s)\sigma_s(g)\;$\;$\;\left(\sum\limits_{t\in T(f)}h(t,s)\sigma_t(f)\right)$
 converges absolutely and uniformly with respect to $t\;\;\bigl(s\bigr)$ \;
(the series are majorized by convergent numerical series). Therefore, one can integrate the double series term by term with respect 
to continuous function of bounded variation $f_{\gc}(t)\;\;\bigl(g_{\gc}(s)\bigr).$

Denote
$$
H_g(t)\doteq \int_c^d h(t,s)\,dg_{\gc}(s),\;\;\;H_f(s)\doteq\int_a^b h(t,s)\,df_ {\gc}(t),
$$
$$
\mathcal H_g(t)\doteq (*)\!\int\limits_c^d h(t,s)\,dg(s)= H_g(t)+\sum\limits_{s\in T}h(t, s)\sigma_s(g),
$$
$$
\mathcal H_f(s)\doteq (*)\!\int\limits_a^b h(t,s)\,df(t)=H_f(s)+\sum\limits_{t\in T}h(t, s)\sigma_t(f).
$$
In this notation, the statement of the theorem has the following form:
\begin{gather}\label{chang1}
(*)\!\int\limits_a^b\mathcal H_g(t)\,df(t)=(*)\!\int\limits_c^d\mathcal H_f(s)\,dg(s).
\end{gather}

Let us prove equality (\ref{chang1}).
\begin{multline*}
(*)\!\int\limits_a^b\mathcal H_g(t)\,df(t)=\int_a^b\mathcal H_g(t)\,df_{\gc}(t)+\sum\limits_{ t\in T(f)}\mathcal H_g(t)\sigma_t(f)=\\=
\int_a^b\left(H_g(t)+\sum\limits_{s\in T(g)}h(t,s)\sigma_s(g)\right)\,df_{\gc}(t)+ \\
+\sum\limits_{t\in T(f)}\left(H_g(t)+\sum\limits_{s\in T(g)}h(t,s)\sigma_s(g)\right)\sigma_t(f)=\\
\int_a^b\int_c^d h(t,s)\,dg_{\gc}(s)\,df_{\gc}(t)
+\int_a^b\left(\sum\limits_{s\in T(g)}h(t,s)\sigma_s(g)\right)\,df_{\gc}(t)+\\+
\sum\limits_ {t\in T(f)}\left(\int\limits_c^d h(t,s)\,dg_{\gc}(s)\right)\sigma_t(f)+ 
\sum\limits_{t\in T(f)}\sum\limits_{s\in T(g)}\sigma_s(g)\sigma_t(f)
\end{multline*}
According to (\ref{chRS}), the first term on the right-hand side can be written as 
$\int_c^d\int_a^b h(t,s)\,df_{\gc}(t)\,dg_{\gc}( s)$. In view of assertion c), the second
(third) term can be written as
$$
\sum\limits_{s\in T(g)}\left(\int\limits_a^b h(t,s)\,df_{\gc}(t)\right)\sigma_s(g)\;\;
\left(\int_c^d\left(\sum\limits_{t\in T(f)}h(t,s)\sigma_t(f)\right)\,dg_{\gc}(s)\right).
$$
Finally, by assertion b), the last term is written as $\sum\limits_{s\in T(g)}\sum\limits_{t\in T(f)}\sigma_t(f)\sigma_s(g). $
Therefore,
\begin{multline*}
(*)\!\int\limits_a^b\mathcal H_g(t)\,df(t)=\int_c^d\int_a^b h(t,s)\,df_{\gc}(t)\,dg_ {\gc}(s)+\\
+\int_c^d\left(\sum\limits_{t\in T(f)}h(t,s)\sigma_t(f)\right)\,dg_{\gc}(s)+
\sum\limits_{s\in T(g)}\left(\int\limits_a^b h(t,s)\,df_{\gc}(t)\right)\sigma_s(g)+\\+
\sum\limits_{s\in T(g)}\sum\limits_{t\in T(f)}\sigma_t(f)\sigma_s(g)=(*)\!\int\limits_c^ d\mathcal H_f(s)\,dg(s).\hfill \square
\end{multline*}


\bigskip

\poonkt
{Limit theorems\;}\phantom{01234567890123456789}\phantom{01234567890123456789}\phantom{01234567890123456789}  \

{\bf{1.}}\;The ``first'' limit theorem (see Theorem \ref{ch4thm3}) holds for the (*)-integral.
\begin{teo}
\label{Fst1}
{\it Let $f_n\in{\bf{N}}[a,b],\,n\in\mathbb N,\;f_n(t)\rightrightarrows f(t)\;(n\to \infty)$,
$g\in{\bf{BV}}[a,b]$. Then}
\begin{gather}
\label{Fst11}
\underset{n\to\infty}{\lim}(*)\!\int\limits_a^b f_n(t)\,dg(t)=(*)\!\int\limits_a^b f(t)\ ,dg(t).
\end{gather}
\end{teo}
\doc\,
Note that, by the assertions of Section 9.4.1, the conditions of the theorem imply that 
$f\in{\bf{N}},\;f_n\to f$ in space
${\bf{N}}[a,b]$ and $\|f_n\|\leqslant M\;\;
(n\in\mathbb N)$ for some $\;M\geqslant 0.$
Therefore, the  integrals in both sides of (\ref{Fst11}) exist. Convergence
\begin{gather} \label{Fst12}
\int_a^bf_n(t)\,dg_{\gc}(t)\to \int_a^bf(t)\,dg_{\gc}(t)\qquad(n\to\infty)
\end{gather}
follows from Theorem \ref{Fst3}. Due to the uniform boundedness of sequence $\{f_n(t)\}_{n=1}^{\infty}$,
the series $\sum\limits_{t\in T(g)\bigcap(a,b)}f_n(t)\sigma_t(g)$ converges uniformly with respect to $n\in\mathbb N,$, 
so one can pass to the limit $n\to\infty$ under the series sign. Consequently,
\begin{gather}
\label{Fst13}
\sum\limits_{t\in T(g)\bigcap(a,b)}f_n(t)\sigma_t(g)\to \sum\limits_{t\in T(g)\bigcap(a,b) }f(t)\sigma_t(g)\qquad(n\to\infty).
\end{gather}
\hfill $\square$
\begin{sle}\label{SFst1}
{\it Let $g_n\in{\bf{N}}[a,b],\,n\in\mathbb N,\;g_n(t)\rightrightarrows g(t)\linebreak
(n\to\infty),
f\in{\bf{BV}}[a,b]$. Then}
\begin{gather} \label{SFst11}
\underset{n\to\infty}{\lim}(*)\!\int\limits_a^b f(t)\,dg_n(t)=(*)\!\int\limits_a^b f(t)\, dg(t).
\end{gather}
\end{sle}
\doc\,We apply successively  (\ref{chang0}), (\ref{Fst11}) and again (\ref{chang0}).
\hfill $\square$

\medskip

There is also an analogue of Theorem \ref{doparcel}.
\begin{teo}
\label{Fst2}
{\it Let $f_n,\,f\in{\bf{N}}[a,b],\;\;\|f_n\|\leqslant M,\;\;n\in\mathbb N,\linebreak
f_n(t)\to f(t),\;t\in [a,b]\;(n\to \infty),\;\;g\in{\bf{BV}}[a,b]$. Then
the limit relation $(\ref{Fst11})$ holds.}
\end{teo}
\doc\,
Theorem \ref{doparcel} implies convergence (\ref{Fst12}).
Since the uniform boundedness of sequence $\{f_n(t)\}_{n=1}^{\infty}$ is now postulated,
convergence (\ref{Fst13}) also takes place.
\hfill $\square$
\begin{sle}
\label{SFst2}
{\it Let
$
g_n,\,g\in{\bf{N}}[a,b],\,n\in\mathbb N,\;g_n(t)\to g(t),
\linebreak
t\in [a,b]\;\;(n\to \infty),\;\;\|g_n\|\leqslant M,\;f\in{\bf{BV}}[a,b]
$.
Then the limit relation} (\ref{SFst11}) holds.
\end{sle}


{\bf{2.}}\;There is an analogue of Helly's theorem for the (*)-integral.
\begin{teo}\label{ChelliM}
{\it Let $f\in{\bf{N}}[a,b]$, 
$$g_n(t)\to g(t)\;\;(n\to\infty) \text{ quasi-uniformly on $[a,b ]$;}
$$
there exists a constant $V>0$ such that $\bigvee\limits_a^b(g_n)\leqslant V$,
$n=1,2,\ldots$;
the series
\begin{equation}
\label{modrsk}
\sum\limits_{t\in T}\,\bigl|\sigma_t(g_n)\bigr|\qquad \Bigl(T\doteq \bigcup\limits_{n=1}^\infty T(g_n)\Bigr)
\end{equation}
converges uniformly with respect to $n\in\mathbb N.$
Then}
\begin{gather}
\label{ch1}
\underset{n\to\infty}{\lim}(*)\!\int\limits_a^b f(t)\,dg_n(t)=(*)\!\int\limits_a^b f(t)\, dg(t).
\end{gather}
\end{teo}
\doc\,
By (\ref{PN}), $\;\;T(g)\subset T.$ Arguing as in the proof of Helly's theorem
(see the proof of Theorem \ref{Chelli2}), one shows that $g\in{\bf{BV}}[a,b].$

By Theorem \ref{osntQU},
\begin{equation}
\label{modrsk+}
\bigl(g_{n}\bigr)_{\gd}(t)\to g_{\gd}(t),\;\;\bigl(g_{n}\bigr)_{\gc}(t )\to g_{\gc}(t)\;(n\to\infty).
\end{equation}
Arguing as in the proof of Theorem \ref{bypart},
we write down the definition of the (*)-integral taking into account (\ref{exd4}).
The boundedness of $f(\cdot)$ and the uniform convergence of series (\ref{modrsk}) give us the uniform with respect to 
$t\in [a,b]$ and $n\in\mathbb N$
convergence of series
$\sum\limits_{t\in T}f(t)\sigma_t(g_n)$.
Therefore, one can pass to the limit under the series sign as $n\to\infty.$
Taking into account Theorem \ref{Chelli2} and Corollary \ref{sochrsk1}, we obtain
\begin{multline*}
(*)\!\int\limits_a^b f(t)\,dg_n(t)=\int_a^b f(t)\,d\bigl(g_{n}\bigr)_{\gc}(t)+ \sum\limits_{t\in T}f(t)\sigma_t(g_n)\to \\
\to \int_a^b f(t)\,dg_{\gc}(t)+\sum\limits_{t\in T}f(t)\sigma_t(g)=(*)\!\int\limits_a^ b f(t)\,dg(t)\;\;\;
(n\to\infty). \hfill \square
\end{multline*}

\medskip

The following modification of Theorem \ref{ChelliM} (following the previous proof) is often more convenient. 

\medskip

{\bf{Theorem\;14.6\,(bis)}}\;
{\it Assume that $f\in{\bf{N}}[a,b],$ there exists constant $V>0$ such that $\bigvee\limits_a^b(g_n)\leqslant V,\;\; n=1,2,\ldots;\;$ convergence
$(\ref{modrsk+})$ takes place and
 series $(\ref{modrsk})$ converges uniformly with respect to $n\in\mathbb N.$ Then the limit relation (\ref{ch1}) holds.}


\medskip

{\bf{3.}}\;Theorems \ref{Fst1} and \ref{ChelliM} allow us to prove the following result,
which finds immediate application in the theory of equations with generalized functions in the coefficients 
(see, for example, \cite{DK06},\,\cite{derr18}).

\begin{teo}
\label{okpredt}
{\it Let $f_n\in{\bf{N}},\,n\in\mathbb N,\;f_n(t)\rightrightarrows f(t),\;g_n(t)\rightrightarrows g(t) \linebreak
(n\to\infty),\;\bigvee\limits_a^b(g_n)\leqslant V,\;V>0, n=1,2,\ldots$. Then} 
$$
\underset{n\to\infty}{\lim}(*)\!\int\limits_a^b f_n(t)\,dg_n(t)=(*)\!\int\limits_a^b f(t)\,dg(t).
$$
\end{teo}
\doc\,Let $\varepsilon>0$ be arbitrarily small.
Since sequence $\{f_n\}_{n=1}^\infty$ converges uniformly, there exists a natural number $N_1$ such that, for $n>N_1$,
$
\|f_n-f\|_{{\bf{N}}}<\dfrac{\varepsilon}{2V}.
$
Hence, by the estimate from 
subsection {\bf{14.1.3}}
$\left|(*)\!\int\limits_a^b f_n(t)\,dg_n(t)-(*)\!\int\limits_a^b f(t)\,dg_n(t)\right|< \frac{\varepsilon}{2}.$
In view of Theorem \ref{ChelliM}, there exists a positive integer $N_2,$ such that, for $n>N_2$,
$$
\left|(*)\!\int\limits_a^b f(t)\,dg_n(t)-(*)\!\int\limits_a^b f(t)\,dg(t)\right|<\frac{\varepsilon}{2}.
$$
So, for
$n>\max\{N_1,\,N_2\}$,
\begin{multline*}
\left|(*)\!\int\limits_a^b f_n(t)\,dg_n(t)-(*)\!\int\limits_a^b f(t)\,dg(t)\right|\ls \\
\ls\left|(*)\!\int\limits_a^b f_n(t)\,dg_n(t)-(*)\!\int\limits_a^b f(t)\,dg_n(t)\right. +\\
+\left|(*)\!\int\limits_a^b f(t)\,dg_n(t)-(*)\!\int\limits_a^b f(t)\,dg(t)\right|<\varepsilon. \hfill \square
\end{multline*}


{\bf{4.}}\;Theorem \ref{Fst1} allows to prove a formula for finding the total variation {\it of indefinite *-integral}
$\bigl($for the definition of $g_{\pi}(\cdot)$, see (\ref{ch25+})$\bigr).$
 \begin{teo}
 \label{teoravvar0}
 {\it Let $x\in{\bf{N}},\,g\in{\bf{BV}},$
 $$
 \Phi(t)=(*)\!\int\limits_{a}^{t}x(s)\,dg(s),\quad \Phi^{(\pi)})(t)=( *)\!\int\limits_{a}^{t}x(s)\,dg_{\pi}(s).
 $$
Then}
\begin{equation}
\label{teoravvar}
\bigvee\limits_a^b(\Phi)=\bigvee\limits_a^b(\Phi^{(\pi)})=(*)\!\int\limits_{a}^{b}|x(t) |\,dg_{\pi}(t)
\end{equation}
\end{teo}
\doc\,
Let $\tau=\{s_j\}_{j=1}^p$~ be an arbitrary partition of interval $[a,b].$
Let us estimate the sums $v_{\tau}\bigl(\Phi\bigr),\;v_{\tau}\bigl(\Phi^{(\pi)}\bigr):$
\begin{multline*}
v_{\tau}\bigl(\Phi\bigr)=\sum\limits_{j=1}^p\Bigl|\Phi(s_j)-\Phi(s_{j-1})\bigr|=\\ =
\sum\limits_{j=1}^p\left|(*)\!\int\limits_{a}^{s_j} x(t)\,dg(t)-(*)\!\int\limits_ {a}^{s_{j-1}} x(t)\,dg(t)\right|=
\sum\limits_{j=1}^p\left|(*)\!\int\limits_{s_{j-1}}^{s_j} x(t)\,dg(t)\right|\leqslant \\
\leqslant \sum\limits_{j=1}^p (*)\!\int\limits_{s_{j-1}}^{s_j} |x(t)|\,dg_{\pi}(t) =
(*)\!\int\limits_{a}^b|x(t)|\,dg_{\pi}(t).
\end{multline*}
Analogously,
\begin{multline*}
v_{\tau}\bigl(\Phi^{(\pi)}\bigr)=\sum\limits_{j=1}^p\Bigl|\Phi^{(\pi)}(s_j)-\Phi ^{(\pi)}(s_{j-1})\bigr|=\\=
\sum\limits_{j=1}^p\left|(*)\!\int\limits_{s_{j-1}}^{s_j} x(t)\,dg_{\pi}(t)\right|\leqslant
\sum\limits_{j=1}^p (*)\!\int\limits_{s_{j-1}}^{s_j} |x(t)|\,dg_{\pi}(t)=
(*)\!\int\limits_{a}^b|x(t)|\,dg_{\pi}(t).
\end{multline*}
Therefore, the following estimates hold:
\begin{multline}
\label{varintocensn}
\bigvee\limits_a^b\bigl(\Phi\bigr)\leqslant (*)\!\int\limits_{a}^b|x(t)|\,dg_{\pi}(t),\quad
\bigvee\limits_a^b\bigl(\Phi^{(\pi)}\bigr)\leqslant (*)\!\int\limits_{a}^b|x(t)|\,dg_{\pi} (t),\\
\bigvee\limits_a^b\bigl(\Phi\bigr)\leqslant\bigvee\limits_a^b\bigl(\Phi^{(\pi)}\bigr).
\end{multline}
Let us represent $x(\cdot)$ as the uniform limit of a sequence of step functions.
To do this, we consider the uniform partition
\begin{equation}
\label{ravnrazb}
\tau_m=\{t_k\}_{k=0}^m,\;t_k=a+\dfrac{k}{m}(b-a),
\end{equation}
and put
$$
x_m(a)=x(a),\;x_m(t)=x\bigl(t_k\bigr)\;\;\text{for}\;\; t_{k-1}<t\leqslant t_k,\;k=1,2,\ldots, m;
$$
then
$x_m(t)\rightrightarrows x(t)$ on $[a,b]$ and, by Theorem \ref{Fst1},
\begin{equation}
\label{varintnach}
(*)\!\int\limits_{a}^b|x(t)|\,dg_{\pi}(t)\stackrel{(A)}=\underset{m\to\infty}{lim} \;(*)\!\int\limits_{a}^b|x_m(t)|\,dg_{\pi}(t).
\end{equation}
Let us also define
\begin{multline}
\label{shXm}
X_m(t)\doteq (*)\!\int\limits_{a}^{t}x_m(s)\;dg(s)\;\;\bigl(\rightrightarrows \Phi(t)\bigr), \\
X_m^{(\pi)}(t)\doteq (*)\!\int\limits_{a}^{t}x_m(s)\;dg_{\pi}(s)\;\;\bigl( \rightrightarrows \Phi^{(\pi)}(t)\bigr).
\end{multline}
The uniform convergence on $[a,b]$ also follows from Theorem \ref{Fst1}.

We have the following chain of statements (see explanation below)
\begin{multline}
\label{varintvikl}
(*)\!\int\limits_{a}^b|x_m(t)|\,dg_{\pi}(t)\stackrel{(B)}=\sum\limits_{k=1}^m( *)\!\int\limits_{t_{k-1}}^{t_k}|x_m(t)|\,dg_{\pi}(t)\stackrel{(C)}=\\=
\sum\limits_{k=1}^m(*)\!\int\limits_{t_{k-1}}^{t_k}|x(t_k)|\,dg_{\pi}(t)\stackrel {(D)}=
\sum\limits_{k=1}^m|x(t_k)|(*)\!\int\limits_{t_{k-1}}^{t_k}\,dg_{\pi}(t)\stackrel {(E)}=\\=
\sum\limits_{k=1}^m\left|x(t_k)(*)\!\int\limits_{t_{k-1}}^{t_k}\,dg_{\pi}(t)\right|\stackrel{(F)}=
\sum\limits_{k=1}^m\left|(*)\!\int\limits_{t_{k-1}}^{t_k}x(t_k)\;dg_{\pi}(t)\right|\stackrel{(G)}=\\=
\sum\limits_{k=1}^m\left|(*)\!\int\limits_{t_{k-1}}^{t_k}x_m(t)\;dg_{\pi}(t)\right|\stackrel{(H)}=v_{\tau_m}\bigl(X_m^{(\pi)}\bigr),
\end{multline}
Explanations: (A) --- Theorem \ref{Fst1}; (B) --- splitting the interval of\linebreak  
integration into $m$ parts and the integral into $m$
terms; (C) --- substitution of the definition of $x_m;$ (D) --- since the factor does not depend on $t$, we take it out of the integral sign;
(E) --- since the integral is positive, we bring it under the absolute value sign; (F) --- since $x(t_k)$ does not depend on $t$, we bring it under the
integral sign; (G) --- within the limits of integration instead of $x(t_k)$ we write $x_m(t);$
(H) --- under the absolute value sign we have the difference $X_m(t_k)-X_m(t_{k-1}),$
that is, it is the standard sum $v_{\tau}$ from the definition of the total variation of function $X_m;$

Let us show that
\begin{equation}
\label{varintravlim}
\underset{m\to\infty}{\lim}\;\Bigl(v_{\tau_m}\bigl(X_m^{(\pi)}\bigr)-v_{\tau_m}\bigl(\Phi^{( \pi)}\bigr)\Bigr)=0.
\end{equation}
To do this, we use inequality $|\alpha|-|\beta|\leqslant |\alpha-\beta|:$
\begin{multline}
\label{raznxmpifi}
\Bigl|v_{\tau_m}(X_m^{(\pi)})-v_{\tau_m}(\Phi^{(\pi)})\Bigr|=\\=
\left|\sum\limits_{k=1}^m\left(\Bigl|(*)\!\int\limits_{a}^{t}x_m(s)\;dg_{\pi}(s) \Bigr|-
\Bigl|(*)\!\int\limits_{a}^{t}x(s)\;dg_{\pi}(s)\Bigr|\right)\right|\leqslant\\
\leqslant\sum\limits_{k=1}^m\left|(*)\!\int\limits_{t_{k-1}}^{t_k}\bigl(x_m(t)-x(t)) \bigr)\,dg_{\pi}(t)\right|
\leqslant\sum\limits_{k=1}^m(*)\!\int\limits_{t_{k-1}}^{t_k}\bigl|x_m(t)-x(t)\bigr|\ ,dg_{\pi}(t)=\\=
(*)\!\int\limits_{a}^{b}\bigl|x_m(t)-x(t)\bigr|\,dg_{\pi}(t)\to 0\;\;(m \to\infty);
\end{multline}
the passage to the limit is done according to Theorem\; \ref{Fst1}.

We obtain (\ref{varintnach})--(\ref{varintravlim}):
\begin{equation}
\label{varintocensv}
(*)\!\int\limits_{a}^b|x(t)|\,dg_{\pi}(t)\leqslant\bigvee\limits_a^b\bigl(\Phi^{(\pi)} \bigr).
\end{equation}
Together with the second estimate in (\ref{varintocensn}) this implies the second equality in (\ref{teoravvar}).

We will prove the first equality in (\ref{teoravvar}). First, let $x(t)\equiv C.$
\begin{multline*}
\Phi(t)\!=\!C(*)\!\int\limits_{a}^{t}\,dg(t)\!=\!
C\Bigl(\int\limits_{a}^{t}\,dg_{\gc}(t)+\sigma_a^+(g)\!+\!\sum\limits_{s\in T(g)\bigcap (a ,b)}\sigma_s(g)+\sigma_t^-(g)\Bigr)\!=\\=
C\Bigl(g_{\gc}(t)-g_{\gc}(a)+g_{\gd}(t)\Bigr)=C\bigl(g(t)-g(a)\bigr) ,\quad
\bigvee\limits_a^b(\Phi)=|C|\bigvee\limits_a^b(g);\\
\Phi^{(\pi)}(t)=C(*)\!\int\limits_{a}^{t}\,dg_{(\pi)}(t)=
C\bigl(g_{\pi}(t)-g_{\pi}(a)\bigr)=Cg_{\pi}(t),\\
\bigvee\limits_a^b\bigl(\Phi^{(\pi)}\bigr)=|C|g_{\pi}(b)=|C|\bigvee\limits_a^b(g),
\end{multline*}
that is, in this case, the first equality in (\ref{teoravvar}) is true.
Consider partition (\ref{ravnrazb}) and sequences
$\{x_m\}_{m=1}^\infty$, $\{X_m\}_{m=1}^\infty$, $\{X_m^{(\pi)}\}_{m =1}^\infty$
(see above).
According to what was just proved,
\begin{multline}
\label{ravvarXmpi}
\bigvee\limits_a^b\bigl(X_m^{(\pi)}\bigr)=\sum\limits_{k=1}^m\bigvee\limits_{t_{k-1}}^{t_k}\bigl (X_m^{(\pi)}\bigr)=
\sum\limits_{k=1}^m |x(t_k)|\bigvee\limits_{t_{k-1}}^{t_k}(g)=\\=
\sum\limits_{k=1}^m\bigvee\limits_{t_{k-1}}^{t_k}(X_m)=\bigvee\limits_a^b\bigl(X_m).
\end{multline}
By the second equality in (\ref{teoravvar}) and Theorem \ref{Fst1}, 
 \begin{equation}
\label{ShvarXmpi}
 \bigvee\limits_a^b\bigl(X_m^{(\pi)}\bigr)=(*)\!\int\limits_{a}^{b}|x_m(t)|\,dg_{\pi} (t)\to
 (*)\!\int\limits_{a}^{b}|x(t)|\,dg_{\pi}(t)=\bigvee\limits_a^b\bigl(\Phi^{(\pi) }\bigr).
 \end{equation}
$(m\to\infty).$ To end the proof of the first equality in (\ref{teoravvar}), we have to show that
\begin{equation}
\label{ShvarXm}
V\doteq\underset{m\to\infty}{lim}\;\bigvee\limits_a^b\bigl(X_m\bigr)=\bigvee\limits_a^b\bigl(\Phi)\doteq F.
\end{equation}
(The existence of limit $V$ follows from (\ref{ravvarXmpi}) and (\ref{ShvarXmpi}).)

If $V\leqslant F,$ then, starting from some $m$, $\bigvee\limits_a^b\bigl(X_m\bigr)\leqslant F$, and,
by Theorem \ref{shvar}, $V=F.$

Suppose that $V>F.$ Then $\alpha\doteq\dfrac{V-F}{2}>0$ and so $F<F+\alpha<V.$

Let $\tau=\{s_j\}_{j=1}^p$~ be an arbitrary partition. Arguing as in (\ref{raznxmpifi}), we show that
 \begin{multline*}
\bigl|v_{\tau}(X_m)-v_{\tau}(\Phi)\bigr|=\\=
\left|\sum\limits_{j=1}^p\left(\left|(*)\!\int\limits_{s_{j-1}}^{s_j}x_m(t)\,dg(t )\right|-
\left|(*)\!\int\limits_{s_{j-1}}^{s_j}x(t)\,dg(t)\right|\right)\right|\leqslant \\
\leqslant\sum\limits_{j=1}^p(*)\!\int\limits_{s_{j-1}}^{s_j}\bigl|x_m(t)-x(t)\bigr|\ ,dg_{\pi}(t)=
(*)\!\int\limits_{a}^{b}\bigl|x_m(t)-x(t)\bigr|\,dg_{\pi}(t)\to 0
\end{multline*}
$(m\to\infty).$ Note that for every $\varepsilon>0$ there is a number $m_0,$ {\it not depending on  partition} $\tau,$
such that
\begin{equation}
\label{raznsum}
\bigl|v_{\tau}(X_m)-v_{\tau}(\Phi)\bigr|<\varepsilon\;\;\;\text{for}\;\;m\geqslant m_0.
\end{equation}
Let us take $\varepsilon<\alpha$ and fix $m_0$  according to (\ref{raznsum}). Let us also fix number $m_1\;(\geqslant m_0)$
such that $\bigl|\bigvee\limits_a^b\bigl(X_{m_1}\bigr)-V\bigr|<\dfrac{\varepsilon}{2};$
we finally fix a partition $\tau^{'}_1$ such that
$v_{\tau^{'}_1}(X_{m_1})>\bigvee\limits_a^b\bigl(X_{m_1}\bigr)-\dfrac{\varepsilon}{2}.$ Then
$v_{\tau^{'}_1}\bigl(X_{m_1}\bigr)>V-\varepsilon.$ Given that
$v_{\tau^{'}_1}(\Phi)\leqslant F,$ we obtain a contradiction with the choice of $\varepsilon:$
$$
\varepsilon>v_{\tau^{'}_1}\bigl(X_{m_1}\bigr)-v_{\tau^{'}_1}(\Phi)>F+\alpha-F=\alpha>\varepsilon.
$$
So, $V=F.$
\hfill $\square$

\begin{sle}
\label{varac}
Let $x\in {\bf{AC}}$ be such that $x^{'}\in {\bf{N}}.$ Then
\begin{equation}
\label{fvarac}
\bigvee\limits_a^b(x)=\int\limits_a^b|x^{'}(t)|\,ds
\end{equation}
\end{sle}
(see Theorem 8 in \cite[IX.4]{nat97}).


\bigskip

\poonkt
{Approximation by absolutely continuous functions\,}\phantom{01234567890123456789}\phantom{01234567890123456789}\phantom{01234567890123456789} \

Often, when dealing with Stieltjes-type integrals, one needs to approximate discontinuous integrable
and integrating functions  by continuous functions
(see, for example, \cite{derr88},\, \cite{fil},\, \cite{Lev}), so that the integral remains close to the original one.
The limit theorems of the previous section cannot be applied directly to the (*)-integral
since they require quasi-uniform convergence. Using Theorem \ref{rlpr} below, we prove a limit theorem of the required type.

It will be more convenient to pass from a discrete parameter ($n\to \infty$) to a continuous one ($\varepsilon\to 0+$). Since this transition does not affect the subject matter, 
the results from the previous sections will be used without any additional explanations.

Let $y\in{\bf{R}}$ and, in addition, let $y(\cdot)$ be right-\;\;$\Bigl($left-$\Bigr)$continuous. 
According to Theorem \ref{rlpr}, we put
\begin{equation}
\label{aprrl}
y_{\varepsilon}(t)\doteq \dfrac{1}{\varepsilon}\int\limits_{t}^{t+\varepsilon}y(s)\,ds\qquad
\left(y_{\varepsilon}(t)\doteq \dfrac{1}{\varepsilon}\int\limits_{t-\varepsilon}^{t}y(s)\,ds\right)\quad (\varepsilon>0).
\end{equation}
In general, if $y\in{\bf{R}}\;\;\;\Bigl(y\in{\bf{BV}},\;\;y\in{\bf{CH}}\Bigr)$, then we put
\begin{equation}
\label{apr}
y(t)=y_{+}(t)+y_{-}(t),\;\;\;
y_{\varepsilon}(t)\doteq \dfrac{1}{\varepsilon}\int\limits_{t}^{t+\varepsilon}y_{+}(s)\,ds+
\dfrac{1}{\varepsilon}\int\limits_{t-\varepsilon}^{t}y_{-}(s)\,ds.
\end{equation}
 $(\varepsilon>0).$ Clearly,
\begin{multline}
\label{aprsv}
y_{\varepsilon}\in{\bf{AC}},\;\;\text{and}\;\;y_{\varepsilon}^{'}\in{\bf{R}}\;\;
\Bigl(y_{\varepsilon}^{'}\in{\bf{BV}},\;\;y_{\varepsilon}^{'}\in{\bf{CH}}\Bigr),\\
y_{\varepsilon}(t)\to y(t)\;\;\text{at}\;\;\varepsilon\to 0.
\end{multline}
Also, {\it if $y\in{\bf{BV}}\;\;\Bigl(y\in{\bf{CH}}\Bigr),$ then} \\
\begin{equation}
\label{shvar+}
\bigvee\limits_{a}^{b}\bigl(y_{\varepsilon}\bigr)\to\bigvee\limits_{a}^{b}(y)\qquad
\Bigl(\bigvee\limits_a^b\bigl(y_{\gd}\bigr)_{\varepsilon}\to\bigvee\limits_{a}^{b}(y_{\gd})\Bigr)\; \;
(\varepsilon\to 0).
\end{equation}
It suffices to prove the first assertion.

First, let $y(\cdot)$ be right continuous. Since $y_{\varepsilon}\in{\bf{AC}},$ we have
\begin{multline*}
\bigvee\limits_{a}^{b}\bigl(y_{\varepsilon}\bigr)=\int\limits_a^b |y^{'}_{\varepsilon}(t)|\,dt=
\dfrac{1}{\varepsilon}\int\limits_a^b |y(t+\varepsilon)-y(t)|\,dt\leqslant
\dfrac{1}{\varepsilon}\int\limits_a^b \bigvee\limits_{t}^{t+\varepsilon}(y)\,dt=\\=
\dfrac{1}{\varepsilon}\int\limits_a^b \bigvee\limits_{a}^{t+\varepsilon}(y)\,dt
-\dfrac{1}{\varepsilon}\int\limits_a^b\bigvee\limits_{a}^{t}(y)\,dt=
\dfrac{1}{\varepsilon}\int\limits_{a+\varepsilon}^{b+\varepsilon} \bigvee\limits_{a}^{t}(y)\,dt
-\dfrac{1}{\varepsilon}\int\limits_a^b\bigvee\limits_{a}^{t}(y)\,dt=\\=\!
\dfrac{1}{\varepsilon}\int\limits_b^{b+\varepsilon}\bigvee\limits_{a}^{t}(y)\,dt\!-\!
\dfrac{1}{\varepsilon}\int\limits_a^{a+\varepsilon}\bigvee\limits_{a}^{t}(y)\,dt=\!
\!\bigvee\limits_{a}^{b}(y)-
\dfrac{1}{\varepsilon}\int\limits_a^{a+\varepsilon}\bigvee\limits_{a}^{t}(y) \,dt\leqslant\bigvee\limits_{a}^{b}(y).
\end{multline*}
We obtain by Theorem \ref{shvar} convergence (\ref{shvar+}). We arrive at the same conclusion for the left-continuous
functions $y(\cdot).$ In the general case, convergence (\ref{shvar+}) follows from definition (\ref{apr}).
\hfill $\square$

\begin{teo}
\label{osneps}
{\it Let $x\in{\bf{R}},\;g\in {\bf{BV}}.$ Define $x_{\varepsilon}(\cdot)$ and $g_{\varepsilon}(\cdot )$
according to $(\ref{apr}),$ so that, due to $(\ref{aprsv})$,
\begin{equation}
\label{property2}
x_{\varepsilon}(t)\to x(t),\;\;g_{\varepsilon}(t)\to g(t)\;\;\;\;(\varepsilon\to 0).
\end{equation}
Then}
\begin{equation}
\label{osnpredperech03}
(*)\!\int\limits_a^b x_{\varepsilon}(t)\,dg_{\varepsilon}(t)\to(*)\!\int\limits_a^b x(t)\,dg(t )\qquad (\varepsilon\to 0),
\end{equation}
\begin{equation}
\label{shvarneoprint}
\bigvee\limits_a^b\bigl(\Phi_{\varepsilon}\bigr)\to \bigvee\limits_a^b(\Phi)\qquad (\varepsilon\to 0),
\end{equation}
{\it where}
$$
\Phi_{\varepsilon}(t)=(*)\!\int\limits_a^t x_{\varepsilon}(s)\,dg_{\varepsilon}(s),\quad
\Phi(t)=(*)\!\int\limits_a^t x(s)\,dg(s).
$$
\end{teo}

The proof of this theorem is based on the following auxiliary result.

Let
$\;
J_{\varepsilon}(x,g)=(*)\!\int\limits_a^b x_{\varepsilon}(t)\,dg_{\varepsilon}(t),\qquad
I_{\varepsilon}(x,g)=(*)\!\int\limits_a^b x_{\varepsilon}(t)\,dg(t).
$
\begin{lem}
\label{osnepslem}
{\it Under the assumptions of the theorem},
\begin{equation}
\label{Vspomogatnpredst04}
J_{\varepsilon}(x,g)=I_{\varepsilon}(x,g)+o(1)\;\;\text{at}\;\;\varepsilon\to 0,
\end{equation}
\end{lem}
\doc\,
Assume first that $x(\cdot),\;g(\cdot)$ are right continuous.
\begin{multline*}
J_{\varepsilon}(x,g)=(*)\!\int\limits_a^b x_{\varepsilon}(t)\,dg_{\varepsilon}(t)=
\int_a^b x_{\varepsilon}(t)\,dg_{\varepsilon}(t)=\\=\int_a^b\left(\dfrac{1}{\varepsilon}\int\limits_{t}^ {t+\varepsilon}x(s)\,ds\right)\,
d\left(\dfrac{1}{\varepsilon}\int\limits_{t}^{t+\varepsilon}g(s)\,ds\right)=\\=
\int_a^b\dfrac{g(t+\varepsilon)-g(t)}{\varepsilon}\,dt\left(\dfrac{1}{\varepsilon}\int\limits_{t}^{t+\varepsilon }x(s)\,ds\right)=\\=
\int_a^b\,\dfrac{g(t+\varepsilon)-g(t)}{\varepsilon}\,\left(\dfrac{1}{\varepsilon}\int\limits_{a}^{a+\varepsilon}x(p+t-a)\,dp\right)\,dt=\\=
\dfrac{1}{\varepsilon}\int\limits_{a}^{a+\varepsilon}\left(\int_a^b x(p+t-a)\,\dfrac{g(t+\varepsilon)-g(t) }{\varepsilon}\,dt\right)\,dp
\end{multline*}
(reversing the order of integration in the iterated Riemann integral is allowed according to \cite[sect. 594]{ficht3}).
\begin{multline*}
I_{\varepsilon}(x,g)=(*)\!\int\limits_a^b x_{\varepsilon}(t)\,dg(t)=
(*)\!\int\limits_a^b\left(\dfrac{1}{\varepsilon}\int\limits_{t}^{t+\varepsilon}x(s)\,ds\right)\,dg( t)=\\=
(*)\!\int\limits_a^b\,\left(\dfrac{1}{\varepsilon}\!\int\limits_{a}^{a+\varepsilon} x(p+t-a)\,dp\right)\,dg(t)=\\=
\dfrac{1}{\varepsilon}\int\limits_{a}^{a+\varepsilon}\left((*)\!\int\limits_{a}^{b} x(p+t-a)\,dg (t)\right)\,dp
\end{multline*}
(the change in the order of integration is carrie out according to Theorem (\ref{chali})).

Consider the difference
$
I_{\varepsilon}(x,g)-J_{\varepsilon}(x,g)=\dfrac{1}{\varepsilon}\int\limits_{a}^{a+\varepsilon}\,\eta(p ,\,\varepsilon)\,dp,
$\;
where
$$
\eta(p,\,\varepsilon)\doteq (*)\!\int\limits_{a}^{b} x(p+t-a)\,dg(t)-\int_a^b x(p+t-a) \,\dfrac{g(t+\varepsilon)-g(t)}{\varepsilon}\,dt
$$
Since
\begin{equation}
\label{chitr}
\dfrac{g(t+\varepsilon)-g(t)}{\varepsilon}\,dt=d\dfrac{1}{\varepsilon}\int\limits_{t}^{t+\varepsilon}g(s) \,ds
\end{equation}
then, denoting
$G(t,\,\varepsilon)\doteq g(t)-g_{\varepsilon}(t),$
we arrive at equality
\begin{multline*}
\eta(p,\,\varepsilon)=(*)\!\int\limits_{a}^{b} x(p+t-a)\,d\left(g(t)-\dfrac{1}{ \varepsilon}\int\limits_{t}^{t+\varepsilon}g(s)\,ds\right)=\\=
(*)\!\int\limits_{a}^{b} x(p+t-a)\,d\bigl(g(t)-g_{\varepsilon}(t)\bigr)=(*)\! \int\limits_{a}^{b} x(p+t-a)\,dG(t,\,\varepsilon)
\end{multline*}

According to Theorem \ref{predstvar01},
\begin{equation}
\label{chitr01}
G(t,\,\varepsilon)\to 0\;\;\;\text{at}\;\;\varepsilon\to 0,\;\;t\in [a,b],
\end{equation}
the jump function $G(\cdot,\,\varepsilon)$ coincides with the jump function $g(\cdot),$ i.e. it does not depend on $\varepsilon$,
so that the relations (\ref{modrsk+}) hold for $G(\cdot,\,\varepsilon)$ and the series of jumps of  function $G(\cdot,\,\varepsilon)$ converges
uniformly with respect to $\varepsilon$ (it does not depend on this parameter at all).
Further (see the definition of norm $\|\cdot\|^{*}$ in (\ref{normbv2})),
\begin{multline}
\label{chitr02}
\bigvee_{a}^{b}\bigl(g_{\varepsilon}\bigr)=\dfrac{1}{\varepsilon}\int_{t}^{t+\varepsilon}|g(s)|\,ds \leqslant
\underset{t\in [a,b]}{sup}\;|g(t)|,\\
\bigvee_{a}^{b}\bigl(G\bigr)\leqslant \bigvee_{a}^{b}(g)+
\bigvee_{a}^{b}\bigl(g_{\varepsilon}\bigr)=\|g\|^{*}_{{\bf{BV}}}.
\end{multline}
Therefore, one can apply Helly's theorem in form {\bf{14.6}}\,(bis),
which yields $\eta(p,\,\varepsilon)\to 0\;\;\text{for}\;\varepsilon\to 0$ and any $p\geqslant a.$ 
Representation
(\ref{Vspomogatnpredst04}) for right-continuous functions is thus proved.
Arguing in a ``dual'' way, we obtain representation (\ref{Vspomogatnpredst04}) for left-continuous functions.

In the general case, according to Theorems \ref{predstvar01} and \ref{rlpr},
$$
(*)\!\int\limits_a^b x(t)\,dg(t)=(*)\!\int\limits_a^b \bigl(x_{+}(t)+x_{-}(t) \bigr)\,d\bigl(g_{+}(t)+g_{-}(t)\bigr).
$$
For integrals
$(*)\!\int\limits_a^b x_{+}(t)\,dg_{+}(t),\;\;(*)\!\int\limits_a^b x_{-}(t )\,dg_{-}(t)$
the validity of presentation (\ref{Vspomogatnpredst04}) was established above.

Consider a pair of functions $\;\;\bigl(x_{+}(\cdot),\,g_{-}(\cdot)\bigr).\;$ The argument for \linebreak
$\bigl(x_{-}(\cdot),\,g_{+}(\cdot)\bigr)\;\;$ is similar.
Let us put, as above,
$$J_{\varepsilon}(x_{+},g_{-})
\doteq (*)\!\int\limits_a^b \bigl(x_{+}\bigr)_{\varepsilon}(t)\,d\bigl(g_{-}\bigr)_{\varepsilon}( t),
\quad I_{\varepsilon}(x_{+},g_{-})\doteq(*)\!\int\limits_a^b \bigl(x_{+}\bigr)_{\varepsilon}(t) \,dg_{-}(t).$$ The difference
$$I_{\varepsilon}(x_{+},g_{-})-J_{\varepsilon}(x_{+},g_{-})=(*)\!\int\limits_a^b \bigl(x_{+}\bigr)_{\varepsilon}(t)\,dg_{-}(t)-
(*)\!\int\limits_a^b \bigl(x_{+}\bigr)_{\varepsilon}(t)\,d\bigl(g_{-}\bigr)_{\varepsilon}(t)$$
is a difference of RS-integrals,\; which can be put in the form
\begin{multline*}
\dfrac{1}{\varepsilon}\int\limits_{a}^{a+\varepsilon}dp\left((*)\!\int_a^b x_{+}(p+t-a)\,dg_{-} (t)\right)-\\-
\dfrac{1}{\varepsilon}\int\limits_{a}^{a+\varepsilon}dp\int_a^b x_{+}(p+t-a)\,\dfrac{g_{-}(t)-g_ {-}(t-\varepsilon)}{\varepsilon}\,dt=\\=
\dfrac{1}{\varepsilon}\int\limits_{a}^{a+\varepsilon}dp
\left((*)\!\int_a^b x_{+}(p+t-a)\,d\bigl(G_{-}\bigr)_{\varepsilon}(t)\right),
\end{multline*}
where $\bigl(G_{-}\bigr)_{\varepsilon}(t)=g_{-}(t)-\bigl(g_{-}\bigr)_{\varepsilon}(t).$
Denote the inner integral by
$$
\Phi(p,\,\varepsilon)\doteq(*)\!\int_a^b x_{+}(p+t-a)\,d\bigl(G_{-}\bigr)_{\varepsilon}( t).
$$
Similarly to (\ref{chitr01}),\,(\ref{chitr02}), we establish that
$$
\bigl(G_{-}\bigr)_{\varepsilon}(t)\to 0\;\;(\varepsilon\to 0,\;t\in [a,b]),\;\;\;
\bigvee_{a}^{b}\bigl(G_{-}\bigr)_{\varepsilon}\leqslant \|g\|^{*}_{{\bf{BV}}},
$$
and also show the independence of the jump function and the jump series $\bigl(G_{-}\bigr)_{\varepsilon}(\cdot)$ from $\varepsilon.$
Therefore, by Helly's theorem in form {\bf{14.6}}\,(bis),
\begin{equation}
\label{chitr03}
\Phi(p,\,\varepsilon)\to 0\;\;\;\text{if}\;\;\varepsilon\to 0\;\;\;\text{and any}\;\; p\geqslant a.
\end{equation}

Let $p\to a+;$ for a fixed $\varepsilon>0$, by Theorem \ref{Fst2},
\begin{equation}
\label{chitr04}
\Phi(p,\,\varepsilon)\to (*)\!\int_a^b x_{+}(t)\,d\bigl(G_{-}\bigr)_{\varepsilon}(t)= \Phi(a,\,\varepsilon),
\end{equation}
So, function $\Phi(\cdot,\,\varepsilon)$ is right continuous at point $a.$

Equality (\ref{Vspomogatnpredst04}) is equivalent to 
\begin{equation}
\label{chitr05}
I_{\varepsilon}(x_{+},g_{-})-J_{\varepsilon}(x_{+},g_{-})=
\dfrac{1}{\varepsilon}\int\limits_{a}^{a+\varepsilon}\Phi(p,\,\varepsilon)\,dp\to 0\;\;\;(\varepsilon\to 0 )
\end{equation}
Let us prove this. To this end, we will prove that
\begin{equation}
\label{chitr06}
\Psi(\varepsilon,\,\delta)=\dfrac{1}{\delta}\int\limits_{a}^{a+\delta}\bigl|\Phi(p,\,\varepsilon)\bigr| \,dp\to 0\quad
\text{at}\;\varepsilon\to 0,\;\delta\to 0
\end{equation}
(meaning the double limit).

Let us suppose the opposite. This means that there are $\gamma_0>0$ and sequences $\varepsilon_n\to 0,\;\delta_n\;\;(n\to\infty)$
such that
$$
\dfrac{1}{\delta_n}\int\limits_{a}^{a+\delta_n}|\Phi(p,\,\varepsilon_n)|\,dp\geqslant\gamma_0.
$$
In view of (\ref{chitr03}), there is a positive integer $N_1$ such that
$|\Phi(a,\,\varepsilon_n)|<\dfrac{\gamma_0}{3}$ for all $n>N_1;$
due to (\ref{chitr04}) there is a positive integer $N_2\;\;(\geqslant N_1)$ such that
$|\Phi(p,\,\varepsilon_n)-\Phi(a,\,\varepsilon_n)|<\dfrac{\gamma_0}{3}$
for $a\leqslant p\leqslant a+\delta_n$ for all $n>N_2;$
\;then for all $n>N_2\;$
$\;|\Phi(p,\,\varepsilon_n)|\leqslant |\Phi(a,\,\varepsilon_n)|+|\Phi(p,\,\varepsilon_n)-\Phi(a,\ ,\varepsilon_n)|<\dfrac{2\gamma_0}{3}$
and
$$
\gamma_0\leqslant \dfrac{1}{\delta_n}\int\limits_{a}^{a+\delta_n}|\Phi(p,\,\varepsilon_n)|\,dp\leqslant\dfrac{2\gamma_0} {3}
$$
We arrive at a contradiction. Therefore, (\ref{chitr06}), (\ref{chitr05}), and equality
 (\ref{Vspomogatnpredst04}) are true.
  \hfill $\square$

Proof of the theorem.\,
 By Lemma \ref{osnepslem}, representation (\ref{Vspomogatnpredst04}) is valid, and by Theorem \ref{Fst2} 
 $$
 I_{\varepsilon}(x,g)=(*)\!\int\limits_a^b x_{\varepsilon}(t)\,dg(t)\to (*)\!\int\limits_a^b x( t)\,dg(t)\;\;(\varepsilon\to 0),
 $$
 that is, (\ref{osnpredperech03}) holds.

Since (see (\ref{shvar+})) $\bigl(g_{\varepsilon}\bigr)_{\pi}(t)\to g_{\pi}(t)\;\;(\varepsilon \to 0),$
convergence (\ref{shvarneoprint}) follows from representation (\ref{teoravvar}) and Theorem \ref{Fst2}.
\hfill $\square$
 
 If $x\in{\bf{BV}},\;\;g\in{\bf{CH}},$ then representation (\ref{apr}) remains valid, we have (\ref{aprsv}),\,
 (\ref{property2}),
 therefore, by formula (\ref{bp1}), we have the following corollary.
 \begin{sle}
 \label{sledstosn}
{\it Let $x\in{\bf{BV}},\;\;g\in{\bf{CH}},\;\;\varepsilon\to 0.$ Then}
$$
(*)\!\int\limits_a^b x_{\varepsilon}(t)\,dg_{\varepsilon}(t)\to(*)\!\int\limits_a^b x(t)\,dg(t )+
 \sum\limits_{t\in T}\Bigl(\sigma^{+}_t(x)\sigma^{+}_t(g)-\sigma^{-}_t(x)\sigma^{-} _t(g)\Bigr).
$$
\end{sle}

\begin{center}
\begin{large}
\section{On differential equations with generalized functions}
\end{large}
\end{center}

\hskip 1mm
\poonkt{Statement of the problem\,}\phantom{0123456789012345678901234567890123456789} \

{\bf{1.}}\;
Ordinary differential equations with generalized functions in the coefficients have been in the focus of researchers' attention in the last few decades of the twentieth century (see, for example, \cite{Anoch},\, \cite {Colomb85}~--- \cite{derr00},\, \cite{derr18},\, \cite{derr96},\, \cite{DK05},\,
\cite{Egorov90},\, \cite{fil},\, \cite{Kurz57}~--- \cite{MisSam},\, \cite{Orlov88},\, \cite{RadNgu90},\, \cite{SamPer},\, \cite{Sar95},\,
\cite{silvavint97},\, \cite{ZavSes91},
and also in this century (see the bibliography in \cite{derr18},\, \cite{Rod21})).
Let us illustrate one of the problems that arises in the study of such equations using the following
example of a scalar ordinary linear differential equation of order $n$ $(n\geqslant 2).$ Cauchy problem for such an equation
is usually written in the following form (see, for example,  \cite{fil},\, \cite{Lev}; see also \cite{derr18}):
\begin{multline}
\label{4eqsk}
\bigl(Lx\bigr)(t)\doteq x^{(n)}+p'_1(t) x^{(n-1)}+\ldots + p'_n(t)x=p'_ {n+1}(t),\;\;x^{(k-1)}(a)=\gamma_{k},\\
k=1,\ldots,n,
\end{multline}
where the coefficients $p'_k(\cdot)$ are the derivatives in the sense of the theory of generalized functions of, 
generally speaking, discontinuous functions.
Expression $p'_1(t) x^{(n-1)}\;$ (and with it the whole expression $Lx$)
is devoid of the classical meaning since, formally, there is a product
of a generalized function $\bigl(p'_1(\cdot)\bigr)$ and a discontinuous function
$\bigl(x^{(n-1)}(\cdot)\bigr).$
Such multiplication is impossible in the theory of generalized functions. Each approach to defining a solution to problem
(\ref{4eqsk}) is at the same time a way to define such a product (see \cite{derr88},\, \cite{fil},\, \cite{Lev} for details).
In \cite{derr88}, the author reduced problem (\ref{4eqsk}) to a Cauchy problem for a quasi-differential 
equation (see below) with usual functions in the coefficients. This reduction is based on the use of LS-integrals $\int\limits_a^t p_1(s)\,dp_k(s)\;(k=2,\ldots,n+1).$
It is shown that if the anti-derivatives of coefficients $p_k(\cdot)$ are absolutely continuous (in this case it is a classical problem),
then the solution of problem (\ref{4eqsk}) coincides with the solution of the quasi-differential problem. If some conditions are satisfied
(see, for example, condition $\mathfrak B$ below),
then the quasi-differential problem does not lose its meaning, and the solution of problem (\ref{4eqsk}) means the solution of
the quasi-differential problem.


\medskip

{\bf{2.}}
We will consider problem (\ref{4eqsk}) under the following \linebreak
assumptions: we say that condition 

$\mathfrak A_k\;\;(k=n,n+1)$ holds  if $p_i\in {\bf{AC}}
\;\;(i=1,\ldots,k);$

$\mathfrak V\;\;$ holds if $\;p_i\in {\bf{BV}}\;\;(i=1,\ldots,n+1);\;$
if at the same time $\;p_i(\cdot)\;\;(i=
2,\ldots,n+1)$ have no common {\it one-sided} discontinuities with $p_1(\cdot)$, then we say that conditions  $\mathfrak V_{\delta}$ hold;

$\mathfrak B$ holds if $p_1$ is B-measurable and bounded, and $p_k\in{\bf{BV}}\;(k=\linebreak= 
2,\ldots,n+1);$

$\mathfrak C$ holds if $p_1\in {\bf{N}},\;\;p_k\in{\bf{BV}}\;\;(k=2,\ldots,n+1); $ we note as obvious 
obvious: \;
conditions $\;\mathfrak C\;$ imply conditions $\mathfrak B$;

$\mathfrak{C}_{\delta}$ holds if $p_1\in{\bf{R}},\;\;p_k\in {\bf{BV}},\;\;(k=2,\ldots ,n+1);$

$\mathfrak D$ holds if $p_1\in {\bf{BV}},\;\;p_k\in{\bf{CH}}\;\;(k=2,\ldots,n+1). $

If conditions $\mathfrak A_{n+1}$ are satisfied, then problem (\ref{4eqsk})~ is classical (the coefficients of the 
equation are summable functions).
Conditions $\mathfrak V,\;\mathfrak V_{\delta}$ were thoroughly studied in \cite{Lev}. There it is established
that conditions $\mathfrak V_{\delta}$ (in the terminology of \cite{Lev} -- ``$\delta$-correctness'')
are necessary and sufficient for the solution of problem (\ref{4eqsk}) to be obtained from the solutions of a sequence of
classical problems using \;{\it one}\; passage to the limit (the necessity is due to the use of the Darboux--Stieltjes integral).
Conditions $\mathfrak B$ were considered by the author in
\cite{derr88} --- \cite{derr95},\, \cite{derr18}. There, in absence of $\delta$-correctness for the aforementioned approximation, one needs 
\;{\it two}\; passages to the limit.
Conditions $\;\;\mathfrak C,\;\mathfrak{C}_{\delta},\;\mathfrak D\;\;$ are considered below. The use of
the *-integral allows to obtain the solution as the limit of a sequence of classical problems using \;{\it one}\; passage to the limit.


\bigskip

\poonkt{Quasi-differential equations\;}\phantom{0123456789012345678901234567890123456789} \

{\bf{1.}}\;
Let
$\mathcal{P}\!=\!\bigl(p_{ik}\bigr)_{i,k=0}^n$~ be a lower triangular matrix,
$p_{ik}:[a,b]\to\mathbb R$, 
$p_{ii}(t)>0\;(i=0,\ldots,n)$ a.e. to $[a,b].\;$
We assume that conditions $\mathbb P$ are satisfied: functions
$$
1/p_{ii}(\cdot)\;(i=1,\ldots,n-1),\;p_{ik}(\cdot)/p_{ii}(\cdot)\;(i=1 ,\ldots,n,\,k=0,\ldots,i-1)
$$
are summable on $\;[a,b].\;$
Let us define the quasi-derivatives of function \linebreak
$x:[a,b]\to \mathbb R$ by formulas (cf. \cite{Shin40}) :
\begin{equation}
\label{eqKdU0}
^0_{\mathcal{P}}x\doteq p_{00}x,\;\;
^k_{\mathcal{P}}x\doteq p_{kk}\dfrac{d}{dt}\left(^{k-1}_{\mathcal{P}}x\right)+
\sum\limits_{i=0}^{k-1}p_{ik}\left(^i_{\mathcal{P}}x\right)\;\;(k=1,\ldots,n).
\end{equation}
The equation
\begin{equation}
\label{eqKdU1}
\bigl(^n_{\mathcal{P}}x\bigr)(t)=f(t),\quad t\in [a,b]\quad (f:\,[a,b]\to \mathbb R)
\end{equation}
is called a linear quasi-differential equation.
A solution of equation (\ref{eqKdU1}) is a function $x:[a,b]\to \mathbb R$ having absolutely continuous quasi-derivatives
$^k_{\mathcal{P}}x,\;k=1,\ldots,n-1$ that satisfy a.e.\,on $[a,b]$ equation (\ref{eqKdU1}).

If function $f(\cdot)/p_{nn(\cdot)}$ is summable on $[a,b]$ (under conditions $\mathbb P$), then there is a unique solution to equation
(\ref{eqKdU1}) satisfying initial conditions
\begin{equation}
\label{eqKdU2}
\bigl(^k_{\mathcal{P}}x\bigr)(a)=\xi_k\quad (\xi_k\in\mathbb{R},\;\;k=0,\ldots,n-1),
\end{equation}
see \cite{derr88},\;\cite{derr18}.
In terms of the column-vectors
$$\xi=(\xi_0,\ldots,\xi_{n-1})^\top, \quad \;_{\mathcal{P}}x\doteq \bigl(^0_{\mathcal{P}} x,\ldots,^{n-1}_{\mathcal{P}}x\bigr)^\top$$
(column of quasi-derivatives), initial conditions (\ref{eqKdU2}) are written as follows: $_{\mathcal{P}}x(a)=\xi.$ Note that the last raw of
matrix $\mathcal{P}$ does not appear in the definition of column $_{\mathcal{P}}x$. 
We will need the column of ordinary derivatives
$_\mathcal{E}x=\bigl(x,x^{'},\ldots,x^{(n-1)}\bigr)^\top=_{\mathcal{P}}x$
for $\mathcal{P}=\mathcal{E}_{n+1}$, where $\mathcal{E}_{n+1}$~ is the identity matrix of order $n+1.$
The above statement from \cite{derr88} follows from the fact that problem (\ref{eqKdU1}), (\ref{eqKdU2}) is equivalent to the following 
Cauchy problem with respect to $\;z(t)\doteq \bigl(_{\mathcal P}x\bigr)(t)$:
\begin{equation}
\label{eqKdUeqvS}
z^{'}=A(t)z+F(t)\;\;(t\in [a,b]),\qquad z(a)=\xi,
\end{equation}
where 
$$
\;F(t)=\left(0,\ldots,0,\dfrac{f(t)}{p_{nn}(t)}\right)^{\top},\;\;\; A(t)=\bigl(a_{ij}(t)\bigr)_{i,j=1}^n,\;\;a_{ij}(t)=0
$$ %
$(i=1,\ldots,n-2,\;\;j=i+2,\ldots,n),\;\;a_{i,i+1}(t)=\dfrac{1} {p_{ii}(t)}\;\;(i=1,\ldots,n-1),\linebreak
a_{ij}(t)=-\dfrac{p_{i,j-1}(t)}{p_{ii}(t)}\;\;$ $\;(i=1,\ldots,n ,\;j=1,\ldots,i).$


\medskip

{\bf{2.}}\;Assume that conditions $\mathfrak C$ are satisfied.
We construct lower triangular matrices $H=\bigl(h_{ik}\bigr)_{i,k=1}^n$ and $\mathcal{P}$ from the primitives of the coefficients of  equation (\ref{4eqsk}) in the following way:
\begin{equation}
\label{strKdU1}
h_{11}(t)\doteq 1,\quad h_{nn}(t)\doteq \exp\bigl(p_1(t)-p_1(a)\bigr)=h_{n+1,n+1}( t),
 \end{equation}
\begin{equation}
\label{strKdU2}
h_{n,n-k+1}(t)\doteq (*)\!\int\limits_{a}^t h_{nn}(s)\,dp_k(s)\;\;\;(k =2,\dots,n,\;t\in I),
 \end{equation}
\begin{equation}
\label{strKdU3}
h_{kk}(t)\doteq exp\int\limits_{a}^t \dfrac{h_{k+1,k}(s)}{h_{k+1,k+1}(s)}\ ,ds,
 \end{equation}
\begin{equation}
\label{strKdU4}
h_{k,k-i+1}(t)\doteq \int\limits_{a}^t \dfrac{h_{kk}(s)h_{k+1,k-i+1}}{h_{ k+1,k+1}(s)}\,ds,
 \end{equation}
$k=n-1,\ldots,2;\;i=2,\ldots,k;\;t\in I\;$
(for $n=2$ equalities (\ref{strKdU3}),\,(\ref{strKdU4}) are absent);

Further, we assume (for brevity, we omit argument $t$)
\begin{equation}
\label{strKdU5}
p_{00}=1,\quad p_{kk}=h_{k+1,k+1}/h_{kk}\;\;\;(k=1,\ldots,n),
\end{equation}
\begin{equation}
\label{strKdU6}
\left.
\begin{array}{lcr}
b^0_{ki}=-h_{k,k-i}\;\;(i=1,\ldots,k-1);\quad b^0_{kk}=0,\\
b^j_{kj}=\dfrac{b^{j-1}_{kj}}{h_{k-j+1,k-j+1}},\;\;b^j_{ki}= b^{j-1}_{ki}-\dfrac{b^{j-1}_{kj}\cdot h_{k-j+1,k-i+1}}{h_{k-j+ 1,k-j+1}}
\end{array}
\right\}
\end{equation}
$(i=j+1,\ldots,k,\;j=1,\ldots,k,$ for $j=k$ the last equality disappears; 
$k=2,\ldots,n$);
\begin{multline}
\label{strKdU7}
p_{10}=h_{21},\;p_{k,k-j}=p_{kk}\cdot b^j_{kj}\;\;(j=1,\ldots,k,\;k=2 ,\ldots,n,\;t\in[a,b]),\\
\mathcal P=\bigl(p_{ik}\bigr)_{i,k=0}^n.
\end{multline}


{\bf{3.}}\;The next properties follow directly from the definitions of matrices $H$ and $\mathcal{P}$
(in parentheses we list those properties that take place when conditions $\mathfrak A_n$ are satisfies):

1)\;$h_{kk}(t)>0,\;\;p_{kk}(t)>0,\;\;\det\,H(t)=\prod\limits_{i=1} ^nh_{ii}(t)>0,\;\;H(a)=\mathcal{E}_n\linebreak
(t\in [a,b],\;\;k=1,\ldots,n);$

2)\;$h_{nn}\in {\bf{N}}\;\;\bigl(h_{nn}\in {\bf{AC}}\bigr),\;\;h_{nk} \in {\bf{BV}}\;\;
\bigl(h_{nk}\in {\bf{AC}}\bigr),\;\;h_{nk}(\cdot)$ are continuous at the points of continuity $p_{n-k+1}\;\; (k=1,\ldots,n-1);$

3)\;$h^{(n-k)}_{kk}\in {\bf{N}}\;\;\bigl(h^{(n-k)}_{kk}\in {\bf{AC }}\bigr),\;
h^{(n-k)}_{ki}\in {\bf{BV}}\;\;\bigl(h^{(n-k)}_{ki}\in {\bf{AC}}\bigr) \linebreak
(k=1,\ldots,n-1,\;i=1,\ldots,k-1);$

4)\;the elements of matrices $H^{-1}_k\;\;\bigl(H_k=(h_{ij})_1^k\bigr)$ have derivatives of order $n-k\;\;(k=1, \ldots,n)$
in ${\bf{N}}\;\;\bigl({\bf{AC}}\bigr);$

5)\;$1/p_{kk}(\cdot),\;\;p_{k,k-j}(\cdot)/p_{kk}(\cdot)\;$ have derivatives of order $n-k-1\linebreak
(j=1,\ldots,k,\;k=1,\ldots,n-1),$
$(k=1,\ldots,n)$ in ${\bf{AC}};$

6)\;$p_{nn}(t)=1,\;\;p_{n,n-i}\in {\bf{N}}\;\;\bigl({\bf{AC}}\bigr )\;\;(i=1,\ldots,n);$

7)\;the properties of the elements of matrix $A$ in (\ref{eqKdUeqvS}), constructed from the elements of matrix $\mathcal P,$ are determined by
properties 5) and 6); thus, the elements in the last row of matrix $A$ are the ``worst'':
they belong to B-space ${\bf{N}}\;\;\bigl({\bf{AC}}\bigr).$

These properties mean that conditions $\mathbb P$ are satisfied for matrix $\mathcal P$. Therefore, Cauchy problem (\ref{eqKdU1}),\,(\ref{eqKdU2})
is uniquely solvable for any $\xi\in\mathbb R^n$ and any $f\in {\bf{N}}.$ Equalities (\ref{strKdU4}),\,(\ref{strKdU5}) imply
\begin{equation}
\label{sledst2}
p_{kk}\cdot h^{'}_{k,k-i+1}=h_{k+1,k-i+1}\;\;(i=1,\ldots,k,\, k=2,\ldots,n-1,\,n>2),
\end{equation}
and the relations (\ref{strKdU6}) are the formulas from Gauss' method for solving linear systems
\begin{equation}
\label{sledst3}
\sum\limits_{\nu=0}^{k-1}p_{k\nu}\cdot h_{\nu+1,1}=0,\;\;\;\;\sum\limits_{\ nu=0}^{k-1}p_{k\nu}\cdot h_{\nu+1,j+1}+p_{kk}\cdot h_{kj}=0
\end{equation}
with respect to  unknowns $p_{k0},\,p_{k,1},\ldots,p_{k,k-1}\;\bigl($here, for every $k\linebreak
(k=2,\ldots,n)$, one has associated a linear system of $k$
equations, the first of which (corresponding to $j=0$) is written separately; $j=1,\ldots,k-1$~ are numbers of the other equations$\bigr)$.
It follows from definitions (\ref{strKdU2}) and (\ref{defali}) that if $p_k(\cdot)\;\;(k=1,\ldots,n)$ are absolutely continuous
(that is, condition $\mathfrak A_n$ is satisfied), then a.e. on $[a,b]$
\begin{equation}
\label{sledst4}
h^{'}_{n,n-k+1}=h_nn\cdot p^{'}_{k}\quad (k=1,\ldots,n).
\end{equation}

\begin{lem}
\label{existkd1}
{\it If function $x:[a,b]\to\mathbb R$ has derivative $x^{(n-1)}(\cdot),$ then it also has quasi-derivatives
$\bigl(^k_{\mathcal{P}}x\bigr)(\cdot)\;\;(k=0,\ldots,n-1)$ and the following identities take place:
$$
^k_{\mathcal{P}}x=\sum\limits_{i=0}^{k}h_{k+1,k-i+1}\cdot x^{(k-1)}\quad ( k=0,\ldots,n-1),\eqno{(\Gamma_k)}
$$
or, in the vector-matrix form}, $_{\mathcal{P}}x=H _{\mathcal{E}}x.$
\end{lem}
\doc\,
Identity $(\Gamma_0)$ is trivial. The quasi-derivative $^1_{\mathcal{P}}x$ exists and identity $(\Gamma_1)$ is true:
$^1_{\mathcal{P}}x=h_{22}\cdot x^{'}+h_{21}\cdot x.$
Suppose that there exist quasi-derivatives up to order $k-1$ and identities $(\Gamma_1),\ldots,(\Gamma_{k-1})\linebreak
(k>1)$ are valid. Then, by definition (\ref{eqKdU0}) and property 3),
$$
^k_{\mathcal{P}}x=p_{kk}\cdot\sum\limits_{i=0}^{k-1}\Bigl(h^{'}_{k,k-i}x^{( k-i+1)}+h_{k,k-i}x^{(k-i)}\Bigr)+
\sum\limits_{\nu=0}^{k-1}p_{k\nu}\sum\limits_{i=0}^{\nu}h_{\nu+1,\nu-i+1} x^{(\nu-i)}.
$$
Applying $(\ref{sledst2})$ and $(\ref{strKdU5})$ successively, we regroup the terms and change the order of summation in the double sum.
As a result, we obtain
\begin{multline*}
^k_{\mathcal{P}}x=\sum\limits_{i=0}^{k}h_{k+1,k-i+1}x^{(k-i)}+x\cdot\sum\limits_{\nu=0}^{k-1}p_{k\nu}h_{\nu+1,1}+\\+
\sum\limits_{j=1}^{k-1}x^{(j)}\left(\sum\limits_{\nu=j}^{k-1}p_{k\nu}h_{\ nu+1,j+1}+p_{kk}h_{kj}\right).
\end{multline*}
Hence, due to linear system $(\ref{sledst3}),$ which corresponds to this $k$, we obtain identity $(\Gamma_k).$
By induction, $^k_{\mathcal{P}}x$
exist and all identities $(\Gamma_k)\;\;(k=1,\ldots,n-1)$ are satisfied.
\hfill $\square$

Properties 1)--- 4) allow us to reverse the last result.

\begin{lem}
\label{existkd2}
{\it If there exist absolutely continuous quasi-derivatives 
$^0_{\mathcal{P}}x,\,\ldots,\,^{(n-1)}_{\mathcal{P}}x,$ then there also exist
ordinary derivatives up to order $n-1,\linebreak
x^{(n-1)}\in{\bf{N}},$ and} $_{\mathcal{E}_n}x=H^{-1}\;_{\mathcal{P}} x.$
\end{lem}
\begin{lem}
\label{existkd3}
{\it If $\mathfrak A_n$ holds, then}
$$
\bigl(^n_{\mathcal{P}}x\bigr)(t)=h_{nn}(t)\bigl(Lx\bigr)(t).
$$
\end{lem}
\doc\,
By definitions \;(\ref{strKdU1}),\;(\ref{strKdU2}),
$h_{nk}\!\in\!{\bf{AC}}$ $(k=1,\ldots,n).$ Therefore, we can substitute identities
$(\Gamma_k)\;(k=1,\ldots,n-1)$ in definition (\ref{eqKdU0}) of a quasi-derivative of order $n.$
Arguing as in the proof of Lemma \ref{existkd1} and using successively identities (\ref{sledst2}),\,(\ref{strKdU5}),\,(\ref{sledst3})
for $k=n,$ we arrive at the assertion of the lemma.
\hfill $\square$

\medskip

Thus, in this section we proved the following result.
\begin{teo}
\label{svedkdu1}
{\it If conditions $\mathfrak A_{n+1}$ are satisfied, then problem 
$(\ref{4eqsk})$ is equivalent to problem $(\ref{eqKdU1}),\;(\ref{eqKdU2})$
with matrix $\mathcal{P}$ constructed in {\bf{2}} and}
$$f(t)=h_{nn}(t)\times p^{'}_{n+1}(t),\;\;\xi=H(a)\gamma\;\;\bigl( \gamma=(\gamma_0,\ldots,\gamma_{n-1}\bigr)^\top.$$
\end{teo}


\bigskip

\poonkt{Definition of solution of problem (\ref{4eqsk})\;} \ \phantom{01}

{\bf{1.}}\;
Consider a semi-homogeneous problem (\ref{4eqsk})
\begin{equation}
\label{exodu1}
\bigl(Lx\bigr)(t)=0\;\;\;(t\in [a,b]),\;\;_{\mathcal{E}_n}x(a)=\gamma.
\end{equation}
As we already noted, under conditions $\mathfrak C$ matrix $\mathcal{P}$ satisfies condition $\mathbb P.$
Therefore, the semi-homogeneous problem
\begin{equation}
\label{exqdu1}
\bigl(^n_{\mathcal{P}}x\bigr)(t)=0\;(t\in [a,b]),\;\;_{\mathcal{P}}x(a) =\xi
\end{equation}
has a unique solution $\widetilde x(\cdot),$ which, in view of Theorem \ref{svedkdu1}, 
can naturally be taken as the solution to problem (\ref{exodu1}).
So, under conditions $\mathfrak C$,
{\it the solution of  problem}  $(\ref{exqdu1})$ is defined as {\it the solution of  problem}  $(\ref{exodu1})$ 
 with matrix $\mathcal{P}$ and $\xi=H(a)\gamma$ defined in 
 {\bf{15.2.2}}.
          
This solution is uniquely determined, it exists and is unique.\linebreak  
According to Lemma \;\;\ref{existkd2}\;\; $$\widetilde x^{(n-1)}\in {\bf{N}}\quad \text{and} \quad \widetilde x^{(n)}$$ is a generalized 
function.

Let \;$X=\bigl(^{i-1}_{\mathcal{P}}x_k\bigr)_{i,k=1}^n$~ be such a fundamental matrix 
of homogeneous equation $\;\;^n_{\mathcal{P}}x=0$
that $X(a)=\mathcal{E}_n.\;$ Then $\;\;\bigl(_{\mathcal{P}}\widetilde x\bigr)(t)=\linebreak= 
X( t)H(a)\gamma,\;$ and, according to Lemma \;\ref{existkd2},
$$
\bigl(_{\mathcal{E}_{n+1}}\widetilde x\bigr)(t)=H^{-1}(t)X(t)X^{-1}(a)H (a)\gamma.
$$
$\bigl($Note that although $H(a),\;X(a)$~ are identity matrices
(see property 1)), we keep this notation in case the initial conditions and the lower limit of integration are different points.$\bigr)$
If $\widetilde X(\cdot)$~ is the fundamental matrix of the equation $\;\;Lx=0\;(\widetilde X(a)=\mathcal E_n),$
then, according to Lemma \ref{existkd2},
\begin{equation}
\label{svjafm}
\widetilde X(t)=H^{-1}(t)X(t).
\end{equation}


{\bf{2.}}\;
Le us return to the original problem (\ref{4eqsk}). Let us first assume that conditions $\mathfrak A_{n+1}$ are satisfied.
Let $M=\bigl(m_{ik}\bigr)_{i,k=1}^n$~ be the matrix
corresponding to equation $Lx=0,$ i.e.\;\;
$$m_{i,i+1}(t)=1\;(i=1,\ldots,n-1), \quad
m_{nk}(t)=-p^{'}_{n-k+1}(t) \;(k=1,\ldots,n),$$ and the remaining elements of this matrix are zero.
Thus $\widetilde X^{'}=M(t)\widetilde X$. Since $X^{'}=A(t)X$
$\bigl($recall: matrix $A$ is constructed from matrix $\mathcal P$
as was shown in subsection 
{\bf{15.2.1}}, 
see (\ref{eqKdUeqvS})$\bigr),$ then, due to equality (\ref{svjafm}) and the fact that the fundamental
matrix determines uniquely (continuous!) matrix of the system, we obtain
\begin{equation}
\label{vostA}
A=(H^{'}+HM)H^{-1}.
\end{equation}
We put
\begin{equation}
\label{pravchast}
h_{n0}(t)\doteq (*)\!\int\limits_{a}^t h_{nn}(s)\,dp_{n+1}(s)\qquad (t\in [a, b]),
\end{equation}
$h_0=(0,\ldots,0,h_{n0})\top;\;\;h_0(t)\in\mathbb R^n,\;\;(t\in [a,b]). $
\begin{teo}
\label{svedkdu2}
{\it If conditions $\mathfrak A_{n+1}$ are satisfied, then problem (\ref{4eqsk}) 
is equivalent to problem
\begin{equation}
\label{sist23}
y^{'}=A(t)y+F(t),\qquad y(a)=\xi
\end{equation}
for $\xi=H(a)\gamma,\;\;\;F=\bigl(h_{n0/h_{nn}}\bigr)(0,\ldots,-h_{n-1,n- 1},\,h_{n,n-1})^\top,$
and
\begin{equation}
\label{prpodst}
_{\mathcal{E}_n}x=H^{-1}(y+h_0),
\end{equation}
where $x$~ is the solution of  problem} (\ref{4eqsk}).
\end{teo}
\doc\,
Let $m=(0,\ldots,0,\,p^{'}_{n+1})^\top.$ Problem (\ref{4eqsk}) can be written as system
\begin{equation}
\label{sist27}
\bigl(_{\mathcal{E}_n}x\bigr)^{'}=M(t)_{\mathcal{E}_n}x+m(t),\qquad _{\mathcal{E} _n}x(a)=\gamma.
\end{equation}
Since now $h^{'}_0=h_{nn}m=Hm,\;\;Ah_0=F,$ in view of equality (\ref{vostA}), we see that substitutions (\ref{prpodst}) and
$y=H(_{\mathcal{E}_n}x-h_0)$ transfer problem (\ref{sist27}) into problem (\ref{sist23}) and vice versa.
\hfill $\square$

Let conditions $\mathfrak C$ be satisfied. We define $h_{n0}(\cdot),\;F(\cdot),\;\xi$ as above. Now, the components of $F(\cdot)$
belong to ${\bf{N}}$. Problem (\ref{sist23}) has unique solution $y(\cdot)$. This solution has absolutely continuous components, 
and the components of $y^{'}(\cdot)$ belong to ${\bf{N}}.$

As above, in subsection 
{\bf{15.3.1}},
{\it function $x(\cdot)$ defined by equality} (\ref{prpodst}) is said to be 
{\it a solution of problem}  $(\ref{4eqsk})$.

Thus defined solution {\it exists and is unique}. Furthermore, its derivative of order $n-1$~ is a function from {\bf{N}},
and the derivative of order $n$~ is a generalized function.

Under conditions $\mathfrak V_{\delta}$, solution of problem $(\ref{4eqsk})$ defined above coincides with  
the solution in the sense of \cite{Lev} since
the *-integral, just like the LS-integral, becomes the Darboux-Stieltjes integral (see, for example, \cite{Lev}).


\medskip
  
{\bf{3}}.\;Our definition of solution of problem (\ref{4eqsk}), defined in subsections 
{\bf{15.3.1,\,15.3.2}}
would be rather formal if it did not prove a continuous dependence of the solution on the antiderivatives of the coefficients of (\ref{4eqsk}). 
This kind of continuity can be obtained using Theorem \ref{okpredt}  (as was done in
\cite{derr88},\, \cite{derr18}). However, in this case the ``approximating'' problems will also be problems 
with generalized functions in the coefficients,
since both the uniform and the quasi-uniform convergence in Theorem \ref{okpredt}  do not allow to approximate discontinuous functions by
absolutely continuous functions. In order for the ``approximating'' problems to be problems with summable coefficients,
we will restrict slightly the class of problems under consideration (\ref{4eqsk})\;(i.e.\,replace conditions 
$\mathfrak{C}$ with $\mathfrak{C}_{\delta}$)
and will use Theorem \ref{osneps}. We will need the following result from article \cite{Lev}, which we present here  without proof.

Assertion $U$: {\it Consider a sequence of Cauchy problems
\begin{equation}
\label{ZCI}
\dot X_m=\mathcal{A}_m(t)X_m+\mathcal{F}_m(t),\qquad X_m(a)=I,\;\;\;m=0,1,2,\ldots,
\end{equation}
$(a\leqslant t\leqslant b),$ where entries of matrix $\mathcal{A}_m(\cdot)$ $\bigl($components of vectors $\mathcal{F}_m(\cdot)\bigr),
\;\;m=0,1,2,\ldots$ are summable functions on $[a,b]$, moreover, entry-wise $($component-wise$)$
$$
\mathcal{A}_m(t)\to \mathcal{A}_0(t),\quad \mathcal{F}_m(t)\to \mathcal{F}_0(t) \quad (m\to\infty)
$$
almost everywhere with respect to the Lebesgue measure. We also assume that all function sequences have
majorants summable in the sense of Lebesgue. Then all problems $(\ref{ZCI})$ are uniquely solvable, the  components of solutions
$X_m(\cdot)\;\;(m=0,1,2,\ldots)$ are absolutely continuous, and} $X_m(t)\;\rightrightarrows\;X(t).$

In order to be able to apply Theorem \ref{osneps}, let us go back to the continuous parameterization.

Recall that if
conditions $\mathfrak{C}_{\delta}$ are  satisfied, then conditions $\mathfrak C$ and $\mathfrak B$ also hold.
Under these conditions, we represent $p_k(\cdot)\;\;(k=1,\ldots,n+1)$
as the sum of right-continuous and left-continuous functions and consider functions
$(p_k)_{\varepsilon}(\cdot)\;\;(\varepsilon>0,\;k=1,\ldots,n+1),$ defined by equalities (\ref{apr}).

Consider a family of classical Cauchy problems
$$
(x)_{\varepsilon}^{(n)}+(p_1)'_{\varepsilon}(t) (x)_{\varepsilon}^{(n-1)}+\ldots +
(p_n)'_{\varepsilon}(t)(x)_{\varepsilon}=(p_{n+1})'_{\varepsilon}(t),\;\;(x)_{\varepsilon }^{(k-1)}(a)=\gamma_{k},
$$
$k=1,\ldots,n.$
\begin{teo}
\label{apprcont}
{\it Let conditions $\mathfrak{C}_{\delta}$ be satisfied. Then}
$$
_{\mathcal{E}_n}(x)_{\varepsilon}(t)\;\rightrightarrows\; _{\mathcal{E}_n}x(t)\quad
(\varepsilon\to 0).
$$
\end{teo}
\doc\,
In view of (\ref{aprsv}), (\ref{shvar+}),
\begin{equation}
\label{aprsv+}
(p_k)_{\varepsilon}\in {\bf{AC}},\;\;\;(p_k)_{\varepsilon}(t)\to p_k(t)\;\;\;(k=1,\ldots ,n+1,\;\;\;\varepsilon\to 0)
\end{equation}
(and also
$\;\;\;\bigvee\limits_a^b\bigl((p_k)_{\varepsilon}\bigr)\to \bigvee\limits_a^b\bigl(p_k\bigr)\;(k=2,\ldots,n+1)).\;$
Substituting 
in definitions (\ref{strKdU1}) --- (\ref{strKdU7}),\;(\ref{pravchast}) functions $(p_k)_{\varepsilon}\;\;(k=1,\ldots,n+1)$  instead of functions $p_k(\cdot)$, we obtain functions $(h_{ik})_{\varepsilon}(\cdot),\;(p_{ik})_{\varepsilon}(\cdot),$ and also functions $(a_{ik})_{ \varepsilon}(\cdot)$
and solution $(x)_{\varepsilon}$ of problem 
(\ref{4eqsk}) according to Theorem \ref{svedkdu2}.
As $\;\varepsilon\to 0\;$, we have by (\ref{aprsv+})
$$
(h_{nn})_{\varepsilon}(t)\to h_{nn}(t),\;\;\;(h_{n,\,n-k+1})_{\varepsilon}(t)\to h_{n,\,n-k+1}(t)\;\;\;(k=2,\ldots,n+1),
$$ 
see Theorem \ref{osneps}),
$(h_{kk})_{\varepsilon}(t)\to h_{kk}(t),\;$
$\;(h_{k,\,k-i+1})_{\varepsilon}(t)\to h_{n,\,n-k+1})(t)\linebreak
k=n -1,\ldots,2;\;i=2,\ldots,k;\;t\in [a,b]\;$
(according to the Lebesgue theorem). Moreover, $(h_{kk})_{\varepsilon}(\cdot)\;\;(k=2,\ldots, n)$ are bounded from above and from below by positive constants. Further,
due to the arithmetic properties of the limit, $$(p_{ik})_{\varepsilon}(t)\to p_{ik}(t)\;\;\;(i=1,\ldots,n ,\;k=
0,1,\ldots,i),
$$
$$
(a_{ik})_{\varepsilon}(t)\to a_{ik}(t)\;\;
(\varepsilon\to 0,\;i=1,\ldots,n,\;k=1,\ldots,i+1,
t\in[a,b]).$$ Moreover,
$(p_{kk})_{\varepsilon}(\cdot)\;\;(k=2,\ldots, n)$
are bounded from above and from below by positive constants (see (\ref{strKdU5})), and
$(a_{ik})_{\varepsilon}(\cdot)$ are bounded. Thus, the conditions of Assertion $U$ stated above are satisfied, so 
$$_{\mathcal{E}_n}(x)_{\varepsilon}(t)\;\rightrightarrows\; _{\mathcal{E}_n}x(t).$$
\hfill $\square$

If conditions of Theorem \ref{apprcont} are satisfied, we will say that solution $x(\cdot)$ 
is $\delta$-{\it correct}. Conditions $\mathfrak{C}_{\delta}$ thus entail  $\delta$-correctness of the solution.


\medskip

{\bf{4}}.\;
Now let conditions $\mathfrak{D}$  be satisfied. The *-integrals in (\ref{strKdU2}),
(\ref{pravchast}), which define the last row of  matrix $H(\cdot)$ and vector $h_0(\cdot)$ in  the right-hand side,  exist (see (\ref{defali2})).
The other definitions in sections 
{\bf{15.2.2,\;15.3.1,\;15.3.2}} 
do not change. The properties and the results in section {\bf{3}} change (slightly) 
in an obvious  way. The definitions of solutions in sections 
{\bf{15.3.1,\;15.3.2}}
do not change. 
However, in this case it is not possible to prove
the $\delta$-correctness of solution (see, for example, Corollary \ref{sledstosn}).


\begin{center}
\begin{large}
\section{General form of continuous linear functional}
\end{large}
\end{center}

\poonkt{General form of continuous linear functionals in spaces ${\bf{R}},\,{\bf{N}}$\; } \
\phantom{1234567890123456789012345678901234567890}

{\bf{1.}}\; Let us consider linear continuous functionals on $B$-spaces ${\bf{BV}}$, ${\bf{R}}$,
${\bf{N}}.$
Let us put
\begin{equation}
\label{flBVdef}
f(x)\doteq (*)\!\int\limits_a^b \mathfrak{f}(t)\,dx(t),\quad \mathfrak{f}\in{\bf{N}},
\end{equation}
and
\begin{equation}
\label{flNdef+}
g(x)\doteq (*)\!\int\limits_a^b x(t)\,d\mathfrak{g}(t),\quad \mathfrak{g}\in{\bf{BV}}.
\end{equation}
These are linear bounded (and, therefore, continuous) functionals on spaces ${\bf{BV}}
\bigl($in the case of (\ref{flBVdef})$\bigr)$ and
${\bf{R}},\;{\bf{N}}\;\;\bigl($in the case of (\ref{flNdef+})$\bigr).$

Linearity follows from statement {\bf{14.1.2}} a),\, boundedness follows from statement {\bf{14.1.3}}:
\begin{equation}
\label{Nozenflsv}
\|f\|\leqslant \|\mathfrak f\|_{\bf{N}}\qquad \text{and}\qquad \|g\|\leqslant \|\mathfrak{g}\|_{\bf{BV}}.
\end{equation}


Let us find the norm of functional (\ref{flBVdef}). To do this, it remains to prove the inequality opposite to the first inequality from (\ref{Nozenflsv}).

Since $\|\mathfrak f\|_{{\bf{N}}}=\sup_{t \in [a,b]}|\mathfrak f(t)|,$ 
for any $\varepsilon>0$ there is a point $t_{\varepsilon}\in[a,b]$ such that 
$|\mathfrak f(t_{\varepsilon})|>\|f\|_{{\bf{N}}}-\varepsilon$. We select
$\;\widehat x(t)=\mathfrak h_{t_{\varepsilon}}(t)$. Then $\|\widehat x\|_{{\bf{BV}}}=1\; $ and
$$
\|f\|=\underset{x\in{\bf{BV}},\,\|x\|=1}{\sup}\;|f(x)|\geqslant |f(\widehat x )|=
\left|(*)\!\int\limits_a^b \mathfrak{f}(t)\,d\widehat x(t)\right|=|\mathfrak f(t_{\varepsilon})|>\| f\|_{{\bf{N}}}-\varepsilon|.
$$
Due to the arbitrariness of $\varepsilon>0,$ this means that $\|f\|\geqslant \|\mathfrak f\|_{{\bf{N}}},$ 
and, together with (\ref{Nozenflsv}), that
$\|f\|=\|f\|_{{\bf{N}}}.$

\medskip

Let us also find the norm of the functional (\ref{flNdef+}).

For any $\varepsilon>0$ there exists partition $\tau_1$ such that
$$v_{\tau_1}(\mathfrak g_c)>
\bigvee\limits_a^b(g_c)-\dfrac{\varepsilon}{2}.$$
Let us expand, if necessary, $\tau_1$ to $\tau\doteq\{t_k\}_{k=0}^m$ so that inequality
$$\;\;\sum\limits_{s\in\tau}\bigl|\sigma_s(\mathfrak g)\bigr|>\bigvee\limits_a^b(\mathfrak g_d)-\dfrac{\varepsilon}{2 },$$
$$\Bigl(\bigvee\limits_a^b(\mathfrak g_d)=
\sum\limits_{s\in\T(\mathfrak g)}\bigl|\sigma_s(\mathfrak g)\bigr|\Bigr)\;$$ 
holds $\bigl($see 
(\ref{predstvar1})
$\bigr)$.
The following inequality will also be satisfied:
\begin{equation}
\label{2rapb}
v_{\tau}(\mathfrak g_c)\;\;\Bigl(\geqslant v_{\tau_1}(\mathfrak g_c)\Bigr)\;>\bigvee\limits_a^b(g_c)-\dfrac{\varepsilon }{2}.
\end{equation}
Further, we put
$$
\widehat x(t)=
\left\{
\begin{array}{lcr}
0\;\;\text{for}\;\;t\in T\doteq T(\mathfrak g)\setminus\tau,\\
{\rm sign}\,\sigma_{t_k}(g)\;\;\text{for}\;\;t\in\tau,\\
{\rm sign}\,\left(\mathfrak{g}(t_k)-\mathfrak{g}(t_{k-1})\right)\;\;\text{for}\;\;t\in(t _ {k-1},\,t_{k})\setminus T.\\
\end{array}
\right.
$$
Obviously, $\widehat x\in{\bf{N}},\;\|\widehat x\|_{{\bf{N}}}=1.$ We have (taking into account (\ref{2rapb}))
\begin{multline*}
\|g\|=\underset{x\in{\bf{N}},\,\|x\|=1}{\sup}\;|g(x)|\geqslant |g(\widehat x )|=\left|(*)\!\int_a^b\widehat x(t)\,d\mathfrak{g}(t)\right|=\\=
\left|\int_a^b\widehat x(t)\,d\mathfrak{g_c}(t)+\sum\limits_{t\in T(\mathfrak g)}\widehat x(t)\sigma_{t }(g)\right|=\\=
\left|\sum\limits_{k=1}^m\bigl|\mathfrak{g}(t_k)-\mathfrak{g}(t_{k-1})\bigr|+\sum\limits_{k= 0}^m\bigl|\sigma_{t_k}(g)\bigr|\right|=\\=
v_{\tau}(\mathfrak g_c)+\sum\limits_{t\in\tau}\bigl|\sigma_t(\mathfrak g)\bigr|>
\bigvee\limits_a^b(g_c)-\dfrac{\varepsilon}{2}+\bigvee\limits_a^b(g_d)-\dfrac{\varepsilon}{2}=\bigvee\limits_a^b(g)-\varepsilon
\end{multline*}
(see the 
equality (\ref{var22})). 
Since $\varepsilon$ is arbitrary, we obtain the inequality opposite to the second inequality from
(\ref{Nozenflsv}), that is, equality $\|g\|=\linebreak=
\|\mathfrak g\|_{{\bf{BV}}}\;$ holds.

\medskip

{\bf{2.}}\;
Let us show that equality (\ref{flNdef+}) is the general form of a linear continuous functional
on spaces ${\bf{R}},\,{\bf{N}}$ 
(cf.\,\cite{Tvrdy96}).

Below ${\bf{X}}_0$ denotes the subspace of a $B$-space ${\bf{X}}$ of functions $x(\cdot)$ having property $x(a) =0.$
\begin{teo}
\label{soprNR}
{\it The space conjugate to the Banach spaces ${\bf{N}}$ and ${\bf{R}}$ is the space ${\bf{BV}}_0$ }
\end{teo}
\doc\,
Equality (\ref{flNdef+}) defines mapping
\begin{equation}
\label{mathfrak G}
\mathfrak G:{\bf{BV}}_0\to{\bf{N}}^{*},\quad \mathfrak G(\mathfrak g)=g\quad (\mathfrak{g}\in{\bf{BV}}_0,\;\;\;g\in{\bf{N}}^{*}).
\end{equation}

Let us present the statement of the theorem in the form of the following chain of simpler statements.

1.\;{\it Mapping $(\ref{mathfrak G})$ is injective.}

Let $\;\mathfrak{g}_1,\,\mathfrak{g}_2\in{\bf{BV}}_0$, $\mathfrak{g}_1\ne\mathfrak{g}_2$.
Let us denote $\mathfrak h(t)\doteq\mathfrak g_1(t)-\mathfrak g_2(t).$

If $\mathfrak h_c(t)\equiv 0,$ then there is $t_0\in[a,b]$ such that $\sigma\doteq 
|\sigma_{t_0}(\mathfrak h)|>0$. Putting
$x_0(t)={\rm sign}\,\sigma_{t_0}(\mathfrak h)$ for $t=t_0,\;\;x_0(t)=0$ for $t\ne t_0,$ we obtain (see ( \ref{flNdef+}))
that $(g_1-g_2)(x_0)=|\sigma_{t_0}(\mathfrak h)|>0.$

Let $\mathfrak h_c(t)$ be not identically zero. Then there is a point $t_1\in[a,b]$ such that $\mathfrak h_c(t_1)\ne 0.$
In this case we assume
$$x_0(t)=\left\{
\begin{array}{lcr}
1\;\;\text{for}\;\;t\in [0,t_1]\setminus T(\mathfrak h),\\
0\;\;\text{for}\;\;t\in [0,t_1]\bigcap T(\mathfrak h),\\
0\;\;\text{for}\;\;t_1<t\leqslant 1,\\
\end{array}
\right.
$$
Then
\begin{multline*}
|(g_1-g_2)(x_0)|=(*)\!\int\limits_a^b x_0(t)\,d\mathfrak h(t)= \\ =
\int\limits_a^b x_0(t)\,d\mathfrak h_c(t)+\sum\limits_{s\in [0,t_0]\bigcap T(x_0)}x_0(s)\sigma_s(\mathfrak h )=\mathfrak h_c(t_0)\ne 0.
\end{multline*}
Thus, $g_1\ne g_2.$

2.\;{\it The mapping $(\ref{mathfrak G})$ is surjective.}

Let $f\in\mathcal L\bigl({\bf{N}},\mathbb R\bigr)\;\bigl(={\bf{N}}^{*}\bigr).$

Put
$u_0(t)\equiv 0.$ For $s\in (0,1],\;$
$$u_s(t)=\left\{
\begin{array}{lcr}
1\;\;\text\;\;{at}\;\;0\leqslant t<s,\\
0\;\;\text\;\;{at}\;\; t\leqslant 1.\\
\end{array}
\right.
$$
For each $s$,\;$u_s(\cdot)\in {\bf{N}}[0,1].$
Let us define function $\mathfrak f:[0,1]\to \mathbb R$ by $\mathfrak f(s)\doteq f(u_s).$ Note that $\mathfrak f(0)=0.$
Let
$$\tau=\{s_k\}_{k=0}^n, \quad 0=s_0<s_1<...<s_n=1,$$ be a partition of interval $[0,1]$ and let $\sigma_k\doteq {\rm sign}\,\bigl(\mathfrak f(s_k)-
\mathfrak f(s_{k-1})\bigr).$
Then
\begin{multline}
\notag
v_{\tau}(\mathfrak f)=\sum\limits_{k=1}^n|\mathfrak f(s_k)-\mathfrak f(s_{k-1})|=
\sum\limits_{k=1}^n\sigma_k(\mathfrak f(s_k)-\mathfrak f(s_{k-1}))=\\=
\sum\limits_{k=1}^n\sigma_k\left(f(u_{s_k})-f(u_{s_{k-1}})\right)=
f\left(\sum\limits_{k=1}^n\sigma_k\bigl(u_{s_k}(t)-u_{s_{k-1}}(t)\bigr)\right)\leqslant \\ \leqslant
 \|f\|\cdot\|y\|_{{\bf{N}}}=\|f\|,
\end{multline}
where $y(t)=\sum\limits_{k=1}^n\sigma_k\bigl(u_{s_k}(t)-u_{s_{k-1}}(t)\bigr),\;y \in{\bf{N}},\;\|y\|_{{\bf{N}}}=1$

Thus, $\mathfrak f\in {\bf{BV}}_0$. We obtain mapping
\begin{equation}
\label{mathfrak H}
\mathfrak H:{\bf{N}}^{*}\to{\bf{BV}}_0,\quad \mathfrak H(f)=\mathfrak f\quad (f\in{\bf{N} }^{*},\;\;\;\mathfrak f\in{\bf{BV}}_0)
\end{equation}
which assigns to each linear continuous functional on ${\bf{N}}$ the function $\mathfrak f\in {\bf{BV}}_0$ according to the  procedure describe above,\; $\mathfrak H(f)=\mathfrak f.$

Let $g\in {\bf{N^{*}}}[0,1]$ be defined by formula (\ref{flNdef+}) using function $\mathfrak f$ constructed from the functional
$f\,\,\;\Bigl(g=\mathfrak G(\mathfrak f)=\mathfrak G\bigl(\mathfrak H(f)\bigr)\Bigr).$ We successively obtain
\begin{multline*}
g(u_s)=(*)\!\int\limits_a^b u_s(t)\,d\mathfrak f(t)=(*)\!\int\limits_a^s\,d\mathfrak f(t) =\\=
\int\limits_a^s\,d\mathfrak f_c(t)+\sum\limits_{t\in T(\mathfrak f)\bigcap [a,s)}\sigma_t(\mathfrak f)+\sigma_s^{ -}(\mathfrak f)=
\mathfrak f_c(s)+\mathfrak f_d(s)=\mathfrak f(s).
\end{multline*}
This means that $g=f,\;$, that is,  mapping $\;\mathfrak G\;$ is surjective.

3.\;${\bf{N}}^{*}={\bf{R}}^{*}={\bf{BV}}_0.$

Due to the linearity of the *-integral and equality
$\|g\|=\|\mathfrak g\|_{{\bf{BV}}}$, $\mathfrak G$ is an isometric isomorphism 
from ${\bf{BV}}_0$ to ${\bf{N}}^{*}.$ This  means that
${\bf{N}}^{*}={\bf{BV}}_0.\;$ Since ${\bf{C}}\subset{\bf{R}}\subset{\bf{N}},$ then 
${\bf{N}}^{*}\subset{\bf{R }}^{*}\subset{\bf{C}}^{*}={\bf{BV}}_0$. Therefore,
${\bf{N}}^{*}={\bf{R}}^{*}={\bf{C}}^{*}={\bf{BV}}_0.$
\hfill $\square$

From the proof of the last theorem, we obtain right away the following statements.
\begin{sle}
\label{vzaimnobr}
{\it The mappings $\mathfrak G$ and $\mathfrak H$ are inverses of each other.}
\end{sle}
\begin{sle}
\label{obsvidlofN}
{\it Equality $(\ref{flNdef+})$ gives the general form of a linear continuous functional on Banach spaces $\,{\bf{N}}$ and ${\bf{R}}.$}
\end{sle}
It is interesting to compare the statement of Corollary \ref{obsvidlofN} with the accuracy of the *-pair $({\bf{N}},\,{\bf{BV}})\;$ (see \cite{derr20}).

Note that $B$-space ${\bf{N}}$ is not reflexive:
$$
{\bf{N}}^{**}={\bf{BV}}_0^{*}\ne{\bf{N}}.
$$

Let $\mathfrak d(t)\doteq\chi_{\mathbb Q}(t)$ be the Dirichlet function,\;$T(d)=[a,b],\;\;\mathfrak{d }\notin{\bf{N}}[a,b].$
However, the LS integral
\begin{equation}
\label{flBVdef!}
d(x)\doteq (LS)\!\int\limits_a^b \mathfrak{d}(t)\,dx(t),
\end{equation}
exists for every $x\in{\bf{BV}}_0,$ hence $d$ is a linear continuous functional
on space ${\bf{BV}}_0$, i.e. $d\in {\bf{BV}}_0^{*}.$ 
At the same time, (*)-integral $(*)\!\int\limits_a^b \mathfrak{d}(t)\,dx(t)$ does not exist.

\newpage

\begin{center}
\begin{large}
\section{Appendix. Proofs of inequalities}
\end{large}
\end{center}

\poonkt{\;H\"{o}lder's inequality\;}\phantom{01234567890123456789} \

Let $p,\,q\in \mathbb R$ satisfy  relations
\begin{equation}
\label{nera1}
1< p<+\infty, \quad 1< q<+\infty, \quad \frac{1}{p}+\frac{1}{q}=1,
\end{equation}
let $T=\{1,2,\ldots,n\}$\;\;or\;\;$T=\mathbb N=\{1,2,\ldots\},$
$x_k,\,y_k\;\;(k\in T)$ are complex numbers. In the case $T=\mathbb N$ it is assumed that the series
$\sum\limits_{k\in T}x_k,\;\sum\limits_{k\in T}y_k$ converge absolutely.
For such $x_k,\, y_k$, the following H\"{o}lder's inequality holds:
\begin{equation}\label{neraH_0}
\sum\limits_{k\in T}|x_k\,y_k)|\leqslant \left(\sum\limits_{k\in T}|x_k|^p\right)^{\frac{1}{p} }\cdot
\left(\sum\limits_{k\in T}|y_k|^q\right)^{\frac{1}{q}}.
\end{equation}
Proof. 1. First, we prove an auxiliary inequality between real numbers
\begin{equation}\label{nera2}
a^{\frac{1}{p}}\cdot b^{\frac{1}{q}}\leqslant \frac{a}{p}+\frac{b}{q}\qquad (a> 0,\;b>0).
\end{equation}

Let $0<\alpha <1.$ Consider function
$$
f:[0,+\infty)\to \mathbb R,\quad f(t)=t^{\alpha}-\alpha t.
$$
Since its derivative $f'(t)=0$ at $t=1$, and when we pass through this point
the derivative changes its sign from $+$ to $-,$ the function $f$ has a maximum at this point equal to $1-\alpha$. Therefore,
$$
t^{\alpha}-\alpha t\leqslant 1-\alpha\qquad (t\geqslant 0,\;0<\alpha <1).
$$
Let $a>0,\;b>0$, put $t\doteq \frac{a}{b},\quad \beta\doteq 1-\alpha\;(0<\beta <1,\;\alpha +\beta =1).$ 
Then from the above inequality we successively obtain
$$
\frac{a^{\alpha}}{b^{\alpha}}-\alpha \frac{a}{b}\leqslant \beta,\quad a^{\alpha}\cdot b^{\beta} \leqslant \alpha a+\beta b.
$$
After substituting $\alpha =\frac{1}{p},\;\beta=\frac{1}{q},$ we arrive at (\ref{nera2}).

2. First let $x_k,\,y_k\;\;(k\in T)$ be such that
\begin{equation}\label{nera3}
\sum\limits_{k\in T}|x_k|^p=\sum\limits_{k\in T}|y_k|^q=1.
\end{equation}
 Under this condition, inequality (\ref{neraH_0}) takes  form
\begin{equation}
\label{nera4}
\sum\limits_{k\in T}|x_k\,y_k|\leqslant 1.
 \end{equation}
Let us prove it.

Let us substitute $|x_k|^p$ and $|y_k|^q$ into inequality (\ref{nera2}) instead of $a$ and $b$, respectively (for every $k\in T$). Then
$$
|x_k\,y_k|=\bigl(|x_k|^p\bigr)^{\frac{1}{p}}\cdot \bigl(|y_k|^q\bigr)^{\frac{1}{ q}}\leqslant
\frac{|x_k|^p}{p}+\frac{|y_k|^q}{q}\quad (k\in T).
$$
Summing up the resulting inequality (for $T=\mathbb N$ we also need to pass to the limit)
taking into account  equalities (\ref{nera3}) and (\ref{nera1}), we obtain (\ref{nera4}).

3. Now let $x_k,\,y_k\;\;(k\in T)$ be arbitrary. We put
$$
\widetilde x_k\doteq \frac{x_k}{\left(\sum\limits_{j\in T}|x_j|^p\right)^{\frac{1}{p}}},\qquad
\widetilde y_k\doteq \frac{y_k}{\left(\sum\limits_{j\in T}|y_j|^q\right)^{\frac{1}{q}}}.
$$
Then, clearly, $\widetilde x_k$ and $\widetilde y_k\;\;(k\in T)$ satisfy conditions (\ref{nera3}). So, according to what we already proved,
$$
\sum\limits_{k\in T}\frac{|\widetilde x_k\widetilde y_k|}{\left(\sum\limits_{j\in T}|x_j|^p\right)^{\frac{1 }{p}}\cdot
\left(\sum\limits_{j\in T}|y_j|^q\right)^{\frac{1}{q}}}\leqslant 1,
$$
from which (\ref{neraH_0}) immediately follows.
\hfill $\square$

A more general inequality 
\begin{equation}
\label{neraH}
(*)\!\int\limits_a^b|x(t)\,y(t)|\,dg(t)\leqslant \left((*)\!\int\limits_a^b|x(t) |^p\,dg(t)\right)^{\frac{1}{p}}\cdot
\left((*)\!\int\limits_a^b|y(t)|^q\,dg(t)\right)^{\frac{1}{q}},
\end{equation}
where $x,\,y\in {\bf{N}},\;g\in {\bf{BV}}$, is proved similarly. Inequality (\ref{neraH_0}) can be obtained from (\ref{neraH}) by taking as $g(\cdot)$
appropriate jump function.


\poonkt{\;Minkowski's inequality\;}\phantom{01234567890123456789} \

Let $x,\,y\in {\bf{N}},\;g\in {\bf{BV}},\;\;1<p<\infty.$ Then the following inequality holds:
\begin{multline}
\label{neraM}
\left((*)\!\int\limits_a^b|x(t)+y(t)|^p\,dg(t)\right)^\frac{1}{p}\leqslant\\
\leqslant \left((*)\!\int\limits_a^b|x(t)|^p\,dg(t)\right)^{\frac{1}{p}}+ \left((* )\!\int\limits_a^b|y(t)|^p\,dg(t)\right)^{\frac{1}{p}}.
\end{multline}
It is called the Minkowski inequality.

Proof. Let $\;\;q=\dfrac{p}{p-1}\;\;$ $\bigl($see (\ref{nera1})$\bigr).\;$ Apply inequality (\ref {neraH}):
\begin{multline*}
\!(*)\!\int\limits_a^b|x(t)+y(t)|^p\,dg(t)\leqslant (*)\!\int\limits_a^b\bigl(|x(t )|+|y(t)|\bigr)^p\,dg(t)\\=
(*)\!\int\limits_a^b|x(t)|\bigl(|x(t)|\!+\!|y(t)|\bigr)^{p-1}\,dg(t)\!
+\!(*)\!\int\limits_a^b|y(t)|\bigl(|x(t)|\!+\!|y(t)|\bigr)^{p-1}\,dg(t)\! \\
\leqslant\left((*)\!\int\limits_a^b|x(t)|^p\,dg(t)\right)^{\frac{1}{p}}\cdot
\left((*)\!\int\limits_a^b\bigl(|x(t)|+|y(t)|\bigr)^{(p-1)q}\,dg(t)\right )^{\frac{1}{q}}\\+
\left((*)\!\int\limits_a^b|y(t)|^p\,dg(t)\right)^{\frac{1}{p}}\cdot
\left((*)\!\int\limits_a^b\bigl(|x(t)|+|y(t)|\bigr)^{(p-1)q}\,dg(t)\right )^{\frac{1}{q}}.
\end{multline*}
Since $(p-1)q=p,$ we can divide both sides of the resulting inequality by the common factor of the right side, obtaining
\begin{multline*}
\left((*)\!\int\limits_a^b|x(t)+y(t)|^p\,dg(t)\right)^{1-\frac{1}{q}} \\
\leqslant\left((*)\!\int\limits_a^b|x(t)|^p\,dg(t)\right)^{\frac{1}{p}}
+\left((*)\!\int\limits_a^b|y(t)|^p\,dg(t)\right)^{\frac{1}{p}},
\end{multline*}
which, in view of equality $1-\frac{1}{q}=\frac{1}{p}$, yields (\ref{neraM}).
\hfill $\square$

Choosing appropriate jump function in (\ref{neraM}) as $g(\cdot)$, we obtain Minkowski's inequality in the form
\begin{equation}
\label{neraM_0}
\left(\sum\limits_{k\in T}|x_k+y_k|^p\right)^\frac{1}{p}\leqslant \left(\sum\limits_{k\in T}|x_k| ^p\right)^{\frac{1}{p}}+ \left(\sum\limits_{k\in T}|y_k|^p\right)^{\frac{1}{p}},
\end{equation}
However, inequality (\ref{neraM_0}) can also be obtained directly from inequality (\ref{neraH_0}) by following the proof
inequalities (\ref{neraM}).




\end{document}